\documentclass[reqno,oneside]{amsart} 
\pdfoutput=1
\usepackage{amsmath}
\usepackage{amsthm,
amssymb,amsmath,enumerate,stackrel,
verbatim,
}
\usepackage{xr-hyper}
\usepackage[all,cmtip]{xy}
\usepackage{enumitem}
\usepackage{soul}

\usepackage[pagebackref,hyperindex,linktocpage=true]{hyperref}
\hypersetup{
    colorlinks,
    linkcolor={red!50!black},
    citecolor={red!50!black}, 
    urlcolor={red!50!black}, 
    filecolor={red!50!black} 
}

\usepackage{tikz}
\usetikzlibrary{matrix,arrows}

\usepackage{tikz-cd}

\usepackage{leftidx}

\usepackage{theoremref}

\makeatletter
\@addtoreset{equation}{section}
\makeatother

\numberwithin{equation}{section}

\theoremstyle{definition}
\newtheorem{Defi}{Definition}[section] \newcommand{\defi}{\begin{Defi}} \newcommand{\xdefi}{\end{Defi}} \newtheorem{DefiLemm}[Defi]{Definition and Lemma} \newcommand{\defilemm}{\begin{DefiLemm}} \newcommand{\xdefilemm}{\end{DefiLemm}} 
\newtheorem{Bsp}[Defi]{Example} \newcommand{\exam}{\begin{Bsp}} \newcommand{\xexam}{\end{Bsp}} 
\newtheorem{Syno}[Defi]{Synopsis} \newcommand{\syno}{\begin{Syno}} \newcommand{\xsyno}{\end{Syno}} 
\newtheorem{Bem}[Defi]{Remark} \newcommand{\rema}{\begin{Bem}} \newcommand{\xrema}{\end{Bem}} 
\newtheorem{Notation}[Defi]{Notation} \newcommand{\nota}{\begin{Notation}} \newcommand{\xnota}{\end{Notation}} 
\newtheorem{Convention}[Defi]{Convention} \newcommand{\conv}{\begin{Convention}} \newcommand{\xconv}{\end{Convention}} 
\newtheorem{Warning}[Defi]{Warning} \newcommand{\warn}{\begin{Warning}} \newcommand{\xwarn}{\end{Warning}} 
\newtheorem{Situation}[Defi]{Situation} \newcommand{\situ}{\begin{Situation}} \newcommand{\xsitu}{\end{Situation}}
\newtheorem{Assumption}[Defi]{Assumption} \newcommand{\assu}{\begin{Assumption}} \newcommand{\xassu}{\end{Assumption}} 

\theoremstyle{plain}
\newtheorem{Theo}[Defi]{Theorem} \newcommand{\theo}{\begin{Theo}} \newcommand{\xtheo}{\end{Theo}} 
\newtheorem{Satz}[Defi]{Proposition} \newcommand{\prop}{\begin{Satz}} \newcommand{\xprop}{\end{Satz}} 
\newtheorem{Lemm}[Defi]{Lemma} \newcommand{\lemm}{\begin{Lemm}} \newcommand{\xlemm}{\end{Lemm}} 
\newtheorem{Coro}[Defi]{Corollary} \newcommand{\coro}{\begin{Coro}} \newcommand{\xcoro}{\end{Coro}}
\newtheorem{Ques}[Defi]{Question} \newcommand{\ques}{\begin{Ques}} \newcommand{\xques}{\end{Ques}}
\newtheorem{Conj}[Defi]{Conjecture} \newcommand{\conj}{\begin{Conj}} \newcommand{\xconj}{\end{Conj}}

\newcommand{\refsect}[1]{Section \ref{sect--#1}}
\newcommand{\refit}[1]{(\ref{item--#1})}
\newcommand{\refeq}[1]{(\ref{eqn--#1})}

\newcommand{\eqn}{\begin{equation}} \newcommand{\xeqn}{\end{equation}}
\newcommand{\eqnarr}{\begin{eqnarray*}} \newcommand{\xeqnarr}{\end{eqnarray*}}
\newcommand{\eqnarra}{\begin{eqnarray}} \newcommand{\xeqnarra}{\end{eqnarray}}

\newcommand{\pf}{\begin{proof}} \newcommand{\xpf}{\end{proof}}

\newif \ifDraft % whether it is a draft or otherwise meant to be sent to other people etc.
\Draftfalse

\ifDraft
  \usepackage[top=1.5cm, bottom=1.5cm, left=.5cm, right=7.5cm, heightrounded, marginparwidth=7cm, marginparsep=.2cm]{geometry}
\else
  \usepackage[top=1.5cm, bottom=1.5cm, 
  left=2cm, right=2cm, % for bug-fixing: before it was 2+2cm
  heightrounded]{geometry}

  \usepackage{letltxmacro}
  \LetLtxMacro\Oldfootnote\footnote
  \renewcommand{\footnote}[2][]{\relax}
\fi

\newif \ifHideSomeComments
\HideSomeCommentsfalse

\ifHideSomeComments
\else
\fi

\newcommand{\nc}{\newcommand}
\nc{\StP}[1]{\cite[Tag~\href{http://stacks.math.columbia.edu/tag/#1}{#1}]{StacksProject}}
\nc{\StPd}[2]{\cite[Tags~\href{http://stacks.math.columbia.edu/tag/#1}{#1}, \href{http://stacks.math.columbia.edu/tag/#2}{#2}]{StacksProject}} 

\nc{\on}{\operatorname}
\nc{\aff}{{\on{aff}}}
\nc{\modi}{{\on{mod}}} 
\nc{\even}{{\on{even}}}
\nc{\odd}{{\on{odd}}}
\nc{\naive}{{\on{naive}}}
\nc{\hofib}{\on{hofib}}
\nc{\Bun}{\on{Bun}}
\nc{\ad}{{\on{ad}}}
\nc{\lft}{{\on{lft}}}
\nc{\modulo}{\on{mod}} 

\nc{\Weil}{{\on{Weil}}} 
\nc{\FWeil}{{\on{FWeil}}} 
\nc{\cons}{{\on{cons}}} 
\nc{\tot}{{\on{Tot}}} 

\nc{\str}{\on{-}}
\nc{\perf}{{\on{perf}}}
\nc{\Rel}{{\on{Pos}}}
\nc{\lan}{\langle}
\nc{\ran}{\rangle}

\nc{\PM}{{\on{PM}}}
\nc{\calH}{{\mathcal H}}
\nc{\calP}{{\mathcal P}} 
\nc{\calC}{{\mathcal C}} 
\nc{\calF}{{\mathcal F}} 
\nc{\calB}{{\mathcal B}} 
\nc{\calI}{{\mathcal I}}
\nc{\calJ}{{\mathcal J}} 
\nc{\calO}{{\mathcal O}}
\nc{\calM}{{\mathcal M}}
\nc{\co}{\colon}

\newcommand{\category}[1]{\mathrm{#1}}

\newcommand\restr[2]{{ 
  \left.\kern-\nulldelimiterspace 
  #1 
  \vphantom{\big|} 
  \right|_{#2} 
  }}

\makeatletter 
\newcommand*{\doublerightarrow}[2]{\mathrel{
  \settowidth{\@tempdima}{$\scriptstyle#1$}
  \settowidth{\@tempdimb}{$\scriptstyle#2$}
  \ifdim\@tempdimb>\@tempdima \@tempdima=\@tempdimb\fi
  \mathop{\vcenter{
    \offinterlineskip\ialign{\hbox to\dimexpr\@tempdima+1em{##}\cr
    \rightarrowfill\cr\noalign{\kern.5ex}
    \rightarrowfill\cr}}}\limits^{\!#1}_{\!#2}}}
\newcommand*{\triplerightarrow}[1]{\mathrel{ 
  \settowidth{\@tempdima}{$\scriptstyle#1$}
  \mathop{\vcenter{
    \offinterlineskip\ialign{\hbox to\dimexpr\@tempdima+1em{##}\cr
    \rightarrowfill\cr\noalign{\kern.5ex}
    \rightarrowfill\cr\noalign{\kern.5ex}
    \rightarrowfill\cr}}}\limits^{\!#1}}}
\makeatother

\newcommand{\Ho}{\category{Ho}} 
\newcommand{\Mod}{\category{Mod}} 
\newcommand{\coMod}{\category{coMod}}

\newcommand{\ModZ}{\Mod_\Z} 
\newcommand{\ModL}{\Mod_\Lambda}

\renewcommand{\Pr}{\category{Pr}}
\newcommand{\PrL}{\Pr^\category{L}} 
\newcommand{\PrSt}{\Pr^{\category{St}}} 
\newcommand{\PrStZ}{\PrSt_\Z}

\newcommand{\Fun}{\category{Fun}} 
 
\newcommand{\Corr}{\category{Corr}} 
\newcommand{\Sch}{\category{Sch}} 
\newcommand{\AffSch}{\category{AffSch}}

\newcommand{\Alg}{\category{Alg}} 
 
\def\ft{\mathrm{ft}} 
 
\newcommand{\dual}{\vee} 
\newcommand{\Zar}{\mathrm{Zar}} 
\newcommand{\Set}{\category{Set}}
\newcommand{\Sat}{\category{Sat}}

\newcommand{\GT}{\hat{G}_1} %Extended dual group

\newcommand{\grAb}{\category{grAb}} 
\newcommand{\grModZ}{\grAb}
 
\def\Gm{\mathbf {G}_\mathrm m} 
\def\Gmmono{\Gm\mathrm{-mono}} 
\def\Ga{\mathbf {G}_\mathrm a}

\def\IC{\mathrm{IC}} 

\newcommand{\colim}{\operatornamewithlimits{colim}} 
\newcommand{\coeq}{\operatornamewithlimits{coeq}} 

\def\id{{\rm id}} 
\def\inv{{\rm inv}} 
\def\ev{{\operatorname {ev}}} 
 
\def\opp{{\rm op}} 
 
\def\To#1#2{\mathop{\count0=#1 \loop\ifnum\count0>0 \smash-\mkern-7mu \advance\count0 -1 \repeat \mathord\rightarrow}\limits^{#2}}

\def\Maps{\mathop{\rm Maps}\nolimits} 
 
\def\Char{\mathop{\rm char}\nolimits} 
\def\CH{\mathop{\rm CH}\nolimits} 
 
\def\av{\mathop{\rm av}\nolimits} 
\def\coav{\mathop{\rm coav}\nolimits} 
\def\Hom{\mathop{\rm Hom}\nolimits} 
\def\Du{\mathrm D} 
\def\PGL{\mathrm {PGL}} 
\def\GL{\mathrm {GL}}

\def\Ab{\mathop{\rm Ab}\nolimits}

\def\Gr{\mathop{\rm Gr}\nolimits} 
\def\Hck{\mathop{\rm Hck}\nolimits} 
\def\sgl{{\{*\}}} 
\def\Fl{\mathop{\rm Fl}\nolimits} 
\def\FlCal{\mathop{\mathcal{F}\ell}\nolimits}

\def\Ext{\mathop{\rm Ext}\nolimits} 
 
\def\IHom{\underline{\Hom}}

\def\Aut{\mathop{\rm Aut}\nolimits} 
 
\def\Map{\Maps} 
\def\Rep{\category{Rep}}

\def\et{\mathrm{\acute et}} 
\def\Nis{\mathrm{Nis}}

\def\coker{\operatorname{coker}} 
 
\def\adj{\mathrm{adj}} 
\def\sico{\mathrm{sc}} 
\def\der{\mathrm{der}} 
\def\CT{\mathrm{CT}} 
\def\cl{\mathrm{cl}} 
\def\sph{\mathrm{sph}} 

\def\Z{{\mathbf Z}} 
\def\Fp{{{\mathbf F}_p}} 
\def\Fq{{{\mathbf F}_q}}

\def\Q{{\mathbf Q}} 
\def\Qp{\Q_p} 
\def\Zp{\Z_p} 
 
\def\Qq{{\overline \Q}} 
 
\def\RR{{\bf R}} 
\def\C{{\mathbf {C}}} 
\def\A{{\bf A}} 
\renewcommand{\P}[1][1]{\mathbf P^{#1}} 
\def\Gm{\mathbf {G}_\mathrm m}

\def\qq{{\mathbf{q}}} 

\def\Aa{{\mathcal{A}}}

\def\Ee{{\mathcal{E}}}
\def\Ff{{\mathcal{F}}}
\def\calH{{\mathcal{H}}}
\def\Jj{{\mathcal{J}}}
\def\Mm{{\mathcal{M}}}
\def\Pp{{\mathcal{P}}}
\def\Qq{{\mathcal{Q}}}
\def\Ss{{\mathcal{S}}}

\def\Uu{{\mathcal{U}}}

\def\H{{\rm H}} 
\def\pe{{{}^\mathrm{p}} \! } 
\def\pH{\pe\H}

\def\im{{\rm im}} 
\def\SH{\category{SH}}

\def\R{\mathrm{R}} 
\def\Betti{\mathrm{B}} 
\def\red{\mathrm{r}} 
\def\redx{(\red)} 
\def\DM{\category{DM}} 
\def\DMr{\DM_\red}

\def\DMrx{\DM_{\redx}}

\def\DTM{\category{DTM}} 
\def\DTMr{\DTM_\red} 
\def\DTMrx{\DTM_{\redx}}

\def\sw{\rm sw} 
 
\def\anti{\mathrm{anti}} 

\def\Vect{\category{Vect}} 
\def\fd{\category{fd}}

\def\Perv{\category{Perv}} 
\def\MTM{\category{MTM}} 
\def\MTMr{\MTM_\red} 
\def\MTMrx{\MTM_{\redx}} 
 
\def\Sat{\category{Sat}} 
\def\Satr{\Sat_\red} 
\def\Satrx{\Sat_{\redx}}

\def\ii{$\infty$}

\def\fil{{\rm fil}}

\def\SL{{\rm SL}}

\def\Spec{\mathop{\rm Spec}}

\newcommand{\M}{\mathrm{M}}

\newcommand{\comp}{\mathrm{c}} 
\newcommand{\locc}{\mathrm{lc}}

\newcommand{\fraka}{\mathfrak{a}} 
\newcommand{\PreStk}{\category{PreStk}}

\newcommand{\an}{\mathrm{an}} 
 
\newcommand{\D}{\category{D}}

\def\sbuildrel#1\over#2{\mathrel{\smash{\mathop{\kern0pt #2}\limits^{#1}}}}

\let\x\times
\let\ol\overline
\renewcommand{\t}{\otimes}

\newcommand{\bx}{\boxtimes}

\renewcommand{\r}{\rightarrow}

\newcommand{\hr}{\hookrightarrow}
\newcommand{\pot}[1]{ [\hspace{-0,5mm}[ {#1} ]\hspace{-0,5mm}] }
\newcommand{\rpot}[1]{ (\hspace{-0,7mm}( {#1} )\hspace{-0,7mm}) }

\mathchardef\mhyphen="2D

\def\matrix#1{\null\,\vcenter{\normalbaselines
    \ialign{\hfil$##$\hfil&&\quad\hfil$##$\hfil\crcr
      \mathstrut\crcr\noalign{\kern-\baselineskip}
      #1\crcr\mathstrut\crcr\noalign{\kern-\baselineskip}}}\,}

\newdimen\harrowsize
\def\mapright#1{\smash{\mathop{\hbox to\harrowsize{\rightarrowfill}}\limits^{#1}}}

{\catcode`@=11
\gdef\cal{\fam\tw@}

\global\let\over\@@over
\global\let\atop\@@atop
\global\let\above\@@above
\global\let\overwithdelims\@@overwithdelims
\global\let\atopwithdelims\@@atopwithdelims
\global\let\abovewithdelims\@@abovewithdelims
\gdef\eqalign#1{\null\,\vcenter{\openup\jot\m@th
\ialign{\strut\hfil$\displaystyle{##}$&$\displaystyle{{}##}$\hfil
      \crcr#1\crcr}}\,}
\newskip\xcentering \global\xcentering=0pt plus 1000pt minus 1000pt
\gdef\eqalignno#1{\displ@y \tabskip\xcentering
  \halign to\displaywidth{\hfil$\@lign\displaystyle{##}$\tabskip\z@skip
    &$\@lign\displaystyle{{}##}$\hfil\tabskip\xcentering
    &\llap{$\@lign##$}\tabskip\z@skip\crcr
    #1\crcr}}
\global\def\cases#1{\left\{\,\vcenter{\normalbaselines\m@th
    \ialign{$##\hfil$&\quad##\hfil\crcr#1\crcr}}\right.}
\gdef\eqlabel#1{\refstepcounter{equation}\label{eqn--#1}\eqno\hbox{\@eqnnum}}
}

\pagecolor{white}

\begin{document}

\author{Robert Cass, Thibaud van den Hove and Jakob Scholbach}

\title{The geometric Satake equivalence\\for integral motives}

\begin{abstract} We prove the geometric Satake equivalence for mixed Tate motives over the integral motivic cohomology spectrum. This refines previous versions of the geometric Satake equivalence for split reductive groups. Our new geometric results include Whitney--Tate stratifications of Beilinson--Drinfeld Grassmannians and cellular decompositions of semi-infinite orbits. With future global applications in mind, we also achieve an equivalence relative to a power of the affine line. 
Finally, we use our equivalence to give Tannakian constructions of Deligne's modification of the dual group and a modified form of Vinberg's monoid over the integers.
\end{abstract}

\maketitle

\tableofcontents

\section{Introduction}

\subsection{Geometry of the Beilinson--Drinfeld affine Grassmannians} 
The affine Grassmannian $\Gr_G$ associated to a split reductive group $G$ lies at the nexus of the Langlands program, representation theory, and the study of moduli of \(G\)-bundles on a curve.
Its topology encodes representation-theoretic information via the geometric Satake equivalence, which roughly says that \(L^+G\)-equivariant perverse sheaves on $\Gr_G$ are equivalent to algebraic representations of the Langlands dual group. 
We refer to \cite[§1.1.2]{BaumannRiche:Satake} for the history of this equivalence, and to \cite[§1]{RicharzScholbach:Motivic} for an overview of how motives fit into the picture. 

This paper contributes to our understanding of $\Gr_G$ and a notable family of higher-dimensional analogues, the Beilinson--Drinfeld affine Grassmannians $\Gr_{G, I}$. 
These are the moduli spaces of $G$-torsors on the curve $X = \A^1$, equipped with a trivialization away from $|I|$ many given points in $X$, where $I$ is a finite non-empty set.
This object interpolates between powers of $\Gr_G$, making it the basis of factorization techniques in the geometric Langlands program. In number theory, its various incarnates are essential in the geometrization of the Langlands program over fields such as $\Fq(t)$ and $\Qp$ \cite{LafforgueV:Chtoucas, FarguesScholze:Geometrization}. It also plays a key role in Bezrukavnikov's equivalence \cite{Bezrukavnikov:Two}, which provides a (tamely ramified) geometric Langlands correspondence over local fields of equicharacteristic (see also \cite{ArkhipovBezrukavnikov:Perverse, Gaitsgory:Central}).

Concerning the geometry of $\Gr_{G,I}$, we prove the following key theorems.

\theo[\thref{BD.WT}] \thlabel{BD.WT.intro}
The stratification of $\Gr_{G,I}$ in \thref{BD.strata}, which combines the Schubert cells in \(\Gr_G\) with a stratification of $X^I$, is admissible Whitney--Tate in the sense of Definitions \ref{Whitney--Tate} and \ref{admissible.stratification}.
\xtheo

This implies that there is a viable subcategory $\DTM(\Gr_{G,I})$ of stratified Tate motives inside the stable \ii-category $\DM(\Gr_{G,I})$ of all motivic sheaves, with integral coefficients, on $\Gr_{G,I}$. 
Objects in this subcategory are precisely those motivic sheaves whose restriction to each stratum is Tate, which means that they can be constructed out of twists $\Z(n)$, $n \in \Z$, of the constant sheaves on the strata. This theorem also implies that if $G$ is defined over a scheme $S$ satisfying the Beilinson--Soulé vanishing condition, such as $S = \Spec \Z$, or a finite field, or a global field, there is an abelian subcategory $\MTM(\Gr_{G,I}) \subset \DTM(\Gr_{G,I})$, akin to the abelian subcategory of perverse sheaves inside the derived category of sheaves.
This theorem is nontrivial since the closures of the Schubert cells are usually singular. The proof builds on the techniques of \cite{RicharzScholbach:Motivic}, which treats the case $I = \sgl$.

A central novelty of this paper is the application of constant term functors in the context of Tate motives.
For a maximal split torus and Borel \(T\subseteq B\subseteq G\) defined 
over \(S\), there is a Jacquet functor \refeq{CT.Def}  \(\CT : \DM_{\Gm}(\Gr_{G,I})\to \DM(\Gr_{T,I})\), for a certain \(\Gm\)-action on \(\Gr_{G,I}\).  
Using this functor for arbitrary sheaves is a standard idea in geometric representation theory and the classical geometric Satake equivalence \cite[§1.10.1]{BaumannRiche:Satake}.
However, it did not appear in the motivic context of \cite{RicharzScholbach:Motivic}, and it is highly delicate to decide whether $\CT$ preserves the much smaller category of stratified Tate motives. 
We achieve this by proving the following geometric result. 

\theo[\thref{cellularity of intersection:torus}]\thlabel{intro.thm cellularity}
The intersection of any Schubert cell with any semi-infinite orbit in \(\Gr_G\) (whose irreducible components are also known as Mirkovi\'c--Vilonen cycles) admits a filtrable decomposition into cellular schemes, i.e., products of copies of $\A^1$ and $\Gm$.
\xtheo

The proof of this theorem is a highly combinatorial refinement of previous results of 
Gaussent--Littelmann \cite{GL:LSGalleries}.
We revisit certain proofs in \cite{GL:LSGalleries}, in order to address several shortcomings in loc.~cit.~in view of what we need.
Namely, \cite[Theorem 4]{GL:LSGalleries} gives a disjoint union, rather than a filtrable decomposition, and only works over an algebraically closed field.
Since there are no buildings available over more general base schemes, we need different arguments, and we explain how to reduce to the case of algebraically closed fields, where we can use \cite{GL:LSGalleries}.
As corollaries of the proof, we moreover obtain triviality results of certain torsors (\thref{Torsor is trivial over decomposition} and \thref{Av.Fiber}), which will be used in Section \ref{subsec:Hopf}.

\subsection{Application: a global integral motivic Satake equivalence}

The above geometric insights about $\Gr_{G,I}$ are the basis of the motivic geometric Satake equivalence.
To state it, we build a full abelian subcategory
$$\Sat^{G,I} \subset \MTM(\Gr_{G,I}).$$
Informally, it consists of those motives which admit a (necessarily unique) $L^+_IG$-equivariant-structure, and which have no subquotients supported away from the open locus $X^\circ \subset X^I$ where all the coordinates are distinct. Here we note that essentially by definition, the action of $L^+_IG$ preserves strata, so it is sensible to consider equivariant stratified Tate motives. See \thref{ULAremark} for a comparison with the definition in \cite{FarguesScholze:Geometrization}.

We prove that the global Satake category $\Sat^{G,I}$ enjoys a fusion product, which makes this into a symmetric monoidal category.
The pushforward along the structural map $\Gr_{G,I} \r X^I$ defines a fiber functor for $\Sat^{G,I}$ by taking total cohomology.
To identify $\Sat^{G,I}$ as a category of representations in a natural manner, it is useful to observe that the Langlands dual group $\hat G$, which we take to be defined over $\Z$, comes naturally with a grading, cf.~\refeq{grading on dual group}. 
We can consider representations of $\hat G$ taking values in graded abelian groups and, therefore, also in $\MTM(S)$ or the equivalent category $\MTM(X^I)$, where $X^I$ carries the trivial stratification.

\theo[\S \ref{sect--Tannakian}, \thref{MTM of local Grassmannian}] \thlabel{Main.Thm}
Suppose $S$ satisfies the Beilinson--Soulé vanishing condition, e.g., $S = \Spec \Z$.
After fixing a suitable pinning of $\hat{G}$, there is a canonical equivalence of symmetric monoidal categories 
\[\Sat^{G,I} \cong \Rep_{ \hat{G}^{I}}(\MTM(S)).\]
\xtheo

Such an equivalence involving sheaves or motives on the Beilinson--Drinfeld affine Grassmannian does not appear in \cite{BaumannRiche:Satake} or \cite{RicharzScholbach:Motivic}, but can instead be compared to \cite[§VI]{FarguesScholze:Geometrization}.
For $I = \sgl$, the above condition on subquotients is vacuous, so that $\Sat^{G, \sgl} = \MTM(L^+G \backslash \Gr_G)$.
The resulting symmetric monoidal equivalence
\[\MTM(L^+G\backslash \Gr_G) \cong \Rep_{ \hat{G}}(\MTM(S))\]
refines the result in \cite{RicharzScholbach:Motivic} to integral coefficients.
The above theorem is the first geometric Satake equivalence that is applicable to the base scheme $S = \Spec \Fp$, and with integral (and therefore also $\Fp$-)coefficients.
That feature distinguishes it from Mirkovi\'c--Vilonen's result \cite{MirkovicVilonen:Geometric} which uses analytic sheaves and therefore requires the base to be $\Spec \C$ and offers no  control over Tate twists, i.e., the extra grading mentioned above. Unlike \cite{RicharzScholbach:Motivic}, the present paper is logically independent of, say, \cite{MirkovicVilonen:Geometric}. We refer to \thref{Betti comparison} for a discussion of the compatibility under Betti realization.

The above result can be regarded as a unification of various Satake-type equivalences such as the one for analytic sheaves in 
\cite{Ginzburg:Perverse, MirkovicVilonen:Geometric}, for \(\ell\)-adic sheaves \cite{Richarz:New,Zhu:Ramified}, and algebraic and arithmetic D-modules \cite{BeilinsonDrinfeld:Quantization, XuZhu:Bessel}. 
Namely, simplified versions of our arguments apply to reprove the geometric Satake equivalence for these sheaf theories, and we have explained the compatibility in the case of complex analytic sheaves in \thref{Betti comparison}.
Moreover, let \(\D\) be a sheaf theory defined over a base field \(k\) over \(S\), which satisfies excision, and admits both a six-functor formalism and a perverse t-structure with heart \(\Perv\).
If it also admits a realization functor \(\DM\to \D\) which is compatible with the six operations and right t-exact when restricted to stratified Tate motives, we expect an equivalence
\[\MTM(L^+G\backslash \Gr_G) \otimes_{\MTM(S)} \Perv(\Spec k) \cong \Perv((L^+G\backslash \Gr_G)\times_S \Spec k),\]
where \(\otimes\) denotes Lurie's tensor product of presentable \(\infty\)-categories.
Although we do not pursue this in the present paper, it would allow us to deduce a Satake equivalence for \(\D\) from \thref{Main.Thm}.
That would give another perspective on how the universality of motivic sheaves allows one to deduce statements of interest in geometric representation theory in various sheaf theories.

We conclude the paper by relating motives to the generic spherical Hecke algebra $\calH^\sph_G(\mathbf q)$, which interpolates between the spherical Hecke algebras $\calH^{\sph}_G:=C_\comp(G(\Fq\pot{t})\backslash G(\Fq\rpot{t})/G(\Fq\pot{t}),\Z)$ for all prime powers $q$ via
$\calH^\sph_G(\mathbf q) \t_{\Z[\mathbf q], \mathbf q \mapsto q} \Z = \calH^\sph_G$.
We consider a variant of the Satake category (for $I = \sgl$), denoted $\Satr^{G, \sgl}$, in which the higher motivic cohomology of $S$ has been suppressed (and therefore $\MTM(S)$ gets replaced by the category of graded abelian groups; among other aspects, this enables us to drop the condition that $S$ satisfies Beilinson--Soulé vanishing).
It contains a subcategory of anti-effective motives in which only \emph{negative} Tate twists are allowed (Definitions \ref{anti-effective} and \ref{Whitney--Tate.smooth}).
Let $\hat G_1$ be Deligne's modification of the Langlands dual group \cite{FrenkelGross:Rigid, Deligne:Letter2007}.
It is contained in the Vinberg monoid $V_{\hat G,\rho_\adj}$ of $\hat G$ \cite{Zhu:Integral} (constructed using the same pinning of \(\hat{G}\) as in \thref{Main.Thm}). If $R(-)$ denotes the representation ring of a symmetric monoidal category we have an isomorphism $R(V_{\hat G,\rho_\adj}) \cong \calH^\sph_G(\mathbf q) $. 
The relationship between these objects is summarized below, where the right vertical arrow is closely related to the classical Satake isomorphism, and the composite of the lower horizontal arrows is obtained by taking the character of a representation (see \thref{Generic:Hecke:Iso} for more details).

\theo We have a commutative diagram as follows, relating the ``automorphic side'' on the left to the ``Galois side'', whose decategorification gives the generic Hecke algebra, which in turn specializes to Hecke algebras over $S = \Spec \Fq$:
$$\xymatrix{
\Satr^{G, \sgl, \anti} \ar@{^{(}->}[d]_{\text{adjoin }\Z(1)} \ar[r]^\cong_(.45){\text{\thref{Reps of Vinberg monoid}}} & \Rep_{V_{\hat G,\rho_\adj}}(\Ab) \ar@{^{(}->}[d] \ar@{|->}[r]^-{R(-)} & \calH^{\sph}_G(\mathbf q) \ar@{^{(}->}[d] \ar[r]^{- \t_{\Z[\mathbf q], \mathbf q \mapsto q} \Z} & \calH_G^\sph \ar@{^{(}->}[d] \\
\Satr^{G, \sgl} \ar[r]^-\cong & \Rep_{\hat G_1}(\Ab) = \Rep_{\hat G}(\grAb) \ar@{|->}[r]^-{R(-)} &  \calH^{\sph}_G(\mathbf q^{\pm 1}) \ar[r]^(.45){- \otimes_{\Z[\mathbf q]} \Z[q^{\pm \frac 12}]} & \Z[q^{\pm \frac 12}][X_*(T)]^{W_0}.
}$$
\xtheo

In particular, the motivic formalism, which allows us to consider anti-effective objects, leads naturally to a categorification of Zhu's integral Satake isomorphism \cite{Zhu:Integral}.

\subsection{Future directions}

This paper paves the way for numerous applications of motives in the Langlands program.

In \cite{CassvdHScholbach:Central}, we will build on the Whitney--Tate stratifications of $\Gr_{G,I}$ and construct a motivic refinement of Gaitsgory's central sheaves functor \cite{Gaitsgory:Central}, which is a first step towards enhancing Bezrukavnikov's work \cite{Bezrukavnikov:Two} to a tamely ramified motivic local Langlands correspondence.

For an application in modular representation theory, we address in \thref{H.Independent} a question asked in the recent work \cite[§1.6.1]{EberhardtScholbach:Integral}, which provides a step toward deducing the Finkelberg–Mirkovi\'c conjecture \cite{FinkelbergMirkovic:Semi-infinite} from its graded version \cite{AcharRiche:Reductive}.

Due to their ability of also handling mod-$p$-coefficients, motivic sheaves with integral coefficients such as the ones in this paper hold the promise of being the cohomology theory of choice in the mod-$p$-Langlands program over $S = \Spec \Fp$.
For example, we aim to establish a (non-formal) link between the category $\DTM(\Gr_{G}, \Fp)$ of Tate motives with $\Fp$-coefficients, and the category $\D_\et(\Gr_G, \Fp)$ of étale $\Fp$-sheaves considered in \cite{Cass:Perverse}. 
This would provide a categorification of the specialization map from $\calH^\sph_G(\mathbf q)$ to the mod $p$ Hecke algebra of $G(\Fq(\!(t)\!))$, 
obtained by sending both \(\qq\) and \(p\) to \(0\). Such specialization maps are a key tool in the mod-$p$-Langlands program, cf. \cite{MFV:I}.

The global Satake equivalence is also instrumental in the ongoing program \cite{RicharzScholbach:Intersection,RicharzScholbach:Motivic} aiming for a motivic approach to V.~Lafforgue's work on the global Langlands parametrization over function fields over $\Fp$ \cite{LafforgueV:Chtoucas}.
All these future ideas are driven by the philosophy, uttered by Langlands \cite{Langlands:Automorphic}, of relating motives (as opposed to Galois representations) to automorphic forms.
Such results will settle the independence of the choice of $\ell$ inherent in the usage of $\ell$-adic cohomology.

\subsection{Outline}

In comparison to proofs of other Satake equivalences in the literature, such as \cite{MirkovicVilonen:Geometric, FarguesScholze:Geometrization}, we need to deal with a number of foundational questions in order to prove \thref{Main.Thm}.
Our ambient sheaf formalism is the one of motivic sheaves with integral coefficients. Constructions of stratified Tate motives, including t-structures and equivariant motives, are recalled and developed in \refsect{motivic sheaves}, see especially \thref{equivariant.MTM}.
This theory satisfies Nisnevich, but not étale descent.
This is an important feature rather than a bug because étale descent would prohibit us from treating mod-$p$-sheaves over $S = \Spec \Fp$  \cite[Corollary A.3.3]{CisinskiDeglise:Etale}.
As a result, we need to prove that all torsors we encounter, in particular the Beilinson--Drinfeld Grassmannians and their related loop groups, are Nisnevich-locally trivial (\thref{BD.Zar}).

In Section \ref{sect--aff Grass} we review some geometric properties of affine Grassmannians, which we throughout consider over an arbitrary base scheme. We then prove \thref{intro.thm cellularity}.

Next, we study geometric properties of the Beilinson--Drinfeld affine Grassmannians in Section \ref{sect--BD Grass}, again over an arbitrary base.
We give various interpretations of the convolution product, both local and global, and relate them to each other.
We then prove \thref{BD.WT.intro}. This shows that there are well-behaved categories of stratified Tate motives on Beilinson--Drinfeld Grassmannians, and we conclude this section by showing that the convolution product preserves these motives.

In Section \ref{sect--Satake category}, we introduce and study the global Satake categories \(\Sat^{G,I}\). An instrumental tool for this is the family of constant term functors \(\CT\) (\thref{CT}), which we show preserve the Satake categories, and in particular stratified Tate motives, using \thref{intro.thm cellularity}.
We then construct the fusion product and show it preserves the Satake category (\thref{fusion.Satake}), which gives \(\Sat^{G,I}\) the structure of a symmetric monoidal category.
To prove \thref{fusion.Satake} we relate the convolution product to motivic nearby cycles, cf.~\thref{prop--convolution-CT}.

In Section \ref{sect--Tannakian} we prove the global integral motivic Satake equivalence as stated above (\thref{Main.Thm}).
We show  in \thref{Sat.Anti} that under the equivalence \(\Sat_\red^{G,\sgl}\cong \Rep_{ \hat{G}}(\grAb)\cong \Rep_{\hat G_1}(\Ab)\), anti-effective motives correspond to representations of the Vinberg monoid \(V_{\hat{G},\rho_\adj}\), and we deduce from this a generic Satake isomorphism (\thref{Generic:Hecke:Iso}).

We note that we cannot deduce \thref{fusion.Satake} from the corresponding statement for sheaves in the complex-analytic topology, since the Betti realization functor is not conservative with integral coefficients. 
Later on, we do implicitly use an $\ell$-adic realization functor when appealing to \cite[Corollary 6.4]{RicharzScholbach:Motivic}, which says that the compact objects in the Satake category with $\Q$-coefficients are semisimple when $I=\{*\}$ and $S = \Spec \Fp$. Ultimately this semisimplicity is deduced from the decomposition theorem, and it is only necessary for the identification of the Langlands dual group \(\hat{G}\) with the Tannakian group of the Satake category. Otherwise, we avoid using realization functors, in order to give a motivic proof of the Satake equivalence.

\subsection*{Acknowledgements}
We thank Esmail Arasteh Rad, Patrick Bieker, St\'ephane Gaussent, Somayeh Habibi, Tom Haines, Shane Kelly, C\'edric P\'epin, Timo Richarz, Markus Spitzweck, Can Yaylali, and Xinwen Zhu for helpful discussions.
We also thank Timo Richarz and the anonymous referees for very useful comments on earlier versions of this paper.
T.v.d.H. thanks the participants of a workshop on geometric Satake organized in Clermont-Ferrand in January 2022, as well as Tom Haines for writing \cite{Haines:Pavings} in response to a question he had.
R.C. and J.S. also thank the Excellence Cluster Mathematics Münster for logistical support. R.C. was supported by the National Science Foundation under Award No. 2103200. T.v.d.H. acknowledges support by the European Research Council (ERC) under the European Union’s Horizon 2020 research and innovation programme (grant agreement 101002592), and by the Deutsche Forschungsgemeinschaft (DFG), through the TRR 326 GAUS (project number 444845124).

\section{Motivic sheaves}
\label{sect--motivic sheaves}

\nota \thlabel{nota-basescheme}
Throughout this paper, we fix a connected base scheme $S$ that is smooth of finite type over a Dedekind ring or a field.
By convention, a scheme is always supposed to be of finite type over $S$.
\xnota

\subsection{Recollections}
\label{sect--recollections}

\subsubsection{Motives}
\label{sect--motives}
The sheaf formalism used in this paper is the category of motivic sheaves defined using the Nisnevich topology, with integral coefficients. 
The point in considering this sheaf theory is that it is universal, i.e., exists over any base scheme, and is independent of the choice of an auxiliary cohomology theory such as $\ell$-adic cohomology.

For a scheme $X$ over a base scheme $S$ as above, the category $\DM(X)$ of \emph{motives over $X$} (or motivic sheaves on $X$) was constructed by Spitzweck \cite{Spitzweck:Commutative} (building on the works of Ayoub, Bloch, Cisinski, Déglise, Geisser, Levine and Morel--Voevodsky; we refer to op.~cit.~for further references).
Very briefly, the stable $\A^1$-homotopy category $\SH(X)$ is the algebro-geometric incarnation of the category of spectra in classical homotopy theory.
Broadly speaking, $\DM(X)$ compares to $\SH(X)$ as the derived category of abelian groups $\D(\Ab)$ compares to spectra.
More precisely, Spitzweck has constructed a motivic ring spectrum $\M\Z \in \SH(\Spec \Z)$ representing motivic cohomology.
Defining $\M\Z_X := \pi^* \M\Z$ for $\pi : X \r \Spec \Z$, the category of motives on $X$ is defined as $\DM(X) := \Mod_{\M\Z_X}(\SH(X))$.

The category $\DM(X)$ shares many properties with the analogously defined category of motives with rational coefficients, as reviewed in, say \cite[Synopsis~2.1.1]{RicharzScholbach:Intersection}.
The abstract categorical properties of $\DM(X)$ as well as the existence of various pullback and pushforward functors (see points (i)--(v), (vii)) hold unchanged. 
For an $S$-scheme $X$ as above, $\SH(X)$ and thus $\DM(X)$ is compactly generated \cite[Théorème~4.5.67]{Ayoub:Six2} by $f_\sharp \Z(k)$, for $f : Y \r X$ smooth, $k \in \Z$. For finite type maps $f$, $f^*$ and $f_!$ preserve compact objects \cite[Proposition~4.2.4, Corollary~4.2.12]{CisinskiDeglise:Triangulated}.
(The preservation of compact objects under $f_*$ and $f^!$, which is proved for arbitrary maps in \cite[Theorem~4.2.48, Corollary~4.2.28]{CisinskiDeglise:Triangulated} under additional hypotheses including rational coefficients is not used in this paper.)
In particular, their right adjoint $f_*$, $f^!$ as well as (trivially) $f^*$, $f_!$ preserve arbitrary colimits.
For two maps $f_1, f_2$, the projection formula (point vii) there) implies an isomorphism \cite[Lemma~2.2.3]{JinYang:Kuenneth}
$$(f_1 \x f_2)_! (- \boxtimes - ) \stackrel \cong \r (f_{1!} - \boxtimes f_{2!} -).$$ 
Trivially, a similar formula holds for *-pullback functors.
Localization triangles, base change, homotopy invariance and relative purity (see points (ix)--(xii) there) hold unchanged.
If $X$ is smooth over a field, then $\Hom_{\DM(X)}(\Z, \Z(n)[m]) = \CH^n(X, 2n-m)$ (a higher Chow group), similar to point (xiii) in \cite[Synopsis 2.1.1]{RicharzScholbach:Intersection}.
Since, over $S = \Spec \Z$ or $\Spec \Fp$, say, resolution of singularities currently requires alterations, as opposed to just blow-ups, Verdier duality is not known to be an involution on $\DM(X)$ (with integral coefficients); for the same reason we do not claim the existence of a weight structure on the categories $\DM(X)$ (cf.~points (viii) and (xv) there).
In contrast to the case of motives with rational coefficients, which form an h-sheaf (point (xiv), \cite[Theorem~2.1.13]{RicharzScholbach:Intersection}), the functor $X \mapsto \DM(X)$ is only a Nisnevich sheaf. This follows from the corresponding sheaf property of the stable homotopy category $\SH$ \cite[Proposition~6.24]{Hoyois:Six}.

\defi
\thlabel{anti-effective}
For a scheme \(X\), the subcategories $$\DM(X)\supset \DTM(X)\supset\DTM(X)^\anti$$ of \emph{Tate motives} (resp.~\emph{anti-effective Tate motives}) are defined to be the presentable stable subcategories generated by $\Z(n)$ with $n \in \Z$ (resp.~$n \le 0$).
\xdefi

The terminology \emph{anti-effective} reflects the fact that we consider the opposite of the usual notion of effective motives in the literature, e.g.,~\cite[Definition~11.1.10]{CisinskiDeglise:Triangulated}.
While anti-effective motives are less immediate to define than effective motives, they necessarily arise in any situation where one glues Tate motives, since $i^! \Z = \Z(-c)[-2c]$ for a closed immersion $i$ of codimension $c$ between smooth $S$-schemes.
On the same note, as was pointed out by T.~Richarz, the homological functor $p_\sharp$ for, say, the projection $p \colon \P_S \r S$ does produce effective motives ($p_\sharp \Z = \Z \oplus \Z(1)[2]$), but we do not use this functor in this paper.
We will eventually relate representations of the Vinberg monoid with a certain category consisting of anti-effective (stratified) Tate motives, cf.~\thref{Reps of Vinberg monoid}.

\subsubsection{Betti realization}
\label{sect--realization}

In order to relate our results with the geometric Satake equivalence in \cite{MirkovicVilonen:Geometric}, we will use the Betti realization functor
$$\rho_\Betti : \DM(X) \r \D(X^\an)$$
taking values in the derived category of sheaves on the analytic space associated to any scheme $X$ (by convention always of finite type) over $S = \Spec \C$.
This functor can be constructed by using Ayoub's Betti realization functor $\SH(X) \r \D(X^\an)$ \cite[Définition~2.1]{Ayoub:Note}, and using that for $S = \Spec \C$ the spectrum $\M\Z$ constructed by Spitzweck (cf.~\refsect{motives}) is isomorphic to the classical Eilenberg--MacLane spectrum, which is mapped to $\Z$ under the above functor \cite[Theorem~5.5]{Levine:Comparison}.
The restriction of $\rho_\Betti$ to the subcategory of constructible motives is compatible with the six functors \cite{Ayoub:Note}.

\rema \thlabel{Betti not conservative}
Betti realization is known to be conservative on compact stratified Tate motives with \emph{rational} coefficients (e.g., \cite[Lemma~3.2.8]{RicharzScholbach:Intersection}), but fails to be conservative for integral coefficients.
This requires us to use methods that are logically entirely independent of, say, \cite{MirkovicVilonen:Geometric}; cf.~ also \thref{Betti comparison}.
\xrema

\subsubsection{Reduced motives}
\label{sect--reduced.motives}
The category of \emph{reduced motives},  introduced in \cite{EberhardtScholbach:Integral}, and its full subcategory of \emph{reduced Tate motives} on a scheme $X / S$ are defined as
$$\eqalign{
\DMr(X) & := \DM(X) \t_{\DTM(S)} \D(\grModZ), \cr
\DTMr(X) & := \DTM(X) \t_{\DTM(S)} \D(\grModZ).}\eqlabel{DTMr}$$ 
Here $\D(\grModZ)$ denotes the derived \ii-category of $\Z$-graded abelian groups and $\t$ denotes Lurie's tensor product (of presentable stable \ii-categories). Referring to op.~cit. for further discussion, we only note that reduced motives therefore behave like motives; i.e., the properties in points (i)--(v), preservation of compact objects under $f^*$ and $f_!$ (as in point (vi there), points (vii), (ix)--(xii) and (xiv) (for the Nisnevich topology) in \cite[Synopsis~2.1.1]{RicharzScholbach:Intersection} hold for $\DMr$. The difference between $\DM$ and $\DMr$ is that the higher motivic cohomology of the base scheme $S$ has been removed.
Indeed, by definition, $\DTMr(S) = \D(\grModZ)$, independently of the choice of the base scheme $S$.
We also refer to \thref{DTMr.independence} for a more general result asserting the independence of the choice of the base scheme $S$ for reduced motives.
Reduced motives will allow us to exhibit the Tannaka dual of the Satake category as a group associated to a $\Z$-graded Hopf algebra, cf.~\thref{Satake.Hopf}.

We will write $\DMrx(X)$ to denote either $\DM(X)$ or $\DMr(X)$, and similarly with $\DTMrx(X)$. We may sometimes omit the subscript ${}_{\redx}$ to ease the notation, but unless specifically noted our results hold for both regular and reduced motives. 

\subsubsection{Functoriality}
\label{sect--functoriality}
The assignment $X \mapsto \DMrx(X)$ can be organized into a lax symmetric monoidal functor
$$\DMrx : \Corr (\Sch_S^{\ft}) \r \PrStZ.$$
The target is the \ii-category of presentable stable $\Z$-linear categories with colimit-preserving functors.
The source denotes the \ii-category of correspondences: its objects are finite type $S$-schemes; morphisms are zig-zags $X \stackrel f \gets Y \stackrel g \r Z$. The functor $\DMrx$ maps this to $g_! f^* : \DMrx(Z) \r \DMrx(X)$.
Thus, this functor encodes the existence of *-pullbacks, !-pushforwards (along maps of finite type $S$-schemes), and the existence of their right adjoints.
The lax monoidality encodes the existence of the exterior product $\boxtimes : \DMrx(X) \t_{\D(\Mod_\Z)} \DMrx(Y) \r \DMrx(X \x Y)$ functors, as well as various projection formulas.
We refer to \cite[§2]{EberhardtScholbach:Integral} for a slightly more detailed survey, including references to the original works where the functor has been constructed.

The functor $\DTM(S) \r \D(\grModZ)$ used in the definition of $\DMr(X)$ gives rise to a natural transformation, called the \emph{reduction functor} 
$$\rho_\red : \DM \r \DMr.$$
It is compatible with the !- and *-pullback and pushforward functors, $\otimes$ and $\boxtimes$ \cite{EberhardtScholbach:Integral}. At least for stratified Tate motives, it is compatible with Verdier duality (\thref{reduction.Hom}). 

\subsubsection{Motives on prestacks}
The above formalism of (reduced) motives on finite type $S$-schemes extends formally (by means of appropriate Kan extensions) to a functor
$$\DMrx^! : (\PreStk_S)^\opp \r \PrSt_\Z.\eqlabel{DM.prestack}$$
The source category is the \ii-category of prestacks, i.e., presheaves of anima on the category of affine, but not necessarily finite type $S$-schemes. 
This construction is parallel to the one for rational coefficients in \cite[\S2]{RicharzScholbach:Intersection};
\refsect{equivariant.motives} unravels this definition in the case of equivariant motives.

\exam
\thlabel{DM.ind-scheme}
If $X = \colim_{k \in K} X_k$ is an ind-scheme, there is an equivalence 
$$\DM(X) = \colim_{\text{!-pushforwards}} \DM(X_k) = \lim_{\text{!-pullbacks}} \DM(X_k).\eqlabel{DM.colim.lim}$$
Here, the (co)limit is taken in $\PrSt_\Z$.
This was shown in \cite[Corollary~2.3.4]{RicharzScholbach:Intersection} for motives satisfying étale descent; since the above claim only concerns (ind-)schemes, it also holds for motives in the Nisnevich topology.
This implies that every object $M \in \DM(X)$ is of the form 
$$M = \colim i_{k!} i_k^! M,$$ 
where $i_k\colon X_k \r X$.
We say that $M$ is \emph{bounded} if it is in the image of the canonical insertion functor $i_{k!}\colon \DM(X_k) \r \DM(X)$ for some $k$.
\xexam

\subsubsection{Hyperbolic localization}
The following statement is needed below in order to decompose the fiber functor into weight functors.

\prop
\thlabel{hyperbolic.localization}
Let $X$ be an (ind-)scheme with an action of $\Gm$ that is respecting some ind-presentation $X = \colim X_i$.
We also assume the $\Gm$-action is Nisnevich-locally linearizable (i.e., 
for each $i$, $X_i$ admits Nisnevich-locally a cover by $\Gm$-stable affine open subschemes).
Consider the fixed points $X^0$, the attractors $X^+$ and repellers $X^-$ of this action: 
$$\xymatrix{
& X^\pm = \IHom_{\Gm}(\A^1_\pm, X) \ar[dr]^{p^\pm} \ar[dl]_{q^\pm} \\
X^0 = \IHom_{\Gm}(S, X) & & X.}$$
Here $\A^1_\pm$ is $\A^1$ with the $\Gm$-action given by $\lambda \cdot t := \lambda^{\pm 1} t$.
These functors are representable by (ind-)schemes.

There is a natural transformation
$$h: q^-_* p^{-!} \r q^+_! p^{+*}.$$
The restriction of this natural transformation to the full subcategory $\DM(X)^{\Gmmono} \subset \DM(X)$ of $\Gm$-monodromic motives (i.e., the subcategory generated under arbitrary colimits by the image of the forgetful functor $\DM(\coeq (\Gm \x X \stackrel[a]{p} \rightrightarrows X)) \r \DM(X)$) is an equivalence.
\xprop

\pf
In the context of étale torsion sheaves on algebraic spaces the statement above is Richarz' version of hyperbolic localization \cite{Richarz:Spaces}. 
The proof in loc.~cit. only uses the *- and !-functors, localization, and homotopy invariance
for étale torsion sheaves, and can be repeated verbatim for motives on schemes.
In contrast to the situation for étale sheaves, both $q^-_* p^{-!}$ and $q^+_! p^{+*}$ preserve arbitrary colimits (cf.~\refsect{motives}), so that the statement above holds for arbitrary (as opposed to finite) colimits of weakly $\Gm$-equivariant sheaves, as stated.
(Another proof in the context of D-modules (over schemes in characteristic 0) due to Drinfeld--Gaitsgory \cite{DrinfeldGaitsgory:Theorem} that again only uses these formal properties of a sheaf context can also be adapted verbatim to motives.)

The statement for $\Gm$-monodromic motives on an ind-scheme $X$ as above follows since again all four functors appearing in the natural transformation $h$ are colimit-preserving and $\DM(\coeq (\Gm \x X \rightrightarrows X)) = \colim_i \DM(\coeq(\Gm \x X_i \rightrightarrows X_i))$.
\xpf

\subsection{Stratified Tate motives}

In the sequel, we will be using standard terminology about stratified (ind-)schemes, as in \cite[\S3]{RicharzScholbach:Intersection}.

\defi
\thlabel{Whitney--Tate}
(\cite[§4]{SoergelWendt:Perverse}, \cite[3.1.11]{RicharzScholbach:Intersection})
Let $\iota : X^\dagger =\bigsqcup_{w \in W} X^w \r X$ be a stratified ind-scheme. 
We say $M$ is a \emph{stratified Tate motive} if $\iota^* M \in \DTM(X^\dagger)$.
The stratification $\iota$ is called \emph{Whitney--Tate} if $\iota^* \iota_* \Z \in \DTM(X^\dagger)$. (Equivalently, $\iota^{v*} \iota^w_* \Z \in \DTM(X^v)$ for all the strata $X^w \stackrel{\iota^w} \r X \stackrel{\iota^v} \gets X^v$.) We similarly define an \emph{anti-effective stratified Tate motive} and an \emph{anti-effective Whitney--Tate stratification} by replacing $\DTM(X^\dagger)$ with $\DTM(X^\dagger)^{\anti}$ everywhere.

If we have an (anti-effective) Whitney--Tate stratification, we denote by $\DTM(X, X^\dagger)^{(\anti)} \subset \DM(X)$ the full subcategory of (anti-effective) stratified Tate motives. 
This category is called the category of \emph{(anti-effective) stratified Tate motives}, and is also denoted by $\DTM(X)^{(\anti)}$ if the choice of $X^\dagger$ is clear from the context.

A Whitney--Tate stratification is called \emph{universally Whitney--Tate} if for any scheme $Y \r S$, the natural map
$$p^* \iota_* \Z \r (\id_Y \x \iota)_* p^{\dagger*} \Z\eqlabel{universally.WT}$$
resulting from the following cartesian diagram is an isomorphism:
$$\xymatrix{
Y \x_S X^\dagger \ar[d]^{p^\dagger} \ar[r]^{\id \x \iota} & Y \x_S X \ar[d]^p \\
X^\dagger \ar[r]^\iota & X.
}$$
\xdefi

\rema
\thlabel{Whitney--Tate.smooth}
Having an (anti-effective) Whitney--Tate stratification ensures that $\iota^*$, $\iota^!$, $\iota_!$ and $\iota_*$ preserve the categories of (anti-effective) stratified Tate motives, cf.~\cite[§4]{SoergelWendt:Perverse}.

Any stratification such that the closures $\ol {X^w}$ are smooth over $S$ is  
anti-effective Whitney--Tate. This follows from excision and relative purity \cite[Remark~4.7]{Wildeshaus:Intermediate}, i.e., the isomorphism $i^! \Z = \Z(-c)[-2c]$ for a codimension $c$ closed immersion of smooth $S$-schemes. 
This also explains why the dual notion of ``effective Whitney--Tate'' stratifications is not a sensible definition. 

If $X$ is universally (anti-effective) Whitney--Tate, then the product stratification $Y \x_S X^\dagger \r Y \x_S X$ is (anti-effective) Whitney--Tate for any scheme $Y / S$. 
\xrema

\lemm
\thlabel{reduction.Hom}
For a Whitney--Tate stratified ind-scheme $X$, the reduction functor $\rho_\red : \DTM(X, X^\dagger) \r \DTMr(X, X^\dagger)$ is compatible with the internal $\Hom$-functor.
In particular, if Verdier duality preserves Tate motives (e.g., $X^\dagger$ is smooth), $\rho_\red$ is compatible with Verdier duality.
\xlemm

\pf
The functor $\iota^!$ is conservative and satisfies $\iota^! \IHom(M, N) = \IHom(\iota^* M, \iota^!N)$.
We can therefore assume $X$ consists of a single stratum. In this case we conclude using $\IHom(\Z(k), N) = N(-k)$.
\xpf

\defi
\thlabel{admissible.stratification}
A map $f : X \r Y$ is called \emph{admissible} if $f$ is smooth and if $f_! f^!$ preserves the subcategory $\DTM(Y)^{\le 0}$, the smallest subcategory of $\DTM(Y)$ containing $\Z(n)$, $n \in \Z$, stable under extensions and colimits (and therefore also shifts $[k]$ for $k > 0$ 
and direct summands).

A stratified map $(X, X^\dagger) \stackrel {(\pi, \pi^\dagger)} \r (Y, Y^\dagger)$ \cite[Definition~3.1.1(ii)]{RicharzScholbach:Intersection} is called admissible if $\pi^\dagger$ is admissible. If $Y=S$, we also say that $X$ is \emph{admissibly stratified}.
\xdefi

\rema
\thlabel{cellular.stratification}
The admissibility of a stratification enforces in particular that the motives of the individual strata $X^w$ are Tate motives over each $Y^w$.
The condition that we only allow positive shifts will be used in order to construct the motivic t-structure.

A Whitney--Tate stratification for which 
the strata $X^w$ are isomorphic to $\Gm^{n_w} \x_S \A^{m_w}$ is admissible. Such stratifications are called \emph{strongly cellular}.

More generally, if each $X^w$ is smooth and admits (in its own right) a strongly cellular Whitney--Tate stratification, then the  stratification of $X$ by the $X^w$ is admissible. We call these stratifications \emph{cellular}, as in \cite[Def 3.1.5]{RicharzScholbach:Intersection}.

If $(X, X^\dagger) \r (Y, Y^\dagger)$ is admissible (e.g., cellular) and $Y$ is admissibly stratified, then $X$ is also admissibly stratified.
\xrema

We introduce admissible (as opposed to cellular) schemes in order to have a t-structure on the open subscheme $U \subset \A^n$ consisting of points in $\A^n$ whose coordinates are pairwise distinct.
The scheme $U$ is not cellular (e.g. for $n=3$, $U$ has no $\mathbf F_2$-points), but admissible by the following lemma.
We will later use this in the context of the Beilinson--Drinfeld Grassmannian (\thref{BD.WT}).

\lemm
\thlabel{admissible.SNC.divisor}
Let $$D := \bigcup_{i \in I} D_i \stackrel i \r X \gets U := X \setminus D$$
be the inclusion of a strict normal crossings divisor (i.e., the $D_J := \bigcap_{j \in J} D_j$ for all $J \subset I$ are smooth over $S$, including $D_\emptyset = X$) and its complement.
Suppose for each $J$, the structural map $\pi_J : D_J \r S$ has the property $\pi_{J!} \pi_J^! \Z \in \DTM(S)^{\le 0}$.
Then $f_! f^! \Z \in \DTM(S)^{\le 0}$, for $f : U := X \setminus D \r S$.
\xlemm

\pf
Let $D^{(n)} = \bigsqcup_{J \subset I, |J| = n} D_J \stackrel {\pi^{(n)}} \r S$ be the disjoint union of the $n$-fold intersections of the individual divisors, so that $D^{(0)} = X$. 
By relative purity (i.e., the isomorphism $g^! = g^* (d)[2d]$ for any smooth map $g$ of relative dimension $d$), it suffices to see that $f_! f^* \Z \in \DTM(S)^{\le 2 d}$, with $d = \dim X$.
By localization (cf.~\cite[Lemme~2.2.31]{Ayoub:Six1}), this object is the homotopy limit of a diagram of the form (with transition maps induced by unit maps, using that $D^{(n)} \r D^{(n-1)}$ is proper)
$$\pi^{(0)}_! \pi^{(0)*} \Z \r \pi^{(1)}_! \pi^{(1)*} \Z \r \pi^{(2)}_! \pi^{(2)*} \Z \r \dots .\eqlabel{localization.SNC divisor}$$
By assumption and relative purity for the smooth maps $\pi^{(n)}$, we have $\pi^{(n)}_! \pi^{(n)^*} \Z \in \DTM(S)^{\le 2(d-n)}$, which implies our claim.
\xpf

The following lemma allows to zig-zag between reduced motives on the Hecke prestack over $\Spec \Q$ and over $\Spec \Fp$ (cf.~\thref{independence.star}).

\lemm 
\thlabel{DTMr.independence}
Consider a cartesian diagram
$$\xymatrix{
X'^\dagger \ar[r]^{\iota'} \ar[d]^{f^\dagger} &
X' \ar[d] \ar[r]^{\pi'} & S' \ar[d]^f \\
X^\dagger \ar[r]^{\iota} & X \ar[r]^\pi & S,
}$$
in which $\iota$ determines a universally admissibly Whitney--Tate stratified (ind-)scheme and $S'$ is an $S$-scheme such that $f^* \pi_* \iota_* \Z_{X^\dagger} \stackrel \cong \r \pi'_* \iota'_* f^{\dagger*} \Z_{X^\dagger}$.
Then $f^* : \DTMr(X) \r \DTMr(X')$ is an equivalence. Here $X' := X \x_S S'$ and reduced motives are taken with respect to the respective base schemes, i.e., $S$ for $X$ and $S'$ for $X'$, so that $\DTMr(X') := \DTM(X') \t_{\DTM(S')} \D(\grAb)$. 
\xlemm

\pf
This is the content of \cite[Proposition~4.25]{EberhardtScholbach:Integral} if the stratification is cellular (as opposed to just admissible).
As in loc.~cit., using the universality, one reduces to the case where $X$ is a single stratum.
We consider the monad $T = \pi_* \pi^*$ associated to the adjunction $\pi^* : \DTMr(S) \rightleftarrows \DTMr(X) : \pi_*$.
By definition of $\DTMr$, 
the image of $\pi^*$ generates $\DTMr(X)$, so that $\pi_*$ is conservative.
It is also colimit-preserving, so that the Barr--Beck--Lurie theorem implies that $\DTMr(X) = \Alg_T(\DTMr(S))$.
In order to establish the equivalence, we first observe that the claim holds true for $S$ in place of $X$ by the definition in \refeq{DTMr}.
By our assumption, $f^*$ commutes with $\pi_*$, so that $f^*$ maps the monad $\pi_* \pi^*$ to the monad $\pi'_* \pi'^*$.
\xpf

\subsection{t-structures}
\label{sect--t-structures}

In this subsection, we summarize some basic properties related to motivic t-structures.
Throughout, we use cohomological conventions concerning t-structures, as in \cite[Définition~1.3.1]{BeilinsonBernsteinDeligne:Faisceaux}, but opposite to \cite[Definition~1.2.1.1]{Lurie:HA}.
The construction works in parallel for reduced and regular motives, except that in the latter case we always (have to) assume (in addition to our running assumption in \thref{nota-basescheme}) that $S$ satisfies the \emph{Beilinson--Soulé vanishing} condition
$$\H^n(S, \Q(k)) := \Hom_{\DM(S)}(\Z, \Q(k)[n]) = 0 \ \ \text{for } n < 0, k \in \Z.\eqlabel{BS.vanishing}$$
This is satisfied, for example, if $S = \Spec \Q$, $\Spec \Z$ or $\Spec \Fp$ by work of Borel and Quillen, cf.~\cite[Lemma~41]{Kahn:AlgebraicKTheory}.
Recall that $\DTMrx$ denotes either $\DTM$ or $\DTMr$.

\lemm
\thlabel{t-structure.heart}
The category $\DTMrx(S)$ carries a right-complete t-structure such that
$$\DTMrx(S)^{\le 0} = \langle \Z(k) \rangle \eqlabel{normalization.t-structure}$$
(i.e., the closure under colimits and extensions of these objects).
The objects $\Z(k)$ are a set of compact generators of the heart of this t-structure.
Both the $\le 0$- and the $\ge 0$-aisle of the t-structure are closed under filtered colimits.
Therefore the truncation functors $\tau^{\ge n}, \tau^{\le n}$ and $\pH^n := \tau^{\ge n} \tau^{\le n}$ preserve filtered colimits.
\xlemm

\pf
The existence of the t-structure and stability under filtered colimits is a generality about t-structures generated by compact objects in cocomplete \ii- (or triangulated) categories. See, e.g.,~\cite[Theorem A.1]{AlonsoEtc}.
The t-structure is right complete since the objects $\Z(k) \in \DTMrx(S)^{\le 0}$ compactly generate $\DTMrx(S)$ under colimits and shifts:
indeed, applying the dual of \cite[1.2.1.19]{Lurie:HA} to $\DTMrx(S)$, 
we need to show $\bigcap \DTMrx(S)^{\ge n} = 0$ (recall that we use cohomological notation). 
For $X \in \bigcap \DTMrx(S)^{\ge n}$ and any $n \in \Z$, we have $\Hom_{\Ho(\DTMrx(S))}(\Z(k)[n], X)=0$ since $\Z(k)[n] \in \DTMrx(S)^{\le -n}$. Thus $X = 0$.
The generators $\Z(k)$ are in the heart since the vanishing in \refeq{BS.vanishing} implies a similar one for $\Z$-coefficients \cite[Lemma~3.4]{Spitzweck:Mixed}. 
The fact that the $\Z(k)$ compactly generate the heart is, e.g., \cite[Ch.~III, Lemma~3.1.(iii)]{BeligiannisReiten:Homological} together with \cite[Theorem~1.11]{AdamekRosicky:Locally}.
\xpf

\rema
\thlabel{compact t-structure}
It is a delicate question whether the t-structure restricts to one on the subcategory of $\DTM(S)$ spanned by compact objects. 
This is known for $\Q$-coefficients by \cite{Levine:Tate} and also for $\Z$-coefficients for $S = \Spec \Z$ (unpublished work of Markus Spitzweck).
It also holds for $\DTMr(S)$.
We will not use such a property in this paper. 
\xrema

\lemm
\thlabel{t-structure.stratified}
\thlabel{IC.motives}
Suppose $(X, X^\dagger = \coprod_w X^w) \r (Y, Y^\dagger = \coprod_u Y^u)$ is an admissibly stratified map between two Whitney--Tate stratified schemes.
Suppose that $\DTMrx(Y^\dagger)$ carries a t-structure whose $\le 0$-aisle is generated (under colimits and extensions) by $\Z_{Y^u}(k)[\dim_S Y^u]$, which we assume lie in the heart.
\begin{enumerate}
  \item 
Then $\DTMrx(X^\dagger)$ also carries a right complete t-structure whose $\le 0$-aisle is generated by $\Z_{X^w}(k)[\dim_S X^w]$.
Again, these objects compactly generate the heart of the t-structure.
\item
\label{item--MTM glued}
The category $\DTMrx(X)$ carries a right complete t-structure glued from the t-structures on the strata, i.e., on $\DTMrx(X^\dagger)$.
Its heart is generated under extensions by the intermediate extensions along the maps $\iota^w : X^w \r X$ \cite[Définition~1.4.22]{BeilinsonBernsteinDeligne:Faisceaux}, i.e., by the objects
$$\IC_{w, L} := \iota^w_{!*} L := \im (\pH^0 \iota^w_! L \r \pH^0 \iota^w_* L) \in \MTMrx(X, X^\dagger)\eqlabel{IC.w.L}$$
for $L \in \MTMrx(X^w)$.
We call $\IC_{w, L}$ the \emph{(reduced) intersection motive}.
\item
\label{item--exactness functors}
The aisle $\DTMrx(X)^{\le 0}$ (resp.~$\DTMrx(X)^{\ge 0}$) is generated by $\iota_! \DTMrx(X^\dagger)^{\le 0}$ (resp.~$\iota_* \DTMrx(X^\dagger)^{\ge 0}$).
Thus, $\iota^*$ is right t-exact and $\iota^!$ is left t-exact.
If $i : \ol {X_w} \r X$ is the inclusion of the closure of a stratum, $i_! = i_*$ is t-exact.
\end{enumerate}
\xlemm

\defi \thlabel{defi--MTM(anti)}
The heart of these t-structures will be denoted by $\MTMrx(X, X^\dagger)$ or just $\MTMrx(X)$ if $X^\dagger$ is clear from the context. Its objects are called \emph{mixed Tate motives}.
Again, there is an obvious variant for anti-effective Whitney--Tate stratifications. The heart of the t-structure on $\DTMrx(X)^\anti$ is denoted $\MTMrx(X)^\anti$.  
\xdefi

\pf
If $\pi^w : X^w \r Y^u$ is the map of a stratum in $X$ to the stratum in $Y$, say of dimension $d$, then $\Z_{X^w} = \pi^{w!} \Z_{Y^u} (-d)[-2d]$. Thus, the orthogonality condition $\Hom(\Z_{X^w}, \Z_{X^w}(k)[n]) = 0$ for $n < 0$ holds by the admissibility condition for $\pi^\dagger$.

\refit{MTM glued} is a generality about glued t-structures, cf.~\cite[Theorem~3.4.2]{Achar:Perverse}.
(In a certain situation where some group acts transitively on $X^w$, we will describe a more narrow set of generators in \thref{MTM.equivariant.ind-scheme}.)
\refit{exactness functors} is also standard, cf.~\cite[§1.4]{BeilinsonBernsteinDeligne:Faisceaux}.
\xpf

From here on we may sometimes write $\dim X$ instead of $\dim_S X$, but we will always mean the dimension relative to $S$ unless otherwise stated.

\lemm \thlabel{Ext.To.Ext2}
Let $X$ be an ind-scheme with an admissible Whitney--Tate stratification. Let $j \colon U \r X$ be the inclusion of an open union of strata and let $i \colon Z \r X$ be the complement. Let $A, B \in \MTMrx(X)$ be such that $A$ has no quotients supported on $Z$ and $B$ has no subobjects supported on $Z$.
Then there is a natural isomorphism
$$\Hom(j_{!*}j^*A, B) \cong \Hom(A, B).$$
\xlemm

\pf
By the assumption on $A$, there is an exact sequence
$$0 \r i_*\pH^{-1}(i^*A)   \r  \pe j_!j^* A \r  A \r 0.$$ We have $\Hom(i_*\pH^{-1}(i^*A) ,B) = 0$ by the assumption on $B$, so $\Hom(A, B) \cong \Hom(\pe j_! j^* A, B)$ 
where $\pe j_! = \pH^0 j_!$. Since $j_{!*}j^*A$ also has no quotients supported on $Z$, replacing $A$ with $j_{!*}j^*A$ in the same argument gives $\Hom(j_{!*}j^* A, B) \cong \Hom(\pe j_! j^* A, B)$.
\xpf

\lemm \thlabel{Anti.Orthogonal}
Let $X$ be a smooth admissible $S$-scheme. Then, for unstratified Tate motives,
$$\eqalign{
\DTMrx(X)^\anti & = \{M \in \DTMrx(X) \: : \: \Map_{\DTMrx(X)}(\Z(p), M) = 0 \text{ for all } p \geq 1\}, \cr
\MTMrx(X)^\anti & = \{M \in \MTMrx(X) \: : \: \Hom_{\MTMrx(X)}(\Z(p)[\dim X], M) = 0 \text{ for all } p \geq 1\}.}$$ 
(Here $\Map$ denotes the mapping complex between the indicated objects; the $n$-th cohomology of this complex is the Hom-group $\Hom(\Z(p), M[n])$ in the triangulated category underlying $\DTMrx(X)$.)
Moreover, $\MTMrx(X)^\anti \subset \MTMrx(X)$ is stable under subquotients.
\xlemm

\pf
We first prove the description of $\DTMrx(X)^\anti$.
To show ``$\subset$'', it suffices to see $\Map_{\DTMrx(X)}(\Z(p), M) = 0$ or equivalently $\Hom_{\DTMrx(X)}(\Z(p), M[s]) = 0$ for all $s \in \Z$. It suffices to consider the case  $M = \Z(n)$ for $n \leq 0$, in which case this group is given by $\H^s(X, \Z(n-p))$, which vanishes if $p \geq 1$ and $n \leq 0$. Indeed, by assumption (and \thref{nota-basescheme}), $X$ is smooth over $\Spec B$, where $B$ is a Dedekind ring or a field. In the latter case, the above Hom-group is a higher Chow group of codimension $n-p$ cycles, which vanishes for $n-p < 0$. We reduce the vanishing in the former case to the latter case by using the distinguished triangle $\bigoplus_{\mathfrak p \subset B} i_{\mathfrak p*} \Z(-1)[-2] \r \Z \r \eta_* \Z$, where $i_{\mathfrak p}$ is the immersion of a closed point and $\eta$ the generic point of $\Spec B$ \cite[§7]{Spitzweck:Commutative}.
For reduced motives, this vanishing still holds by \cite[(3.4)]{EberhardtScholbach:Integral}. 

To show ``$\supset$'', let us write $\mathcal C$ for the right hand category.
The inclusion $i \colon \DTMrx(X)^\anti \r \mathcal{C}$ admits a right adjoint $R$ 
by the adjoint functor theorem.
In order to show $iR = \id$, it suffices to show $\Map(\Z(e), iRM) \r \Map(\Z(e), M)$ is an equivalence for all $M \in \mathcal{C}$ and $e \in \Z$. For $e \geq 1$ this is immediate because the cofiber of $iRM \r M$ is in $\mathcal{C}$. For $e \leq 0$, by adjunction and $Ri = \id$, we have equivalences $\Map(\Z(e), iRM) = \Map(\Z(e), RM) = \Map(\Z(e), M)$. 

The proof for $\MTMrx(X)^\anti$   is analogous. 
This also implies that $\MTMrx(X)^\anti$ is stable under subquotients. 
Indeed, if $N \r M$ is an injection in $\MTMrx(X)$, so is $\Hom(\Z(p)[\dim X], N) \r \Hom(\Z(p)[\dim X], M)$.
Stability under quotients then follows as well.
\xpf

\subsection{A motivic computation}
On several occasions, including our computation of constant term functors (cf.~\thref{first stratification of intersection}), we will encounter the following geometric situation.

\defi
\thlabel{Defi filtrable} \cite[Definition 2]{BB:Properties}
A decomposition of a scheme \(X\) into locally closed subschemes \((X_\alpha)_{\alpha\in A}\) is called \emph{filtrable}, if there exists a finite decreasing sequence \(X=X_0\supset X_1\supset \ldots \supset X_m=\varnothing\) of closed subschemes of \(X\), such that for each \(j=1,\ldots,m\), the complement \(X_{j-1}\setminus X_j\) is one of the \(X_\alpha\)'s.
\xdefi

In particular, every stratification is filtrable. While not every filtrable decomposition is a stratification, it is enough for the purposes of inductively applying localization, as in the following proof.

\lemm \thlabel{lemm--cellular-coh}
Let $f \colon Y \to S$ be a scheme with filtrable decomposition by cells isomorphic to products of $\A^1$ and $\Gm$. If $\dim Y \leq d$ then the following holds:
\begin{itemize}
   \item $f_! f^*$ maps $\DTMrx(S)^{\le 0}$ to $\DTMrx(S)^{\leq 2d}$, and in particular, we have $f_! \Z \in \DTMrx(S)$.
   \item \label{item--top coho} We have $\pe H^{2d} (f_!\Z) = \bigoplus_{C_d \subset Y} \Z(-d)$, i.e., one summand for each cell $C_d \subset Y$ of dimension $d$.
   \item For $L \in \MTMrx(S)$  there is a canonical isomorphism $\pe H^{2d} (f_!f^*L) \cong \pe H^{2d} (f_!\Z) \otimes L$.
 \end{itemize} 
\xlemm

\pf
By excision, i.e., the distinguished triangle $j_! j^* \r \id \r i_! i^*$ for a closed embedding $i$ with complement $j$, $f_!\Z \in \DTMrx(S)$. 
We prove the remaining statements by induction on the number of cells in $Y$. 
The statement is clear if $Y$ has one cell by the K\"unneth formula for $f_!$ (cf.~\refsect{motives}), since $\A^1$ has cohomology $\Z(-1)[-2]$ and $\Gm$ has cohomology $\Z[-1] \oplus \Z(-1)[-2]$. For the inductive step, let $Z \subset Y$ be a closed cell and let $U$ be the complement. Let $f_U$ and $f_{Z}$ be the structure maps of $U$ and $Z$ to $S$. If $\dim Z < d$ then by excision we have $\pe H^{2d}(f_{U!} \Z) \cong \pe H^{2d}(f_! \Z)$, so we are done by induction. If $\dim Z = d$, we have an exact sequence 
$$\pe H^{2d-1}(f_{Z!} \Z) \to \pe H^{2d}(f_{U!} \Z) \to \pe H^{2d}(f_{!} \Z) \to \pe H^{2d}(f_{Z!} \Z) \to 0.\eqlabel{exact.sequence.blah}$$ 
The left term is a free graded abelian group in which all Tate twists lie between $0$ and $-d + 1$, and by induction the next term is a free abelian group with Tate twist $-d$. Hence the left map is zero, so that $H^{2d}(f_{!} \Z)$ is an extension of $H^{2d}(f_{U!} \Z)$ by $\Z(-d)$. We have $\Ext^1_{\DTMrx(S)}(\Z(-d), \Z(-d)) = 0$: for reduced motives this holds since $\DTMr(S) = \D(\grModZ)$. For regular motives, if $S$ is smooth over a field, the group equals $\CH^{-1}(S, 1)=0$. If $S$ is smooth over a Dedekind ring, we reduce to the field case by the same argument as in \thref{Anti.Orthogonal}.
Using this vanishing, we are done by induction. To prove the last statement, by the projection formula we have $f_!f^* L \cong f_! \Z \otimes L$. By the previous excision argument, $f_!f^* L \in \DTM(S)^{\leq 2d}$. Since $(-) \otimes L$ preserves $\DTMrx(S)^{\leq 0}$, applying $\pe H^{2d}$ to the projection formula gives $\pe H^{2d}(f_!f^* L) \cong \pe H^{2d}(f_! \Z) \otimes L$.
\xpf

\subsection{Equivariant motives}
\label{sect--equivariant.motives}

\subsubsection{Basic definitions and averaging functors}
The functor in \refeq{DM.prestack} gives a category $\DMrx(Y)$ of (reduced) motives on any prestack $Y$ over $S$, and a !-pullback functor between such categories, for \emph{any} map of prestacks.
An important example of a prestack is a quotient prestack
$$G \setminus X := \colim \left (\dots G \x_S G \x_S X \triplerightarrow{}  G \x_S X \doublerightarrow{}{} X \right),$$
where $X$ is any prestack acted upon by a group prestack $G$. 
(An example coming up below is the quotient $LG/L^+G$ of the loop group, which is an ind-scheme, by the positive loop group, which is a group scheme, although not of finite type.)
For such quotients, the definition gives
$$\DMrx(G \setminus X) = \lim \left ( \DMrx(X) \doublerightarrow{}{} \DMrx(G \x_S X) \triplerightarrow{} \DMrx(G \x_S G \x_S X) \r \dots \right ),\eqlabel{equivariant.DM},$$
where the limit is formed using !-pullback (along the various action and projection maps). 
\rema
In colloquial terms this means that an object $M \in \DM(G \setminus X)$ is a collection of motives $M_n \in \DM(G^{\x n} \x X)$ together with isomorphisms $a^! M_0 \cong M_1 \cong p^! M_0$ and likewise for the higher order terms, subject to compatibility conditions.
The category is equivalent to its full subcategory spanned by objects for which $M_1 = p^! M_0$ (and the right hand isomorphism is the identity), and similarly for $M_n$, $n \ge 2$. In other words, one may forget the $M_n$ for $n \ge 1$, and only keep $M_0$ and $a^! M_0 \cong p^! M_0$ and higher-order isomorphisms. 
We may also complete the semi-cosimplicial diagram of \ii-categories in \refeq{equivariant.DM} into a cosimplicial diagram, by adding !-pullbacks along the various unit maps $G^n \r G^{n+1}$. 
By the finality of $\Delta_{\text{inj}}$ in $\Delta$ \cite[Lemma~6.5.3.7]{Lurie:Higher}, this yields an equivalent \ii-category \cite[Proposition~4.1.1.8]{Lurie:Higher}; cf.~also similar a similar discussion in \cite[§III.15, p.~187]{KiehlWeissauer}.
In other words the limit in \refeq{equivariant.DM} encodes the usual idea of $G$-equivariant sheaves on $X$.
Denoting the natural map $X \r G \setminus X$ by $u$, the map $u^! : \DM(G \setminus X) \r \DM(X)$ is simply the functor forgetting everything but the $M_0$ above.
\xrema

\rema
If $G$ is smooth over $S$, then one can replace !-pullbacks by *-pullbacks in the discussion above, giving an equivalent \ii-category.
The usage of !-pullbacks (as opposed to *-pullbacks) has the advantage of giving a uniform theory applicable to all prestacks (including arbitrary quotients $G \setminus X$, and also ind-schemes), as is explained in \cite[Remark~2.2.2, iv)]{RicharzScholbach:Intersection}.
\xrema

\rema
\thlabel{equivariant.functoriality}
(Functoriality for equivariant motives)
Suppose $f : X \r Y$ is a $G$-equivariant map of prestacks, and write $\ol f : G \setminus X \r G \setminus Y$ for the induced map. If $f^!$ admits a left adjoint $f_!$, then the adjoint functor theorem (cf.~\cite[Lemma~2.2.9]{RicharzScholbach:Intersection}) guarantees the existence of a left adjoint, denoted $\ol f_!$, of $\ol f^!$. 
If $G$ is a pro-smooth group scheme over $S$, then the !-pullback along the (pro-smooth) action and projection maps commute with $(\id_{G^n} \x f)_!$. Therefore $\ol f_!$ can be computed naively, i.e., $\ol f_!$
is compatible (via the forgetful functors, i.e., !-pullbacks along $X \r G \setminus X$ etc.) with $f_!$.

This construction of adjoints can be iterated.
For example, if $f$ is a proper schematic map, there are adjoints $(\ol f^*, \ol f_!, \ol f^!)$ between the categories $\DMrx(G \setminus X)$ and $\DMrx(G \setminus Y)$.
Again, if $G$ is pro-smooth, then these functors can be computed as $f^*$ etc. on the level of the underlying motives.
Similarly, if $f$ is, say, a $G$-equivariant map of finite type $S$-schemes (so that $f^*$ exists) and $G$ is pro-smooth, then there is an adjunction  $\ol f^* : \DMrx(Y/G) \rightleftarrows \DMrx(X/G) : \ol f_*$, whose adjoints reduce to the usual $f^*$ and $f_*$ on the level of non-equivariant motives.
\xrema

Recall that the Verdier duality functor is defined as $\Du : \DM(X)^{\opp} \r \DM(X)$, $\Du(\calF) := \IHom(\calF, \omega_X)$, where $\omega_X := p^! \Z$ is the dualizing sheaf, with $p : X \r S$ being the structural map.
Note that the usage of the terminology ``duality'' is an abuse: in the generality of motives with integral coefficients over $\Spec \Z$, say, $\Du$ need not be an equivalence.

\lemm
\thlabel{averaging.functor}
Let a smooth algebraic group $G$ act on a scheme $X$.
\begin{enumerate}
\item 
\label{item--av.DM.exists}
The forgetful functor $u^! : \DMrx(G \setminus X) \r \DMrx(X)$ admits a left adjoint $\coav := \coav_G$ and a right adjoint $\av := \av_G$,
called \emph{(co)averaging functors}.
\item 
 \label{item--av.DM.explicit}
The composite $u^! \av_G$ can be computed as $a_* p^* = a_* p^! (-d)[-2d]$, where $d := \dim G / S$ and $G \x X \stackrel[p]{a} \rightrightarrows X$ are the action and projection map.
Likewise, $u^! \coav_G = a_! p^!$.

\item
\label{item--av.duality}
Verdier duality (denoted by $\Du$) exchanges averaging and coaveraging functors in the sense that there is a natural isomorphism of functors $\Du_{X/G} \coav = \av \Du_X$.

\item 
\label{item--av.DM.red}
The reduction functor $\rho_\red$ commutes with $\av_G$ and $\coav_G$ (for $\DM$, respectively $\DMr$).

\item
\label{item--av.base.change}
If $f : Y \r X$ is a map of $G$-schemes, then $u^!$ and $\av$ commute with $f^!$. 
Moreover, $\coav$ commutes with $f^*$ (and therefore with $f^!$ if $f$ is smooth).

\item
\label{item--av.box}
Given another such pair $(G', X')$, there is an isomorphism
$$\coav_{G \x_S G'}(- \boxtimes -)  \stackrel \cong \r \coav_{G}(-) \boxtimes \coav_{G'}(-).$$
The same holds for the averaging functors if $S = \Spec k$ is a field of characteristic zero.
\end{enumerate}
\xlemm

\pf
\refit{av.DM.exists}: This follows from the adjoint functor theorem, but also from the following explicit description, which proves \refit{av.DM.explicit}.
We describe equivariant motives via the limit description as in \refeq{equivariant.DM}, using the isomorphism $X \cong G \setminus (G \x X)$.
The following diagram displays the low degrees of the simplicial diagrams whose colimits are shown in the bottom line:
$$\xymatrix{
X \ar@{=}[d] \ar@<.5ex>[r]^(.4)e & G \x (G \x X) \ar@<.5ex>[l]^(.6){p_X}  \ar@<.5ex>[d]^{p_{G \x X}} \ar@<-.5ex>[d]_{a_G \x \id_X} \ar[r]^(.6){\id_G \x a} & G \x X \ar@<.5ex>[d]^{p_X} \ar@<-.5ex>[d]_{a}
\\
X \ar@{=}[d] \ar@<.5ex>[r]^(.4)e & G \x X \ar@<.5ex>[l]^(.6){p_X} \ar[d] \ar[r]^(.6)a & X \ar[d] \\
X \ar@<.5ex>[r] & G \setminus (G \x X) \ar@<.5ex>[l]^(.6){\cong} \ar[r]^(.6){G \setminus a} & G \setminus X.
}$$
The left horizontal bottom maps are isomorphisms of prestacks.
Under this equivalence, the functor $u^!$ is induced levelwise by !-pullback along the right horizontal maps, i.e., $(\id_{G^n} \x a)^!$ in degree $n \ge 0$.
These functors admit a left (resp.~a right) adjoint, namely 
$(\id_{G^n} \x a)_!$ (resp.~$(\id_{G^n} \x a)_* (-d)[-2d]$), which both commute with !-pullback along the vertical maps in the right and middle diagram since the squares are cartesian and $G$ is smooth.
Thus, they assemble to the asserted adjoints of $u^!$.

Under the isomorphism $X \cong G \setminus (G \x X)$, the composite $u^! \av_G$ is given by $a^! a_* (-d)[-2d]$.
In terms of motives on $X$, this means we have to evaluate
$$e^! a^! a_* p^! (-d)[-2d],$$
which is isomorphic to $a_* p^*$.
Likewise, $u^! \coav_G$ corresponds to the endofunctor $a^! a_!$ on $\DM(G \x X)$, and $e^! (a^! a_!) p^! = a_! p^!$.

\refit{av.duality}: This follows from \refit{av.DM.explicit} since for any map $f : Y \r Z$ of finite type $S$ schemes, we have $\Du f_! = f_* \Du$, and for any smooth map (such as $f = a$ or $p$) we have $f^* \Du_Z = \IHom(f^* -, f^* \omega_Z) = \IHom(f^! -, \omega_Y) = \Du_Y f^!$, as consequences of the projection formula, resp.~the projection formula for $f_\sharp$ vs.~$\otimes$ and relative purity \cite[Theorem~2.4.50, \S 1.1.33]{CisinskiDeglise:Triangulated}.

\refit{av.DM.red}: By definition, $u^!$ commutes with $\rho_\red$, so there is a natural map $\rho_\red \av_{\DM} \r \av_{\DMr} \rho_\red$ (and analogously for $\coav$). 
To check it is an isomorphism it suffices to append the conservative functor $u^!$, so the claim follows from \refit{av.DM.explicit} since $\rho_\red$ commutes with *-functors.

\refit{av.base.change}: This is similarly reduced to the observation that $a_* p^*$ (resp.~$a_! p^!$) commutes with $f^!$ (resp.~$f^*$) as asserted, by base-change and relative purity (applied to the smooth map $p$).

\refit{av.box}: By the definition of $\boxtimes$ on prestacks of the form $X / G$ \cite[Appendix~A]{RicharzScholbach:Motivic}, the forgetful functor is compatible with exterior products.
This gives a map 
$\av_G (-) \boxtimes \av_{G'}(-) \r \av_{G \x G'}(- \boxtimes -)$. In order to check it is an isomorphism, we apply $(u \x u')^! = u^! \boxtimes u'^!$, so we need to prove $(a_* p^* -) \boxtimes (a'_* p'^* -) = (a \x a')_* (p \x p')^* (- \boxtimes -)$.
This holds by the Künneth formula for *-pushforwards \cite[Theorem~2.4.6]{JinYang:Kuenneth} (this needs the assumption on $S$, which is used to apply resolution of singularities).
The argument for $\coav$ is similar but only uses the compatibility of $\boxtimes$ with !-pullbacks along smooth maps (such as the action and projection maps for the smooth group scheme $G$) and !-pushforwards (i.e., the projection formula, cf.~\cite[Lemma~2.2.3]{JinYang:Some}). 
\xpf

\subsubsection{Equivariant Tate motives}

\defi
\thlabel{equivariant.DTM}
For a stratified (ind-)scheme $X$ with an action of a group scheme $G$ (not necessarily of finite type over $S$), the category $\DTMrx(G \setminus X)$ of \emph{equivariant (reduced) Tate motives} on the prestack quotient $G \setminus X$ is defined as the full subcategory of $\DMrx(G \setminus X)$ whose underlying object in $\DMrx(X)$ is in $\DTMrx(X)$, cf.~\cite[\S 3.1]{RicharzScholbach:Intersection}.
Note that this is a slight abuse of notation in that it depends not only on the prestack $G \setminus X$, but rather on $X$ and $G$.
\xdefi

\lemm
\thlabel{equivariant.MTM}
\thlabel{MTM.equivariant.ind-scheme}
We continue to assume that $S$ satisfies the conditions around \refeq{BS.vanishing}.
Let $(X, X^\dagger) \stackrel {(\pi, \pi^\dagger)} \r (Y, Y^\dagger) \r S$ be a Whitney--Tate stratified ind-scheme over an admissibly stratified Whitney--Tate scheme $Y$. We also suppose that $\pi^\dagger$ is admissible.
Finally, let a pro-smooth group $G = \lim G_k / Y$ act on $X = \colim X_k$ such that the following conditions are met:
\begin{itemize}
  \item The $G$-action preserves $X_k$ and factors there over $G_k$. 
  In addition the $G_k$-action preserves the stratification, cf.~\cite[Definition~3.1.26]{RicharzScholbach:Intersection}.
  \item $\ker (G \r G_k)$ is a split pro-unipotent $Y$-group \cite[Definition~A.4.5]{RicharzScholbach:Intersection}.
  \item $Y^\dagger \x_Y G_k$ is cellular (over $Y^\dagger = \bigsqcup_w Y^w$, i.e., over each stratum $Y^w \subset Y$).
\end{itemize}
Then $\DTMrx(G \setminus X) = \colim \DTMrx(G_k \setminus X_k)$ (\thref{equivariant.DTM}) carries a (unique) right-complete t-structure such that the forgetful functor to $\DTMrx(X)$ is t-exact.
Its heart denoted by $\MTMrx(G \setminus X)$ is equivalent to 
$$\colim \MTMrx(G \setminus X_k) = \colim \MTMrx (G_k \setminus X_k).$$
Each $M \in \MTMrx(G \setminus X)$ is a filtered colimit of bounded subobjects, namely (for $i_k : X_k \r X$)
$$M = \colim_k i_{k!} \pe i_k^! M.\eqlabel{obj MTM colimit}$$

The objects $\IC_{w, L} := \iota^w_{!*} \pi^{w*} (L[\dim \pi^w])$,  for $L \in \MTMrx(Y^w)$ naturally lie in $\MTMrx(G \setminus X)$.
If, in addition to the above, each stratum $X^w_k \subset X_k$ is of the form $X^w_k = (G_k / H^w_k)_\Nis$ for some subgroup scheme $H^w_k$ that is  cellular and fiberwise (over $Y^w$) connected, then the objects $\IC_{w, L}$ generate $\MTMrx(G \setminus X)$ under filtered colimits and extensions.
\xlemm

\pf 
The t-structure on $\DTMrx(X)$ afforded by \thref{t-structure.stratified} yields a t-structure on $\DTMrx(G \setminus X)$, which is shown (using the cellularity of $G|_{Y^\dagger}$) as in \cite[Proposition~3.2.15]{RicharzScholbach:Intersection}.
Using the notation of \thref{equivariant.functoriality}, the object $\im (\pH^0 \ol \iota_{!}^w L \r \pH^0 \ol \iota_{*}^w L) \in \MTMrx(G \setminus X)$ maps to $\IC_{w, L}$ under the forgetful functor.

As in \refeq{DM.colim.lim}, one has an equivalence
$\colim \MTMrx(G \setminus X_k) = \lim \MTMrx(G \setminus X_k)$, where the limit is formed using the right adjoints of $(i_{kk'})_!$ for $k' > k$, i.e., the truncated functors $\pe i_{kk'}^!$.
This formally implies the isomorphism \refeq{obj MTM colimit}.
In addition, since the $(i_{kk'})_!$ are t-exact and therefore $\pe i_{kk'}^!$ are left t-exact, $i_{k!} \pe i_k^! M$ is indeed a subobject of $M$. 

In particular, to see the final claim, we may replace $X$ by some $X_k$ and $G$ by $G_k$, which is smooth.
As in \thref{t-structure.stratified}\refit{MTM glued}, we may then replace $Y$ by a single stratum.
By assumption, the stratification $X^\dagger$ on $X = (G / H)_\Nis$ is admissible over $Y$.
The admissibility of $\pi^\dagger: X^\dagger \r Y$ implies that $\pi^{w!} [-\dim_Y X^w]$ is t-exact.
Given these two properties of this pullback functor, the proof in \cite[Proposition~3.2.23]{RicharzScholbach:Intersection} carries over: the cellularity and fiberwise connectedness of $H$ implies an equivalence of categories $\MTMrx(G \setminus X) = \MTMrx(Y)$.
\xpf

\rema
In contrast to the case of rational coefficients in \cite[Proposition~3.2.23]{RicharzScholbach:Intersection}, for integral coefficients it is not necessarily true that $\MTM(S)$ (or the other $\MTM$-categories considered above) is generated under extensions by the motives $\Z(k)$ and $\Z/n(k)$.
In fact, in unpublished work, Spitzweck shows that the cofiber of the natural map (in $\DTM(\Spec \Z)$), $\Z/2 \to \Z/2(1)$ corresponding to the element $-1 \in K_1(\Z)$ lies in the heart of the motivic t-structure.
\xrema

\subsection{Rational and modular coefficients}
\label{subsection:coefficients}

In order to compute the dual group of the Satake category, we will also need to work with rational and $\Fp$-coefficients.
For $\Lambda = \Q$ or $\Fp$ and a prestack $X$, the category $\DMrx(X, \Lambda)$ of \emph{(reduced) motives on $X$ with $\Lambda$-coefficients} is defined as $\DMrx(X) \t_{\D(\ModZ)} \D(\ModL)$. By construction, $\DM(X, \Q)$ is the category of Beilinson motives \cite[\S 14]{CisinskiDeglise:Triangulated}. The existence and properties of the six functors holds without any change, as does the definition and properties of (stratified) Tate motives. Thus, the category $\DTM(X, \Q)$ of stratified Tate motives with rational coefficients is exactly the one considered in \cite{RicharzScholbach:Intersection,RicharzScholbach:Motivic}.

\lemm
\thlabel{DTM.Q}
Let $\Lambda = \Q$ or $\Fp$.
The category $\DTMrx(S, \Lambda)$ carries a t-structure whose $\le 0$-aisle is generated under colimits and extensions by $\Lambda(k)$, $k \in \Z$.
These objects lie in the heart of the t-structure, which is denoted by $\MTMrx(S, \Lambda)$. Moreover, they form a family of compact generators of the heart.

The forgetful functor $U: \DTMrx(S, \Lambda) \r \DTMrx(S)$ is t-exact. For $\Lambda = \Q$ it is fully faithful with essential image being the objects $\calF$ such that $n \cdot \id_\calF$ is an isomorphism for all nonzero integers $n$.
For $\Lambda = \Fp$ its restriction to $\MTMrx(S, \Fp)$ is fully faithful, with essential image being the objects $\calF \in \MTMrx(S)$ such that $p \cdot \id_\calF = 0$.

For $S = \Spec \Fq$, the reduction functor $\rho_\red : \DTM(S, \Q) \r \DTMr(S, \Q)$ is an equivalence.
\xlemm

\pf
The functor $- \t \Lambda : \DTM(S) \r \DTM(S, \Lambda)$ is right t-exact by definition of the t-structures.
Being the right adjoint of a right t-exact functor, $U$ is left t-exact.
The full faithfulness of $U$ for $\Lambda = \Q$ holds by $\Q \t_\Z \Q = \Q$.

In order to check $\Lambda(k) \in \MTM(S, \Lambda)$ it suffices to see $\Hom_{\DTM(S, \Lambda)}(\Lambda, \Lambda(k)[-n]) = \Hom_{\DTM(S)} (\Z, \Lambda(k)[-n]) = 0$ for $n > 0$. This is clear for $\Lambda = \Q = \colim \frac 1n \Z$. To see this vanishing for $\Lambda = \Fp$, use the exact sequence
$$\dots \r \underbrace{\Hom(\Z, \Z(k)[-n])}_{=0} \r \Hom(\Z, \Fp(k)[-n]) \r \Hom(\Z, \Z(k)[-n+1]) \stackrel p \r \Hom(\Z, \Z(k)[-n+1]).$$
The right-hand group vanishes for $n > 1$. For $n = 1$, it vanishes unless $k = 0$. In this case the group equals $\Z$ (cf.~\thref{nota-basescheme} for our standing assumptions on $S$), on which the $p$-multiplication is injective.
This also means that $U(\Fp(k))$ and $U(\Q(k))$ lie in $\MTM(S)$. Thus $U$ is right t-exact.

To check the full faithfulness of $U|_{\MTM(S, \Fp)}$, it suffices to observe that 
$\pH^0(- \t \Fp)$ is its left adjoint and that for $A \in \MTM(S, \Fp)$, $\pH^0(U(A) \t \Fp)$ is the cokernel of the $p$-multiplication on $U(A)$, which is therefore isomorphic to $A$.

The final statement is \cite[Proposition~5.3]{EberhardtScholbach:Integral}.
\xpf

\rema
Once the t-structure for motives with $\Lambda$-coefficients exists on the base scheme $S$, it extends mutatis mutandis to admissibly stratified schemes, as in \thref{t-structure.stratified}.
The above t-exactness and full faithfulness properties of $U$ carry over to that situation.
\xrema

\rema
\thlabel{functor i}
The reduction functor $\rho_\red : \DTM(S) \r \DTMr(S)$ admits a section $i : \DTMr(S) = \D(\grAb) \r \DTM(S)$; this functor is the unique colimit-preserving functor sending $\Z(k)$ (i.e., $\Z$ in graded degree $-k$) to $\Z(k)$ (i.e., the $k$-fold Tate twist).
Given that $\Z/n(k) \in \MTM(S)$, it restricts to a functor
$$i : \MTMr(S) = \grAb \r \MTM(S).$$
This functor is faithful, since $\Hom_{\MTMr(S)}(A, B(k)) = 0$ for any two abelian groups $A, B$ and $k \ne 0$.
Note it is not \emph{fully} faithful; e.g.~$\Hom_{\MTM(\Spec \Z)} (\Z/2, \Z/2(1)) \ne 0$.
However, the restriction of $i$ to flat (or, ind-free) graded abelian groups is fully faithful.
\xrema

\section{Affine Grassmannians}
\label{sect--aff Grass}

\subsection{Definitions and basic Whitney--Tateness properties}
Throughout this paper, $G$ denotes the base change to $S$ of a split reductive group over $\Z$ (all reductive groups are assumed to be connected). We fix a split maximal torus and a Borel $T \subset B \subset G$, also defined over $\Z$. Let $X^*(T)$ (resp. $X_*(T)$) be the group of (co)characters, and let $X_*(T)^+$ be the monoid of dominant cocharacters with respect to \(B\).
By a parabolic subgroup of $G$, we mean a subgroup $P \subset G$ containing $B$ associated to a subset of the simple roots. 

In order to give some constructions uniformly, we let $\mathcal{G}$ be a smooth affine group scheme over $S$. The loop group (resp.~ positive loop group) is the functor $L\mathcal{G} : \AffSch_S^\opp \r \Set$, $\Spec R \mapsto \mathcal{G}(R\rpot{t})$ (resp.~ $L^+\mathcal{G}(R) = \mathcal{G}(R\pot{t})$). 
The affine Grassmannian of $\mathcal{G}$ is the étale sheafification of the presheaf quotient 
$$\Gr_{\mathcal{G}} := (L\mathcal{G}/L^+\mathcal{G})_\et.\eqlabel{definition.Gr.G}$$  
If \(S=\Spec A\) is affine, we can more generally associate a positive loop group with a smooth affine group \(\mathcal{G}\) over \(A\pot{t}\) via \(L^+\mathcal{G}(R)=\mathcal{G}(R\pot{t})\), as in the case of parahoric groups below, and similarly for the loop group and affine Grassmannian.
The above choices determine a standard apartment of $G$ with origin $0$ and a standard alcove $\mathbf{a}$. 
(For most of this paper, a rudimentary understanding of the standard apartment is sufficient, e.g.,~as reviewed in \cite[§4.1]{RicharzScholbach:Intersection}. The theory of buildings will be used in more detail in Section \ref{subsec--GL}, and the necessary definitions will be introduced there.)
For any facet $\mathbf{f}$ in the closure of $\mathbf{a}$, there is an associated parahoric subgroup $\mathcal{P} \subset LG$. 
If $\mathbf{f}$ contains the origin in its interior, then $\mathcal{P} = L^+G$, and if $\mathbf{f} = \mathbf{a}$ we write $\mathcal{I}$ for the corresponding Iwahori subgroup. 

\defilemm
\thlabel{Sheafifications of quotient agree}
For a parahoric subgroup $\calP \subset LG$, the Zariski, Nisnevich, and étale sheafifications of the quotient $LG/\calP$ agree.
We denote this common quotient by $\Fl_\calP$, and call it the \emph{(partial) affine flag variety} of \(\calP\).  The sheaf $\Fl_\calP$ is represented by an ind-projective scheme over $S$. As usual, we denote the \emph{affine Grassmannian} by $\Gr_G:= \Fl_{L^+G}$,  and the \emph{full affine flag variety} by $\Fl := \Fl_\calI$. 
\xdefilemm

\pf Because sheafification commutes with base change we may assume $S =\Z$. Then the agreement of the different sheafifications of $LG/\calP$ follows from \cite[Def. 5 ff.]{Faltings:Loops} (see also \cite[Theorem 2.3.1]{dHL:Frobenius}). 
The representability of $\Fl_\calP$ by a separated ind-finite type ind-scheme is a consequence of \cite[Corollary 11.7]{PappasZhu:Kottwitz} and \cite[Corollary 3.11 (i)]{HainesRicharz:TestFunctionsWeil}, which also show that $\Fl_\calP$ is ind-projective as soon as it is ind-proper. The ind-properness of $\Gr_G$ follows from \cite[Corollary 3.11 (iii)]{HainesRicharz:TestFunctionsWeil} (cf.~also \cite[Theorem 4.5.1 (iv)]{BeilinsonDrinfeld:Quantization}), and hence $\Fl$ is also ind-proper since the projection $\Fl \r \Gr_G$ is a $G/B$-torsor. Now by considering the surjection $\Fl \r \Fl_{\calP}$ we conclude that the target is also ind-proper by \StP{03GN}.
\xpf

The following statement will be used to classify $L^+G$-equivariant motives on $L^+G$-orbits, which by the Cartan decomposition are indexed by $X_*(T)^+$.
For \(\mu\in X_*(T)\), we denote by \(t^\mu\colon S \xrightarrow{t} L\Gm \xrightarrow{\mu} LT \to LG \to \Gr_G\) the corresponding point in the affine Grassmannian.
Fix two parahorics \(\Pp,\Qq\subseteq LG\) corresponding to facets \(\mathbf{f},\mathbf{f}'\) in the closure of $\mathbf{a}$. 
Then the \(\Qq\)-orbits on \(\Fl_{\Pp}\) are indexed by \(W_{\mathbf{f}'} \backslash \widetilde{W} / W_{\mathbf{f}}\), where \(\widetilde{W}\) is the extended affine Weyl group (or Iwahori-Weyl group), and \(W_{\mathbf{f}},W_{\mathbf{f}'}\) are the stabilizing subgroups as in \cite[(4.2.10)]{RicharzScholbach:Intersection}.
For \(w\in \widetilde{W}\), we denote by \(\dot{w}\in LG(\Z)\) a representative of \(w\).

\prop\thlabel{Stab.Split}
The stabilizer \(\Qq_w\subseteq \Qq\) of \(\dot{w}\cdot e\in \Fl_{\Pp}(\Z)\) is represented by a closed subgroup, which is an extension of a split reductive \(\Z\)-group by a split pro-unipotent \(\Z\)-group in the sense of \cite[Definition A.4.5]{RicharzScholbach:Intersection}.
The étale sheaf-theoretic image \[\Fl_{\Pp}^w := \Qq\cdot \dot{w}\cdot e\subseteq \Fl_{\Pp}^{\leq w}\]
agrees with \((\Qq/\Qq_w)_{\et}\), where \(\Fl_{\Pp}^{\leq w}\) is the scheme-theoretic image. Using a superscript \(n\) to denote jet groups, this quotient agrees with \((\Qq^n/\Qq_w^n)_{\Zar}\) for \(n\gg 0\).
For such \(n\), \(\Qq_w^n\) is an extension of a split reductive \(\Z\)-group by a split unipotent \(\Z\)-group.
\xprop
We note that \(\Qq\) arises as the positive loop group of a \(\Z\pot{t}\)-group scheme, so that we can indeed consider \(\Qq^n\). When $\mathcal{P} = \mathcal{Q} = L^+G$, we write the orbits and their closures as $\Gr_G^\mu$ and $\Gr_G^{\leq \mu}$ for $\mu \in X_*(T)^+$.
\pf
The representability of \(\Qq_w\) follows from \cite[Lemma 4.3.7]{RicharzScholbach:Intersection}, where it is shown that \(\Qq_w\) agrees with the subgroup scheme of \(\Qq\) corresponding to the subset \(\mathbf{f}'\cup w\mathbf{f}\) of the standard apartment; we refer to loc.~cit.~for details.
Then the desired description of \(\Qq_w\) was shown in \cite[Lemma 4.2.7, Remark 4.2.8]{RicharzScholbach:Intersection}.
Moreover, we have \(\Fl_{\Pp}^w\cong (\Qq/\Qq_w)_{\et} \cong (\Qq^n/\Qq_w^n)_{\et}\) for \(n\gg 0\) by \cite[Lemma 4.3.7 (ii)]{RicharzScholbach:Intersection}.
It remains to show that \(\Qq^n\to (\Qq^n/\Qq_w^n)_{\et}\) admits sections Zariski-locally.

If \(\Qq=\mathcal{I}\) is an Iwahori, then \(\Fl^w_{\Pp}\) is an affine space, and \(\Qq^n\to \Fl^w_{\Pp}\) admits a global section by the proof of \cite[Proposition 4.3.9 (i)]{RicharzScholbach:Intersection}.
For a general parahoric \(\Qq\supseteq \mathcal{I}\), the \(\Qq\)-orbit \(\Fl^w_{\Pp}\) contains a unique open dense \(\mathcal{I}\)-orbit, which yields a section of \(\mathcal{I}^n\to \Qq^n\to \Fl^w_{\Pp}\) over this orbit.
It then suffices to find enough \(\Z\)-points of \(\Qq^n\), to translate this local section to sections over a Zariski cover of \(\Fl^w_{\Pp}\).
For this, we claim that translates by representatives of $W_{\mathbf{f}'}$ suffice. 
Indeed, if \(\mathsf{Q}\) is the maximal reductive quotient of \(\Qq^0\), then  \(\Fl^w_{\Pp}\) is \(\mathsf{Q}\)-equivariantly an affine bundle over a partial flag variety $X$ for \(\mathsf{Q}\) (this is similar to \cite[Proposition 8.7.(b)]{PappasRapoport:Twisted}). Moreover, the image of $\calI$ in \(\mathsf{Q}\) is a Borel, and $W_{\mathbf{f}'}$ is identified with the Weyl group of \(\mathsf{Q}\). Now the claim follows since $X$ is covered by the Weyl-translates of its open Borel orbit \cite[II.1.10 (5)]{Jantzen:Representations}.
\xpf

\rema
By \cite[Proposition 4.4.3]{RicharzScholbach:Intersection} the formation of \(\Fl_{\Pp}^{\leq w}\) commutes with base change along $S \r \Spec \Z$ up to a nil-thickening. Since nil-thickenings induce equivalences on categories of motives, we can safely view \(\Fl_{\Pp}^{\leq w}\) (and the orbits \(\Fl_{\Pp}^{w}\)) over an arbitrary base $S$ as a base change from $S = \Spec \Z$; cf.~\cite[\S 4.4]{RicharzScholbach:Intersection} for more details.
\xrema

\prop \thlabel{Orbit.Eq}
Let  \(w \in W_{\mathbf{f}'} \backslash \widetilde{W} / W_{\mathbf{f}}\) and let $e \colon S \r \mathcal{Q} \backslash \Fl_\calP^w$ be the inclusion of the point corresponding to $\dot{w}$.  Then we have an equivalence
$$e^![\dim \Fl_{\Pp}^{w}] \colon \MTMrx(\mathcal{Q} \backslash \Fl_\calP^w) \cong \MTMrx(S).$$
\xprop

\pf 
The orbit $\Fl_\calP^w$ is affinely stratified by Iwahori orbits \cite[Proposition 4.3.9]{RicharzScholbach:Intersection}, so it is cellular and hence admissible by \thref{cellular.stratification}. We will verify the assumptions of  \thref{MTM.equivariant.ind-scheme}, where we take $Y = S$ and $X = \Fl_\calP^w$ with the trivial stratification, acted upon by $\mathcal{Q}$. The jet group construction \cite[(4.2.5)]{RicharzScholbach:Intersection} expresses $\mathcal{Q}$ as a pro-algebraic $S$-group, with the action on $\Fl_\calP^w$ factoring through some jet quotient by \cite[Lemma A.3.5]{RicharzScholbach:Intersection} (or by \thref{Stab.Split} above). The pro-smoothness of $\mathcal{Q}$, cellularity of the jet quotients, and split pro-unipotence of the kernels are all shown in \cite[Lemma 4.2.7]{RicharzScholbach:Intersection}. Now thanks to \thref{Stab.Split}, the final part of the proof \thref{MTM.equivariant.ind-scheme} shows the desired equivalence,  cf.~also \cite[Theorem~5.3.4]{RicharzScholbach:Intersection} for a related statement.
\xpf

\lemm \thlabel{GrP.LC} 
For any parabolic subgroup $P \subset G$, $\Gr_P$ is represented by an ind-scheme of ind-finite type over $S$. The natural morphism $\Gr_P \r \Gr_G$ identifies $\Gr_P$ with the attracting locus of a $\Gm$-action on $\Gr_G$, and it restricts to a locally closed embedding on connected components of $\Gr_P$. If $k$ is a field, then $\Gr_P(k) \r \Gr_G(k)$ is a bijection.
\xlemm

\pf
Choose a cocharacter $\eta \in X_*(T)^+$ which is orthogonal to the simple roots associated to $P$, but not orthogonal to any other simple root. Then the action of $\Gm$ on $G$ via conjugation by $\eta$ extends to a $\Gm$-action on $\Gr_G$. By \cite[Theorem 3.17]{HainesRicharz:TestFunctionsWeil}, $\Gr_P$ identifies with the attractor $(\Gr_G)^+$ for this $\Gm$-action. 
The representability of $\Gr_P$ then follows from \cite[Theorem 2.1 (iii)]{HainesRicharz:TestFunctionsWeil} (cf.~also \cite[Theorem 4.5.1 (i)]{BeilinsonDrinfeld:Quantization}). To show that we have a locally closed embedding, we proceed as in the field case \cite[Lemma 3.7]{HainesRicharz:TestFunctions}, but using statements valid over a general base in \cite{HainesRicharz:TestFunctionsWeil} 
(recall that a locally closed embedding of ind-schemes is a map such that after base change by an affine scheme, we have a locally closed embedding of schemes in the usual sense). First note that there exists a closed embedding $G \r \mathbf{G}\mathbf{L}_n$ for some $n$. 
Then the proof of \cite[Lemma 3.16]{HainesRicharz:TestFunctionsWeil} shows the $\Gm$-action on $\Gr_G$ is Zariski-locally linearizable. Thus, we may write $\Gr_G = \colim X_i$, where each $X_i$ is projective over $S$ and $\Gm$-stable, and there is a $\Gm$-equivariant Zariski cover $U_i \r X_i$ which is affine.
Then $X^+ = \colim (X_i)^+$, and by \cite[Lemma 1.11]{Richarz:Spaces} (cf.~also \cite[Lemma 1.4.9 (i)]{DrinfeldGaitsgory:Theorem}) we have $(U_i)^+ = (U_i)^0 \times_{(X_i)^0} (X_i)^+$. 
Since $(U_i)^+$ is representable by a closed subscheme of $U_i$ \cite[Lemma 1.9]{Richarz:Spaces} (cf.~also \cite[Proposition 1.6.2 (ii)]{DrinfeldGaitsgory:Theorem}), we conclude that $\Gr_P \r \Gr_G$ is Zariski-locally on $\Gr_P$ a locally closed immersion. The final claim about points over an algebraically closed field $k$ will then imply it is a locally closed immersion on each connected component.
For this, 
we note that $\Gr_G(k) = LG(k)/L^+G(k)$ and $\Gr_P(k) = LP(k)/L^+P(k)$.
It follows that $\Gr_P(k) \r \Gr_G(k)$ is injective, and it is surjective by the Iwasawa decomposition of $G(k(\!(t)\!))$. The final claim for an arbitrary field then follows as well.
\xpf

In particular, \(\Gr_P\) induces a decomposition of \(\Gr_G\) into locally closed sub-ind-schemes.
At least when \(P=B\) is a Borel, we will see in \thref{Semi-infinite:Strat} that this decomposition is even a stratification, i.e., that the closure relations hold.
By \cite[Corollary 1.12]{Richarz:Spaces}, there are bijections \(\pi_0(\Gr_B)\cong \pi_0(\Gr_T)\cong X_*(T)\). The resulting connected components of \(\Gr_B\) are denoted by \(\Ss_\nu^+\) for \(\nu\in X_*(T)\), and called the \emph{semi-infinite orbits}. 
If $B^-$ is the opposite Borel subgroup, we denote by $\Ss_\nu^-$ the corresponding connected component of $\Gr_{B^-}$.
By \thref{GrP.LC}, we can and will view the \(\Ss_\nu^\pm\) as locally closed sub-ind-schemes of \(\Gr_G\).
We note that \(\Gr_B = \coprod_{\nu\in X_*(T)} \Ss_\nu^+\) is the attractor for the \(\Gm\)-action on \(\Gr_G\) induced by a regular dominant cocharacter, while \(\Gr_{B^-} = \coprod_{\nu\in X_*(T)} \Ss_\nu^-\) is the repeller for the same action.

The next proposition is an extension (but not strictly speaking a corollary) of the Whitney--Tateness of partial affine flag varieties \cite[Theorem~5.1.1]{RicharzScholbach:Intersection}.
It will be used in order to show the Whitney--Tateness of the Beilinson--Drinfeld Grassmannian.
The following lemma serves to show the anti-effectivity. 

\lemm \thlabel{A.Anti.Effective}
Let $X$ and $Y$ be ind-schemes, each having a Whitney--Tate stratification by affine spaces. Let $\pi \colon X \r Y$ be a smooth map which sends strata onto strata, and such that for each stratum $X_w \subset X$, the induced map on strata $X_w \r \pi(X_w)$ is a relative affine space. Then the functors $\pi_!$ and $\pi^* \pi_!$ preserve anti-effective stratified Tate motives.
\xlemm

\pf
By excision and base change, it suffices to consider a motive $\iota_{w !} \Z$, where $\iota_w \colon X_w \r X$ is a stratum. Then the lemma follows from the fact that the structure map $f \colon \A^n_S \r S$ satisfies $f_!(\Z) \cong \Z(-n)[-2n]$.
\xpf

\prop
\thlabel{Fl.universally.WT}
For any parahoric subgroups $\calP$, $\calP' \subset LG$, the stratification of the partial affine flag variety $\Fl_\calP$ by $\calP'$-orbits is anti-effective universally Whitney--Tate.
\xprop

\pf
This follows by revisiting the proof in \cite[Theorem~5.1.1]{RicharzScholbach:Intersection}: Beginning with the case where $\calP = \calP' = \calI$ is the Iwahori subgroup, let $\iota : \Fl^\dagger \r \Fl$ be the stratification map and, for any $S$-ind-scheme $Y$, let $\iota' : \Fl^\dagger \x_S Y \r \Fl \x_S Y$ be the product stratification.
For an element $w$ of the extended affine Weyl group, one shows by induction on the length $l(w)$ that $\iota'^! \iota'_{w!} \Z \in \DTM(\Fl^\dagger \x Y)$: this is clear if $l(w) = 0$.
Inductively, for $w = vs$ for a simple reflection $s$ and an element with $l(v) = l(w)-1$, there is a cartesian diagram, where $\calP_s$ is the parahoric subgroup associated to $s$ and the map $\pi$ arises from the inclusion $\calI \subset \calP_s$:
$$\xymatrix{
(\Fl^v \sqcup \Fl^{w}) \x Y \ar[r] \ar[dr]_{\pi^\dagger} & \pi^{-1}(\Fl_{{\calP_s}}^v \x Y) \ar[d]^{\tilde \pi} \ar[r] & \Fl \x Y \ar[d]^{\pi} \\
& \Fl_{{\calP_s}}^v \x Y \ar[r] & \Fl_{{\calP_s}} \x Y.
}$$
The map $\pi^\dagger$ is isomorphic to the disjoint union of $\id_{\Fl^v}$ and the projection $p : \A^1_{\Fl^v_{\calP_s}} \r \Fl^v_{\calP_s}$. More generally, the map $\pi$ is smooth and proper, and the induced map from each stratum of $\Fl \times Y$ onto its image in $\Fl_{{\calP_s}} \times Y$ is either an isomorphism or an affine space of relative dimension one.
Applying Verdier duality $\Du$ to the localization sequence \cite[(5.1.2)]{RicharzScholbach:Intersection} and noting that $\omega_{\A^n} = \Z(n)[2n]$ gives a fiber sequence
$$\iota'_{v*} \Z(-1)[-2] \r \pi^! \pi_! \iota'_{v*} \Z(-1)[-2] \r \iota'_{w*} \Z. \eqlabel{Fl.WT}.$$
Thus, the fact $p_! \Z = \Z(-1)[-2]$ shows that we have a Whitney--Tate stratification, and the condition of \refeq{universally.WT} being an isomorphism on the summand corresponding to $v$ implies the same for the summand of $w$.

The Whitney--Tateness for general $\calP, \calP' \subset LG$ is then treated identically as in loc.~cit. 
We prove the anti-effectivity of the stratification using the same reduction steps, as follows.

\textit{First case: $\calI =  \calP = \calP'$.} We apply \thref{A.Anti.Effective} to $\pi$ (we may assume $Y = S$), and conclude the claim using \refeq{Fl.WT} and $\pi^! = \pi^*(1)[2]$.

\textit{Second case: $\calI = \calP' \subset \calP$.}
By \cite[Lemma~4.3.13]{RicharzScholbach:Intersection}, the projection $\Fl \r \Fl_{\calP}$ satisfies the hypotheses of \thref{A.Anti.Effective}, where both ind-schemes are stratified by $\calI$-orbits. For a stratum $\iota_w \colon \Fl_{\calP}^w \r \Fl_{\calP}$, let $s \colon \Fl_{\calP}^w \r \Fl$ be a section which includes $ \Fl_{\calP}^w$ as a stratum of $\Fl$. Then $\iota_{w*} \Z = \pi_* s_*\Z$, so we conclude by the first case and \thref{A.Anti.Effective}.

\textit{Third case: $\calP$, $\calP'$ arbitrary.}
Let $\iota \colon \Fl_{\calP}^\dagger \r \Fl_{\calP}$ be the stratification by $\calP'$-orbits, and let $\iota'$ be the stratification of $\Fl_{\calP}^\dagger$ by $\calI$-orbits. To check if $\iota^* \iota_* \Z$ is anti-effective, we apply \thref{Anti.Orthogonal}. By localization, the condition $\Maps(\Z(p), \iota^* \iota_* \Z) = 0$ is equivalent to $\Map(\iota'_! \iota'^* \Z(p), \iota^* \iota_* \Z) = \Map(\iota'^* \Z(p), \iota'^! \iota^*\iota_* \Z) = 0$. Since the motive $\Z$ on $\Fl^\dagger_{\calP}$ is anti-effective with respect to $\iota'$,  we conclude by the second case.
\xpf

\lemm
\thlabel{base.change.Fl}
Consider some schematic map $f : X' \r X''$ of (ind-)schemes over $S$, some $M \in \DM(X')$ and some stratified Tate motive $N \in \DTM(\Fl_\calP)$,
where the stratification on $\Fl_\calP$ is by any $\calP'$-orbits.
Then the following natural maps are isomorphisms:
$$\eqalign{f_* M \boxtimes N \r & (f \x \id)_* (M \boxtimes N), \cr
(f \x \id)^! (M \boxtimes N) \r & f^! M \boxtimes N.}$$
\xlemm

\rema
Resolution of singularities implies that *-pushforward functors are compatible with exterior products: for a field $k$, and two maps $X' \stackrel f \r X''$, $Y' \stackrel g \r Y''$, the natural map 
$$f_* M \boxtimes g_* N \r (f \x g)_* (M \boxtimes N)$$ 
is an isomorphism for any $M \in \DM(X')$, $N \in \DM(Y')$ if $k$ is of characteristic 0 or if $\Char k$ is invertible in the ring of coefficients \cite[Theorem~2.4.6]{JinYang:Kuenneth}.
Since below we are interested in motives with integral coefficients, and work over $\Spec \Z$, we need to supply a more specific argument.
\xrema

\pf
We may refine our stratification and replace $\calP'$ by the Iwahori subgroup $\calI$.
The proof of \thref{Fl.universally.WT} implies that (cf.~\cite[Proposition~5.2.2]{RicharzScholbach:Intersection}) $\DTM(\Fl_\calP)$ is generated (under colimits) by $\pi_! \DTM(\Fl)$, where $\pi: \Fl \r \Fl_\calP$ is the quotient map.
This map is proper, so the projection formula and the fact that *-pushforwards along schematic maps (as well as any !-pullback) preserve colimits
reduce the claim for $\calP$ to the one for $\calI$.
In this case, again by loc.~cit., the category $\DTM(\Fl)$ is the smallest presentable subcategory containing the subcategories $\tau_* \DTM(S)$, where $\tau : S \r \Fl$ ranges over the closed embeddings of the base points of connected components, and stable under $\pi_s^* \pi_{s*}$, where $\pi_s: \Fl \r \Fl_{\calP_s}$ is as in the proof above.
Now, $\boxtimes$ commutes with $\tau_*$ and also with $\pi_s^*$ and $\pi_{s*}$, since this map is smooth and proper.
\xpf

The following corollary will be used in order to show that Beilinson--Drinfeld affine Grassmannians are Whitney--Tate stratified.

\coro
\thlabel{WT.prod}
Let $X$ be any stratified Whitney--Tate (ind-)scheme $X$.
Then the product stratification on $X \x_S \Fl_\calP$ (i.e., strata are products of $X^w$ times $\calP'$-orbits, for an arbitrary fixed parahoric subgroup $\calP'$) is again Whitney--Tate.
\xcoro

\pf
Abbreviate $\Fl := \Fl_\calP$ and write $\iota_X$ and $\iota_{\Fl}$ for the stratification maps.
It suffices to have an isomorphism $(\iota_X \x \iota_{\Fl})_* \Z = \iota_{X*} \Z \boxtimes \iota_{\Fl*} \Z$, since in any case $*$-pullbacks commute with exterior products.
We have
$$(\iota_X \x \iota_{\Fl})_* \Z = (\iota_X \x \id)_* (\id \x \iota_{\Fl})_* \Z.$$
Since $\Fl$ is universally Whitney--Tate, $(\id \x \iota_{\Fl})_* \Z = p^* \iota_{\Fl*} \Z$, where $p : X^\dagger \x {\Fl} \r {\Fl}$ is the projection.
This motive can also be written as $\Z_{X^\dagger} \boxtimes \iota_{\Fl*} \Z$.
Applying $(\iota_X \x \id)_*$ to it gives, by \thref{base.change.Fl}, $\iota_{X*} \Z \boxtimes \iota_{\Fl*} \Z$.
\xpf

\subsection{Semi-infinite orbits over an algebraically closed field}
\label{subsection:semi-infinite orbits} In this section we assume that $S = \Spec k$ is the spectrum of an algebraically closed field and that $G$ is simple and simply connected. 
See \thref{reduction to simply connected case} for how to generalize the argument to general reductive groups.
Let $\mathcal{I} \subset L^+G$ be the Iwahori subgroup. The $\mathcal{I}$-orbits in $\Gr_G$ are parametrized by $X_*(T)$, see \cite[Example 4.2.12]{RicharzScholbach:Intersection}. For $\lambda,\nu \in X_*(T)$, let $\Gr_{G}^{\mathcal{I} \lambda} \subset \Gr_G$ be the corresponding $\mathcal{I}$-orbit, and
consider the semi-infinite orbit \(\Ss_\nu^+\). View the intersection $\Ss_\nu^+ \cap \Gr_G^{\mathcal{I} \lambda}$ as a reduced subscheme of $\Gr_G^{\mathcal{I} \lambda}$. 

Our purpose in this section is to give a short proof that $\Ss_\nu^+ \cap \Gr_G^{\mathcal{I} \lambda}$ has Tate cohomology. This is logically independent of the rest of the paper, since in the case of $L^+G$-orbits $\Gr_G^{\lambda}$, we will prove in \thref{cellularity of intersection:torus} the stronger result that $\Ss_\nu^+ \cap \Gr_G^{\lambda}$ has a filtrable cellular decomposition over $\Z$ for any split $G$. While the latter result is needed in several places, we have included this section to give a flavor of the more complicated proof of \thref{cellularity of intersection:torus}, and because it can be used to give short proof, without delving into the detailed combinatorics in \refsect{Galleries.Review}, that the constant term functors for $B$ preserve Tate motives over an algebraically closed field. 

\prop
\thlabel{CT.Tate}
If $p \colon \Ss_\nu^+ \cap \Gr_G^{\mathcal{I} \lambda} \rightarrow S$ is the structure morphism, $p_!(\Z) \in \DTM(S)$.
\xprop

\pf We can assume that $\Ss_\nu^+ \cap \Gr_G^{\mathcal{I} \lambda} \neq \emptyset$.
Pick a regular cocharacter $\Gm \rightarrow T$, so that the reduced locus of $\Gr_T$ is the set of fixed points for the resulting $\Gm$-action on $\Gr_G$. Because  $G$ is simple and simply connected, by \cite[Theorem 2.5.3]{Zhu:Introduction} we may identify $\Gr_G$ with the flag variety of an affine Kac--Moody group.
Choose a Bott--Samelson resolution $m \colon X \rightarrow \overline{\Gr_{G}^{\mathcal{I} \lambda}}$ such that $m$ is an isomorphism over $\Gr_{G}^{\mathcal{I} \lambda}$. 
Details about the construction of $X$ can be found in \cite[\S 4]{Juteau:ParitySheaves2014}.

Let $f = \restr{m}{m^{-1}(\Ss_\nu^+ \cap \overline{\Gr_{G}^{\mathcal{I} \lambda}})}$. Let $i \colon (\Ss_\nu^+ \cap \overline{\Gr_{G}^{\mathcal{I} \lambda}}) \setminus (\Ss_\nu^+ \cap \Gr_{G}^{\mathcal{I} \lambda}) \rightarrow \Ss_\nu^+ \cap \overline{\Gr_{G}^{\mathcal{I} \lambda}}$ be the closed immersion, and let $\overline{p} \colon \Ss_\nu^+ \cap \overline{\Gr_{G}^{\mathcal{I} \lambda}} \r S$ be the structure map. By applying localization to $f_!\Z$ we have an exact triangle
$$p_! \Z \rightarrow \overline{p}_!  f_! \Z \rightarrow \overline{p}_! i_! i^* f_! \Z.$$
We will prove the middle and right terms lie in $\DTM(S)$.

Note that the variety $X$ is smooth, projective, and it has
a $\Gm$-action such that $m$ is $\Gm$-equivariant. 
As $X$ embeds equivariantly into a product of affine flag varieties, this $\Gm$-action has isolated fixed points. 
Since $k$ is algebraically closed, the attractors for this action of $\Gm$ as in \cite{BB:Theorems} then give a decomposition of $X$ into affine spaces. 
By the existence of $T$-equivariant ample line bundles on affine flag varieties, cf.~\cite[\S 1.5]{Zhu:Introduction}, $X$ also embeds $\Gm$-equivariantly into some projective space.
By \cite[Th.~3]{BB:Properties}, this implies the Bia\l ynicki-Birula decomposition of $X$ is a filtrable decomposition into affine cells. Since $\Ss_\nu^+ \cap \overline{\Gr_{G}^{\mathcal{I} \lambda}}$ is an attractor, the fiber $m^{-1}(\Ss_\nu^+ \cap \overline{\Gr_{G}^{\mathcal{I} \lambda}})$ is a union of attractors.
By repeatedly applying localization and noting that affine spaces have Tate cohomology it follows that $\overline{p}_!  f_! \Z \in \DTM(S).$

The proof of \cite[Prop. 4.11]{Juteau:ParitySheaves2014} shows that $m_!(\Z)$ is isomorphic to a composition of $*$-pullbacks and $!$-pushforwards along stratified maps between Kac--Moody flag varieties which are stratified by (affine) $\mathcal{I}$-orbits.
In particular, $m_!(\Z) \in \DTM(\overline{\Gr_{G}^{\mathcal{I} \lambda}}),$ where $\overline{\Gr_{G}^{\mathcal{I} \lambda}}$ is stratified by $\mathcal{I}$-orbits.
By proper base change, the restriction of $f_!(\Z)$ to each $\Ss_\nu^+ \cap \Gr_G^{\mathcal{I} \lambda'} \subset \overline{\Gr_{G}^{\mathcal{I} \lambda}}$ is Tate. 
Thus $\overline{p}_! i_! i^* f_! \Z \in \DTM(S)$ by induction on $\lambda$ and localization with respect to the stratification by intersections of $\Ss_\nu^+$ with $\mathcal{I}$-orbits. 
\xpf

\subsection{Intersections of Schubert cells and semi-infinite orbits}\label{subsec--GL}

In the following subsection, we prove that the intersections of the Schubert cells and semi-infinite 
orbits admit a filtrable cellular decomposition, following \cite{GL:LSGalleries}. This will later allow us to show the constant term functors preserves Tate motives. 
Over an algebraically closed field, the latter statement can also be shown using \thref{CT.Tate}. 
The proof in this section, while longer and more combinatorial, works over any base. The stronger cellularity result will moreover be used, among other things, to show that the Hopf algebra arising from the Tannakian formalism is flat.
For the rest of this section, we will work over \(S=\Spec \Z\) for simplicity; the general case then follows by base change. 

Before diving into the necessary combinatorics, let us give a brief overview of the picture; certain notions below will be made precise later on.
Since the main goal is to describe intersections of Schubert cells and semi-infinite orbits, we fix a Schubert variety \(\Gr_G^{\leq \mu}\).
By choosing a so-called minimal combinatorial gallery \(\gamma_\mu\) in the standard apartment joining \(0\) with \(\mu\), we obtain a Bott--Samelson resolution \(\Sigma(\gamma_\mu)\to \Gr_G^{\leq \mu}\), given by an iterated fibration with partial flag varieties as fibers.
This resolution is equivariant for the \(\Gm\)-action induced by a regular anti-dominant cocharacter, and hence induces a map on attractors.
On \(\Gr_G^{\leq \mu}\), the attractors are exactly the intersections with the (negative) semi-infinite orbits.
Since the Bott--Samelson resolution restricts to an isomorphism over the open Schubert cell, we are thus led to studying the intersection of the attractors of \(\Sigma(\gamma_\mu)\) with the preimage of \(\Gr_G^\mu\).
The advantage of working with \(\Sigma(\gamma_\mu)\) is that it is smooth projective; in particular its attractor locus is smooth (which fails for \(\Gr_G^{\leq \mu}\)) and induces a decomposition of \(\Sigma(\gamma_\mu)\) (which fails for \(\Gr_G^\mu\) by \cite[Lemma 1.4.9]{DrinfeldGaitsgory:Theorem}). 

It turns out that the fixed points of the \(\Gm\)-action on \(\Sigma(\gamma_\mu)\) can be indexed by certain combinatorial galleries in the standard apartment (which is in particular independent of the base scheme), \thref{decomposition of bott samelson}, whereas over algebraically closed fields, the points of \(\Sigma(\gamma_\mu)\) correspond to more general galleries (of a fixed type) in the Bruhat--Tits building \cite[Definition-Proposition 1]{GL:LSGalleries}.
Using this interpretation, the attractors correspond to the fibers of a certain retraction map from the Bruhat--Tits building onto the standard apartment, and the preimage of \(\Gr_G^\mu\) corresponds to the minimal galleries \cite[Lemma 10]{GL:LSGalleries}.
The fibers of this retraction map can then be identified with the products of affine root groups \cite[Lemma 13]{GL:LSGalleries}, and at least for regular \(\mu\), the open subset of minimal galleries corresponds to a product of (possibly punctured) affine root groups.
For general \(\mu\), the situation is more complicated, and the varieties in question will only be decomposed into cells.
Given that all this is combinatorial in nature, it extends from algebraically closed fields to \(\Spec \Z\), cf.~Propositions \ref{regular case} and \ref{deodhar in the infinite setting}.
It also turns out that the resulting decompositions are filtrable.
Finally, the fact that the intersection of a Schubert cell and a semi-infinite orbits can be decomposed into intersections of the Schubert cell with attractors for \(\Sigma(\gamma_\mu)\) plays a crucial role in \thref{Torsor is trivial over decomposition}, which will allow us to show the motivic Hopf algebra governing the Satake category is a reduced motive in \thref{H.Independent}.

\nota
Consider \(T\subseteq B\subseteq G\) as before. Recall that if \(\ev:L^+G\to G\) is the evaluation of \(t\) at \(0\), then 
the Iwahori \(\mathcal{I}\subseteq L^+G\subseteq LG\) is the inverse image 
\(\ev^{-1}(B)\). Let \(N:=N_G(T)\) be the normalizer of \(T\) in \(G\), which 
gives us the finite Weyl group \(W=N/T\)
of \(G\). We also have the affine Weyl group \(W^{\fraka}\), which agrees with 
the extended affine Weyl group \(N(\Z\rpot{t})/T(\Z\pot{t})\) when \(G\) is simply connected.
We denote by \(\Phi\) the roots of \(G\), and \(\Phi^+\subset \Phi\) the positive roots with respect to \(B\). 
Similarly, we denote by \(R^+\subseteq R\) the (positive) affine roots.
Finally, we fix a Chevalley system of $G$. In particular, for any root \(\beta\in \Phi\), we get an isomorphism 
\(x_\beta\colon \Ga\r U_\beta\), where \(U_\beta\) is the root group associated to \(\beta\).
\xnota

\subsubsection{Review of combinatorial galleries} \label{sect--Galleries.Review}
For the rest of this section, we will assume \(G\) is semisimple and simply connected,
unless mentioned otherwise (which will be the case for the main theorems in this section). Under this assumption, a weaker version of what we want was already proved by Gaussent--Littelmann
in \cite{GL:LSGalleries} for complex groups. However, their arguments work almost verbatim over any algebraically closed field, and we will also use this in the sequel.
Our aim will be to generalize (and strengthen) their work to groups over \(\Spec \Z\). 
We first recall the most important notation and terminology from \cite{GL:LSGalleries}, and refer to loc.~cit.~ for details. 
By the following remark, we can work over an arbitrary algebraically closed field \(k\) for now.

\rema
\thlabel{building.independence.k.Z}
For any field \(k\), the natural maps \(N(\Z\rpot{t})/T(\Z\rpot{t}) \xrightarrow{\cong} N(k\rpot{t})/T(k\rpot{t})\)
and \(T(\Z\rpot{t})/T(\Z\pot{t})\xrightarrow{\cong} T(k\rpot{t})/T(k\pot{t})\) are isomorphisms.  
In particular, 
there is a natural isomorphism \(N(\Z\rpot{t})/T(\Z\pot{t})\xrightarrow{\cong} N(k\rpot{t})/T(k\pot{t})\), i.e., an isomorphism between the 
(extended) affine Weyl groups of \(G\) and \(G_k\). Using this, we get an isomorphism of apartments
\(\Aa(G,T)\cong \Aa(G_k,T_k)\), equivariant for the actions of the affine Weyl groups of \(G\) and \(G_k\).
This isomorphism is moreover compatible with the identifications \(\Aa(G,T)\cong X_*(T)\otimes_{\Z}\RR\) and
\(\Aa(G_k,T_k)\cong X_*(T)\otimes_{\Z} \RR\), if we choose the canonical Chevalley valuations on \(G\otimes_{\Z} \Z\pot{t}\) 
and \(G_k \otimes_k k\pot{t}\) as basepoints of the apartments.
\xrema

Consider the group \(X_*(T)\) of cocharacters of \(G\), and let \(\Aa_G:=X_*(T)\otimes_{\Z} \RR\).
The affine Weyl group 
\(W^{\fraka}\) acts on \(\Aa_G\) by affine reflections. The reflection hyperplanes (also called \emph{walls}) in 
\(\Aa_G\) for this action are all of the form \(\mathbf{H}_{\beta,m}=\{a\in \Aa_G\mid \langle 
a,\beta\rangle=m\}\), for some positive root \(\beta\in \Phi^+\) and \(m\in \Z\). Let \(s_{\beta,m}\in W^{\fraka}\) 
denote the corresponding affine reflection, and \(\mathbf{H}^+_{\beta,m}=\{a\in \Aa_G\mid \langle 
a,\beta\rangle\geq m\}\) and \(\mathbf{H}^-_{\beta,m}=\{a\in \Aa_G\mid \langle a,\beta\rangle\leq m\}\) the associated 
closed half-spaces. 

\rema
\thlabel{Remark:Tits convention}
We use the Tits convention for the action of \(W^{\fraka}\) on \(\Aa_G\), i.e., we let \(\lambda\in T(k\rpot{t})/T(k\pot{t})=X_*(T)\) 
act via the translation by \(-\lambda\) \cite[§1.1]{Tits:Corvallis}. 
The advantage of this convention is that if we let \(W^{\fraka}\) act on the set 
\(R\) of affine roots as usual, via \(w(\alpha)(x):=\alpha(w^{-1}x)\) for \(\alpha\in R\), 
\(w\in W^{\fraka}\) and \(x\in \Aa_G\), then \(wU_{\alpha} w^{-1}=U_{w(\alpha)}\).
\xrema

Let \(\mathbf{H}^{\fraka}:=\bigcup_{\beta\in \Phi^+,m\in \Z} \mathbf{H}_{\beta,m}\)
denote the union of the reflection 
hyperplanes. Then the connected components of \(\Aa_G\setminus \mathbf{H}^{\fraka}\) are called \emph{open alcoves}, and 
their closures simply \emph{alcoves}. More generally, a \emph{face}
of \(\Aa_G\) is a subset \(F\) that can be 
obtained by intersecting closed affine half-spaces and reflection hyperplanes, one for each \(\beta\in \Phi^+\) and 
\(m\in \Z\). One example of an alcove is the \emph{fundamental alcove} \(\Delta_f=\{a\in \Aa_G\mid 0\leq 
\langle a, \beta\rangle \leq 1, \forall \beta\in \Phi^+\}\).

It is well-known that \(W^{\fraka}\) is generated by the affine reflections \(S^{\fraka}\),
consisting of those 
\(s_{\beta,m}\in W^{\fraka}\) such that \(\Delta_f\) contains a codimension 1 face lying in \(\mathbf{H}_{\beta,m}\). 
For a face \(F\subseteq \Delta_f\), we define the \emph{type} 
\(S^{\fraka}(F)\) as the subset of \(S^{\fraka}\) consisting of those \(s_{\beta,m}\in S^{\fraka}\) such that \(F\) is contained in the hyperplane 
\(\mathbf{H}_{\beta,m}\). 
In particular, \(S^{\fraka}(0)=S\) consists of those \(s_{\beta,0}\) such that \(H_{\beta,0}\) contains a codimension 1 face of \(\Delta_f\), and \(S^{\fraka}(\Delta_f)=\varnothing\). 
Since \(W^{\fraka}\) acts simply transitively on the set of alcoves, we can 
translate any face of \(\Aa_G\) to a face of \(\Delta_f\) using this action, and use this to define the type of arbitrary faces 
in \(\Aa_G\). Moreover, to any proper subset \(t_\bullet\subset S^{\fraka}\), we can associate a parahoric subgroup \(\Pp_{t_\bullet}:=\bigcup_{w\in 
W^{\fraka}_{t_\bullet}} \mathcal{I} w\mathcal{I}\subseteq LG\), where \(W^{\fraka}_{t_\bullet}\) is the subgroup of \(W^{\fraka}\) generated by \(t_{\bullet}\). We call this the \emph{standard parahoric of type \(t_{\bullet}\).} 
For example, we have \(\Pp_{\varnothing}=\mathcal{I}\), and \(\Pp_S=L^+G\). Conversely, any parahoric subgroup \(\Pp\) 
containing \(\mathcal{I}\) arises uniquely in this way; we denote the associated subgroup of \(W^{\fraka}\) by \(W^{\fraka}_\Pp\).

Now, consider the isomorphisms \(x_\beta:\Ga\r U_\beta\) arising from the Chevalley system.
Let \(v\) denote the canonical valuation on \(k\rpot{t}\), and define, for any \(r\in \RR\), the subgroup
\[U_{\beta,r}:=1\cup \{x_\beta(f)\mid f\in k\rpot{t}, v(f)\ge r\}\subseteq G(k\rpot{t}).\]
Letting \(\ell_\beta(\Omega):=-\inf_{x\in \Omega} \langle x,\beta \rangle \) 
for any \(\varnothing \neq \Omega\subseteq \Aa_G\), we can define the subgroups 
\(U_{\Omega}:=\langle U_{\beta,\ell_\beta(\Omega)}\mid \beta\in \Phi\rangle\) of \(G(k\rpot{t})\), and use this to define the
affine building of \(G\).

\defi
The \emph{affine building} \(\Jj^{\fraka}\) of \(G\) is the quotient \(G(k\rpot{t})\times \Aa_G/\sim\), where two pairs 
\((g,x)\) and \((h,y)\) are equivalent if there is some \(n\in N(k\rpot{t})\) such that \(nx=y\) and \(g^{-1}hn\in U_x\).

There is an \(N(k\rpot{t})\)-equivariant injection \(\Aa_G\hr \Jj^{\fraka}\colon x\mapsto (1,x)\). An \emph{apartment} of 
\(\Jj^{\fraka}\) is a subset of the form \(g\Aa_G\), for some \(g\in G(k\rpot{t})\).
In particular, \(\Aa_G\) is an apartment, called the \emph{standard} apartment. 

Finally, we define the faces in \(\Jj^{\fraka}\) as the \(G(k\rpot{t})\)-translates of the faces in \(\Aa_G\). The type of a face in \(\Jj^{\fraka}\) is defined similarly, by translating to a face in \(\Aa_G\); this is a well-defined notion.
\xdefi
The following definition is a central topic in \cite{GL:LSGalleries}; this and all of the remaining definitions in \refsect{Galleries.Review} are taken from op.~cit.

\defi
\thlabel{gallery}
A \emph{gallery} in \(\Jj^{\fraka}\) is a sequence of faces 
  \[\gamma=(\Gamma'_0\subset \Gamma_0\supset \Gamma'_1\subset \ldots \supset \Gamma'_p\subset \Gamma_p\supset \Gamma'_{p+1})\]
  in \(\Jj^{\fraka}\), such that 
  \begin{itemize}
    \item \(\Gamma'_0\) and \(\Gamma'_{p+1}\),  called the \emph{source} and \emph{target} of \(\gamma\), are vertices,
    \item the \(\Gamma_j\)'s are all faces of the same dimension, and
    \item for \(1\leq j\leq p\), the face \(\Gamma'_j\) is a codimension one face of both \(\Gamma_{j-1}\) and \(\Gamma_j\).
  \end{itemize}
The \emph{gallery of types} of such a gallery \(\gamma\) is the sequence of types of the faces of \(\gamma\):
\[t_\gamma=(t'_0\supset t_0\subset t'_1\supset \ldots t'_p\subset t_p\supset t'_{p+1}),\]
where the \(t'_j\) and \(t_j\) are the types of \(\Gamma'_j\) and \(\Gamma_j\) respectively.
\xdefi

We will be especially interested in galleries of a more combinatorial nature, depending on a fixed minimal gallery. To define these, we need the following extension of cocharacters.

\rema
\thlabel{generalized cocharacters}
If \(G_{\adj}=G/Z_G\) is the adjoint quotient of \(G\) and \(T_\adj\subseteq G_\adj\) the adjoint torus, then we have natural inclusions \(X_*(T)\subseteq X_*(T_\adj)\subseteq X_*(T)\otimes_\Z \R\). The lattice \(X_*(T_\adj)\subseteq X_*(T)\otimes_\Z \R\) consists exactly of those vertices that are a face of \(\Aa_G\). In what follows, we will often consider all elements of \(X_*(T_\adj)\subseteq X_*(T)\otimes_\Z \R\) instead of just the cocharacters; this will help us when applying the results of this section to non-simply connected groups. Note that notions such as dominance and regularity extend to \(X_*(T_\adj)\).
\xrema

For \(\mu\in X_*(T_\adj)\), let 
\(\mathbf{H}_\mu=\bigcap_{\langle\mu,\alpha\rangle =0} \mathbf{H}_{\alpha,0}\) be the intersection of those 
reflection hyperplanes corresponding to the roots orthogonal to \(\mu\) 
(so that \(\mathbf{H}_\mu=\Aa_G\) for regular \(\mu\)).

\defi
\thlabel{Definition joining}
Fix some \(\mu\in X_*(T_\adj)\), let \(F_f\) be the face corresponding to \(0\in \Aa_G\), and \(F_\mu\) the face 
corresponding to \(\mu\in \Aa_G\); note that both are vertices. A gallery \(\gamma=(\Gamma'_0\subset 
\Gamma_0\supset \Gamma'_1\subset \ldots \supset \Gamma'_p\subset \Gamma_p\supset \Gamma'_{p+1})\) contained in \(\Aa_G\) is 
said to \emph{join 0 with \(\mu\)} if its source is \(F_f\), its target \(F_\mu\), and if the dimension of 
the large faces \(\Gamma_j\) is equal to the dimension of \(\mathbf{H}_\mu\).
\xdefi

In fact, we are interested in those galleries that are minimal in a precise sense. While one can define what it 
means for an arbitrary gallery in \(\Jj^{\fraka}\) to be minimal, we will content ourselves to give an equivalent 
definition for galleries joining 0 with some dominant \(\mu\in X_*(T_\adj)\), cf.~ \cite[Lemma 4]{GL:LSGalleries}. We say that a 
reflection hyperplane \(\mathbf{H}\) \emph{separates} a subset \(\Omega\) and a face \(F\) of \(\Aa_G\), if 
\(\Omega\) lies in a closed half-space defined by \(\mathbf{H}\), and \(F\) is contained in the opposite open 
half-space. For two faces \(E,F\) in \(\Aa_G\), let \(\Mm_{\Aa_G}(E,F)\) be the set of such 
hyperplanes separating \(E\) and \(F\). It is known that this set is finite.

\defi
Let \(\gamma_\mu=(F_f\subset \Gamma_0\supset \Gamma'_1\subset \ldots \supset \Gamma'_p\subset 
\Gamma_p\supset F_\mu)\) be a gallery joining 0 with \(\mu\). For each \(0\leq j\leq p\), let \(\calH_j\) 
be the set of reflection hyperplanes \(\mathbf{H}\) such that \(\Gamma'_j\subset \mathbf{H}\) and 
\(\Gamma_j\nsubseteq \mathbf{H}\). We say \(\gamma_\mu\) is \emph{minimal} if all the faces of \(\gamma_\mu\) are contained 
in \(\mathbf{H}_\mu\), and if there is a disjoint union \(\bigsqcup_{0\leq j<p} 
\calH_j=\Mm_{\Aa_G}(F_f,F_\mu)\). 
\xdefi

\rema\thlabel{First parahoric of minimal gallery}
By minimality, \(\Gamma_0\) is contained in \(\mathbf{H}_\mu\), and also in \(\Delta_f\) since \(\mu\) is dominant.
In particular, there is a unique choice for \(\Gamma_0\), namely the facet corresponding to the parahoric contained in \(L^+G\), whose reduction mod \(t\) is the parabolic generated by the root groups of those roots \(\alpha\) for which \(\langle \alpha,\mu\rangle \geq 0\).
\xrema

Let us fix a minimal gallery \(\gamma_\mu\) joining 0 with \(\mu\), with associated gallery of types \(t_{\gamma_\mu} = (S\supset t_0\subset t'_1 \supset \ldots \supset t'_p\subset t_p\supset t_\mu)\). We denote by \(\Gamma(\gamma_\mu)\) the set of 
all galleries of type \(t_{\gamma_\mu}\) and of source \(F_f\) contained in the standard apartment \(\Aa_G\). 
Such galleries are called \emph{combinatorial of type \(t_{\gamma_\mu}\)}.
One can describe \(\Gamma(\gamma_\mu)\) quite 
explicitly; recall that for some type \(t\), we defined \(W_t\) as the subgroup of \(W^{\fraka}\) generated by \(t\). For 
simplicity, we will also write \(W_\mu:=W_{t_\mu}\), and similarly \(W_i:=W_{t_i}\) and \(W'_i:=W_{t'_i}\). Then, by 
\cite[Proposition 2]{GL:LSGalleries}, there is a bijection
\[W\times^{W_0} W_1'\times^{W_1}\ldots \times^{W_{p-1}} W'_p/W_p\to \Gamma(\gamma_\mu),\]
sending an equivalence class \([\delta_0,\delta_1,\ldots, \delta_p]\) to the combinatorial gallery of type \(t_{\gamma_\mu}\) given
by 
\((F_f\subset \Sigma_0\supset \Sigma'_1\subset \ldots \supset \Sigma'_p\subset \Sigma_p\supset F_\nu)\), where 
\(\Sigma_j=\delta_0\delta_1\ldots \delta_j F_{t_j}\).

Most important for us will be the subset of \emph{positively folded combinatorial galleries}. Before we can explain their 
definition, note that we can (and will) assume that for any \([\delta_0,\delta_1,\ldots, \delta_p]\in W\times^{W_0} W_1'\times^{W_1}\ldots \times^{W_{p-1}} W'_p/W_p\), each \(\delta_j\in W'_j\) is the minimal length representative of 
its class in \(W'_j/W_j\). Moreover, the minimal gallery \(\gamma_\mu\) is represented by \([1, \tau_1^{\min}, 
\ldots,\tau_p^{\min}]\), where each \(\tau_j^{\min}\in W'_j\) is the minimal length representative of the longest class in 
\(W'_j/W_j\).
The positively folded galleries are now defined as in \cite[Definition 16]{GL:LSGalleries}.

\defi
Let \(\delta = (F_f\subset \Sigma_0\supset \Sigma'_1\subset \ldots \supset \Sigma'_p\subset \Sigma_p\supset F_\nu)\) be 
a combinatorial gallery in \(\Gamma(\gamma_\mu)\) corresponding to \([\delta_0,\delta_1,\ldots, \delta_p]\), 
where each \(\delta_j\) is a minimal representative of its class in \(W'_j/W_j\).
\begin{enumerate}
  \item If \(j\geq 1\) and \(\delta_j\neq \tau_j^{\min}\), we say \(\delta\) is \emph{folded} around \(\Sigma'_j\).
\end{enumerate}
Now, consider for each \(j\geq 1\), the combinatorial galleries
\[\gamma^{j-1} = [\delta_0,\ldots,\delta_{j-1}, \tau_{j}^{\min},\ldots,\tau_p^{\min}] = (F_f\subset \ldots \subset 
\Sigma_{j-1}\supset \Sigma'_j\subset \Omega_j\supset \Omega'_{j+1}\subset \ldots)\] and
\[\gamma^{j} = [\delta_0,\ldots,\delta_j, \tau_{j+1}^{\min},\ldots,\tau_p^{\min}] = (F_f\subset \ldots \subset 
\Sigma_{j-1}\supset \Sigma'_j\subset \Sigma_j\supset \Sigma'_{j+1}\subset \ldots).\]
Then, by \cite[Lemma 5]{GL:LSGalleries}, there exist positive roots \(\beta_1,\ldots,\beta_q\) and integers 
\(m_1,\ldots,m_q\) such that the small face \(\Sigma'_j\) is contained in \(\bigcap_{i=1}^q \mathbf{H}_{\beta_i,m_i}\), 
and where \(\Sigma_j=s_{\beta_q,m_q}\ldots s_{\beta_1,m_1}(\Omega_j)\).
\begin{enumerate}
  \setcounter{enumi}{1}
  \item If \(\delta\) is folded around \(\Sigma'_j\), we say this folding is \emph{positive} if \(\Sigma_j\subset 
  \bigcap_{i=1}^q \mathbf{H}_{\beta_i,m_i}^+\).
  \item The combinatorial gallery \(\delta\) is \emph{positively folded}, if all of its folds are positive.
\end{enumerate}
  We denote the subset of \(\Gamma(\gamma_\mu)\) consisting of the positively folded combinatorial 
  galleries by \(\Gamma^+(\gamma_\mu)\).
\xdefi

As an example, \(\gamma_\mu\) does not have any folds, so that it is automatically positively folded. 
We will also need the following definition, which is equivalent to the definition given in \cite[Top of p.~58]{GL:LSGalleries}.

\defi
Let \(\delta=(F_f\subset \Sigma_0\supset \Sigma'_1\subset \ldots \supset \Sigma'_p\subset \Sigma_p\supset F_\nu)\in 
\Gamma^+(\gamma_\mu)\) be a positively folded combinatorial gallery. A \emph{load-bearing wall} for \(\delta\) at 
\(\Sigma_j\) is a reflection hyperplane \(\mathbf{H}=\mathbf{H}_{\beta,m}\) with \(\beta\in \Phi^+\) and \(m\in \Z\), such that \(\Sigma'_j\subset \mathbf{H}\) and 
\(\Sigma_j\nsubseteq \mathbf{H}\), and such that \(\Sigma_j\subseteq \mathbf{H}^+_{\beta,m}\).
\xdefi

As \(\delta\) was assumed positively folded, it follows from the definitions that any folding hyperplane is a load-bearing wall.

Finally, to each \(\delta=[\delta_0,\ldots,\delta_p]\in \Gamma^+(\gamma_\mu)\), we will need to attach two sets of 
indices. 
For any affine root \(\alpha\) of \(G\), we denote by \(\Uu_\alpha\) the corresponding 
root subgroup of \(LG\).

\defi
\thlabel{notation:J-inf}
\begin{enumerate}
  \item For any parahorics \(\Qq\subseteq \Pp\) and \(w\in W_\Pp/W_\Qq\), we define the subsets of affine roots \(R^+(w):=\{\alpha>0\mid \Uu_{w^{-1}(\alpha)}\nsubseteq \Qq\}\) and \(R^-(w):=\{\alpha<0\mid w(\alpha)<0, \Uu_\alpha\subseteq \Pp,\Uu_\alpha\nsubseteq \Qq\}\).
  \item Using the same notation as above, we define \(\Uu^+(w):=\prod_{\eta\in R^+(w)}\Uu_\eta\) and \(\Uu^-(w):=\prod_{\theta\in R^-(w)}\Uu_\theta\).
  Note that there is a natural locally closed immersion \(\Uu^+(w)w\Uu^-(w) \subseteq \Pp/\Qq\).
  \item Let \(\delta=[\delta_0,\ldots,\delta_p]\in \Gamma^+(\gamma_\mu)\) be a positively folded combinatorial gallery of type \(t_{\gamma_\mu}\). For any \(j\), let \(\Pp_j\) and \(\Qq_j\) be the parahoric subgroups containing \(\mathcal{I}\) of type \(t'_j\) and \(t_j\) respectively. Consider the set of walls in \(\Aa_G\) that contain \(\Ff_{\Pp_j}\), but not \(\delta_jF_{\Qq_j}\). If we index this set by \(I_j\), then \(I_j\) can be decomposed in as \(I_j=I^+_j\sqcup I^-_j\), such that 
  \(R^+(\delta_j)=\{\alpha_i\mid i\in I^+_j\}\), where \(\alpha_i\) is the positive root corresponding to the wall 
  \(\mathbf{H}_i\), and there is a similar description for \(R^-(\delta_j)\); cf.~ \cite[§10]{GL:LSGalleries}. Then, we define  \(J_{-\infty}(\delta)\subseteq \bigsqcup_{j=0}^p I_j\) as the subset corresponding to those walls that are 
  load-bearing. We also define \(J_{-\infty}^{\pm}(\delta):=J_{-\infty}(\delta) \cap 
  (\bigsqcup_{i=0}^p I_j^{\pm})\).
\end{enumerate}
\xdefi

\subsubsection{Cellular stratifications of Bott--Samelson schemes}

In this subsection, we will apply the methods from \cite{GL:LSGalleries}, and explain how to generalize them to 
more general bases. By \thref{building.independence.k.Z}, we can use notions that depend only on the standard apartment (instead of the whole affine building) over any base, and independently of the base. 
Examples of this are the standard apartment \(\Aa_G=\Aa(G,T)\), and combinatorial galleries of a fixed type, possibly minimal or positively folded.

\nota\thlabel{notation generalized schubert varieties}
Fix some dominant \(\mu\in X_*(T_\adj)\), a minimal gallery \(\gamma_\mu\) in \(\Aa_G\) joining \(0\) with \(\mu\),
and let 
\[t_{\gamma_\mu}=(t'_0\supset t_0\subset t'_1 \supset \ldots \supset t_{j-1} \subset t'_j \supset t_j \subset \ldots t'_p\supset t_p \subset t_\mu)\] 
be its gallery of types. For \(0\le j\le p\), let \(\Pp_j\subseteq LG\) (resp.~ \(\Qq_j\), resp.~ \(\Pp_\mu\)) 
be the parahoric subgroup of type \(t'_j\) (resp.~ \(t_j\), resp.~ \(t_\mu\)) containing \(\mathcal{I}\). 
Note that if \(\mu\) is an actual cocharacter of \(T\)
(rather than \(T_{\adj}\), cf.~\thref{generalized cocharacters}), then \(\Pp_\mu=L^+G\). In general, \(\Pp_\mu\) is exactly the parahoric subgroup such 
that the sheaf quotient \(LG/\Pp_\mu\) from \thref{Sheafifications of quotient agree} is \(LG\)-equivariantly universally homeomorphic to the connected component of \(\Gr_{G_{\adj}}\) corresponding 
to the image of \(\mu\) under \(X_*(T_\adj)\to \pi_1(G_\adj)\cong \pi_0(\Gr_{G_\adj})\). This will allow us to reduce the combinatorics needed to the case of simply connected groups, cf.~ \thref{reduction to simply connected case}.
For simplicity, we will denote the \(L^+G\)-orbit in \(LG/\Pp_\mu\) corresponding to \(\mu\) by \(\Gr_{G}^\mu\), 
and its closure by \(\Gr_{G}^{\le \mu}\), but we emphasize that these are only subschemes of \(\Gr_G\) when \(\mu\) is 
an actual cocharacter of \(T\), as opposed to a cocharacter of \(T_{\adj}\).
\xnota

\defi\thlabel{Defi--BSresolution}
The \emph{Bott--Samelson} scheme \(\Sigma(\gamma_\mu)\) is the contracted product
\[\Pp_0\times^{\Qq_0} \Pp_1 \times^{\Qq_1} \ldots \times^{\Qq_{p-1}} \Pp_p/\Qq_p.\]
\xdefi

\prop
\thlabel{bott samelson are desingularizations}
The multiplication morphism (defined since \(\Qq_p\subseteq \Pp_\mu\))
$$\psi:\Sigma(\gamma_\mu)\to \Gr_G^{\le \mu}$$
is an isomorphism over the open subscheme \(\Gr_G^\mu\).
\xprop
\pf

Denote the restricted morphism by \(\phi\colon\Sigma^\circ(\gamma_\mu)\to \Gr_{G}^\mu\), i.e., 
\(\Sigma^\circ(\gamma_\mu)=\psi^{-1}(\Gr_{G}^\mu)\). By \cite[Lemma 10]{GL:LSGalleries}, $\phi$ is bijective over algebraically closed fields.
Note that \(\Sigma(\gamma_\mu)\) is an iterated sequence of 
Zariski-locally trivial fibrations with smooth connected fibers, and in particular smooth and integral itself. 
As a Schubert cell, the target of $\phi$ is also smooth and integral.
Bott--Samelson resolutions are always birational over a field, and in particular over $\Q$, so that $\phi$ is birational. 
Thus $\phi$ is an isomorphism by Zariski's main theorem \StP{05K0}.
\xpf

Since these Bott-Samelson schemes are smooth projective, they are much better behaved with respect to \(\Gm\)-actions. So, we will first study the decomposition on \(\Sigma(\gamma_\mu)\) induced by a certain action, and later restrict this decomposition to \(\Gr_G^\mu\). 
Recall the notion of filtrable decompositions from \thref{Defi filtrable}.

\defilemm
\thlabel{decomposition of bott samelson}
Consider some regular anti-dominant cocharacter \(\lambda\in X_*(T)\), and the induced \(\Gm\)-action on \(\Sigma(\gamma_\mu)\) via left-multiplication on the first factor. Then the connected components of the attractor locus are indexed by \(\Gamma(\gamma_\mu)\), and we denote them by \(C_\delta\). These \(C_\delta\) form a filtrable decomposition of \(\Sigma(\gamma_\mu)\).
\xdefilemm
\pf
First, note that \(\Sigma(\gamma_\mu)\) admits a \(\Gm\)-equivariant embedding into some projective space with linear \(\Gm\)-action: this holds for any partial flag variety, and \(\Sigma(\gamma_\mu)\) is an iterated Zariski-locally trivial fibration of such.
Then the fixed-point locus for this \(\Gm\)-action is a closed subscheme of \(\Sigma(\gamma_\mu)\) \cite[Theorem 1.8 (i)]{Richarz:Spaces}.
Since on geometric fibers (over \(\Spec \Z\)), this fixed-point locus consists of the points \(\Gamma(\gamma_\mu)\subseteq \Sigma(\gamma_\mu)\) by \cite[Proposition 6]{GL:LSGalleries}, \cite[Corollary 1.16]{Richarz:Spaces} tells us that is already true over \(\Spec \Z\).
The claim about connected components of the attractor locus \(\Sigma(\gamma_\mu)^+\) (which is representable and smooth by \cite[Theorem 1.8 (iii)]{Richarz:Spaces}) then follows from \cite[Proposition 1.17 (ii)]{Richarz:Spaces}; we denote the resulting decomposition into connected components by \(\Sigma(\mu)^+=\coprod_{\delta\in \Gamma(\gamma_\mu)} C_\delta\).

Since the \(\Gm\)-action on \(\Sigma(\gamma_\mu)\) is Zariski-locally linearizable by the above paragraph, the same argument as \thref{GrP.LC} shows that the map \(\Sigma(\gamma_\mu)^+\to \Sigma(\gamma_\mu)\) is a locally closed immersion on each connected component.
To see that it is also a bijection, we may restrict to the geometric fibers by \cite[Corollary 1.16]{Richarz:Spaces}.
In that case, it follows from \cite[Proposition 6]{GL:LSGalleries}, since the Bia\l ynicki-Birula decomposition from \cite[Theorem 4.4]{BB:Theorems} agrees with the attractor locus when working over an algebraically closed base field.

It remains to show that the decomposition \(\Sigma(\gamma_\mu)=\bigsqcup_{\delta\in \Gamma(\gamma_\mu)} C_\delta\) is filtrable. 
But this follows from the proof of \cite[Theorem 3]{BB:Properties} (additionally using \cite[Lemma 1.4.9 (ii)]{DrinfeldGaitsgory:Theorem}), since \(\Sigma(\gamma_\mu)\) is smooth projective and admits a \(\Gm\)-equivariant embedding into some projective space with linear \(\Gm\)-action.
\xpf

By \cite[Proposition 6]{GL:LSGalleries}, the geometric fibers of \(C_\delta\) are affine spaces. In fact, the \(C_\delta\) are already affine spaces over \(\Z\), but we omit details as we will not need this.

By \thref{bott samelson are desingularizations} and \thref{decomposition of bott samelson}, we get an induced filtrable decomposition of \(\Gr_{G}^\mu=\bigsqcup_{\delta \in \Gamma(\gamma_\mu)} X_\delta\), with 
\(X_\delta:=C_\delta\cap \Gr_{G}^\mu\). This is related to the semi-infinite orbits as in the proposition below. Recall that \(w_0\in W\) denotes the longest element in the finite Weyl group of \(G\).

\rema\thlabel{mistake in GL}
In \cite[Remark 14 ff.~and Theorem 3]{GL:LSGalleries}, it is claimed that \(\psi\colon \Sigma(\gamma_\mu)\to \Gr_G^{\leq \mu}\) sends the combinatorial gallery \(\gamma_\mu\) to \(t^\mu\in \Gr_G^{\leq \mu}\).
This is incorrect, as can already be seen in the case \(G=\SL_2\). 
Indeed, let \(T\subset B\subset \SL_2\) denote the torus of diagonal matrices and Borel of upper-triangular matrices, and let \(\mu\) be the unique nonzero dominant quasiminuscule cocharacter.
Then there is a unique choice of \(\gamma_\mu\), given by \([1,\tau]\), where a representative of \(\tau\) is given by
\[\begin{pmatrix}
  0 & t^{-1} \\ -t & 0
\end{pmatrix}\in LG.\]
Let \(\Pp\supset \mathcal{I}\) be the unique standard hyperspecial parahoric of \(G\) different from \(L^+G\). The Bott-Samelson resolution
\[\Sigma(\gamma_\mu) = L^+G \overset{\mathcal{I}}{\times} \Pp /\mathcal{I} \to LG/L^+G\]
sends \(\gamma_\mu\) to 
\[\begin{pmatrix}
  0 & t^{-1} \\ -t & 0
\end{pmatrix} \cdot L^+G = \begin{pmatrix}
0 & t^{-1} \\ -t & 0
\end{pmatrix} \begin{pmatrix}
0 & -1 \\ 1 & 0
\end{pmatrix} \cdot L^+G = \begin{pmatrix}
t^{-1} & 0 \\ 0 & t
\end{pmatrix} \cdot L^+G,\]
i.e., to \(t^{w_0(\mu)}\).
More generally, for any group \(G\), the resolution \(\psi\) sends a combinatorial gallery with target \(F_{\nu}\) to \(t^{w_0(\nu)}\), for \(\nu\in X_*(T_{\adj})\).
\xrema

This mistake does not affect the results of \cite{GL:LSGalleries}, as long as one adds a \(w_0\) in certain places (namely, when connecting the building-theoretic results to the geometry of affine Grassmannians).
This is for example the case in \cite[Theorem 3]{GL:LSGalleries}, as a combination of \cite[Example 8, Theorem 3, and Corollary 4]{GL:LSGalleries} yields the wrong dimensions for the intersections \(\Gr_G^\mu\cap \Ss_{\nu}^-\) (compared to e.g.~\cite[Theorem 3.2]{MirkovicVilonen:Geometric}).
Instead, the corrected statement is as follows:

\prop
\thlabel{first stratification of intersection}
For some \(\nu\in X_*(T_\adj)\), let \(\Gamma(\gamma_\mu,\nu)\) denote those 
galleries in \(\Gamma(\gamma_\mu)\) with target \(F_\nu\). 
Then we have a filtrable decomposition \(\Ss_{w_0(\nu)}^- \cap \Gr_{G}^\mu=\bigsqcup_{\delta\in \Gamma(\gamma_\mu,\nu)} X_\delta\).
\xprop
\pf
The Bott-Samelson resolution \(\psi\colon \Sigma(\gamma_\mu)\to \Gr_G^{\leq \mu}\) is clearly \(L^+G\)-equivariant.
Recall that the \(C_{\delta}\subseteq \Sigma(\gamma_\mu)\) are the attractors for a \(\Gm\)-action induced by a regular anti-dominant cocharacter, and similarly for \(\Ss_{\nu}^-\cap \Gr_G^{\leq \mu}\subseteq \Gr_G^{\leq \mu}\).
As both the source and target of \(\psi\) are proper, \(\psi\) preserves attractors.
It remains to observe that for \(\delta\in \Gamma(\gamma_\mu,\nu)\), the attractor \(C_{\delta}\) gets mapped into \(\Ss_{w_0(\nu)}^-\), cf.~\thref{mistake in GL}, and to use \thref{bott samelson are desingularizations}.
\xpf

\lemm
\thlabel{can assume positively folded}
If \(X_\delta\neq \varnothing\), then \(\delta\in \Gamma^+(\gamma_\mu)\), i.e., \(\delta\) is positively folded.
\xlemm
\pf
Again, this follows from \cite[Lemma 11]{GL:LSGalleries} applied to the geometric points of \(\Spec \Z\).
\xpf

The following lemma and corollary will not be necessary for the main theorem of this section, but will be used later on, such as in the proofs of Propositions \ref{DTM.Hck.X.adjunction} and \ref{Y.Strat}.

\lemm\thlabel{Torsor is trivial over decomposition}
Let \(n\gg 0\) be such that the \(L^+G\)-action on \(\Gr_G^\mu\) factors through \(L^nG\), and let \(\Pp_{w_0(\mu)}^n\subset L^nG\) be the stabilizer of \(t^{w_0(\mu)}\).
Then the \(\Pp_{w_0(\mu)}^n\)-torsor \(L^nG\to \Gr_G^\mu\) is trivial over \(X_\delta\), for any \(\delta\in \Gamma^+(\gamma_\mu)\).
\xlemm
\pf
Consider the following diagram with cartesian square, where the vertical maps are given by the reduction mod \(t\), using the isomorphism \(\Gr_G^\mu\cong L^nG/\Pp_{w_0(\mu)}^n\):
\[\begin{tikzcd}
  L^nG \arrow[r, "g"] & (L^nG/L^{>0}G\cap \Pp_{w_0(\mu)}^n)_{\et} \arrow[r, "f"] \arrow[d] & \Gr_G^\mu \arrow[d]\\
  & G\arrow[r, "q'"] & G/\Pp_{w_0(\mu)}^0.
\end{tikzcd}\]
Let \(U^-\) be the unipotent radical of the opposite Borel of \(B\).
We claim that the lower horizontal map, which is a \(\Pp_{w_0(\mu)}^0\)-torsor, is trivial over \(U^-w\Pp_{w_0(\mu)}^0/\Pp_{w_0(\mu)}^0\), for any element \(w\) in the (finite) Weyl group of \(G\); note that these are affine spaces, and are exactly the attractors for the \(\Gm\)-action on \(G/\Pp_{w_0(\mu)}^0\) induced by a regular anti-dominant cocharacter.
Indeed, there is a vector subgroup $U_{w,w_0(\mu)} \subset U$ which maps isomorphically onto its image in the Bruhat cell $U(w_0 w) \calP_{w_0(\mu)}^0 / \calP_{w_0(\mu)}^0$ \cite[\S 13.8]{Jantzen:Representations}. Thus, there is a section of $q'$ over $U^-w \calP_{w_0(\mu)}^0 / \calP_{w_0(\mu)}^0 = w_0 U (w_0 w) \calP_{w_0(\mu)}^0 / \calP_{w_0(\mu)}^0$ with image $w_0 U_{w,w_0(\mu)} \subset G$. See also \cite[Lemme 6.2]{NgoPolo:Resolutions} in the minuscule case.

Since the diagram is cartesian, the left vertical map is an affine bundle, so that \((L^nG/L^{>0}G\cap \Pp_{w_0(\mu)}^n)_{\et}\) is an affine scheme.
Thus \(g\), which is a torsor under a split unipotent group by \cite[Remark 4.2.8]{RicharzScholbach:Intersection}, is a trivial torsor by \cite[Proposition A.6]{RicharzScholbach:Intersection}.
On the other hand, the previous paragraph implies that the \(\Pp_{w_0(\mu)}^0\)-torsor \(f\) is trivial over the preimage of any \(U^-w\Pp_{w_0(\mu)}^0/\Pp_{w_0(\mu)}^0\).
So it remains to show any \(X_\delta\) is contained in such a preimage.

For this, consider the Bott-Samelson scheme \(\Sigma(\gamma_\mu)\) from \thref{Defi--BSresolution}, and its projection \(\Sigma(\gamma_\mu)\to \Pp_0/\Qq_0\) onto the first factor.
This projection is clearly \(L^+G\)-equivariant, and its target can be identified with \(G/\Pp_{w_0(\mu)}^0\).
Indeed, it is well-known that \(\Pp_{w_0(\mu)}^0\) is the parabolic generated by the root groups of those roots \(\alpha\) such that \(\langle \alpha,\mu\rangle \geq 0\), cf.~\cite[Corollary 2.1.11 ff.]{Zhu:Introduction}, so that this follows from \thref{First parahoric of minimal gallery}.
Moreover, the base point of \(\Sigma(\gamma_\mu)\) gets sent to the base point of \(G/\Pp_{w_0(\mu)}^0\) by \thref{mistake in GL},
so that \(\Sigma(\gamma_\mu)\to \Pp_0/\Qq_0\) restricts to the usual projection \(\Gr_G^\mu\to G/\Pp_{w_0(\mu)}^0\).
Since the source and target of \(\Sigma(\gamma_\mu)\to \Pp_0/\Qq_0\) are proper, this map preserves attractors for the \(\Gm\)-action induced by a regular anti-dominant cocharacter.
We conclude by recalling that the attractors of \(\Sigma(\gamma_\mu)\) are exactly the \(C_\delta\) for \(\delta\in \Gamma(\gamma_\mu)\), and that \(X_\delta = C_\delta\cap \Gr_G^\mu\).
\xpf

\coro \thlabel{Av.Fiber} In the notation of \thref{Torsor is trivial over decomposition}, let 
$a \colon L^nG \times X_\delta \r \Gr_G^{\mu}$ be the action map. Then $$a^{-1}(t^{w_0(\mu)}) \cong \calP_{w_0(\mu)}^n \times X_\delta.$$
\xcoro

\pf
The projection $a^{-1}(t^{w_0(\mu)}) \r X_\delta$ is a $\calP_{w_0(\mu)}^n$-torsor. Let $r \colon X_\delta \r L^nG$ be a section, as per \thref{Torsor is trivial over decomposition}. Then
there is an isomorphism $\calP_{w_0(\mu)}^n \times X_\delta \r a^{-1}(t^{w_0(\mu)}) \colon (p, x) \mapsto (p \cdot r(x)^{-1}, x)$.
\xpf

In the proof of our main result, \thref{cellularity of intersection:torus}, the existence of a cellular decomposition will be reduced to considering galleries with only three faces. The following proposition is a step towards this case. We consider a \emph{triple gallery of 
types} \((t_{j-1}\subset t'_j\supset t_j)\) for some \(0< j\leq p\), and the parahoric subgroups \(\Pp\) and \(\Qq\) 
containing \(\mathcal{I}\) corresponding to \(t'_j\) and \(t_j\) respectively. Moreover, we let \(w\) be (the shortest 
representative of) some element in \(W_\Pp/W_\Qq\), and \(\tau^{\min}\) the shortest representative of the longest 
class in \(W_\Pp/W_\Qq\). Recall also the groups \(\Uu^+(w)\) and \(\Uu^-(w)\) from \thref{notation:J-inf}.

\prop
\thlabel{deodhar in the infinite setting}
The intersection \(\Uu^+(w)w\Uu^-(w)\cap \Uu^+(\tau^{\min})\tau^{\min} \subseteq \Pp/\Qq\) admits a filtrable decomposition into products of \(\A^1\)'s and \(\Gm\)'s.
\xprop
\pf
As in \cite[Proposition 9]{GL:LSGalleries}, we are reduced to showing that for some \(v\in W^{\fraka}_\Pp\), the 
intersection \((\mathcal{I} v^{-1}\Qq/\Qq)\cap( \mathcal{I}^-\Qq/\Qq)\subseteq \Pp/\Qq\) admits a filtrable decomposition as in the statement of 
the proposition. 
Note that \(\Pp/\Qq\) is isomorphic to a partial affine flag variety of the maximal reductive quotient \(G'\) of \(\Pp\) (which is automatically split).

Hence it suffices to show that for any Borel and standard parabolic \(B'\subseteq P'\subseteq G'\) and any elements \(v',w'\in W'\) in the finite Weyl group of \(G'\), the intersection \(B'v'\cdot P'\cap B'^-w'\cdot P'\subseteq G'/P'\) admits the desired filtrable decomposition.
If \(v',w'\) are the minimal length representatives of their classes in \(W'/W'_{P'}\), then \(G'/B'\to G'/P'\) maps \(B'v'\cdot B'\cap B'^-w'\cdot B'\) isomorphically to \(B'v'\cdot P'\cap B'^-w'\cdot P'\) \cite[§5]{Rietsch:Closure}, hence we may assume \(B'=P'\).
In that case we can consider the Deodhar decomposition of \cite[Theorem~1.1, Corollary~1.2]{Deodhar:GeometricDecomposition}; although loc.~cit.~is stated over an algebraically closed field, the cells appearing arise from (punctured) root groups, and hence the proof works over  the base ring \(\Z\).
It remains to see the Deodhar decomposition is filtrable, which was shown in \cite[Lemma 2.5]{Dudas:DLrestriction} (again the proof is purely combinatorial, and hence works over \(\Z\)).
\xpf

We now have all the ingredients to deduce 
the main result of this section.
We begin with the simpler case where \(\mu\) is regular. Then, for any \(\delta=[\delta_0,\ldots,\delta_p]\in \Gamma(\gamma_\mu)\) and any \(1\leq i\leq p\), \(\delta_i\) is either trivial, or a simple reflection \cite[p.80]{GL:LSGalleries}. There is a unique possibility for this simple reflection, and we denote the corresponding affine root by \(\alpha_i\).

\prop
\thlabel{regular case}
Assume \(\mu\) is regular, and let \(\delta\in \Gamma^+(\gamma_\mu)\). Then \(X_\delta \cong \A^k\times \Gm^l\) for some \(k,l\ge 0\).
\xprop
\pf
Since \(\Sigma(\gamma_\mu)\) was defined as an iterated Zariski-locally trivial \(\Pp_i/\Qq_i\)-fibration for some parahoric subgroups \(\Qq_i\subset \Pp_i\), the Bruhat stratification for each such quotient shows that \(\Sigma(\gamma_\mu)\) is stratified by affine spaces, indexed by \(W\times^{W_0} W_1'\times^{W_1}\ldots \times^{W_{p-1}} W'_p/W_p\), where the locally closed strata can be described via root groups. (We note that, although this index set is naturally in bijection with \(\Gamma(\gamma_\mu)\), the resulting decomposition does \emph{not} agree with the one from \thref{decomposition of bott samelson}.)

Then, \cite[Proposition 10]{GL:LSGalleries}, tells us that over any geometric point of \(\Spec \Z\) and any \(\delta=[\delta_0,\delta_1,\ldots,\delta_p]\in \Gamma(\gamma_\mu)\) (with each \(\delta_i\) a minimal representative of its class in \(W'_i/W_i\)), \(C_\delta\) is exactly the locally closed subscheme of \(\Sigma(\gamma_\mu)\) given by 
\[\delta_0\cdot \prod_{\beta<0,\delta_0(\beta)<0} \Uu_\beta \cdot \prod_{i=1}^p \Uu^\bullet_{\alpha_i}\cdot \delta_i\]
where \(\alpha_i\) is as above, \(\Uu^\bullet_{\alpha_i}\) is defined as \(\Uu_{\alpha_i}\) 
(resp.~ \(\Uu^\times_{-\alpha_i}\), resp.~ \(\{0\}\)) when \(i\in J^+_{-\infty}(\delta)\) (resp.~ 
\(i\in J^-_{-\infty}(\delta)\), resp.~ \(i\notin J_{-\infty}(\delta)\)), and 
\(J_{-\infty}(\delta)=J^+_{-\infty}(\delta)\sqcup J^-_{-\infty}(\delta)\) is as in 
\thref{notation:J-inf}. In particular, the same description in terms of (punctured) root groups can be given over \(\Spec \Z\), showing the proposition.
\xpf

In general, \(X_\delta\) will not be as simple, and we have to decompose it further using \thref{deodhar in the infinite setting}. 
Before we state our main result, let us compare Schubert cells in different affine Grassmannians. This will allow us to remove the assumption that \(G\) is semisimple or simply connected. In particular, we omit this assumption for the lemma. We will denote by \(G_\adj\) the adjoint quotient of \(G\), and by \(G_\sico\) the simply connected cover of \(G_\adj\).
Recall from \thref{generalized cocharacters} and \thref{notation generalized schubert varieties} that \(\Gr_{G_\sico}^{\leq \mu}\) makes sense for any \(\mu\in X_*(T_{\adj})^+\supseteq X_*(T_{\sico})^+\), where \(T_{\adj}\subseteq G_{\adj}\) and \(T_{\sico}\subseteq G_{\sico}\) are the maximal tori induced by \(T\subseteq G\).

\lemm
\thlabel{reduction to simply connected case}
Let \(\mu\in X_*(T)^+\) be a dominant cocharacter of \(G\), and denote the induced dominant cocharacter of \(G_\adj\) the same way.
Then there exist universal homeomorphisms \(\Gr_{G}^{\leq \mu}\to \Gr_{G_\adj}^{\leq \mu}\gets \Gr_{G_\sico}^{\leq \mu}\), which are equivariant for the \(L^+G\)-, respectively the \(L^+G_\sico\)-action. Moreover, these morphisms restrict to isomorphisms on the open Schubert cells, and are compatible for the intersections with the semi-infinite orbits.
\xlemm
\pf
First, consider the morphism \(\Gr_G\to \Gr_{G_\adj}\) induced by the quotient \(G\to G_\adj\). 
It is clearly \(LG\)-equivariant, and restricts to a morphism \(\Gr_{G}^{\leq\mu}\to \Gr_{G_\adj}^{\leq \mu}\), as both subschemes are defined as orbit closures. 
This latter morphism is proper, and a universal homeomorphism when restricted to the geometric points of \(\Spec \Z\) by \cite[Proposition 3.5]{HainesRicharz:Cohen}. 
Thus, it is universally bijective and universally closed, hence a universal homeomorphism.
It moreover restricts to an isomorphism over \(\Gr_{G}^\mu \to \Gr_{G_\adj}^\mu\), as both source and target are smooth over \(\Spec \Z\).
The compatibility for the intersections with the semi-infinite orbits follows from the \(LG\)-equivariance.

Next, we note that \(LG_\sico\) acts on \(\Gr_{G_\adj}\) via the natural morphism \(LG_\sico\to LG_\adj\). This realizes any connected component of \(\Gr_{G_\adj}\) as a quotient, up to universal homeomorphism, of \(LG_\sico\) by a hyperspecial parahoric subgroup; 
this follows as in the previous paragraph for the neutral connected component, and in general by conjugating \(L^+G_\sico\) by a suitable element of \(LG_\adj(\Z)\). In particular, this identification is \(LG_\sico\)-equivariant. Consequently, it restricts to a universal homeomorphism \(\Gr_{G_\sico}^{\leq \mu}\to \Gr_{G_\adj}^{\leq \mu}\), and to an isomorphism \(\Gr_{G_\sico}^\mu\to \Gr_{G_\adj}^\mu\) by a similar argument as above.
\xpf

For the following theorem, we again omit the assumption that \(G\) is semisimple or simply connected.

\theo
\thlabel{cellularity of intersection:torus}
Let \(\mu\) and \(\nu\) be cocharacters of \(G\), with \(\mu\) dominant. Then the intersection \(\Gr_{G}^\mu \cap \Ss_\nu^-\) admits a filtrable decomposition into cellular schemes.
\xtheo
\pf
By \thref{reduction to simply connected case}, we may assume \(G\) simply connected, so that the results of this 
subsection apply. By \thref{first stratification of intersection}, it is enough to prove that for each 
\(\delta\in \Gamma(\gamma_\mu)\), the intersection \(X_\delta=C_\delta\cap \Gr_{G}^\mu\) is cellular. Moreover, by 
\thref{can assume positively folded}, we may assume that \(\delta\in \Gamma^+(\gamma_\mu)\).

Now, consider the immersion \(C_\delta\to \Sigma(\gamma_\mu)\). As in the regular case (\thref{regular case}), we need to 
determine the preimage \(X_\delta\) of \(\Gr_{G}^\mu\) under this map. Over algebraically closed fields, this is done in 
\cite[Theorem 4]{GL:LSGalleries}, by first identifying the points of \(\Sigma(\gamma_\mu)\) with certain galleries in 
\(\Jj^{\fraka}\), 
and showing that the galleries corresponding to points in \(\Gr_G^\mu\) are exactly the minimal galleries. 
Moreover, such a gallery is minimal exactly when each triple gallery contained in it is minimal \cite[Remarks 6 and 8]{GL:LSGalleries}.
Using the interpretation of galleries of triples (with a fixed source) as certain points of \(\Uu^+(\delta_i)\delta_i \Uu^-(\delta_i)\subseteq \Pp_i/\Qq_i\) \cite[§8]{GL:LSGalleries}, the minimal galleries correspond exactly to \((\Uu^+(\delta_i)\delta_i \Uu^-(\delta_i) \cap \Uu^+(\tau_i^{\min})\tau_i^{\min}) \subseteq \Pp_i/\Qq_i\) \cite[Lemma 12, Proposition 8]{GL:LSGalleries} (where we used the same notation as in \thref{deodhar in the infinite setting}).
Thus, splitting up \(\delta\) into galleries of triples, the description of \(\Sigma(\gamma_\mu)\) as an iterated Zariski-locally trivial fibration with fibers \(\Pp_i/\Qq_i\) shows that \(X_\delta \subseteq \Sigma(\gamma_\mu)\) is a sub-iterated fibration, with iterated fibers given by \(\Uu^+(\delta_i) \delta_i \Uu^+(\delta_i) \cap \Uu^+(\tau_i^{\min})\tau_i^{\min}\).
Since this description works uniformly over all geometric points of \(\Spec \Z\), it already holds over \(\Spec \Z\).
Hence the theorem follows by inductively applying \thref{deodhar in the infinite setting}.
\xpf

\rema
\begin{enumerate}
  \item By replacing the Borel \(B\) by its opposite Borel, the theorem also holds for the positive semi-infinite orbits \(\Ss_\nu^+\).
  \item Let \(n\gg 0\) be such that the \(L^+G\)-action on \(\Gr_G^\mu\) factors through \(L^nG\), and let \(\Pp_\mu^n\subset L^nG\) be the stabilizer of \(t^\mu\).
  Then \thref{Torsor is trivial over decomposition} implies that the \(\Pp_\mu^n\)-torsor \(L^nG\to \Gr_G^\mu\) is trivial over each cell in the above decomposition of \(\Gr_G^\mu\cap \Ss_\nu^+\).
\end{enumerate}
\xrema

\exam
\thlabel{Example PGL2}
Let us work out what happens for the (non-simply connected) group \(\PGL_2\), whose simply connected cover is \(\SL_2\). 
Considering the natural identification \(X_*(\PGL_2)\cong \Z\), let \(\mu\in \Z_{\geq 0}\) be a dominant cocharacter 
of \(\PGL_2\), which lives in the standard apartment 
\(\Aa_{\SL_2}=X_*(\SL_2)\otimes_{\Z} \RR\cong \RR\) of \(\SL_2\). Note that \(\mu\) is (induced by) a cocharacter of 
\(\SL_2\) exactly when \(\mu\) is even. Consider the unique minimal gallery
\[\gamma_\mu:=\left(\{0\}\subset [0,1] \supset \{1\} \subset \ldots \supset \{\mu-1\} \subset [\mu-1,\mu]\supset \{\mu\}\right)\]
joining \(0\) with \(\mu\). Noting that the affine Weyl group of \(\SL_2\) is generated by two reflections 
\(s_0:x\mapsto -x\) and \(s_1:x\mapsto 2-x\), the gallery of types of \(\gamma_\mu\) is given by
\[\left( \{s_0\} \supset \varnothing \subset \{s_1\} \supset \ldots \subset \{s_{\mu-1 (\modulo 2)}\} \supset \varnothing \subset \{s_{\mu (\modulo 2)}\}\right),\]
and there is a bijection
\[\Gamma(\gamma_\mu)\cong \langle s_0 \rangle \times \langle s_1\rangle \times \ldots \times \langle s_{\mu -2 (\modulo 2)} \rangle \times \langle s_{\mu-1 (\modulo 2)} \rangle.\]

As the face corresponding to the empty type is the fundamental alcove, we see that under this bijection, 
\(\gamma_\mu\) corresponds to \((1,s_1,s_0,\ldots)\), which is just a straight path in \(\Aa_{\SL_2}\) from \(\{0\}\) to 
\(\{\mu\}\). The other combinatorial galleries \((\delta_0,\delta_1,\delta_2,\ldots)\) can be described as 
follows: they begin at \(0\), and if \(\delta_0=1\), then they start in the positive direction (towards \(\mu)\), 
otherwise they start in the opposite direction. After this, if \(\delta_i=1\), the path turns around, so that the 
\(i-1\)th and \(i\)th large faces in \(\delta\) agree, otherwise the path continues in the same direction. However, at the 
points where the path turns around, there is a fold, which is positive exactly when the path is going into the 
negative direction, and turns to the positive direction. 

In particular, there is a unique positively folded 
combinatorial gallery in \(\Gamma^+(\gamma_\mu)\), with source \(0\) and target \(\nu\), where \(\nu\in X_*(T_\adj)\) corresponds to an integer congruent to \(\mu\) modulo 2, such that \(-\mu\leq \nu\leq \mu\). We 
note that in this case, the index set \(J_{-\infty}\) consists of those indices for which a path in \(\Aa_{\SL_2}\) 
moves in the positive direction, and \(J_{-\infty}^-\) those indices which correspond to a fold (necessarily positive). 
So, using the proof of \thref{regular case}, we see that \(\Gr_G^\mu\cap \Ss_{\mu}^-=\Spec \Z\), that 
\(\Gr_G^\mu\cap \Ss^-_{-\mu}\cong \A^{\mu}\), and that \(\Gr_G^{\mu}\cap \Ss^-_{\nu}\cong \Gm\times 
\A^{\frac{\mu-\nu}{2}-1}\) when \(-\mu<\nu<\mu\). 
These dimensions agree with \cite[Theorem 3.2]{MirkovicVilonen:Geometric}, cf.~also \cite[(3.2.2)]{XiaoZhu:Cycles}.
\xexam

\exam
In the previous example, all galleries in \(\Gamma^+(\gamma_\mu)\) had different targets, so that each \(\Gr_G^\mu\cap \Ss_\nu^-\) was already isomorphic to a product of \(\A^1\)'s and \(\Gm\)'s.
This does not hold in general, even for regular cocharacters.
For example, let \(\mu\in \Aa_G\) be the simple reflection of the origin, over the codimension 1 face of \(\Delta_f\) opposite to the origin, and \(\gamma_\mu\) the natural corresponding minimal combinatorial gallery.
For each alcove \(\mathfrak{a}\) adjacent to the origin, let \(\mathbf{H}_{\beta,m}\) be the reflection hyperplane containing the codimension 1 face \(F\) of \(\mathfrak{a}\) opposite to the origin.
Then, if \(\mathfrak{a}\subseteq \mathbf{H}_{\beta,m}^+\), there is a positively folded combinatorial gallery
\[(F_f\subset \mathfrak{a}\supset F \subset \mathfrak{a}\supset F_f)\in \Gamma^+(\gamma_\mu).\]
So if there are multiple alcoves \(\mathfrak{a}\) as above, which is the case as soon as \(G\) has semisimple rank \(>1\), there are multiple galleries in \(\Gamma^+(\gamma_\mu)\) with the origin as target.
We refer to \cite[Examples 7 and 10]{GL:LSGalleries} for more details in the case \(G=\SL_3\).
\xexam

\section{Beilinson--Drinfeld Grassmannians and convolution}
\label{sect--BD Grass}
\subsection{Beilinson--Drinfeld Grassmannians} 
In this section, we collect basic geometric information about the Beilinson--Drinfeld affine Grassmannians, following \cite[\S 3]{Zhu:Introduction} and \cite[\S 3]{HainesRicharz:TestFunctionsWeil} in the case of constant group schemes. 

Let $\cal{G}$ be a smooth affine group scheme over $S$ and let $X:= \A^1_S$. Then we have a distinguished $S$-point $\{0\} \in X(S)$. To define the Beilinson--Drinfeld Grassmannians, we first introduce some notation for working with \'etale $\cal{G}$-torsors.
For $\Spec(R) \in \AffSch_S$, let $X_R := X \times_S \Spec(R)$. If $x \colon \Spec(R) \rightarrow X$ is a morphism, we denote the graph of $x$ by $\Gamma_x \subset X_R$. We fix a trivial $\cal{G}$-torsor $\mathcal{E}_0$ on $X$, and for any $\Spec(R) \in \AffSch_S$ we also denote its base-change to $X_R$ by $\mathcal{E}_0$.

Note that any $\Spec(R) \in \AffSch_S$ can be viewed as an \(X\)-scheme lying over $0$ by composing with inclusion $S \xrightarrow{0}$ X. In the following proposition we make this identification. A theorem of Beauville and Laszlo \cite{BeauvilleLaszlo:Descent} implies that $\Gr_{\mathcal{G}}$ (defined in \refeq{definition.Gr.G}) has the following moduli interpretation, cf.~\cite[Example 3.1 (i)]{HainesRicharz:TestFunctionsWeil}. 
\prop \thlabel{moduliprop}
There is a canonical isomorphism of \'etale sheaves
$$\Gr_{\mathcal{G}}(R) \cong \left\{(\mathcal{E}, \beta) \: : \: \mathcal{E} \text{ is a } \mathcal{G}\text{-torsor on } X_R, \: \beta \colon \restr{\mathcal{E}}{X_R - \Gamma_0} \cong \restr{\mathcal{E}_0}{X_R - \Gamma_0}\right\}.$$
\xprop
 
Let $I$ be a nonempty finite set. For a point $x = (x_i) \in X^I(R)$, let
$$\Gamma_x = \bigcup_{i \in I} \Gamma_{x_i} \subset X_R.$$

\defi \thlabel{BD.Def}
The Beilinson--Drinfeld Grassmannian  for $\mathcal{G}$ over $X^I$ is the functor
$$\Gr_{\mathcal{G}, I}(R) = \left\{(x, \mathcal{E}, \beta) \: : \: x \in X^I(R), \: \mathcal{E} \text{ is a } \mathcal{G}\text{-torsor on } X_R, \: \beta \colon \restr{\mathcal{E}}{X_R - \Gamma_x} \cong \restr{\mathcal{E}_0}{X_R - \Gamma_x}\right\}.$$ 
\xdefi

\rema \thlabel{Other.curves}
The definition of $\Gr_{\mathcal{G}, I}$ makes sense for general smooth curves $X$, but the existence of the necessary t-structure on $\DTM(X^I)$ is not known in general. In future work we hope to extend the results in this paper to other curves for which the t-structure is known to exist, such as $X= \P$.
\xrema

\rema \thlabel{HR.Translate}
Our definition of $\Gr_{\mathcal{G}, I}$ is the following specialization of \cite[Eqn. (3.1)]{HainesRicharz:TestFunctionsWeil}. In the notation of loc.~cit., let $\Spec(\mathcal{O}) = \A_S^I$ and $X = \A^1_{\mathcal{O}}$. If $\Spec(\mathcal{O})$ has affine coordinate functions $x_i$ for $i \in I$ and $X$ has the affine coordinate function $t$, then let $D$ be the divisor on $X$ defined locally by the ideal $\prod_{i \in I} (t-x_i)$. After identifying pairs $(\Spec(R) \in \AffSch_S, x \in \A_S^I(R))$ with objects in $\AffSch_{\Spec(\mathcal{O})}$, \thref{BD.Def} agrees with $\Gr_{(X, \mathcal{G} \times X, D)}$ as defined in \cite[Eqn. (3.1)]{HainesRicharz:TestFunctionsWeil}.
\xrema

\lemm \thlabel{lemm.BDrep}
If $\mathcal{G} = G$ is a split reductive group, $\Gr_{G,I}$ is represented by an ind-projective scheme over $X^I$. If $\mathcal{G} = P$ is a standard parabolic subgroup of $G$, $\Gr_{P,I}$ is represented by an ind-scheme of ind-finite type over $X^I$.
\xlemm

\pf
For any closed immersion of group schemes $G \r \text{GL}_n$, the quotient $\text{GL}_n /G$ is an affine $S$-scheme, cf.~\cite[Corollary 9.7.7]{Alper:Adequate}.  Hence $\Gr_{G,I}$ is ind-projective by \cite[Corollary 3.11 (i)]{HainesRicharz:TestFunctionsWeil}. By \cite[Theorem 2.1 (iii) and Theorem 3.17]{HainesRicharz:TestFunctionsWeil}, this also implies the claim for $\Gr_{P,I}$.
\xpf

In the special case $I =\{*\}$ is a singleton, \thref{moduliprop} implies there is a canonical isomorphism
\begin{equation} \label{Gr.Pt} \Gr_{\mathcal{G}, \{*\}} \cong \Gr_{\mathcal{G}} \times X. \end{equation} For general $I$, the fiber of $\Gr_{\mathcal{G},I}$ over a point in $(x_i) \in X^I$ depends on the partition of $(x_i)$ into pairwise distinct coordinates. More precisely, let $$\phi \colon I \twoheadrightarrow J$$ be a surjection of nonempty finite sets. This induces a partition $$I = \bigcup_{j \in J} I_j, \quad I_j:=\phi^{-1}(j).$$ Let
\begin{equation} \label{Loc.Sub} X^{\phi} := \{(x_i) \in X^I \: : \: x_i = x_{i'}  \text{ if and only if } \phi(i) = \phi(i')\}.\end{equation} This is a locally closed subscheme of $X^I$. In the special case $\phi = \id$ we write $X^\circ = X^{\id}$, which is the locus with pairwise distinct coordinates. For later use, we also define the open subscheme
\begin{equation} \label{Open.Sub} X^{(\phi)} := \{ (x_i) \in X^I \: : \: x_i \neq x_{i'} \text{ if } \phi(i) \neq \phi(i')\} \subset X^I.\end{equation}

\prop\thlabel{prop.factorization of Gr}
There is a canonical isomorphism
$$\restr{\Gr_{\mathcal{G},I}}{X^{\phi}} \cong \prod_{j \in J} \Gr_{\mathcal{G}}  \times X^{\phi}\eqlabel{ZhuFactor}.$$
\xprop

\pf
The arguments in \cite[Proposition 3.1.13, Theorem 3.2.1]{Zhu:Introduction} generalize to an arbitrary base. 
\xpf

We now define the global loop groups. For $x \in X^I(R)$, let $\hat{\Gamma}_x$ be the formal completion of $\Gamma_x$ in $X_R$. Locally on $S$, $\hat{\Gamma}_x$ is the formal spectrum of a topological ring, so by forgetting the topology we can view $\hat{\Gamma}_x$ as an object in $\AffSch_S$, cf.~\cite[\S 3.1.1]{HainesRicharz:TestFunctionsWeil}. 
Following the discussion in loc.~cit., one can view ${\Gamma}_x$ as a Cartier divisor in $\hat{\Gamma}_x$, so $\hat{\Gamma}_x^{\circ} := \hat{\Gamma}_x - {\Gamma}_x \in \AffSch_S$. 

\defi 
The \emph{global positive loop group} $L_I^+\mathcal{G}$ (resp.~the \emph{global loop group} $L_I\mathcal{G}$) is the functor
$$L_I^+\mathcal{G}(R) = \{(x, g) \: : \: x \in X^I(R), \: g \in \mathcal{G}(\hat{\Gamma}_x) \}.$$
$$L_I\mathcal{G}(R) = \{(x, g) \: : \: x \in X^I(R), \: g \in \mathcal{G}(\hat{\Gamma}^\circ_x) \}.$$ 
\xdefi

By \cite[Lemma 3.2]{HainesRicharz:TestFunctionsWeil}, $L^+_I\mathcal{G}$ is represented by a pro-smooth affine group scheme over $X^I$ and $L_I \mathcal{G}$ is represented by an ind-affine ind-scheme over $S$.
By arguments very similar to \thref{prop.factorization of Gr}, both of these groups 
satisfy a factorization property as in \refeq{ZhuFactor}. 
More precisely, for a surjection $\phi \colon I \twoheadrightarrow J$ we have
$$\restr{L_I\mathcal{G}}{X^{\phi}} \cong \prod_{j \in J} L\mathcal{G}  \times X^{\phi}, \eqlabel{ZhuGroupFactor}$$ and the same holds for $L^+_I\mathcal{G}$.
We emphasize that if $I = \{*\}$ is a singleton, there is a canonical isomorphism
$L_{\{*\}}^+\mathcal{G} \cong L^+\mathcal{G} \times X.$ 
The proof of the following lemma was explained to us by T. Richarz.

\lemm \thlabel{BD.Zar}
If $\mathcal{G} = G$ is a split reductive group, the Beilinson--Drinfeld Grassmannian $\Gr_{G,I}$ can be identified with the Zariski, Nisnevich, and \'etale sheafifications of $L_IG/L^+_IG$.
\xlemm

\pf We show that the arguments in \thref{Sheafifications of quotient agree} globalize. Since sheafification commutes with base change, we may assume $S = \Spec \Z$. By \cite[Lemma 3.4]{HainesRicharz:TestFunctionsWeil}, we have a right $L_I^+G$-torsor $L_IG \r \Gr_{G,I}$, which we must show is Zariski-locally trivial. 
The big open cell in \cite[Lemma 3.15]{HainesRicharz:TestFunctionsWeil} is an open sub-ind-scheme of $\Gr_{G,I}$ over which $L_IG \r \Gr_{G,I}$ admits a section. 
Under the factorization isomorphism \refeq{ZhuFactor}, it restricts to products of the big open cell $L^{--}G = \ker(L^{-}G \r G)$ in $\Gr_G$, where $L^{-}G(R) = G(R[t^{-1}])$ and the map is $t^{-1} \mapsto 0$.
By \cite[Definition 5 ff.]{Faltings:Loops}, 
 $\Gr_G$ is covered by left translates of $L^{--}G$ by the points in $LT(\Z)$ given by evaluation of cocharacters in $X_*(T)$ at $t$. 
Fix a surjection $\phi \colon I \r J$. Then $\restr{L_IT}{X^{\phi}} \cong (LT)^J \times X^{\phi}$. 
By the previous discussion, it suffices to show that a point in $((LT)^J \times X^{\phi})(X^{\phi})$
corresponding to a tuple in $X_*(T)^J$ lifts to an $X^I$-point of $L_IT$.
For this we use the explicit description of $L_IT(X^I)$ in \cite[\S 3.1.1]{HainesRicharz:TestFunctionsWeil}. 
As in \thref{HR.Translate}, let $D = \prod_{i \in I}(t-x_i)$, where $|I| = n$, and let $X^I(\!(D)\!)$ be the ring of functions on the complement of $D$ in the completion of $\Z[x_1, \ldots, x_n][t]$ at $D$. If we identify $X^I$ with $\Spec \Z[x_1, \ldots, x_n]$,  then $L_IT(X^I) = T(X^I(\!(D)\!))$. Note that since $X = \mathbb{A}^1$, we have the global invertible functions $(t-x_i)$ on $X^I(\!(D)\!)$.
Choosing $|J|$ of the coordinates $x_1, \ldots, x_n$ to represent the distinct $J$ coordinates over $X^\phi$, it follows that the lift we need exists, and is already defined before passing to the completion.
\xpf

By \thref{BD.Zar},  the group $L_IG$ acts on $\Gr_{G,I}$ on the left. We denote the result of the action of $g \in L_IG(R)$ on $(x, \mathcal{E}, \beta) \in \Gr_{G,I}(R)$ by $(x, g \mathcal{E}, g \beta)$. 
We refer to \cite[\S 3.1]{Richarz:New} for more details on these group actions.
\defi

\thlabel{Hecke.stack}
For any nonempty finite set $I$, we define the \emph{Hecke prestack} by $\Hck_{G,I} := L^+_I G \backslash \Gr_{G,I}$.
We will denote the canonical quotient map by $u \colon \Gr_{G,I} \r \Hck_{G,I}$.
\xdefi

The following result should be well-known, but we include it as we were unable to find a reference.

\prop
\thlabel{fiber product of affine grassmannians}
Let \(K,L,M\) be smooth affine \(S\)-group schemes, and \(K\to M\) and \(L\to M\) group homomorphisms, with \(L\to M\) surjective. Suppose $K\times_M L$ is represented by a smooth affine \(S\)-group scheme, so that $\Gr_{K\times_M L}$ is well-defined. Then the natural morphisms \(\Gr_{K\times_M L}\r \Gr_K\times_{\Gr_M} \Gr_L\), and 
\(\Gr_{K\times_M L,I}\r \Gr_{K,I}\times_{\Gr_{M,I}} \Gr_{L,I}\) are isomorphisms for any finite set \(I\). 
\xprop
\pf
We will only show the first assertion; the case of Beilinson--Drinfeld Grassmannians can be handled analogously.
To construct the inverse, let \(R\) be a scheme over \(S\). Let $x \in X(R)$ be the point $\Spec(R) \r X$ corresponding to the origin.
An element of \((\Gr_K\times_{\Gr_M} \Gr_L)(R)\) can be represented by a pair \((\Ee_K,\beta_K)\), with \(\Ee_K\) a \(K\)-torsor on $\hat{\Gamma}_x$ and  \(\beta_K:\restr{\Ee_K}{\hat{\Gamma}_x^\circ} \cong \restr{\Ee_{K,0}}{\hat{\Gamma}_x^\circ}\), a similar pair \((\Ee_L,\beta_L)\) for \(L\),
and an isomorphism \(\alpha:\Ee_K\times^{K} M\cong \Ee_L \times^{L} M\), 
commuting with \(\beta_K\) and \(\beta_L\) under the natural identifications \(\Ee_{0,K}\times^{K} M\cong \Ee_{0,M} \cong \Ee_{0,L}\times^{L} M\). 
Let us denote \(\Ee_K\times^{K} M\cong \Ee_L \times^{L} M\) by \(\Ee_M\).
Using the natural morphisms \(\Ee_K\cong \Ee_K \times^{K} K \to \Ee_K \times^{K} M \cong \Ee_M\) and \(\Ee_L \to \Ee_M\), 
we can consider the fiber product \(\Ee_K \times_{\Ee_M} \Ee_L\); the surjectivity of \(L\to M\) then ensures this is nonempty.
Moreover, the isomorphisms 
\(\beta_K\) and \(\beta_L\) induce an isomorphism \(\beta_{K\times_M L}\colon \restr{\Ee_K\times_{\Ee_M} \Ee_L}{\hat{\Gamma}^{\circ}_x} \cong \restr{\Ee_{0,K\times_M L}}{\hat{\Gamma}^{\circ}_x}\). We leave it to the reader to verify that the inverse is
$$(\Ee_K,\beta_K, \Ee_L,\beta_L, \alpha) \mapsto (\Ee_K \times_{\Ee_M} \Ee_L, \beta_{K\times_M L}).$$
\xpf

\subsection{Convolution Grassmannians}
\label{Section:Convolution}
\subsubsection{Local case} Recall that $\calP \subset LG$ denotes a parahoric subgroup.

\defilemm
\thlabel{DTM.convolution}
\thlabel{Convolution.Anti}
The \emph{convolution product} is defined to be the functor
$$\star : \DMrx(\calP \backslash LG / \calP) \x \DMrx(\calP \backslash LG / \calP) \r \DMrx(\calP \backslash LG / \calP)$$
$$\calF_1 \star \calF_2 := m_! p^! (\calF_1 \boxtimes \calF_2),\eqlabel{star}$$ 
where the maps are the natural quotient and multiplication maps (which are maps of prestacks):
$$\xymatrix{ \calP \backslash LG / \calP \x \calP \backslash LG / \calP & \ar[l]_-p \calP \backslash LG \x^{\calP} LG / \calP  \ar[r]^-m & \calP \backslash LG / \calP. }\eqlabel{Loc.Conv}$$
The functor $\star$ preserves (anti-effective) stratified Tate motives (\thref{equivariant.DTM}). It endows (at least) the homotopy category $\Ho(\DTMrx(\calP \backslash LG / \calP)^{(\anti)})$ with the structure of a monoidal category.
\xdefilemm

\rema
\thlabel{placid prestacks}
Recall the existence of the functors in the right hand side of \refeq{star}.
For any map of prestacks $f : Y \r Z$, we have a !-pullback functor by construction, cf.~\refeq{DM.prestack}.
If $f$, or its Zariski or Nisnevich sheafification, is an ind-schematic map (as is the case for $m$), then $f^!$ admits a left adjoint $f_!$, cf.~\cite[Lemma~2.2.9, Proposition~2.3.3]{RicharzScholbach:Intersection}, together with \thref{equivariant.functoriality}.
Finally, the exterior product functor for motives on placid prestacks has been constructed in \cite[Corollary~A.15]{RicharzScholbach:Motivic}.
The prestacks $\calP \backslash LG / \calP$ and also the Hecke prestacks $\Hck_{G, I}$ are placid, the point being that the quotient is formed with respect to a pro-smooth group $\calP$ (resp.~$L^+_I G$ in the global case below).
For later use, we note that the !-pullback along a pro-smooth quotient map (such as $\Gr_G \r \calP \backslash \Gr_G$) is compatible with $\boxtimes$.
\xrema

\pf
The proofs of~\cite[Lemma~3.7, Theorem~3.17]{RicharzScholbach:Motivic} carry over to show that Tate motives are preserved under convolution, and that convolution defines a monoidal structure (at least up to homotopy).
Indeed, the associativity of the convolution product proved there for motives with rational coefficients (which satisfy étale descent) also holds for motives only satisfying Nisnevich descent, so \cite[Propositions~2.3.3, 2.4.4]{RicharzScholbach:Intersection} holds in the present setting.
Concerning the preservation of Tate motives under convolution, the proof of \cite[Theorem~3.17]{RicharzScholbach:Motivic} and \cite[Theorem~5.3.4]{RicharzScholbach:Intersection} applies to the present setting since the characterization of the full subcategory $\DTMrx(\calP \backslash LG / \calP) \subset \DM(\calP \backslash LG / \calP)$ in Part ii) there holds for motives satisfying only Nisnevich descent. 
This is proved by reducing to the situation discussed in \cite[Proposition~3.2.23, Corollary~3.2.24]{RicharzScholbach:Intersection} (with ``étale'' being replaced by ``Nisnevich''), which in its turn reduces to the computation of $\DTM_G((G/H)^\Nis)$ considered in \cite[Proposition~3.1.23]{RicharzScholbach:Intersection} and \cite[Proposition~1.1]{RicharzScholbach:IntersectionCorrigendum}.

What is more, these proofs also show more generally that the convolution product defined in \cite[Definition~3.1]{RicharzScholbach:Intersection} (for a triple of parahorics $\calP', \calP, \calP''$) preserves anti-effective Tate motives. More precisely, the argument where all parahorics equal $\calI$ in \cite[Proposition~3.19]{RicharzScholbach:Intersection} applies verbatim. The next argument where $\calP$ is arbitrary in \cite[Proposition~3.26]{RicharzScholbach:Intersection} needs to be modified by using the functors $a^*$ and $b^*$ instead of $a^!$ and $b^!$ in the diagram in loc.~cit. so as to avoid introducing positive twists. The general case follows by reduction to $\calP' = \calP'' = \calI$ as in the final part of \cite[Theorem~3.17]{RicharzScholbach:Intersection}, and preservation of anti-effectivity follows as in the third case of \thref{Fl.universally.WT}. We note that to apply these arguments we must replace all \'etale quotients of loop groups by Zariski or Nisnevich quotients, which we may do by \thref{Sheafifications of quotient agree} and \thref{Stab.Split} (this is relevant in \cite[Proposition~3.1.23, Theorem~5.3.4]{RicharzScholbach:Intersection}).
\xpf

\subsubsection{Global case: Type I}
By replacing $LG$ by $L_IG$ and $\calP$ by $L^+_IG$ in \refeq{Loc.Conv}, and using any number of factors of $\Hck_{G,I}$, we obtain a convolution product 
$$\star : \DMrx(\Hck_{G,I}) \x \cdots \x \DMrx(\Hck_{G,I} ) \r \DMrx(\Hck_{G,I} ).$$ 
$$\calF_1 \star \cdots \star \calF_n := m_! p^! (\calF_1 \boxtimes \ldots \boxtimes \calF_n)[-(n-1)|I|]. \eqlabel{GrI.Conv}$$ This uses \thref{BD.Zar}.
The box products are formed with respect to $X^I$. As above, this functor turns the homotopy category of $\DMrx(\Hck_{G,I})$ into a monoidal category. It will be clear from the context if we mean the local or global version of $\star$. The shift by $-(n-1)|I|$ ensures the box product will be right exact for a t-structure introduced later, cf.~\thref{fusion.Satake} and \thref{fiber.functor.monoidal}.

\subsubsection{Beilinson--Drinfeld convolution Grassmannians}
If $J$ is an ordered finite set we identify $J$ with $\{1, \ldots, |J|\}$. For the rest of this subsection, fix a surjection of nonempty finite sets $\phi \colon I \twoheadrightarrow J$, where $J$ is ordered. If $(x_i) \in X^I$, let $x_{I_j} \in X^{I_j}$ be the corresponding component.

\defi 
The \emph{convolution Beilinson--Drinfeld Grassmannian} over $X^I$ is the functor
$$\widetilde{\Gr}_{G,\phi}(R) = \left\{(x_j, \mathcal{E}_j, \beta_j)_{j =1, \ldots, |J|} \: : \: x_j \in X^{I_j}(R), \: \mathcal{E}_j \text{ is a } G\text{-torsor on } X_R, \: \beta_j \colon \restr{\mathcal{E}_j}{X_R - \Gamma_{x_{j}}} \cong \restr{\mathcal{E}_{j-1}}{X_R - \Gamma_{x_{j}}}\right\}, \eqlabel{Definition convolution Grassmannian}$$
where $\mathcal{E}_0$ is the trivial \(G\)-torsor.
\xdefi

 Since $\widetilde{\Gr}_{G,\phi}$ can be written as a twisted product of the $\Gr_{G, I_j}$ (see \cite[Eq.~(3.1.21) ff.]{Zhu:Introduction} or \refsect{onemoreConv}), and the relevant torsors are Zariski-locally trivial by \thref{BD.Zar}, $\widetilde{\Gr}_{G,\phi}$ is represented by an ind-proper ind-scheme over $X^I$. There is a convolution morphism
$$m_{\phi} \colon \widetilde{\Gr}_{G,\phi} \rightarrow \Gr_{G,I}, \quad (x, \mathcal{E}_i, \beta_i) \mapsto (x, \mathcal{E}_{|J|}, \beta_1 \circ \dots \circ \beta_{|J|})\eqlabel{convolution.morphism}$$ which restricts to an isomorphism over the locus of $X^I$ with pairwise distinct coordinates. The group $L_I^+G$ acts on $\widetilde{\Gr}_{G,\phi}$ on the left by $(x, \mathcal{E}_j, \beta_j) \mapsto (x, g\mathcal{E}_j, g\beta_jg^{-1})$, cf.~\cite[Corollary 3.10 ff.]{Richarz:New}, and $m_{\phi} $ is equivariant for this action. There is also a factorization property for $\widetilde{\Gr}_{G,\phi}$ similar to \refeq{ZhuFactor}. More precisely, let $\phi' \colon I \twoheadrightarrow K$ be a surjection. Then
$$\restr{m_\phi}{X^{\phi'}} \colon \prod_{k \in K} Y_k  \times X^{\phi'} \r \prod_{k \in K} \Gr_{G}  \times X^{\phi'}$$ splits into a product of local convolution maps $Y_k \r \Gr_G$, where $Y_k = \Gr_G \tilde{\times} \cdots \tilde{\times} \Gr_G$ and the number of factors in this twisted product equals the number of $j \in J$ such that $\phi^{-1}(j) \cap (\phi')^{-1}(k) \neq \emptyset$. 
Finally, if $\phi = \id$ we let $\widetilde{\Gr}_{G,I}: = \widetilde{\Gr}_{G,\phi}$ and $m_I := m_{\phi}$.

\subsubsection{Global case: Type II} \label{sect--subsub-TypeII}
As in previous approaches to geometric Satake, we will relate the convolution product $\star$ on $\DMrx(\Hck_{G,I})$ to a fusion product. In order to prove facts about the fusion product, we will use another type of global convolution product described below.

Let $I = I_1 \sqcup I_2$ be a partition into two nonempty finite sets associated to a surjection $\phi \colon I \r \{1,2\}$.
We define a functor on $\AffSch_S$ by
$$\eqalign{
L_{I_1 I_2}G(R) = \{ & ((x_i, \mathcal{E}_i, \beta_i)_{i=1,2}, \sigma) \colon x_i \in X^{I_i}(R), \mathcal{E}_i \text{ is a } G\text{-torsor on } X_R, \cr
& \beta_i \colon \restr{\mathcal{E}_i}{X_R - \Gamma_{x_i}} \cong \restr{\mathcal{E}_{0}}{X_R - \Gamma_{x_i}}, \sigma \colon \restr{\mathcal{E}_0}{\hat{\Gamma}_{x_2}} \cong \restr{\mathcal{E}_1}{\hat{\Gamma}_{x_2}} \}.}$$
There is a commutative diagram, which we explain below.
$$\xymatrix{
\Gr_{G, I_1} \x_{S} \Gr_{G,I_2} \ar[d]^u & L_{I_1I_2}G \ar[d]_q \ar[l]_-{p} 
\\
\Gr_{G, I_1} \x_{S} \Hck_{G,I_2} \ar[d]^v &  \widetilde{\Gr}_{G, \phi} \ar[d]^w \ar[l]_-{\tilde p} \ar[r]^{m_\phi} & \Gr_{G, I} \ar[d]^w & 
\\
\Hck_{G, I_1} \x_{S} \Hck_{G,I_2} & L^+_I G \setminus  \widetilde{\Gr}_{G, \phi} \ar[l]_-{\ol p} \ar[r]^-{\ol m_\phi} & \Hck_{G,I}} \eqlabel{TypeII}$$

\lemm
The projection $p \colon L_{I_1 I_2}G \r \Gr_{G, I_1} \times \Gr_{G, I_2}$ which forgets $\sigma$ is a Zariski-locally trivial $L_{I_2}^+G \times_{X^{I_2}} X^I$-torsor. Consequently, $L_{I_1 I_2}G$ is represented by an ind-scheme over $X^I$.
\xlemm

\pf 
The group $L_{I_2}^+G \times_{X^{I_2}} X^I$ acts on $L_{I_1I_2}G$ by changing $\sigma$, and $p$ is a torsor for this group.
For every point in $\Gr_{G, I_1}(R)$, $\mathcal{E}_1$ is trivializable on $X_R - \Gamma_{x_1 \cup x_2}$, so by \thref{BD.Zar}, $\mathcal{E}_1$ is trivializable on $\hat{\Gamma}_{x_1 \cup x_2}$ Zariski-locally with respect to $R$.
By pulling back along the map $\hat{\Gamma}_{x_2} \r \hat{\Gamma}_{x_1 \cup x_2}$ this shows that a trivialization $\sigma$ as in the definition of $L_{I_1 I_2}G(R)$ exists Zariski-locally on $R$.
\xpf

For a point in $LG_{I_1I_2}(R)$, we can construct a point $(x, \mathcal{E}_j', \beta_j')_{j=1,2} \in \widetilde{\Gr}_{G, \phi}(R)$ as follows. Let $x = x_1 \times x_2$, $\mathcal{E}_1' = \mathcal{E}_1$ and $\beta_1' = \beta_1$. 
Let $\mathcal{E}_2'$ be the bundle obtained using \cite{BeauvilleLaszlo:Descent} and \cite[\S 2.12]{BeilinsonDrinfeld:Quantization} to glue $\restr{\mathcal{E}_1}{X_R - \Gamma_{x_2}}$ to $\restr{\mathcal{E}_2}{\hat{\Gamma}_{x_2}}$ along $\restr{\sigma \circ \beta_2}{\hat{\Gamma}_{x_2}^\circ}$. 
By construction there is an isomorphism $\beta_2' \colon \restr{\mathcal{E}_2'}{X_R - \Gamma_{x_1}} \cong \restr{\mathcal{E}_1}{X_R - \Gamma_{x_1}}$. 
This defines the map $q$, which is a torsor for the action of $L_{I_2}^+G \times_{X^{I_2}} X^I$, that fixes $(\mathcal{E}_1, \beta_1)$, and sends $(\mathcal{E}_2, \beta_2, \sigma)$ to $(g\mathcal{E}_2, g \beta_2, \sigma g^{-1})$, cf.~\cite[\S 1.7.4]{BaumannRiche:Satake}.

\lemm
The map $q$ is a Zariski-locally trivial $L_{I_2}^+G \times_{X^{I_2}} X^I$-torsor.
\xlemm
\pf 
Fix a point $(x_i, \mathcal{E}_i, \beta_i)_{i=1,2} \in \widetilde{\Gr}_{G, \phi}(R)$. 
By \thref{BD.Zar}, after an affine Zariski cover $\Spec(R') \r \Spec(R)$, both bundles $\mathcal{E}_i$ become trivial on $\hat{\Gamma}_{x_1 \cup x_2}$. 
Then the isomorphism $\beta_2 \colon \restr{\mathcal{E}_2}{X_{R'} - \Gamma_{x_2}} \cong \restr{\mathcal{E}_1}{{X_{R'} - \Gamma_{x_2}}}$ shows that $\mathcal{E}_2$ is obtained by gluing $\restr{\mathcal{E}_1}{{X_{R'} - \Gamma_{x_2}}}$ to the trivial bundle $\restr{\mathcal{E}_0}{\hat{\Gamma}_{x_2}}$. 
We can then construct a section over $R'$ by taking $\mathcal{E}_2' = \mathcal{E}_0$ to be the trivial bundle on $X_{R'}$.
\xpf

The maps labelled with $u$, $v$, and $w$ are the natural maps to the prestack quotients by the left $L_{I_2}^+G \times_{X^{I_2}} X^I$, $L_{I_1}^+G \times_{X^{I_1}} X^I$, and $L_I^+G$-actions. 
The map $\tilde p$ exists because $p$ is equivariant for the action of $L_{I_2}^+G \times_{X^{I_2}} X^I$ on the right factor of $\Gr_{G, I_1} \x_S \Gr_{G, I_2}$ and the action via $q$ on $L_{I_1I_2}G$.
Technically, $\tilde p$ only exists after taking the Zariski-sheafification of the target. 
However, for any prestack $Z$ the !-pullback along the sheafification map $Z \r Z_\Zar$ induces an equivalence $\DMrx(Z_\Zar) \stackrel \cong \r \DMrx(Z)$ (and likewise with the Nisnevich topology), so we can safely pretend $\tilde p$ exists as stated.
Likewise, the map $\ol p$ exists because we only take pre-stack quotients on the bottom row.

For $\calF_i \in \DMrx(\Hck_{G,I_i})$, we can form the Type II convolution product $$\ol m_{\phi !} \ol p^! (\calF_1 \boxtimes \calF_2) \in \DMrx(\Hck_{G,I}).$$

\lemm \thlabel{Twisted.Prod.Global}
The motive on \(\widetilde{\Gr}_{G,\phi}\) underlying $ \ol p^! (\calF_1 \boxtimes \calF_2)$ agrees with what is in the literature often denoted by $\calF_1 \widetilde \boxtimes \calF_2$, cf.~\cite[A.1.2]{Zhu:Introduction}. Similarly, the motive on \(\Gr_{G,I}\) underlying $\ol m_{\phi !} \ol p^! (\calF_1 \boxtimes \calF_2)$ agrees with $m_{\phi !}(\calF_1 \widetilde \boxtimes \calF_2)$. This specializes to the construction in \cite[Eqn. (5.6)]{MirkovicVilonen:Geometric} (without the perverse truncation of the box product) under Betti realization $\rho_\Betti$ for compact motives over $S = \Spec \C$. 
\xlemm 

\pf
The left most functor $v^!$ is, by definition, compatible with $\boxtimes$.
The map $m_{\phi}$ is ind-schematic and proper, so that the natural map $m_{\phi !} w^! \r w^! \ol m_{\phi !}$ is an isomorphism.
Thus, the underlying motive of $\ol m_{\phi !} \ol p^! (\calF_1 \boxtimes \calF_2)$ in $\DMrx(\Gr_{G,I})$ is given by $m_{\phi !} \tilde p^! (v^! (\calF_1 \boxtimes \calF_2))$.
By descent for Zariski torsors, $\ol p^! (\calF_1 \boxtimes \calF_2)$ is the unique object in $\DMrx(L^+_I G \backslash  \widetilde{\Gr}_{G, \phi} )$ whose image under $q^!w^! $ in $\DMrx(L_{I_1I_2} G)$ is equivariantly isomorphic $p^! (u \circ v)^! (\calF_1 \boxtimes \calF_2)$.
Again by definition, $u^!$ and $v^!$ are compatible with $\boxtimes$.
In other words, the motive $\tilde p^! (v^!(\calF_1 \boxtimes \calF_2))$ is a twisted external product in the sense of \cite[A.1.2]{Zhu:Introduction}.
The compatibility with \cite[Eqn. (5.6)]{MirkovicVilonen:Geometric} then follows from the fact that Betti realization commutes with the six functors for compact motives \cite{Ayoub:Note}.
\xpf

\rema
There is an analogue of $L_{I_1I_2}G$ for $n$ factors $I_1$, \ldots, $I_n$, cf.~\cite[Definition~3.11]{Richarz:New}. \thref{Twisted.Prod.Global} and \thref{TypeII.Relation.iso}  below generalize to $n$-fold convolution products.
\xrema

\subsubsection{One further convolution Grassmannian}
\label{sect--onemoreConv} 
We will need one further object similar to $L_{I_1I_2}G$ in order to show admissibility of a certain stratification of $\Gr_{G,I}$. Specifically, for $n > 1$ there is an \'{e}tale $L_{\{n\}}^+G$-torsor
$$\mathbf{E} \r \widetilde{\Gr}_{G, \{1, \ldots, n-1\}} \times X$$
whose functor of points records $(x, \mathcal{E}_i, \beta_i)_{i=1, \ldots, n-1} \in \widetilde{\Gr}_{G, \{1, \ldots, n-1\}}(R)$, $x_n \in X(R)$, and a trivialization of $\mathcal{E}_{n-1}$ on $\hat{\Gamma}_{x_n}$. There is an $L_{\{n\}}^+ G$-action on $\mathbf E$ given by changing the trivialization over $\hat{\Gamma}_{x_n}$.
Then $\widetilde{\Gr}_{G,I} \cong \mathbf{E} \times_{X}^{L_{\{n\}}^+G} \Gr_{G,\{n\}}$, where $I = \{1, \ldots, n\}$ and $L_{\{n\}}^+G$ acts diagonally. 
Given an $R$-point of $\mathbf{E}$, the last bundle $\mathcal{E}_{n-1}$ is trivial on $X_R - \Gamma_{x_1 \cup \cdots \cup x_{n-1}}$, and therefore also on $X_R - \Gamma_{x_1 \cup \cdots \cup x_{n-1} \cup x_n}$. 
Hence $\mathcal{E}_{n-1}$ becomes trivial on $\hat{\Gamma}_{x_1 \cup \cdots \cup x_n}$ after passing to some Zariski cover of $R$ by \thref{BD.Zar}, so $\mathbf{E} \r \widetilde{\Gr}_{G, \{1, \ldots, n-1\}} \times X$ is Zariski-locally trivial and $\mathbf{E}$ is represented by an ind-scheme over $X^I$.
Iterating this procedure shows that $\widetilde{\Gr}_{G,I}$ can be written as a twisted product 
$$\widetilde{\Gr}_{G,I} \cong \Gr_{G,\{1\}} \tilde{\times}  \cdots \tilde{\times} \Gr_{G,\{n\}}. \eqlabel{Twisted.Prod}$$
See \cite[Eq. (3.1.22)]{Zhu:Introduction} for more details. 

\subsection{Stratifications of Beilinson--Drinfeld Grassmannians}

In this section we show that a certain stratification of the Beilinson--Drinfeld affine Grassmannians is Whitney--Tate.

\defi \thlabel{BD.strata}
For a surjection of nonempty finite sets $\phi \colon I \twoheadrightarrow J$ and $\mu = (\mu_j) \in (X_*(T)^+)^J$, the corresponding stratum of $\Gr_{G,I}$ is
$$\Gr_{G,I}^{\phi, \mu} := \prod_{j \in J} \Gr_G^{\mu_j} \times X^\phi \subset \Gr_{G,I},$$ where the inclusion is induced by the factorization isomorphism \refeq{ZhuFactor}. Let
$$\iota \colon \Gr_{G,I}^\dagger \rightarrow  \Gr_{G,I} \eqlabel{stratification.BD.Gr}$$ be the inclusion of the disjoint union of the strata. 
\xdefi

\exam
\thlabel{All.Zero} 
As a preparation for the proof below, we consider the case of the trivial group $G$, where $\Gr_{G,I} \cong X^I = \A^I$. 
For $\phi \colon I \twoheadrightarrow J$, let $$j_\phi \colon X^\phi \rightarrow X^I \eqlabel{stratification.XI}$$ 
be the inclusion of the corresponding stratum. For example, 
if $I = \{1,2\}$ and $\phi = \id$, we have $X^{\circ} = X^2 \setminus \Delta_X$, the complement of the diagonal. In general, the closure of $X^\phi$ is the diagonal
$$\overline{X^\phi} := \{(x_i) \in X^I \: : \: x_i = x_i' \text{ if } \phi(i) = \phi(i')\}.$$
In particular, this is smooth and it is a union of strata, so that we have a universal Whitney--Tate stratification, cf.~\thref{Whitney--Tate.smooth}.
\xexam

\lemm
The closure of $\Gr_{G, I}^{\phi, \mu} $ is a union of strata.
\xlemm

\pf 
The closure is supported over $\overline{X^\phi}$, so by using \refeq{ZhuFactor}, we reduce to the case where $\phi$ is a bijection. 
Thus we can assume our initial stratum is of the form $\Gr_{G, I}^{\circ, \mu}$. In this case, using \refeq{Twisted.Prod} we can form the closed subscheme $$Y:= (\Gr_{G}^{\leq \mu_1} \times X) \: \tilde{\times} \cdots \tilde{\times} \: (\Gr_{G}^{\leq \mu_n} \times X) \subset \widetilde{\Gr}_{G,I}.$$ Note that $m_I$ (cf.~\refeq{convolution.morphism}) is an isomorphism over $X^\circ$, and that by computing smooth-locally it follows that $Y$ is the closure of $\Gr_{G, I}^{\phi, \mu}$ in $\widetilde{\Gr}_{G,I}$. Hence $m_I(Y) \subseteq \overline{\Gr_{G, I}^{\circ, \mu}} \subset \Gr_{G,I}$. As $m_I(Y)$ is closed
and it contains $\Gr_{G, I}^{\circ, \mu}$, this containment is an equality. 

It remains to see that $m_I(Y)$ is a union of strata. For this, let $\phi' \colon I \twoheadrightarrow K$ be an arbitrary surjection with $K$ nonempty and put $\lambda_k = \sum_{i \in (\phi')^{-1}(k)} \mu_i$ for $k \in K$. Using that the image of a local convolution morphism $\Gr_G^{\leq \mu_1} \tilde{\times} \cdots \tilde{\times} \Gr_G^{\leq \mu_n} \rightarrow \Gr_G$ is $\Gr_G^{\leq \mu_1+ \cdots + \mu_n}$,  it follows from factorization that
$$\restr{m_I(Y)}{X^{\phi'}} \cong \Gr_G^{\leq \lambda_1} \times \cdots \times \Gr_G^{\leq \lambda_{|K|}} \times X^{{\phi'}}.$$ 
\xpf

\theo \thlabel{BD.WT}
The stratification of $\Gr_{G,I}$ in \refeq{stratification.BD.Gr} is universally admissible Whitney--Tate. 
\xtheo
The key idea for the proof is to interpret restrictions along partial diagonals in \(\Gr_{G,I}\) as convolution. Then we can apply \thref{DTM.convolution}.
\pf
Let $\iota^{\phi, \mu} \colon \Gr_{G, I}^{\phi, \mu} \rightarrow \Gr_{G, I}$ be the inclusion of a stratum as in \thref{BD.strata}, where $\phi \colon I \twoheadrightarrow J$. Then we have to prove $\iota^* \circ \iota^{\phi, \mu}_*( \Z) \in \DTMrx(\Gr_{G, I}^\dagger).$

Fix a bijection $I \xrightarrow{\sim} \{1, \ldots, n\}$. We will induct on $n$. If $n=1$, the proposition follows from smooth base change applied to the projection $\Gr_G \times X \rightarrow \Gr_G$ and the universal Whitney--Tate property of $\Gr_G$ (\thref{Fl.universally.WT}). 
In fact, the Iwahori-orbit stratification on $\Gr_G$ by affine spaces is Whitney--Tate, so we get a cellular Whitney--Tate stratification in the case $n = 1$, also using that $X = \A^1$ is cellular itself.

Now assume $n > 1$. 
By \thref{admissible.SNC.divisor}, the schemes $X^\phi$ and therefore also the strata $\Gr_{G,I}^{\phi, \mu}$ are admissible $S$-schemes in the sense of \thref{admissible.stratification}.
Thus, it remains to show the stratification is universally Whitney--Tate.
If $\phi$ is not injective, there exist $i$, $i' \in I$ such that $\overline{X^{\phi}}$ is contained in the diagonal $$X^{x_i=x_{i'}} = \{(x_i) \in X^I \: : \: x_i = x_{i'} \}.$$ 
This diagonal is a union of strata, so under the obvious identification $X^{x_i=x_{i'}} \cong X^{I - \{i'\}}$ and the factorization property \refeq{ZhuFactor} this case is covered by induction.
 
Now assume $\phi$ is injective. We can assume $I=J$ and $\phi = \id$. Write $X^\circ := X^\id$. To ease notation we let $\phi \colon I \r K$ be a new surjection, and we will compute the fiber of $\iota^{\circ, \mu}_* \Z$ over $X^{\phi}$. We consider two further cases.

If there exists $i \in I$ such that $\phi^{-1}\phi(i)$ is a singleton, let $X_i^\circ \subset X^{I - \{i\}}$ be the locus with pairwise distinct coordinates. 
Consider the open subset \(X^{(\phi)}\subseteq X^I\) from \eqref{Open.Sub}, for which $X^\circ$, $X^{\phi} \subset X^{(\phi)}$. 
Then 
$$\restr{\Gr_{G,I}}{X^\circ} \cong (\restr{\Gr_{G, I-\{i\}}}{X_i^\circ } \times \Gr_{G,\{i\}}) \times_{X^I} X^{(\phi)}$$ (note that the right hand side already lives over \(X^\circ\)) and
$$\restr{\Gr_{G,I}}{X^{(\phi)}} \cong (\Gr_{G,I-\{i\}} \times \Gr_{G,\{i\}}) \times_{X^I} X^{(\phi)}.$$

Hence this case follows by induction applied to the factor $\Gr_{G, I-\{i\}}$ and \thref{WT.prod}, along with universality of the Whitney--Tate stratification of $\Gr_{G,\{i\}}$.

For the final case, suppose the fibers of $\phi$ have at least two elements each. Write $j^\circ \colon X^\circ \rightarrow X$ for the embedding.
Let $\mu = (\mu_i) \in (X_*(T)^+)^I$, and let $j^{\mu_i} \colon \Gr_G^{\mu_i} \rightarrow \Gr_G$ be the inclusion of the corresponding $L^+G$-orbit.
After possibly relabelling $I$ we may assume $\phi$ is order-preserving for the induced order on $K$.

The proof is based on the following diagram (where $\mathbf{E}$ was defined in \refsect{onemoreConv}): 
$$\xymatrix{ 
& \prod_{i \in I} \Gr_G \times X^\circ \times L^+G \ar[d]^{q^\circ} \ar[r]^-{j^2} & 
\mathbf{E} \times_X \Gr_{G,\{n\}} \ar[d]^q & \ar[l]_-{i_2} Y_{<m} \times Y_m \ar[d]^{q^{\phi}} \\
& \prod_{i \in I} \Gr_G \times X^\circ \ar[d]^{\sim} \ar[r]^-{j^1} & \widetilde{\Gr}_{G,I} \ar[d]^{m_I} & \ar[l]_-{i_1} \prod_{k \in K} \Gr_G^{I_k, \text{conv}} \times X^{\phi} \ar[d]^{m_I^{\phi}} \\ 
\prod_{i \in I} \Gr_G^{\mu_i} \x X^\circ \ar[r]^-{j^{\mu_i}} & \prod_{i \in I} \Gr_G \times X^\circ \ar[d] \ar[r]^-{j^0} & \Gr_{G,I} \ar[d] & \ar[l]_-{i_0} \prod_{k \in K} \Gr_G \times X^{\phi} \ar[d] \\
& X^\circ \ar[r]^-{j^\circ} & X & X^{\phi} \ar[l]_-{j_{\phi}}.
}$$
The bottom squares are cartesian by the factorization property of $\Gr_G$. The remaining squares, all of which are also cartesian, will be discussed below.

It suffices to show
$$i_0^* \circ j_*^0 \left(\underset{i \in I}\boxtimes j^{\mu_i}_* \Z \boxtimes \restr{\Z}{X^\circ}\right) \cong \underset{k \in K} \boxtimes \left( \underset{i \in I_k}{\star} j^{\mu_i}_* \Z \right) \boxtimes \left ( j_{\phi}^* j^\circ_*  (\Z ) \right ). \eqlabel{BD.main.eqn}$$ 
Indeed, the convolution product $\underset{i \in I_k}{\star} j^{\mu_i}_* \Z$ is an object of $\DTMrx(\Gr_G)$ by \thref{DTM.convolution}.
Furthermore, this convolution product over an arbitrary base $S$ is obtained by pullback from $S = \Spec \Z$ by proper base change and \thref{WT.prod}. 
For the second factor, we note that $j_{\phi}^* j^\circ_* (\Z ) \in \DTMrx(X^{\phi})$ by \thref{All.Zero}.
Again, this is independent of the base scheme $S$ because the computation behind \thref{Whitney--Tate.smooth} only uses relative purity.

We will prove this formula by induction on the number of cocharacters such that $\mu_i=0$, starting with the case where all $\mu_i = 0$. In this case the above formula holds because the relevant geometry is supported on the image of the trivial section $X^I \rightarrow \Gr_{G,I}$. 

Now suppose $\mu_i \neq 0$ for some $i$. After possibly relabelling we can assume $\mu_n \neq 0$.
The middle part of the above diagram is cartesian, where $I_k = \phi^{-1}(k)$ and we write, for a nonempty finite ordered set $J$,
$$\Gr_G^{J, \text{conv}} := \underbrace{LG \times^{L^+G} \cdots \times^{L^+G} LG \times^{L^+G} \Gr_G}_{|J| \text{ factors}}.$$ 
The morphism $m_I^{\phi}$ is a product of $|K|$ local convolution morphisms, times the identity morphism on $X^{\phi}$. 
The map $q$ at the top is the $L^+_{\{n\}}G$-torsor introduced in \refeq{Twisted.Prod}.
Let $m \in K$ be the largest element. 
The fiber of $\mathbf E \x_X \Gr_{G,\{n\}}$ over $X^{\phi}$ is the product of 
$$Y_{<m}:= \prod_{k < m} \Gr_G^{I_k, \text{conv}} \times X^{\phi}$$ and
$$
Y_{m} :=  \underbrace{LG \times^{L^+G} \cdots \times^{L^+G} LG \times \Gr_G}_{|I_m| \text{ factors}
}.
$$
Here we set $Y_{<m} = S$ if $|K|=1$. The map $q^\circ$ is a trivial $L^+G$-torsor, cf.~\cite[(1.7.5) ff.]{BaumannRiche:Satake}.
The map $q^{\phi}$ is the quotient by the diagonal action of $L^+G$ on the last two factors $LG \times \Gr_G$ in $Y_m$, which is possible by our assumption that \(|I_m|>1\). The point is that in the top row we have split off the factor $\Gr_G$ which supports $j^{\mu_n}_* \Z$.

By proper base change, $$i_0^* \circ j_*^0 \left(\underset{i \in I}\boxtimes j^{\mu_i}_* \Z \boxtimes \restr{\Z}{X^\circ}\right) \cong m_{I*}^{\phi} \circ i_1^* \circ j_*^1 \left(\underset{i \in I}\boxtimes j^{\mu_i}_* \Z \boxtimes \restr{\Z}{X^\circ}\right).$$ Let
$$\mathcal{E} = \underset{1 \leq i \leq n-1}{\boxtimes} j_*^{\mu_i} \Z \boxtimes \restr{\Z}{\Gr_G^{0}} \boxtimes \restr{\Z}{X^\circ},$$ i.e., the object we would like to push-pull but where we set the last cocharacter $\mu_n = 0$.
By smooth base change and \thref{base.change.Fl}
we have
$$q^{\phi*} \circ i_1^* \circ  j^1_{*} (\boxtimes_{i \in I} j^{\mu_i}_* \Z \boxtimes \Z) \cong q^{\phi*} \circ i_1^* \circ  j^1_{*} (\mathcal{E})  \boxtimes j^{\mu_n}_*\Z,$$ cf.~\cite[Prop. 7.4(ii)]{Zhu:Coherence}. Here we view $q^{\phi*} \circ i_1^* \circ  j^1_{*} (\mathcal{E})$ as supported on the closed subscheme of $Y_{< m} \times Y_m$ obtained by replacing $\Gr_G$ with its basepoint in the last factor of $Y_m$, so that the external product makes sense.

 This isomorphism is equivariant for the diagonal action of $L^+G$, so by descent we have 
$$i_1^* \circ  j^1_{*} (\boxtimes_{i \in I} j^{\mu_i}_* \Z \boxtimes \Z) \cong i_1^* \circ  j^1_{*} (\mathcal{E}) \widetilde \boxtimes j^{\mu_n}_*\Z.$$
Here the twisted product is formed on
$$ \prod_{k \in K} \Gr_G^{I_k, \text{conv}} \times X^{\phi} \cong \prod_{k<m} \Gr_G^{I_k, \text{conv}} \times \Gr_G^{I_m - \{n\}, \text{conv}} \times X^{\phi} \tilde{\times} \Gr_G,$$ which comes from the identification (obtained by restricting \eqref{eqn--Twisted.Prod} to the diagonal) 
$$\Gr_G^{I_m, \text{conv}} \cong \Gr_G^{I_m - \{n\}, \text{conv}} \tilde{\times} \Gr_G.$$
The construction of the twisted product of motivic sheaves is analogous to \thref{Twisted.Prod.Global} and is purely local, cf.~\cite[Lemma 3.11 ff.]{RicharzScholbach:Motivic}.

Now we factor $m_I^{\phi}$ by first convolving the left $n-1$ factors, and then convolving with the final factor supporting $j^{\mu_n}_* \Z$, as below. 
$$\xymatrix{ \prod_{k \in K}  \Gr_G^{I_k, \text{conv}} \times X^{\phi}  \ar[r]^-{m_1} \ar[rd]^{m_I^{\phi}} & \prod_{k < m}  \Gr_G \times (\Gr_G \tilde{\times} \Gr_G) \times X^{\phi} \ar[d]^{m_2} \\ & \prod_{k \in K} \Gr_G \times X^{\phi}.}$$
By proper base change and \thref{base.change.Fl}, $$m_{1*} ( i_1^* \circ  j^1_{*} (\mathcal{E}) \widetilde \boxtimes j^{\mu_n}_*\Z) \cong (m_{I*}^{\phi} \circ i_1^* \circ j^1_* (\mathcal{E}))  \widetilde \boxtimes  j^{\mu_n}_*\Z.$$ By induction the above is isomorphic to
$$\underset{k <m} \boxtimes \left( \underset{i \in I_k}{\star} j^{\mu_i}_* \Z\right) \boxtimes \left(\left(\underset{i \in I_m -\{n\}}{\star} j^{\mu_i}_{*} \Z \right) \widetilde \boxtimes  j^{\mu_n}_* \Z  \right) \boxtimes j_{\phi}^* j^\circ_*  (\Z ).$$ Applying the projection formula to $m_{2*} = m_{2!}$, we deduce \refeq{BD.main.eqn}, completing the proof.
\xpf

\rema
The proof of \thref{BD.WT} implies that $\underset{i \in I} \star j^{\mu_i}_* \Z$ is independent of the ordering on $I$, modulo the factor $ j_{\phi}^* j^\circ_*  (\Z )$. This is similar to Gaitsgory's construction of the commutativity constraint \cite{Gaitsgory:Central, Gaitsgory:Braiding}, where this factor disappears when one uses nearby cycles; see also \cite{Zhu:Ramified}.
\xrema

\thref{BD.WT} entitles us to the following definitions, cf.~\thref{Whitney--Tate} and \thref{equivariant.DTM}.

\defilemm
\thlabel{Defi:DTMHck}
We denote by $\DTMrx(\Gr_{G,I})$ the category of \emph{stratified Tate motives} on the Beilinson--Drinfeld Grassmannian. 
The category $\DTMrx(\Hck_{G,I})$ of \emph{Tate motives on the Hecke prestack} is defined to be the full subcategory of $\DMrx(\Hck_{G,I})$ of objects whose underlying motive in $\DMrx(\Gr_{G,I})$ lies in $\DTMrx(\Gr_{G,I})$. 

If the base scheme $S$ satisfies the condition around \refeq{BS.vanishing}, both categories carry a natural t-structure such that the forgetful functor $u^! : \DTMrx(\Hck_{G,I}) \r \DTMrx(\Gr_{G,I})$ is t-exact.
The heart of these t-structures is denoted by $\MTMrx$. We refer to these motives as \emph{mixed (stratified) Tate motives}.
\xdefilemm

\pf 
As noted in the proof of \thref{BD.WT}, it follows from \thref{admissible.SNC.divisor} that the Whitney--Tate stratification of $X^I$ by the $X^\phi$ is admissible. Hence \thref{t-structure.stratified} applies to the stratified map $X^I \r S$, where $S$ has the trivial stratification, so that we have a t-structure on $\DTM(X^I)$ glued from t-structures on the strata. By factorization \refeq{ZhuFactor}, the fiber of the stratified map $\pi_G \colon \Gr_{G,I} \r X^I$ over a stratum $X^\phi$ is a product of copies of $\Gr_G$ over $X^\phi$. Since the stratification of $\Gr_G$ by $L^+G$-orbits is cellular \cite[Corollary 4.3.12]{RicharzScholbach:Intersection}, it follows by \thref{cellular.stratification} that $\pi_G$ is admissibly stratified. Hence we can now apply \thref{t-structure.stratified} to $\pi_G$, giving the t-structure on $\DTMrx(\Gr_{G,I})$.

In order to construct the t-structure on $\DTMrx(\Hck_{G,I})$, we need to verify the remaining assumptions of \thref{MTM.equivariant.ind-scheme} for the action of $L^+_IG$ on $\Gr_{G,I}$. By \cite[Example A.4.12]{RicharzScholbach:Intersection}, the pro-smooth group scheme $L^+_IG$ is a particular case of the pro-algebraic groups considered in \cite[Appendix A.4]{RicharzScholbach:Intersection} (see also \cite[Example 3.1(iv), Lemma 3.2]{HainesRicharz:TestFunctionsWeil}). In particular, $\Gr_{G,I}$ is a colimit of $L^+_IG$-stable closed $X^I$-subschemes of finite type, on each of which the action of $L^+_IG$ factors through a finite type jet quotient with split pro-unipotent kernel \cite[Lemma A.3.5, Proposition A.4.9]{RicharzScholbach:Intersection}. It is not difficult to check that under the factorization isomorphism \refeq{ZhuGroupFactor}, the fiber of a jet quotient of $L_I^+G$ is a product of jet quotients of $L^+G$, cf.~\cite[Example A.4.12]{RicharzScholbach:Intersection}. Now by also using the factorization isomorphism \refeq{ZhuFactor}, the remaining assumptions of \thref{MTM.equivariant.ind-scheme}, which can be checked over strata of $X^I$, have already been verified in the proof of \thref{Orbit.Eq}.
\xpf

\coro
\thlabel{Gr.I.Tate.map}
The natural map \(\pi_G\colon \Gr_{G,I}\to X^I\) is Whitney--Tate, i.e., $\pi_{G*} = \pi_{G!}$ preserves Tate motives, with respect to the stratification of $X^I$ in \thref{All.Zero}.
\xcoro

\pf 
By ind-properness of $\pi_G$ (\thref{lemm.BDrep}), we have $\pi_{G*} = \pi_{G!}$.
It then suffices to show that $\pi_{G!}$
maps the generators of $\DTMrx(\Gr_{G,I})$, namely $\iota_!^{\phi, \mu} \Z$, to an object in $\DTMrx(X^I)$, for $\phi : I \twoheadrightarrow J$ and $\mu = (\mu_j) \in (X_*(T)^+)^J$.
Using the factorization property \refeq{ZhuFactor}, we have
$$\pi_{G!} \iota_!^{\phi, \mu} \Z \cong j_{\phi !} \left (\underset{j \in J} \boxtimes \restr{\pi^{(j)}_! \Z}{X^\phi} \right),$$
where $j_\phi \colon X^\phi \r X^I$ is as in \refeq{stratification.XI}, and $\pi^{(j)} : \Gr_G^{\mu_j} \x_S X \r X$ is the projection.
The stratification of $\Gr_G^{\mu_j}$ by Iwahori-orbits is a stratification by affine spaces, so that pushforward along the structural map $\Gr_G^{\mu_i} \r S$ preserves Tate motives \cite[Lemma~3.1.19]{RicharzScholbach:Intersection}.
Thus, the above expression is an object of $\DTMrx(X^I)$.
\xpf

\subsection{Tate motives on Beilinson--Drinfeld Grassmannians}
\nota
When working on \(\Gr_{G,I}\), for some finite index set \(I\), we will often need to shift by \(|I|\). We will denote \([|I|]\) and \([-|I|]\) by \([I]\) and \([-I]\), and similarly for Tate twists.
\xnota
\subsubsection{Convolution}
In this subsection, we prove that the various convolution functors preserve Tate motives, and we relate Type I and II convolution products of Tate motives. 

In the context of the Type II convolution diagram, consider the surjection $I \sqcup I \r \{1,2\}$ sending the first copy of $I$ to $1$, and the second copy to $2$. The fiber of part of \refeq{TypeII} over the diagonal embedding $X^{I} \r X^{I \sqcup I}$ is the global version of the classical convolution diagram \cite[(Eqn. (4.1)]{MirkovicVilonen:Geometric},
\begin{equation}\begin{tikzcd} \label{TypeI}
  \Gr_{G,I} \times_{X^{I}} \Gr_{G,I} & \arrow{l} \arrow{r} L_{I}G \times_{X^{I}} \Gr_{G,I}  & (L_{I}G \times_{X^{I}}^{L^+_{I}G} \Gr_{G,I})_{\Zar} \arrow{r}{m} & \Gr_{G,I}.
\end{tikzcd}\end{equation}
The left map is the quotient on the left factor and the right map is the quotient by the diagonal action of $L_{I}^+G$.
Denoting the base of a box product with a subscript, for $\calF_1, \calF_2 \in \DMrx(\Hck_{G,I})$,  we can form the twisted product $\calF_1 \widetilde \boxtimes_{X^{I}} \calF_2 \in \DMrx(L_{I}G \times_{X^{I}}^{L^+_{I}G} \Gr_{G,{I}})$ by applying the same arguments as in \thref{Twisted.Prod.Global}.
Let $f_{I} \colon \Gr_{G,I} \r \Hck_{G,I}$ be the quotient map.
Base change along the quotient of $m$ by the left action of $L_{I}^+G$ shows that $$m_!(\calF_1 \widetilde \boxtimes_{X^{I}} \calF_2) \cong f_{I}^!(\calF_1 \star \calF_2)[I]. \eqlabel{TypeI.Relation}$$ 
This result generalizes to an $n$-fold convolution product.

\prop
\thlabel{Convolution:Global:Tate} The convolution product $\star$ in \refeq{GrI.Conv} preserves Tate motives, i.e., it restricts to a functor
$$\star : \DTMrx(\Hck_{G,I}) \x \cdots \x \DTMrx(\Hck_{G,I} ) \r \DTMrx(\Hck_{G,I} ).$$ Likewise, for a surjection \(\phi \colon I\twoheadrightarrow J\) and $I_j =\phi^{-1}(j)$, the convolution product $m_{\phi!} (- \widetilde \boxtimes \cdots \widetilde \boxtimes - )$ from \thref{Twisted.Prod.Global} (for any number of factors) restricts to a functor
$$\DTMrx(\Hck_{G, I_1}) \times \cdots \times \DTMrx(\Hck_{G, I_{|J|}}) \r \DTMrx(\Gr_{G, I}).$$
\xprop
\pf
By continuity we may restrict to bounded objects (\thref{DM.ind-scheme}).
Then we may replace the torsors used to construct the twisted products in both types of convolution by finite type quotients. 
In this case, the twisted products can be formed using either $!$- or $*$-pullback. 
The fibers of \eqref{TypeI} and \refeq{TypeII} over the strata of $X^I$ are products of local convolution diagrams.
Hence, by base change and the compatibility between box products and the two operations of $*$-pullback and $!$-pushforward \cite[Theorem 2.4.6]{JinYang:Kuenneth}, the claim follows from the local case (\thref{DTM.convolution}).
\xpf

As was recalled in \refsect{functoriality}, the functor $Z \mapsto \DM(Z)$ is a lax symmetric monoidal functor out of the category of correspondences (on finite type $S$-schemes). This implies the existence of maps $\DM(Z) \t \DM(Z') \r \DM(Z \x_S Z')$, which moreover are compatible with *-pullbacks and !-pushforwards (and therefore also !-pullbacks along smooth maps and *-pushforwards along proper maps). 
That latter compatibility implies by adjunction the existence of maps
\begin{equation}
  f^! \Z \t f^* M \r f^! M.\label{f! f* tensor}
\end{equation}
Also recall the (formal) extension of these constructions to placid prestacks from \thref{placid prestacks}.

\prop \thlabel{TypeII.Relation}
Let $$i_1 \colon \Gr_{G, I} \r \Gr_{G,I \sqcup I}, \quad i \colon \Hck_{G,I} \times_{X^{I}} \Hck_{G,I} \r \Hck_{G,I} \times_S \Hck_{G,I}$$ be the diagonal embeddings. 
For $\calF_1, \calF_2 \in \DTMrx(\Hck_{G,I})$, there is a canonical map
\begin{equation}
f_{I}^!(\calF_1 \star \calF_2) (-I)[-I] \r i_1^!m_{\phi !}(\calF_1 \widetilde \boxtimes_S \calF_2 ). \label{TypeII.Relation.map}
\end{equation}
\xprop

\rema
The map \eqref{f! f* tensor} and also the map \eqref{TypeII.Relation.map} are not in general isomorphisms.
However, \eqref{f! f* tensor} is an isomorphism if $M = \Z(k)$ is twisted-constant (more generally, if $M$ is dualizable). We will use this in \thref{TypeII.Relation.iso} to provide an instance where \eqref{TypeII.Relation.map} is in fact an isomorphism. 
\xrema

\pf
Note that the definition of $\star$ includes a shift by $-|I|$. 
By base change and \refeq{TypeI.Relation}, it suffices to map $\calF_1 \widetilde \boxtimes_{X^{I}} \calF_2 (-I)[-2|I|]$ to the corestriction of $\calF_1 \widetilde \boxtimes_S \calF_2$ to the diagonal. 
This amounts to mapping $\calF_1 \boxtimes_{X^{I}} \calF_2 (-I)[-2|I|]$ to $i^!(\calF_1 \boxtimes_{S} \calF_2)$. 
To get a map 
$$(\calF_1 \boxtimes_{X^{I}} \calF_2) (-I)[-2|I|] \r i^!(\calF_1 \boxtimes_{S} \calF_2),$$
rewrite $i$ as $i' \times \id \colon \Hck_{G,I} \times_{X^{I}} \Hck_{G,I} \r (\Hck_{G,I} \times_S X^{I}) \times_{X^{I}} \Hck_{G,I}$. 
Here $i'$ is the product of the identity map of $\Hck_{G, I}$ and the structure map to $X^{I}$. 
The map we seek is the canonical map in \cite[Theorem 2.4.6]{JinYang:Kuenneth},
$$i'^!(\calF_1 \boxtimes_S \restr{\Z}{X^{I}}) \boxtimes_{X^{I}}\calF_2 \r (i' \times \id)^! (\calF_1 \boxtimes_S \restr{\Z}{X^{I}} \boxtimes_{X^{I}} \calF_2).\eqlabel{dumdideldum}$$ 
Indeed, by using relative purity to rewrite $\calF_1 \boxtimes_S \restr{\Z}{X^{I}}$ as a $!$-pullback from $\Hck_{G,I}$, we have $i'^!(\calF_1 \boxtimes_S \restr{\Z}{X^{I}}) \cong \calF_1 (-I) [-2|I|]$.
\xpf

\subsubsection{Independence of the base}
The following statement, which is false for non-reduced motives, allows to connect Satake categories over various base schemes.

\lemm
\thlabel{independence.star}
Let $\pi : S' \r S$ be a scheme over $S$.
Let $G' = G \x_S S'$, and denote by $\Gr_{G', I}$ the associated Beilinson--Drinfeld affine Grassmannian over $S'$. 
Then the functor
$$\pi^* : \DTMrx(\Hck_{G, I}) \r \DTMrx(\Hck_{G', I})$$
is a monoidal functor with respect to the convolution product $\star$ and also with respect to $m_{\phi} (- \widetilde \boxtimes \ldots \widetilde \boxtimes - )$.
For reduced motives (but not for non-reduced ones), this functor is an equivalence, where reduced motives are taken with respect to the respective base schemes, i.e.~$S$ for $\Hck_{G, I}$ and $S'$ for $\Hck_{G', I}$. 
\xlemm

\pf
The functor $\pi^*$ exists since $\Gr_{G', I} \r \Gr_{G, I}$ is schematic.
Up to taking reduced subschemes, the stratification on $\Gr_{G', I}$ is just the preimage stratification of the one on $\Gr_{G,I}$.
This follows from \cite[Proposition 4.4.3]{RicharzScholbach:Intersection}.
Thus $\pi^*$ preserves Tate motives.
The functor $\pi^*$ is compatible with these convolution functors since it commutes with $\boxtimes$, the !-pullback functors along pro-smooth maps (the maps $p$ in \refeq{Loc.Conv}, resp.~in \refeq{TypeII}) and !-pushforward (along the maps $m$, resp.~$m_\phi$).
The functor $\pi^*$ is an equivalence for reduced motives by \thref{DTMr.independence}, which is applicable by \thref{BD.WT}.
\xpf

\subsubsection{Forgetting the equivariance}
The following result (\thref{FF.Stab}) will be used implicitly in several places, cf.~\thref{Sub.Quot}; the proceeding result (\thref{lemm-XI-stable}) will be used in the proof of \thref{Sat.Props}.

\prop \thlabel{FF.Stab}
The pullback functor $u^! \colon \MTMrx(\Hck_{G,I}) \rightarrow \MTMrx(\Gr_{G,I})$ is fully faithful, and the image is stable under subquotients.
\xprop

\pf
Since $u^!$ preserves colimits it suffices to prove this for the respective subcategories of bounded objects.
Being supported on a finite type $L^+_IG$-stable closed subscheme $Y \subset \Gr_{G,I}$, they admit a finite filtration by IC motives of mixed Tate motives on the strata.
The action of $L_I^+G$ on $Y$ factors through a smooth quotient  $L^+_IG \twoheadrightarrow H$ such that the kernel of this quotient map is split pro-unipotent by \cite[Lemma A.3.5, Proposition A.4.9]{RicharzScholbach:Intersection}. 
By \cite[Proposition 3.1.27]{RicharzScholbach:Intersection}, we reduce to considering $H$-equivariant motives on $Y$.

Let $a, p \colon H \times_{X^I} Y \r Y$ be the action and projection maps.
Since $H$ is smooth and the fiber of $H$ over each stratum $X^\phi \subset X^I$ is cellular, the preimage stratification on $H \times_{X^I} Y$ is admissible. By smooth base change, we therefore get an exact functor $p^![-d] = p^*(d)[d] \colon \MTMrx(Y) \r \MTMrx(H \times_{X^I} Y)$, where $H$ has relative dimension $d$. The restriction of $p^![-d]$ to the fiber over each stratum $X^\phi \subset X^I$ is fully faithful by \cite[Proposition 3.2.12]{RicharzScholbach:Intersection}. The functor $p^![-d]$ also preserves IC motives. Since homomorphisms between IC motives associated to the same stratum are determined by their restriction to the stratum, and there are no nonzero homomorphisms between IC motives associated to different strata, this implies $p^![-d]$ is fully faithful when restricted to IC motives. Full faithfulness in general then follows by induction on the lengths of filtrations by IC motives as in \cite[Lemma 3.4]{Cass:Perverse}. By the proof of \cite[Proposition 3.2.20]{RicharzScholbach:Intersection} this implies that $u^!$ is fully faithful, with image consisting of mixed Tate motives such that there exists an isomorphism $a^! \calF \cong p^! \cal F$.

To show stability under subquotients, by a standard argument as in \cite[Proposition 4.2.13]{Letellier:Fourier} it suffices to show that the image of the exact functor $p^![-d] \colon \MTMrx(Y) \r \MTMrx(H \times_{X^I} Y)$ is stable under subquotients. 
Here we give $H \times_{X^I} Y$ the preimage stratification. 
The functor $p^![-d]$ has a left adjoint $\pe p_![d]$ and a right adjoint $\pe p_*(-d)[-d]$. 
Then by \cite[\S 4.2.6]{BeilinsonBernsteinDeligne:Faisceaux} it suffices to verify that for all $\calF \in \MTMrx(H \times_{X^I} Y)$, the natural map $\calF \rightarrow \pe p^! p_! \calF$ is an epimorphism. (Note that there some typos in \emph{loc.~cit.}; in particular, condition (b') should state that $B \r u^*u_!B$ is an epimorphism, so that $B$ has a maximal quotient in the image of $u^*$.) For this we may instead take $\calF \in \DTMrx^{\leq 0}(H \times_{X^I} Y)$ and we must show that the homotopy fiber of $\calF \r p^! p_! \calF$ lies in $\DTMrx^{\leq 0}(H \times_{X^I} Y)$. It suffices to verify this condition on generators of $\DTMrx^{\leq 0}(H \times_{X^I} Y)$, so we may take a stratum $j \colon H \times_{X^I} Y_w \rightarrow H \times_{X^I} Y$, and let $\calF = j_! \Z[d + \dim Y_w]$. 
By base change, we have $p^! p_! \calF \cong f^*f_!(\Z[d])(d)[2d] \boxtimes \Z[\dim Y_w]_{Y_w}$, where $f \colon \restr{H}{X^\phi} \rightarrow X^{\phi}$ and $X^\phi \subset X^I$ is the stratum lying below $Y_w$. Since $\restr{H}{X^\phi}$ is cellular we have $f^*f_!(\Z[d])(d)[2d] \in \DTMrx^{\leq 0}(\restr{H}{X^\phi})$. Furthermore, since $\restr{H}{X^\phi}$ has a unique top dimensional cell (and it is open), it follows from excision that $\pe H^0(f^*f_!(\Z[d])(d)[2d]) \cong \Z[d]$. By right exactness of $\boxtimes$, the natural map $\calF \r p^! p_! \calF$ is therefore an isomorphism after applying $\pH^0(-)$, so before applying $\pH^0(-)$ the fiber lies in $\DTMrx^{\leq 0}(H \times_{X^I} Y)$. 
\xpf

\lemm \thlabel{lemm-XI-stable}
Let $p \colon X^I \r S$ be the structure map. Then the image of the exact functor $p^![-d] \colon \MTMrx(S) \rightarrow \MTMrx(X^I)$ is stable under subquotients, where the target consists of stratified mixed Tate motives.
\xlemm

\pf
As in the proof of \thref{FF.Stab}, we must show that for $\calF \in \DTMrx^{\leq 0}(X^I)$, the homotopy fiber of the unit map $\eta: \calF \r p^!p_! \calF$ lies in $\DTMrx^{\leq 0}(X^I)$. We may check this on the generators $j_!\Z [\dim X^\phi]$, where $j \colon X^\phi \rightarrow X^I$. If $X^\phi \neq X^\circ$ then by induction on $|I|$ we have $p^!p_! \calF \in \DTMrx^{\leq -1}(X^I)$ for dimension reasons (this uses excision and the compactly supported cohomology of $\A^n$ for $n \leq |I|$), so it remains to consider $\calF = j_!\Z [I]$ where $X^\phi = X^\circ$. In this case, we argue as in \thref{FF.Stab}. 
We apply \thref{admissible.SNC.divisor} to $X^\circ \subset X^I$. Writing $q : X^\circ \r S$,
we have $\pH^0 (p^! p_! j_! \Z[I]) = \pH^0 (p^*[I] q_! q^! \Z) = p^* [I] \pH^0 (q_! q^! \Z)$. 
By the computation in \refeq{localization.SNC divisor}, this object is isomorphic to $p^*[I] p_! p^! \Z = \Z[I]$.
Therefore $\pH^0(\eta)$ is the natural surjection $\pe j_! \Z[I] \r \Z[I]$ to the IC motive.
\xpf

\section{The global Satake category}
\label{sect--Satake category}

In this section, we construct and study the Satake category.
We do this in a global situation, i.e., as certain motives on the Beilinson--Drinfeld affine Grassmannians.
Throughout this section, we fix a nonempty finite set \(I\).

\label{sect--CTsect}
\subsection{Constant terms}

Given a cocharacter \(\chi\in X_*(T)\), consider the induced conjugation action \(\Gm\times G\to G\colon (t,g)\mapsto \chi(t)\cdot g \cdot \chi(t)^{-1}\). 
The attractor and repeller for this action are opposite parabolics \(P^+\) and \(P^-\) of \(G\), and the fixed points are given by the Levi subgroup \(M=P^+\cap P^-\). 
We will often abbreviate $P := P^+$.
If \(\chi\) is dominant regular, then \(P=B\) is the Borel, \(P^-=B^-\) the opposite Borel, and \(M=T\) is the maximal torus.

Now, \(\Gm\) also acts on \(\Gr_{G,I}\) via \(\Gm\to L^+_I\Gm \xrightarrow{L^+_I\chi} L^+_IT\to L^+_IG\).
This \(\Gm\)-action is Zariski-locally linearizable if $G = \GL_n$ by the proof of \cite[Lemma 3.16]{HainesRicharz:TestFunctionsWeil}, and it follows for general groups as well since $\Gr_{G,I}$ admits a \(\Gm\)-equivariant closed embedding into some $\Gr_{\GL_n,I}$ by \cite[Proposition 3.10]{HainesRicharz:TestFunctionsWeil}.
The fixed points, attractor and repeller are given by \(\Gr_{M,I}\), \(\Gr_{P^+,I}\) and 
\(\Gr_{P^-,I}\) respectively, compatibly with the natural morphisms between them, cf.~ \cite[Theorem 3.17]{HainesRicharz:TestFunctionsWeil}. In particular, the natural projections and inclusions corresponding to these affine Grassmannians only depend on the parabolic \(P^+\), not on \(\chi\). 
We obtain the corresponding hyperbolic localization diagram as follows.
Since the top horizontal maps are $L^+_I M$-equivariant, we get the corresponding diagram of prestacks underneath:
$$\xymatrix{
\Gr_{M, I} \ar[d] & \Gr_{P^\pm, I} \ar[l]_{q^\pm_P} \ar[d] \ar[r]^{p^\pm_P} & \Gr_{G, I} \ar[d] \\
L^+_IM\backslash\Gr_{M, I} & L^+_IM\backslash\Gr_{P^\pm, I} \ar[l]_{\overline{q}^\pm_P} \ar[r]^{\overline{p}^\pm_P} & L^+_IM\backslash\Gr_{G, I}.}\eqlabel{hyperbolic.diagram}$$
  
By a proof analogous to that of \thref{GrP.LC}, the morphisms $p^{\pm}_P$ are locally closed immersions on connected components. Recall that these are indexed by \(\pi_1(M)\), and are exactly the preimages of the connected components of \(\Gr_{M,I}\). 

If $\chi$ is dominant regular, so that $P = B$ is a Borel, the connected components of \(\Gr_{B^{\pm},I}\) are denoted \(\Ss_{\nu,I}^\pm\), and called the semi-infinite orbits as for the usual affine Grassmannian.

\prop
\thlabel{Semi-infinite:Strat}
The semi-infinite orbits determine a stratification of \(\Gr_{G,I}\).
\xprop
\pf We claim that \(\overline{\Ss_{\nu,I}^+} = \bigcup_{\nu'\leq \nu} \Ss_{\nu',I}^+\). 
To prove this we can assume \(S\) is the spectrum of an algebraically closed field.
In this case, consider the usual affine Grassmannians \(\Gr_{T} \leftarrow \Gr_B \to \Gr_G\), 
and let \(\Ss_\nu^+ \subseteq \Gr_B\) be the preimage of the connected component \([\nu]\in \pi_0(\Gr_T)\). 
Taking the closure inside \(\Gr_G\), we get \(\overline{\Ss_\nu^+} = \bigcup_{\nu'\leq \nu} \Ss_{\nu'}^+\) by \cite[Proposition 5.3.6]{Zhu:Introduction}.
From this, we can immediately conclude our lemma in the case \(I=\{*\}\). 
The case of general \(I\) is a straightforward generalization of arguments in the proof of \cite[Proposition 1.8.3]{BaumannRiche:Satake} for $I=\{1,2\}$, which involve the factorization property \refeq{ZhuFactor} and the identification of the \(\Ss_{\nu,I}^+\) with the attractors for a $\Gm$-action.
\xpf

The previous proposition also shows that the semi-infinite orbits \(\Ss_{\nu,I}^-\subseteq \Gr_{B^-,I}\) for the opposite Borel determine a stratification.

\exam\thlabel{examples of intersections}
Let \((\mu_i)_i\in (X_*(T)^+)^I\). Then the restriction of the reduced intersection \(\Gr_{G, I}^{\circ,(\mu_i)_i}\cap \Ss_{(\mu_i)_i,I}^- \subseteq \Gr_{G, I}\) to \(X^\circ \subseteq X^I\) is canonically isomorphic to \(X^\circ\).
Specifically, it consists of the image of the section \(X^{\circ}\to \restr{\Gr_{G,I}}{X^\circ}\) corresponding to \((w_0(\mu_i))_i\), where \(w_0\) is the longest element of the finite Weyl group.
Indeed, it suffices to check this claim after base change to an algebraically closed field, and by \thref{prop.factorization of Gr} we may consider the local affine Grassmannian \(\Gr_G\).
In this case, the claim follows from \cite[(3.6)]{MirkovicVilonen:Geometric} (whose proof works over arbitrary algebraically closed fields).
Similarly, using \cite[Theorem 3.2]{MirkovicVilonen:Geometric}, we see that if \((\mu'_i)_i\in (X_*(T)^+)^I\) is such that the intersection \(\Gr_{G, I}^{\circ,(\mu'_i)_i}\cap \Ss_{(\mu_i)_i,I}^-\subseteq \Gr_{G, I}\) is nonempty, then \(\mu_i\leq \mu_i'\) for each \(i\in I\).
\xexam

We consider again a general parabolic.
If \((\Gr_{P^\pm,I})_\nu\) and \((\Gr_{M,I})_\nu\) are the connected components corresponding to \(\nu\in \pi_1(M)\), we denote the restriction of \(p^\pm_P\) and \(q^\pm_P\) by \[(\Gr_{M, I})_\nu \stackrel{q^\pm_\nu} \gets (\Gr_{P^\pm, I})_\nu \stackrel{p^\pm_\nu} \r \Gr_{G, I}.\eqlabel{hyperbolic.diagram.restricted}.\]
The map $q_P^\pm$ is map between ind-schemes, so that also the functor $q_{P*}$
exists \cite[Theorem~2.4.2]{RicharzScholbach:Intersection}.
The geometry of hyperbolic localization, see e.g.~\cite[Construction~2.2]{Richarz:Spaces}, induces a map
$q^-_{P*} p^{-!}_P \r q^+_{P!} p^{+*}_P$. 
Since $L^+_I M$ is pro-smooth, we obtain functors $\ol q^-_{P*}$ etc.~which are compatible with forgetting the $L^+_I M$-action (cf.~\thref{equivariant.functoriality}), and therefore a 
natural transformation \[(\overline{q}^-_P)_*(\overline{p}^-_P)^! \to (\overline{q}^+_P)_!(\overline{p}^+_P)^*\eqlabel{hyploc.nattrans}\] of functors \(\DMrx(\Hck_{G, I}) \to \DMrx(\Hck_{M,I})\). By \thref{hyperbolic.localization}, this is an equivalence after forgetting the \(L^+_IM\)-equivariance, for \(\Gm\)-monodromic objects. However, as forgetting the equivariance is conservative, and the \(\Gm\)-action factors through \(L^+_IM\), we see that \(\eqref{eqn--hyploc.nattrans}\) is already an equivalence. Hence, the following definition makes sense.

\defi
\thlabel{CT} 
Using the maximal torus quotient \(M/M_{\der}\) of \(M\), we define the degree map as the locally constant function 
\[\deg_P:\Gr_{M,I} \to \Gr_{M/M_{\der},I} \to X_*(M/M_{\der}) \xrightarrow{\langle 2\rho_G-2\rho_M,- \rangle} \Z,\]
where the middle map is given by summing the relative positions, and \(\rho_-\) indicates in which group we take the half-sum of the positive roots. 
If \(P=B\), we will usually write \(\deg:=\deg_B\).

The \emph{constant term functor} associated to \(P\) is $$\CT_P^I:=(\overline{q}^+_P)_!(\overline{p}^+_P)^*[\deg_P] \cong (\overline{q}^-_P)_*(\overline{p}^-_P)^![\deg_P]\colon \DMrx(\Hck_{G,I}) \to \DMrx(\Hck_{M,I}). \eqlabel{CT.Def}$$
Implicit in this definition is the functor \(\DMrx(\Hck_{G, I}) \to \DMrx
(L^+_IM \backslash \Gr_{G,I})\) forgetting part of the equivariance.
\xdefi

\rema
\thlabel{Remark:CTforGr}
\begin{enumerate}
  \item Since pullback and pushforward between quotient stacks are compatible with forgetting the equivariance by \thref{equivariant.functoriality}, the constant term functors satisfy
  \[u^!\CT_P^I \cong (q^+_P)_!(p^+_P)^*u^![\deg_P] \cong (q^-_P)_*(p^-_P)^!u^![\deg_P],\]
  where $u$ denotes both quotient maps \(\Gr_{G,I}\to \Hck_{G,I}\) and \(\Gr_{M,I}\to \Hck_{M,I}\).
  \item Since \(\CT_P^I\) admits a description in terms of both left and right adjoints, it preserves all limits and colimits.
\end{enumerate}
\xrema

\rema
\thlabel{CT:local}
The same discussion as above also works in the setting of usual affine Grassmannians, so that we can define constant term functors
\[\CT_P:=(\overline{q}^+_P)_!(\overline{p}^+_P)^*[\deg_P] \cong (\overline{q}^-_P)_*(\overline{p}^-_P)^![\deg_P]\colon\DMrx(L^+G\backslash \Gr_G) \to \DMrx(L^+M\backslash \Gr_M).\]
Although we will not mention this explicitly, all properties we prove for \(\CT_P^I\) also hold for \(\CT_P\).
\xrema

The following lemma can be compared to \cite[Proposition VI.7.13]{FarguesScholze:Geometrization} and the proof of \cite[5.3.29]{BeilinsonDrinfeld:Quantization}.

\lemm
\thlabel{CT.composition}
Let \(P'\subseteq P\subseteq G\) be parabolic subgroups with Levi quotients \(M'\subseteq M\), and let \(Q:=\im(P'\to M)\) be the parabolic of \(M\) with Levi quotient \(M'\). Then there is a natural equivalence \(\CT_{P'}^I\cong \CT_Q^I\circ \CT_P^I\) of functors \(\DMrx(\Hck_{G, I}) \to \DMrx(\Hck_{M',I})\).
\xlemm
\pf
This follows from base change and \thref{fiber product of affine grassmannians}.
\xpf

The following result is crucial in order to prove the t-exactness of the fiber functor.
Recall the category of (mixed) Tate motives on the Hecke prestack, cf.~ Definitions \ref{Hecke.stack} and \ref{Defi:DTMHck}.
Recall also the notion of a bounded motive on an ind-scheme $Y$ from \thref{DM.ind-scheme}.
We say \(\Ff\in \DMrx(\Hck_{G, I})\) is bounded if \(u^!\Ff\in \DMrx(\Gr_{G,I})\) is bounded.

\prop
\thlabel{CT.DTM}
For any parabolic \(P\subseteq G\) with Levi \(M\) as above, the constant term functor \(\CT_P^I\) preserves stratified Tate motives. Moreover, when restricted to bounded motives \(\Ff\), \(\CT_P^I\) also reflects Tateness, i.e., \(\Ff\in \DTMrx(\Hck_{G, I})\) if and only if \(\CT_P^I(\Ff)\in \DTMrx(\Hck_{M, I})\).
\xprop

\pf
Since $\CT_P^I$ commutes with restriction along the maps $j_\phi$ \refeq{stratification.XI}, we may replace $X^I$ by $X^\phi$.
By using factorization properties, we may assume $\phi$ is injective, and then have to consider $\boxtimes_{i \in I} \CT_P^{\{i\}}$. By the Künneth formula, i.e., the compatibility of !-pushforwards and *-pullbacks with exterior products (\refsect{motives}), and preservation of Tate motives by \(\boxtimes\), it suffices to consider the individual \(\CT_P^{\{i\}}\)'s, i.e., we may assume \(I=\sgl\).
Preservation of Tateness can be checked on bounded objects, so we may assume any \(\Ff\in \DTMrx(\Hck_{G,\sgl})\) below to be bounded.
Moreover, since constant term functors preserve bounded motives and filtered colimits, we can use \thref{CT.composition} to reduce to the case \(P=B\). 
Indeed, in the notation of \thref{CT.composition}, if \(P'\) is a Borel, then so is \(Q\).
But if \(\CT_{P'}^{\sgl}\) preserves Tateness, and \(\CT_Q^{\sgl}\) reflects Tateness (for bounded motives), then \(\CT_{P'}^{\sgl} \cong \CT_Q^{\sgl}\circ \CT_P^{\sgl}\) implies that \(\CT_P^{\sgl}\) preserves Tateness.
A similar argument will show that \(\CT_P^{\sgl}\) reflects Tateness for bounded motives.

To see that \(\CT_B^{\sgl}\) preserves Tateness, it suffices to check this for \((q^+_B)_!(p^+_B)^*\), where the arrows are as in \refeq{hyperbolic.diagram}, cf.~\thref{Remark:CTforGr}.
This is computed via pullback to the semi-infinite orbits \(\Ss_{\nu,\sgl}^+\), and then pushforward to the curve \(X\).
Now, \(\DTM(\Gr_{G, \sgl})\) is generated under colimits, extensions, and twists, by the !-pushforwards of constant motives on the Schubert cells \(\Gr_{G,\sgl}^\mu\), for \(\mu\in X_*(T)^+\).
These are thus sent to the !-pushforwards of the constant motives on the intersections \(\Gr_{G, \sgl}^\mu\cap \Ss_{\nu,\sgl}^+\), which are Tate by \thref{cellularity of intersection:torus} (and using \eqref{Gr.Pt} to see that the maps \(\Gr_{T, \sgl}\leftarrow \Gr_{B,\sgl} \to \Gr_{G, \sgl}\) arise from \(\Gr_T\leftarrow \Gr_B\to \Gr_G\) by taking the product with \(X\)).

It remains to show that \(\CT_B^{\sgl}\) reflects Tateness, so let \(\Ff\in \DMrx(\Hck_{G, \sgl})\) be bounded and assume \(\CT_B^{\sgl}(\Ff)\in \DTMrx(\Hck_{T,\sgl})\). 
For a finite subset \(W\subseteq X_*(T)^+\) closed under the Bruhat order let $i_W \colon \Gr_{G, \sgl}^{W} \rightarrow \Gr_{G, \sgl}$ be the inclusion of the corresponding closed union of Schubert cells. Now fix \(W\subseteq X_*(T)^+\) for which \(i_{W,*}i_W^!\Ff\cong \Ff\) and let \(\mu\in W\) be a maximal element. 
Then \(j : \Gr_{G,\sgl}^\mu\subseteq \Gr_{G,\sgl}^W\) is an open immersion; let $i$ be the complementary closed inclusion of $\Gr^{W \setminus \{\mu\}}_{G,\sgl}$. 
Consider the semi-infinite orbit \(\Ss_{\mu,\sgl}^-\), for which \(\Ss_{\mu,\sgl}^-\cap \Gr_{G,\sgl}^W = \Ss_{\mu,\sgl}^-\cap \Gr_{G,\sgl}^\mu \cong X\) by maximality of \(\mu\in W\) (see \thref{examples of intersections}). Then, as \(X\to \Gr_{T,\sgl}^\mu\) is an isomorphism on reduced loci, \(\CT_B^{\sgl}(\Ff)\in \DTMrx(\Hck_{T,\sgl})\) implies that the !-restriction of \(\Ff\) to \(\Ss_{\mu,\sgl}^-\cap \Gr_{G,\sgl}^\mu\) is Tate. 
The \(L^+_{\sgl}G\)-equivariance of $\Ff$ implies that \(j^! \Ff\) is Tate by
\cite[Propositions 1.1 and 1.3]{RicharzScholbach:IntersectionCorrigendum}. In the exact triangle \(j_!j^!\Ff\to \Ff\to i_*i^*\Ff\) of motives on \(\Hck_{G, \sgl}\) the first two terms are mapped to Tate motives under $\CT_B^{\sgl}$, hence so does the third term.
An induction on $W$ then shows that \(\Ff\) is stratified Tate, finishing the proof.
\xpf

The following result will allow us to reduce many proofs to the case of tori, which is easier to handle by e.g.~ \thref{piT}. 

\lemm
\thlabel{CT.conservative}
The restriction of the constant term functor $\CT_P^I\colon \DMrx(\Hck_{G, I})\to \DMrx(\Hck_{M, I})$ to the subcategory of bounded motives is conservative. 
\xlemm

\pf
As the property of being bounded is preserved by the constant terms, we can assume \(P=B\) is the Borel by \thref{CT.composition}.
Given some bounded \(\Ff\in \DMrx(\Hck_{G,I})\) that satisfies \(\CT_B^I(\Ff)=0\), we will prove \(\Ff=0\).
Recall that \(\Ff\) being trivial can be checked on the strata of \(X^I\),
as hyperbolic localization commutes with the restriction functors. 
So assume \(\Ff\neq 0\), and let \(X^\phi\subseteq X^I\) be a stratum, on which the restriction of \(\Ff\) does not vanish.
By the factorization property \cite[Theorem 3.2.1]{Zhu:Introduction}, we can assume \(\phi\) is bijective, so that \(X^\phi=X^\circ\).
Now, let \((\mu_i)_i\in (X_*(T)^+)^I\) correspond to a maximal stratum 
\(L^+_IG\backslash\Gr_{G,I}^{\circ,(\mu_i)_i}\) on which \(\Ff\) is supported. 
This stratum is isomorphic to the prestack quotient
\(\left(\prod_{i\in I} (L^+G)_{\mu_i}\times X^\circ \right)\backslash X^\circ\), 
where \((L^+G)_{\mu_i}\) is the stabilizer of \(t^{\mu_i}\) in \(L^+G\).
Consider \(\Ss_{(\mu_i)_i}^-:=\restr{\left(\prod_{i\in I}\Ss_{\mu_i,\{i\}}^-\right)}{X^\circ}\), a connected component of the restriction of the semi-infinite orbit \(\Ss_{\sum_i \mu_i,I}^-\).
As \((\Gr_{G,I}^{\circ,(\mu_i)_i}\cap \Ss_{(\mu_i)_i,I}^-) \times_{X^I} X^\circ\cong X^\circ\) by \thref{examples of intersections},
we see that after forgetting the \(L^+_IT\)-equivariance, the restriction of \(\CT_B^I(\Ff)\) to the fiber over $X^\circ$ of the irreducible component of $\Gr_{T,I}$ indexed by $(\mu_i)_i$ (under \refeq{ZhuFactor}) is given by a shift of the !-pullback of 
\(\Ff\) along \(X^\circ\to \left(\prod_{i\in I} (L^+G)_{\mu_i}\times X^\circ\right)\backslash X^\circ\).
Since forgetting the equivariance is conservative, the restriction of \(\Ff\) to the stratum \(L^+_IG\backslash \Gr_{G, I}^{\circ,(\mu_i)_i}\) must vanish, and we get a contradiction. 
\xpf

\lemm
\thlabel{piT}
The pushforward $\pi_{T!} : \DTMrx(\Gr_{T,I}) \r \DTMrx(X^I)$ is t-exact and conservative.
\xlemm

\pf
Since $\pi_{T!} = \pi_{T*}$ it suffices to check this after replacing $X^I$ by $X^{\phi}$.
Over \(X^\phi\), the reduced subschemes of the connected components of \(\Gr_{T,I}\) are just \(X^\phi\), so the claim is immediate from the definitions.
\xpf

As \(\CT_P^I\) preserves Tate motives, we can now prove it is moreover t-exact.

\prop
\thlabel{CT.t-exact}
The constant term functor \(\CT_P^I\colon \DTMrx(\Hck_{G,I}) \to \DTMrx(\Hck_{M,I})\) is t-exact.
In particular, if $\Ff \in \DTMrx(\Hck_{G,I})$ is bounded, then $\Ff$ lies in positive (resp.~negative) degrees if and only if this is true for $\CT_P^I(\Ff)$. 
(The stratifications are those of \thref{BD.strata} and \thref{All.Zero}.)
\xprop

\pf 
It is enough to show \(\CT_P^I\) is t-exact, the second statement follows from \thref{CT.conservative}. By \thref{CT.composition}, we can then also assume \(P=B\), as t-exactness can be checked on bounded objects.

We will  show that \((q^+_B)_! (p^+_B)^*[\deg]\) is right t-exact, while \((q^-_B)_*(p^-_B)^![\deg]\) is left t-exact.
For the right t-exactness, note that \(\DTMrx^{\leq 0}(\Gr_{G,I})\) is generated by \(\iota_! \DTMrx^{\leq 0}(\coprod_{\phi,\mu} \Gr_{G,I}^{\phi,\mu})\) (\thref{t-structure.stratified}\refit{exactness functors}).
So consider some \(\phi\colon I\twoheadrightarrow J\) and \(\mu=(\mu_j)_j \in (X_*(T)^+)^J\), and
let us denote \(\langle 2\rho,\mu\rangle :=\sum_{j\in J}\langle 2\rho,\mu_j\rangle\). 
Now, let \(\nu=(\nu_j)_j\in X_*(T)^J\), denote by \(\Gr_{B,I}^{\phi,\mu,\nu}\) the intersection of the preimages of \(\Gr_{G,I}^{\phi,\mu}\) and \(\Gr_{T,I}^{\phi,\nu}\) in \(\Gr_{B,I}\), and consider the diagram
\[\begin{tikzcd}
  \Gr_{T,I}^{\phi,\nu}\arrow[d, "\iota_T := \iota_T^{\phi,\nu}"'] & & \Gr_{B,I}^{\phi,\mu,\nu} \arrow[ll,"\ol q := q^+_{\phi,\mu,\nu}"'] \arrow[rr, "\ol p := p^+_{\phi,\mu,\nu}"] \arrow[d] & & \Gr_{G,I}^{\phi,\mu} \arrow[d, "\iota_G := \iota_G^{\phi,\mu}"]\\
  \Gr_{T,I} & & \Gr_{B,I} \arrow[ll, "q := q^+_B"'] \arrow[rr, "p := p^+_B"] & & \Gr_{G,I},
\end{tikzcd}\]
We can assume \(\Gr_{B,I}^{\phi,\mu,\nu}\neq \varnothing\).
Using the product description of Beilinson--Drinfeld Grassmannians over \(X^\phi\) given in \thref{prop.factorization of Gr}, we see via \thref{cellularity of intersection:torus} that \(\Gr_{B,I}^{\phi,\mu,\nu}\) admits a filtrable cellular decomposition relative to \(\Gr_{T,I}^{\phi,\nu} \cong X^{\phi}\), of relative dimension \(\langle \rho,\mu+\nu\rangle\).
We have 
\begin{align*}
\deg & = \langle2\rho,\nu \rangle \\
& =  - \langle 2\rho,\mu\rangle + 2\langle \rho,\mu+\nu \rangle \\
& = -\dim_{X^\phi} (\Gr_{G,I}^{\phi, \mu}) + 2 \dim \ol q.
\end{align*}

By base change and localization, as \(\Gr_{B,I}\times_{\Gr_{G,I}} \Gr_{G,I}^{\phi,\mu}\) is a disjoint union of (finitely many) \(\Gr_{B,I}^{\phi,\mu,\nu}\)'s, it suffices to show 
that \(\ol q_! \ol p^* [\deg]\) is right t-exact.

Let $\pi_G : \Gr_{G,I}^{\phi, \mu} \r X^\phi$ indicate the structural map. 
By the smoothness of $\pi_G$, the definition of the motivic t-structure and \thref{t-structure.stratified}\refit{exactness functors}, $\DTMrx^{\le0}(\Gr_{G,I}^{\phi, \mu})$ is generated under colimits by $\pi_G^*[\dim_{X^\phi} (\Gr_{G,I}^{\phi, \mu})] (\DTMrx^{\le 0}(X^\phi))$.
Similarly, the $\le 0$-aisle of the t-structure on $\Gr_{T,I}^{\phi, \nu}$ is generated by $\pi_T^* (\DTMrx^{\le 0}(X^\phi))$, where $\pi_T \colon \Gr_{T,I}^{\phi,\nu} \cong X^{\phi}$.

Therefore
\begin{align*}
\ol q_! \ol p^* \pi_G^* [\dim_{X^\phi} (\Gr_{G,I}^{\phi, \mu}) + \deg] & = \ol q_! \ol q^* \pi_T^* [\dim_{X^\phi} (\Gr_{G,I}^{\phi, \mu}) + \deg] \\
& = \ol q_! \ol q^* \pi_T^* [2 \dim \ol q]
\end{align*} 
is right t-exact, since $\ol q_! \ol q^*[2 \dim \ol q]$ has that property (\thref{lemm--cellular-coh}).

To prove the asserted left t-exactness, we use the same diagram as above, except that $B$ is replaced by the opposite Borel $B^-$.
The map $\ol q$ still has a filtrable cellular decomposition, of relative dimension $\langle \rho, \mu - \nu \rangle$.
In particular, $\deg_{B^-} = - \deg_B$.
The aisle \(\DTMrx^{\geq 0}(\Hck_{G,I})\) is generated by \(\iota_* \DTMrx^{\geq 0}(\coprod_{\phi,\mu} \Hck_{G,I}^{\phi,\mu})\) (\thref{t-structure.stratified}\refit{exactness functors} joint with \thref{equivariant.MTM}).
By base change, we have to show the left t-exactness of $\iota_{T*} \ol q_* \ol p^! [\deg]$.
Since all t-structures are accessible and right complete, and since the functor preserves colimits, it suffices to see that objects in $\iota_* \MTMrx(\coprod_{\phi,\mu} \Hck_{G,I}^{\phi,\mu})$ are mapped to $\DTM(\Hck_{T,I}^{\phi, \nu})^{\ge 0}$.
We now use the equivalence $\pi_G^![-\dim_{X^{\phi}} \Gr_{G,I}^{\phi,\mu}] : \MTM(X^\phi) \stackrel \cong \r \MTM(\Hck_{G,I}^{\phi, \mu})$, which follows from factorization \refeq{ZhuFactor} and \thref{Orbit.Eq} in the case $S = X^\phi$ (more generally \thref{MTM.equivariant.ind-scheme} applies to $\DTMrx(\Hck_{G,I})$; see \thref{justification Hecke} for justification).
We have
$$\ol q_* \ol p^! \pi_G^![-\dim_{X^\phi} (\Gr_{G,I}^{\phi, \mu}) + \deg_{B}] = \ol q_* \ol q^! \pi_T^! [-\dim_{X^\phi} (\Gr_{G,I}^{\phi, \mu}) + \deg_{B}].$$ 
The functor $\ol q_* \ol q^![-\dim_{X^\phi} (\Gr_{G,I}^{\phi, \mu}) + \deg_{B}]$ is right adjoint to the right t-exact $\ol q_! \ol q^*[\dim_{X^\phi} (\Gr_{G,I}^{\phi, \mu}) + \deg_{B^-}]$, and therefore is left t-exact.
\xpf

We denote the natural projection \(\Gr_{G,I}\to X^I\) by \(\pi_G\), and similarly for other groups.

\prop
\thlabel{CT.Reflect.ULA}
If \(\Ff\in \DTM(\Hck_{G, I})\) satisfies \(\pi_{T!}u^!\CT_B^I(\Ff) \in \DTMrx(X^I)\), then it also satisfies \(\pi_{G!}u^!(\Ff)\in \DTMrx(X^I)\), where Tate motives on $X^I$ are defined with respect to the trivial stratification.
\xprop

\pf
We use the notation of \refeq{hyperbolic.diagram.restricted}.
By assumption,
\[\pi_{T!}u^! \CT_B^I(\Ff)\cong \bigoplus_{\nu\in X_*(T)} \pi_{T!} (q^+_\nu)_! (p^+_\nu)^*(u^!\Ff)[\langle 2\rho,\nu\rangle]\] 
lies in $\DTMrx(X^I)$, which is idempotent-closed, so that each of the summands is also contained in $\DTMrx(X^I)$.
On the other hand, the stratification of \(\Gr_{G,I}\) into the semi-infinite orbits \(\Ss_{\nu,I}^+\) (\thref{Semi-infinite:Strat}) gives a filtration on 
\(\pi_{G!}(u^!\Ff)\) with graded pieces \[\pi_{G!}(p^+_\nu)_! 
(p^+_\nu)^*(u^!\Ff) \cong \pi_{T!}(q^+_\nu)_! (p^+_\nu)^*(u^!\Ff).\] 
So \(\pi_{G!}(u^!\Ff)\) is a colimit of extensions of Tate motives on \(X^I\), and hence Tate itself.
\xpf

The following proposition gives two ways to describe the fiber functor we will use later on, similar to \cite[Theorem 1.5.9]{BaumannRiche:Satake}.

\prop \thlabel{Fiber.Decomposition}
There is a natural equivalence
\[\bigoplus_{n\in \Z} \pH^n\pi_{G!}u^!\cong \pi_{T!}u^! \CT_B^I\]
of functors \(\MTMrx(\Hck_{G, I})\to \MTMrx(X^I)\).
\xprop

\pf
More precisely, we will construct a natural equivalence \[\pH^n\pi_{G!}u^!\cong \bigoplus_{\nu\in X_*(T)\colon \langle2\rho,\nu \rangle=n} \pi_{T!}(q^+_\nu)_!(p^+_\nu)^*u^![n]\] of functors \(\MTMrx(\Hck_{G,I})\to \MTMrx(X^I)\), for each \(n\in \Z\).

For any \(n\in \Z\), let \(\Ss_n^+ = \bigsqcup_{\nu\in X_*(T)\colon \langle 2\rho,\nu\rangle=n} \Ss_{\nu,I}^+ \subseteq \Gr_{G,I}\). 
These \(i_n\colon \Ss_n^+\to \Gr_{G,I}\) determine a decomposition of \(\Gr_{G,I}\). Then we have a natural equivalence 
\(\pH^n\pi_{G!} i_{n!}i_n^*u^!\cong \bigoplus_{\langle 2\rho,\nu\rangle=n} \pi_{T!}(q^+_\nu)_!(p^+_\nu)^*u^![n]\) 
of functors \(\MTM(\Hck_{G,I})\to \MTMrx(X^I)\), while \(\pH^k\pi_{G!} i_{n!}i_n^*u^!\) is trivial if \(k\neq n\); this follows from the proof of \thref{CT.t-exact} and \thref{piT}.

As the Bruhat ordering can only compare tuples of cocharacters \(\mu_i\) for which \(\sum_i \langle 2\rho,\mu_i\rangle\) have the same parity, we can decompose \(\Gr_{G,I} = \Gr_{G,I}^{\mathrm{even}}\coprod \Gr_{G,I}^{\mathrm{odd}}\) 
into two clopen sub-ind-schemes, each containing the Schubert cells corresponding to \((\mu_i)_i\in (X_*(T)^+)^I\) for which \(\sum_i \langle 2\rho,\mu_i\rangle\) is even, respectively odd.
Additionally using that \(\overline{\Ss_{\nu,I}^+}=\bigsqcup_{\nu'\leq \nu} \Ss_{\nu',I}^+\) in \(\Gr_{G,I}\), cf.~ the proof of \thref{Semi-infinite:Strat}, we get the closure relations \(\overline{\Ss_n^+}=\Ss_n^+ \sqcup \Ss_{n-2}^+ \sqcup \Ss_{n-4}^+ \sqcup \ldots = \Ss_n^+\sqcup \overline{\Ss_{n-2}^+}\). 
Denote the corresponding inclusion by \(\overline i_n\colon \overline{\Ss_n^+}\to \Gr_{G,I}\). 
We claim that the two natural morphisms 
$$\pH^n\pi_{G!}i_{n!}i_n^*u^! \to \pH^n\pi_{G!}\overline i_{n!}\overline i_n^* u^!\gets \pH^n\pi_{G!}u^!$$
of functors are equivalences, which will finish the proof. 
Since all functors commute with filtered colimits, it suffices to check this for bounded motives $\Ff\in \MTMrx(\Hck_{G, I})$, even after forgetting the equivariance.  
Moreover, using the decomposition of \(\Gr_{G,I}\) into clopen sub-ind-schemes as above, we can assume the support of \(\Ff\) is contained in \(\Gr_{G,I}^{\mathrm{even}}\); 
the case where its support is contained in \(\Gr_{G,I}^{\mathrm{odd}}\) can be handled analogously.

Consider the closed immersion with complementary open immersion \(\overline{\Ss_{n-2}^+}\xrightarrow{i} \overline{\Ss_n^+} \xleftarrow{j} \Ss_n^+\).
Applying \(\overline i_{n!}\) to the localization sequence \(j_!j^*\overline i_n^*\Ff\to \overline i_n^*\Ff \to i_!i^*\overline i_n^*\Ff\) 
gives the exact triangle \(i_{n!}i_n^* \Ff\to \overline i_{n!}\overline i_n^*\Ff\to \overline i_{n-2, !}\overline i_{n-2}^*\Ff\), which in turns gives a long exact sequence
\[\ldots 
\to \pH^{k}(\pi_{G!}i_{n!}i_n^* \Ff) \to \pH^k(\pi_{G!}\overline i_{n!}\overline i_n^*\Ff) \to \pH^k(\pi_{G!}\overline i_{n-2, !}\overline i_{n-2}^*\Ff) \to \pH^{k+1}(\pi_{G!}i_{n!}i_n^* \Ff)\to \ldots\]
We claim that
$$\pH^{k}(\pi_{G!}\overline i_{n!}\overline i_n^* \Ff) = 0 \text{  if  } k> n \text{  or  } k \text{  is odd}.$$ This can be proved by induction on $n$, starting with the observation that \(\pi_{G!}\overline i_{n!}\overline i_n^*\Ff=0\) for \(n\ll 0\) as the support of \(\Ff\) is bounded. Moreover, a bounded motive has support in finitely many $\Ss_n^+$, so that the vanishing of $\pe \H^k \pi_{G!} i_{n!} i_n^* $ for all $k \ne n$ shows that $\pi_{G!}\overline i_{n!}\overline i_n^* \Ff$ lives in finitely many cohomological degrees. Thus, the asserted vanishing holds for $k \ll 0$.
Now the claim follows in general by induction on $k$, using the long exact sequence and the fact that $\pH^k(\pi_{G!} i_{n!}i_n^*\Ff) =0$ if \(k\neq n\).

We further claim that the natural localization morphisms give isomorphisms
\[\pH^n(\pi_{G!}i_{n!}i_n^* \Ff) \xrightarrow{\cong} \pH^n(\pi_{G!}\overline i_{n!}\overline i_n^*\Ff) \xleftarrow{\cong} \pH^n(\pi_{G!}\overline i_{m!}\overline i_m^*\Ff)\] for all \(m\geq n\) such that \(m\equiv n\mod 2\). The first isomorphism is immediate from the previous claim; the second follows from induction on $m$. As the support of \(\Ff\) is bounded, we conclude by noting that \(\overline i_{m!}\overline i_m^*\Ff \cong \Ff\) for \(m\gg0\).
\xpf

\coro
\thlabel{fiber.functor.properties}
The functor \(\pi_{T!}u^!\CT_B^I\cong \bigoplus_{n\in \Z} \pH^n \pi_{G!}u^! \colon \MTMrx(\Hck_{G,I}) \r \MTM(X^I) \) is exact, conservative, and faithful.
\xcoro

\pf
The exactness combines \thref{Fiber.Decomposition} and \thref{CT.t-exact}.
The functor is conservative on bounded motives by \thref{CT.conservative} and \thref{piT}.
The conservativity of an exact functor is equivalent to its faithfulness.
The faithfulness for bounded objects implies the one for unbounded objects.
\xpf

\subsection{The global Satake category}
For a surjection $\phi \colon I \twoheadrightarrow J$ of nonempty finite sets, recall that we have defined the locally closed subscheme $X^{\phi} \subset X^I$ \eqref{Loc.Sub} and the open subscheme \eqref{Open.Sub} $X^{(\phi)} \subset X^I$. We denote the corresponding open immersion and complementary closed immersion into $\Hck_{G,I}$ by
$$j^{(\phi)} \colon \restr{\Hck_{G,I}}{X^{(\phi)}} \rightarrow \Hck_{G,I}, \quad i_{(\phi)} \colon \restr{\Hck_{G,I}}{X - X^{(\phi)}} \rightarrow \Hck_{G,I}.$$ If $\phi = \id$, then $X^{\id} \subset X^I$ is the locus with distinct coordinates. We denote the base change of an $X^I$-scheme to $X^{\id}$ with the symbol $\circ$. For example, $X^\circ = X^{\id} $ and $\Hck_{G,I}^\circ = \Hck_{G,I} \times_{X^I} X^{\circ}$.

\defi
\thlabel{IC.X.circ}
Consider some $\mu \in (X_*(T)^+)^I$ and let $j^{\circ, \mu} \colon \Hck_{G,I}^{\circ, \mu} \rightarrow \Hck_{G,I}$ be the inclusion of the corresponding stratum, as defined in \thref{BD.strata} when $\phi = \id$. By \thref{equivariant.MTM}, we may for $L \in \MTMrx(S)$ define
$$\IC_{\mu, L} \in \MTMrx(\Hck_{G, I})$$ to be the (reduced) intersection motive of this stratum.
\xdefi

\rema
\thlabel{justification Hecke}
Note that the three conditions in \thref{MTM.equivariant.ind-scheme} are satisfied for the following reasons: The quotients $G_k$ are constructed as in \cite[Corollary 3.11]{HainesRicharz:TestFunctionsWeil}, cf.~also \cite[Example A.12(ii)]{RicharzScholbach:Intersection} for the particular case relevant here. The split pro-unipotence of $\ker(G_k \r G)$ is proved in \cite[Proposition A.9]{RicharzScholbach:Intersection}, and cellularity of the fibers of $G_k$ over the $X^\phi$ follows from factorization as in \refeq{ZhuFactor} and the Bruhat decomposition applied to the maximal reductive quotients.
\xrema

\subsubsection{Definition and first properties of the bounded Satake category}
\defi \thlabel{Sat.Def}
Let $W \subset (X_*(T)^+)^I$ be a finite subset closed under the Bruhat order, and let
$\Hck_{G,I}^W$ be the closure of the union of the strata $\Hck_{G,I}^{\circ, \mu}$ for $\mu \in W$.
The \emph{bounded Satake category} $\Sat_{\redx, W}^{G,I}$ is defined as the full subcategory of $\MTMrx(\Hck_{G, I}^W)$ consisting of objects $\mathcal{F}$ that admit a finite filtration 
$$0 = \mathcal{F}_0 \subset \cdots \subset \mathcal{F}_k = \mathcal{F}\eqlabel{filtration}$$ 
for some integer $k$ such that $\mathcal{F}_i \in \MTMrx(\Hck_{G,I}^W)$ and
$$\mathcal{F}_i/\mathcal{F}_{i-1} \cong \IC_{\mu_i, L_i}$$ for some $L_i \in \MTMrx(S)$ and $\mu_i \in (X_*(T)^+)^I$ for all $1 \leq i \leq k$. 
\xdefi

In \refsect{unbounded.objects}, we will take an appropriate colimit over these bounded Satake categories to obtain the global (unbounded) Satake category.

\rema \thlabel{Sub.Quot}
By \thref{MTM.equivariant.ind-scheme}, every object in $\MTMrx(\Hck_{G,I}^W)$  admits a finite filtration by intersection motives of objects in $\MTMrx(X^\phi)$, and the point of \thref{Sat.Def} is to require that only strata supported over $X^\circ$ and objects in $\MTMrx(X^\circ)$ pulled back from $\MTMrx(S)$ appear. 
Moreover, by \thref{FF.Stab} every filtration  by objects in $\MTMrx(\Gr_{G,I}^W)$ of the underlying (non-equivariant) motive of an object in $\MTMrx(\Hck_{G,I}^W)$ comes from a filtration defined in $\MTMrx(\Hck_{G,I}^W)$.
\xrema

\rema \thlabel{ULAremark} 
In \cite[Definition VI.9.1]{FarguesScholze:Geometrization}, the Satake category is defined by restricting to sheaves which are universally locally acyclic (ULA) relative to $X^I$, which ensures that after applying constant term functors, the result is essentially a collection of local systems on $X^I$. 
Instead of developing a notion of ULA Nisnevich motives, we adopt a more ad-hoc approach. It is motivated by \thref{Sat.T} and \thref{CT.Sat} below, which show that constant terms produce collections of objects in $\MTMrx(X^I) \cong \MTMrx(S)$, i.e.~unstratified as opposed to stratified motives.
\xrema

\prop \thlabel{Sat.1} Let $p \colon \Gr_{G,\{*\}} \r \Gr_G$ be the projection coming from the identification (\ref{Gr.Pt}). Then $p^![-1]$ induces a t-exact equivalence
$$p^![-1] \colon \DTMrx(\Gr_G) \xrightarrow{\sim}  \DTMrx(\Gr_{G,\{*\}})$$ with quasi-inverse $p_![1]$.
\xprop

\pf
As $p^!$ commutes with both types of pushforwards and pullbacks between strata, the argument in \cite[Proposition~4.25]{EberhardtScholbach:Integral} reduces us to the case of a single stratum.
Here the result follows from $\A^1$-homotopy invariance.
\xpf

\coro
\thlabel{MTM of local Grassmannian} If $I = \{*\}$, let $\Gr_G^W \subset \Gr_G$ be the union of the $\Gr_G^{\mu}$ for $\mu \in W$.
Then the functor $p^!(-1)[-1]\cong p^*[1]$ induces an equivalence
$$\MTMrx(L^+G \backslash \Gr_G^W) \xrightarrow{\sim}  \Sat_{\redx, W}^{G,\{*\}},$$ which identifies the IC motives in these two categories.
\xcoro

\pf
First note that $\Sat_{\redx, W}^{G,\{*\}} = \MTMrx(L_{\{*\}}^+G \backslash \Gr^W_{G,\{*\}})$ by \thref{Orbit.Eq}, since $X^{\circ} = X$ if $I = \{*\}$. 
By smooth base change and the isomorphism $L_{\{*\}}^+G \cong L^+G \times X$, it  follows that $p^!$ preserves equivariance.
Hence $p^![-1]$ induces a functor as in the proposition. Since $p$ is smooth of relative dimension one and $\Gr_{G,\{*\}}$ has the preimage stratification, it also follows that $p^![-1]$ is t-exact (recall the convention on the normalization of the t-structure from \refeq{normalization.t-structure}).
Now we may conclude by the same argument as in \thref{Sat.1}.
\xpf

\prop \thlabel{Sat.Int.Ext}
For any surjection $\phi \colon I \twoheadrightarrow J$ and $\mathcal{F} \in \Sat_{\redx, W}^{G,I}$, there is a canonical isomorphism 
$$j_{!*}^{(\phi)} (j^{(\phi),*} (\mathcal{F})) \cong \mathcal{F}.$$
where $j^{(\phi)} \colon \restr{\Hck_{G,I}}{X^{(\phi)}} \rightarrow \Hck_{G,I}$.
\xprop

\pf
Both objects are perverse and canonically identified over $X^{(\phi)}$, so it suffices to show $\pe i^*_{(\phi)} \mathcal{F} = 0$ and $\pe i^!_{(\phi)} \mathcal{F} = 0$.
By induction on the length of a filtration of $\calF$ as in \refeq{filtration}, it suffices to consider $\mathcal{F} = \IC_{\mu, L}$ for some $\mu \in (X_*(T)^+)^I$ and $L \in \MTMrx(S)$. This case is immediate because $\IC_{\mu, L}$ is an intermediate extension from the open subset $X^\circ \subset X^{(\phi)}$. 
\xpf

\exam \thlabel{L.Perv.XI}
For the trivial group $G = 1$, we have $\Gr_{1,I} = X^I$.
Since the projection $p : X^I \r S$ and also the structural map of all strata are smooth, the functor $p^*[I]$ is t-exact, so that the intersection motives (in \refeq{IC.w.L}, with respect to the trivial stratification of $X^I$) are just given by $p^* L[I]$ for $L \in \MTMrx(S)$.
\xexam

Recall that $j^{(\phi)}_{!*}$ is fully faithful. This is a generality about recollement of t-structures, cf.~\cite[Remarque 1.4.14.1]{BeilinsonBernsteinDeligne:Faisceaux}, so we arrive at the following corollary of \thref{Sat.Int.Ext}. 

\coro
\thlabel{Sat.Restr.FF}
The restriction functor $$j^{(\phi),*} \colon \Sat_{\redx, W}^{G,I} \rightarrow \MTMrx(\restr{\Hck_{G,I}^W}{X^{(\phi)}}) $$ is fully faithful. 
\xcoro

\subsubsection{The bounded Satake category of a torus}
We now give a complete description of $\Sat^{T,I}_{\redx, W}$, i.e., the case of a torus.
We can enlarge $W \subset (X_*(T)^+)^I$ to a subset of $\widetilde W \subset X_*(T)^I$ by taking orbits of all elements under the coordinate-wise action of $|I|$-copies of the Weyl group. Then $\CT_B^I$
sends motives supported on $\Hck_{G, I}^W$ to motives supported on $\Gr_{T,I}^{\widetilde W}$. To ease the notation we will usually write $\Gr_{T,I}^{W}$ instead of $\Gr_{T,I}^{\widetilde W}$. With this convention, we have $\Gr_{T,I}^{W} = \Gr_{G,I}^W \cap \Gr_{T,I}$.

Note that $$(\Gr_{T,I}^\circ)_{\text{red}} = \coprod_{\mu \in X_*(T)^I} X^\circ.$$ It follows that 
$$\MTMrx(\Gr_{T,I}^\circ) \cong \Fun(X_*(T)^I, \MTMrx(X^\circ)). \eqlabel{Gr.T.circ}$$
Similarly, the connected components of $\Gr_{T,I}^{\circ, W}$ are in bijection with $\widetilde W$.
Additionally, the connected components of $\Gr_{T,I}$ are in bijection with $X_*(T)$. The reduced closure of the connected component of $\Gr_{T,I}^{\circ}$ associated to $(\mu_i)_i \in X_*(T)^I$ is isomorphic to $X^I$, and it is an irreducible component of the connected component of $\Gr_{T,I}$ associated to $\sum_i \mu_i$.  For $W \leq W'$, at the level of irreducible components the inclusion $i_{W, W'} \colon \Gr_{T,I}^W \r \Gr_{T,I}^{W'}$ is identified with the obvious inclusion $\sqcup_{\widetilde W} X^I \r \sqcup_{\widetilde W'} X^I$.

\prop \thlabel{Sat.T}
Let $j_I \colon X^\circ \rightarrow X^I$ and $j^{\circ} \colon \Hck_{T,I}^{\circ, W} \rightarrow \Hck_{T,I}^{W}$ be the inclusions.
We have the following commutative diagram, where the composite $\Theta_T$ is fully faithful and induces an equivalence as indicated:
$$\xymatrix@C=20mm{
\Fun(\widetilde W, \MTMrx(X^I)) \ar[r]^{\prod_{\widetilde W} j_I^{*}} \ar[dd]_\cong^{\Theta_T}  & \Fun(\widetilde W, \MTMrx(X^\circ)) \ar[d]_\cong^{\refeq{Gr.T.circ}} \\ & \MTMrx(\Hck_{T,I}^{\circ, W}) \ar[d]^{j^{\circ}_{!*}} \\
\Sat^{T,I}_{\redx, W} \ar@{^{(}->}[r] & \MTMrx(\Hck_{T,I}^W).
}$$
Here mixed Tate motives on $X^I$ are with respect to the \emph{trivial stratification} (so that $\MTMrx(X^I) \cong \MTMrx(S)$).
In particular, $\Sat^{T,I}_{\redx, W}$ is a compactly generated category. Moreover, for $W \leq W'$ and the associated inclusion $i_{W, W'} \colon \Gr_{T,I}^W \r \Gr_{T,I}^{W'}$, the functors $(i_{W,W'})_!$ and $\pe i_{W,W'}^! $ preserve $\Sat^{T,I}_{\redx, W}$, and with respect to the identification $\Theta_T$ they are given by extension by zero and restriction along $\widetilde W \subset \widetilde W'$. 
\xprop

\pf
Let $L \in \MTMrx(X^I)$. The object of $\Fun(\widetilde W, \MTMrx(X^I))$ supported at $\mu$ with value $L$ is mapped under the above functor to $\IC_{\mu, L} \in \Sat_{\redx, W}^{T,I}$.
Thus $\Theta_T$ exists as shown in the diagram. 

By \thref{FF.Stab}, we may work with motives on $\Gr_{T,I}^W$ instead of $\Hck_{T,I}^W$.
Note that by \thref{L.Perv.XI}, $\IC_{\mu,L}$ is the pullback of $L[I]$ along the structure morphism $X^I \rightarrow S$ of the corresponding irreducible component of $\Gr_{T,I}^W$. Hence $\Theta_T$ is exact. This description of $\IC_{\mu,L}$, along with the fact that $j_{!*}^{\circ}$ is fully faithful, implies that $\Theta_T$ is also fully faithful. To see that $\Theta_T$ is an equivalence, it suffices to show that if $\mu \neq \lambda$, 
then
 $\Ext^1_{\DTMrx(\Gr_{T,I})}(\IC_{\mu,L_1}, \IC_{\lambda,L_2}) = 0$ for any $L_1$, $L_2 \in \MTMrx(S)$.
 
To prove this, let $Y$ be the (reduced) union of the supports of $\IC_{\mu,L_1}$ and $\IC_{\lambda,L_2}$. We can assume $Y$ is connected. Let $i \colon Z \rightarrow Y$ be the inclusion of the intersection of the supports, and let $j \colon U \rightarrow Z$ be the inclusion of the complement. 
By localization, it suffices to prove
$$\Hom_{\DTMrx(Y)}(\IC_{\mu,L_1}, j_*j^*\IC_{\lambda,L_2}[1]) = 0, \quad \Hom_{\DTMrx(Y)}(\IC_{\mu,L_1}, i_*i^!\IC_{\lambda,L_2}[1]) = 0.$$ Using adjunction, the left group is a $\Hom$ group between motives on $U$ supported on disjoint opens, so it vanishes. Now apply adjunction to identify the right group with a $\Hom$ group on $Z$. The scheme $Z$ is stratified by a union of cells $X^\phi$ in the support $X^I$ of $\IC_{\lambda,L_2}$, so by localization we can assume $Z = X^\phi$, where $X^\phi \subset X^I$ has codimension $c > 0$ since $\mu \neq \lambda$.
Then by relative purity the group on the right is $$\Hom_{\DTMrx(X^\phi)}(L_1, L_2(-|I|+c)[1 -2|I|+2c]).$$ Since $1-2|I|+2c < 0$ and we have a t-structure on $\DTMrx(X^\phi)$ with the $L_i$ lying in the same perverse degree by \thref{admissible.SNC.divisor} and \thref{MTM.equivariant.ind-scheme}, this group vanishes.

The category $\MTMrx(X^I)$ is compactly generated (\thref{t-structure.heart}), hence so is the above functor category. The final statement about $(i_{W,W'})_!$ and $\pe i_{W,W'}^! $ follows from the geometry of $\Gr_{T,I}$ and the description of $\Sat^{T,I}_{\redx, W}$ provided by $\Theta_T$.
\xpf 

\subsubsection{Behaviour under constant term functors}
\prop \thlabel{CT.Sat}
The constant term functor restricts to a functor
$$\CT^I_B \colon \Sat_{\redx, W}^{G,I} \rightarrow \Sat_{\redx, W}^{T,I}.$$
\xprop

\pf  
It suffices to show that $\CT_B^I(\IC_{\mu, L}) \in \Sat_{\redx, W}^{T,I}$ for $L \in \MTMrx(S)$ and $\mu \in W$.
To prove this, note that \thref{CT.t-exact} implies $\CT_B^I(\IC_{\mu, L}) \cong j^\circ_{!*} j^{\circ *}(\CT_B^I(\IC_{\mu, L})).$ Hence by \thref{Sat.T}, it suffices to identify $j^{\circ *}(\CT_B^I(\IC_{\mu, L})) \cong \CT_B^I (j^{\circ *}(\IC_{\mu, L})) $ with an object in the image of $\Fun(\widetilde W, \MTMrx(X^I)) \r \Fun(\widetilde W, \MTMrx(X^\circ))$. 
Let $\IC_{\mu,L}^I$ be the $\IC$-motive associated to $L$ and the open embedding $\prod_{i \in I}(\Gr_G^{\mu_i} \times X) \r \prod_{i \in I} (\Gr_G^{\leq \mu_i} \times X)$.
By admissibility of $X^\circ$, pullback along $(\Gr_G)^I \times X^\circ \r (\Gr_G)^I$ is t-exact up to a shift for the product stratifications, cf.~the proof of \thref{MTM.equivariant.ind-scheme}. Hence by factorization \refeq{ZhuFactor}, $j^{\circ *}(\IC_{\mu, L})$ is isomorphic to the restriction of $\IC_{\mu,L}^I$ to $X^\circ$, as both are the pullback of an IC-motive on $(\Gr_G)^I$ 
over the base $S$. Now we conclude, since constant terms commute with pullback along the smooth map $X^\circ \r S$, and since constant terms preserve Tateness by \thref{CT.DTM} (applied to the base $S$, or equivalently, by homotopy equivalence, to the case $I = \{*\}$). In slightly more detail, $j^{\circ *}(\CT_B^I(\IC_{\mu, L}))$ is isomorphic to the restriction of $(\prod_i q_B^+)_! (\prod_{i} p_B^+)^* \IC_{\mu,L}^I$ to $X^\circ$. 
Analogous reasoning as in the proof of \thref{Sat.1} yields a homotopy equivalence $\DTM(\Gr_{G^I} \times X^I) \cong \DTM(\Gr_{G^I, \{*\}})$ which extends the equivalence $\DTM(X^I) \cong \DTM(X)$. 
Using the identification $\Gr_{G^I} \cong (\Gr_G)^I$ of \thref{fiber product of affine grassmannians},  $(\prod_i q_B^+)_! (\prod_{i} p_B^+)^*$ corresponds under these homotopy equivalences to a constant term functor $\MTMrx(\Hck_{G^I, \{*\}}) \r \MTMrx(\Hck_{T^I, \{*\}})$. 
It follows that the map \(\widetilde{W}\to \MTMrx(X^\circ)\) corresponding to \(j^{\circ*}(\CT_B^I(\IC_{\mu, L}))\) takes values in unstratified Tate motives, i.e., pulled back along \(\MTMrx(X)\cong \MTMrx(S)\to \MTMrx(X^I)\).
\xpf

\theo \thlabel{Sat.Props}
The following conditions on $\mathcal{F} \in \MTMrx(\Hck_{G, I}^W)$ are equivalent.
\begin{enumerate}
\item $\mathcal{F}$ belongs to $\Sat_{\redx, W}^{G,I}$.
\item $\CT_B^I(\mathcal{F})$ belongs to $\Sat_{\redx, W}^{T,I}$.
\item For every finite filtration of $\calF$ with subquotients isomorphic to IC motives of objects $L \in \MTMrx(\Hck_{G,I}^{\phi, \mu})$ with $\mu \in W$, \(\phi\) is a bijection (so the stratum is supported over $X^I$) and $L$ is pulled back from $\MTMrx(S)$.
\end{enumerate}
The category $\Sat_{\redx, W}^{G,I} \subset \MTMrx(\Hck_{G, I}^W)$ is stable under subquotients and extensions, and in particular it is abelian.
Moreover, it is compactly generated.
\xtheo

\pf 
Property (1) implies (2) by \thref{CT.Sat}. To show that (2) implies (1), take a filtration of $\calF$ as in (3). The minimal nonzero object in this filtration is the IC motive of some $L \in \MTMrx(\Hck_{G,I}^{\phi, \mu})$. By \thref{Orbit.Eq} (applied to the base $S= X^\phi$), $L$ is pulled back from $\MTMrx(X^\phi)$.
By \thref{Sat.Int.Ext} and the fact that $\CT_B^I$ is conservative and exact, \(\phi\) must be a bijection, so we can write this IC motive as $\IC_{\mu, L}$ for $L \in \MTMrx(X^\circ)$. 
Letting \(\iota_T^{\phi,\mu}\colon \Gr_{T, I}^{\phi,\mu}\hookrightarrow \Gr_{T,I}\) denote the open inclusion, $\iota_{T!*}^{\phi,\mu}\iota_T^{\phi,\mu,*}\CT^I_B(\IC_{\mu, L})$ is a subobject of $\iota_{T!*}^{\phi,\mu}\iota_T^{\phi,\mu,*}\CT_B^I(\calF)$ since the first two functors are t-exact and $\iota_{T!*}^{\phi,\mu}$ preserves injections. The latter motive is an unstratified object of $\MTMrx(X^I)$ by assumption. Thus, by \thref{lemm-XI-stable}, $\iota_{T!*}^{\phi,\mu}\iota_T^{\phi,\mu,*}\CT^I_B(\IC_{\mu, L})$ is pulled back from $\MTMrx(S)$. By the proof of \thref{CT.conservative}, the restriction $\iota_T^{\phi,\mu,*}\CT_B^I(\IC_{\mu, L})\in \MTMrx(X^\circ)$ is isomorphic to a twist of $L$, 
so that $L$ is pulled back from $\MTMrx(S)$. Hence  
$\IC_{\mu, L} \in \Sat_{\redx, W}^{G,I}$, and $\CT_B^I(\calF / \IC_{\mu, L}) \in \Sat_{\redx, W}^{T,I}$ since the latter category is abelian.
By induction on the length of a filtration of $\calF$, it follows that (2) implies (1). It is clear that (3) implies (1), and the previous inductive argument also shows that (2) implies (3). 

Since $\Sat_{\redx, W}^{T,I}$ is abelian, then (2) implies that $\Sat_{\redx, W}^{G,I}$ is abelian. Furthermore, $\Sat_{\redx, W}^{G,I}$ is stable under extensions by definition.  Given that it is abelian, to see that is stable under subquotients it suffices to show it is stable under subobjects. For this, it suffices to show that if the IC motive of some $L \in \MTMrx(\Hck_{G,I}^{\phi, \mu})$ is a subobject of some object in $\Sat_{\redx, W}^{G,I}$, then $\phi = \id$ and $L$ is pulled back from $\MTMrx(S)$, so that the IC motive also lies in $\Sat_{\redx, W}^{G,I}$.
But we already showed this above. 

Since $\CT^I_B$ preserves (and $\Sat^{T,I}_{\redx, W}$ has) all (small) colimits 
by \thref{Remark:CTforGr}, $\Sat^{G,I}_{\redx, W}$ is stable under colimits.
Since $\CT^I_B$ is conservative on bounded objects and preserves colimits, its left adjoint maps a set of compact generators of $\Sat_{\redx, W}^{T,I}$ to a set of compact generators of $\Sat^{G,I}_{\redx, W}$. 
\xpf

\subsubsection{Extension to unbounded objects}
\label{sect--unbounded.objects}
We are now in a position to bind all the bounded Satake categories together.

\defilemm
For $W \le W'$, both adjoints 
$$i_! := (i_{W,W'})_! : \MTMrx(\Hck_{G,I}^W) \leftrightarrows \MTMrx(\Hck_{G,I}^{W'}) : \pe i_{W,W'}^! =: \pe i^!$$ 
preserve the Satake categories.
Therefore 
$$\Satrx^{G, I} := \colim_{i_!} \Sat_{\redx, W}^{G,I} = \lim_{\pe i^!} \Sat_{\redx, W}^{G, I}\eqlabel{eqn.Sat.colim.lim}$$
is an abelian full subcategory of $\MTM(\Hck_{G,I})$.
An object $M \in \MTMrx(\Hck_{G,I})$ lies in this subcategory iff $\pe i_W^! M \in \Sat_{\redx, W}^{G,I}$ for all $W$.
We call this category the \emph{global (unbounded) Satake category}.
\xdefilemm

\pf 
We write $i$ for the transition maps $i_{W,W'}$ and also drop $\redx$, $G$ and $I$ from the notation.
The functors $i_!$ are exact.
They preserve the IC-motives in \thref{IC.X.circ}, and therefore preserve the Satake category.
We now show that its right adjoint (on the level of $\MTM(\Hck)$) $\pe i^!$ also preserves the Satake category.
This will imply that the above limit is well-defined and is, by the adjoint functor theorem, equivalent to the colimit (which is formed in $\Pr$, the category of presentable (ordinary) categories).
To see this claim, note first that $\pe i^! \IC_{\mu, L} = 0$ if $\mu \notin W$, since $\IC_{\mu, L}$ has no perverse subsheaves supported on $\Gr^{\circ, \leq \mu} \setminus \Gr^{\circ, \mu}$. 
If $\mu \in W$, then this object is just $\IC_{\mu, L}$ again.
Given an extension $0 \r A \r B \r C \r 0$ in $\Sat_{W'}$, such that $\pe i^!$ maps the outer terms to $\Sat_W$, we have an exact sequence $0 \r \pe i^! A \r \pe i^! B \r C' \r 0$, where $C' \subset \pe i^! C$ is a subobject.
This subobject, which a priori lies in $\MTM(\Hck^W)$, is an object of $\Sat_W$ by \thref{Sat.Props}.
Since $\Sat_W$ is also stable under extensions, $\pe i^! B \in \Sat_W$.
Therefore $\pe i^!$ preserves the Satake category.
In addition, $\pe i^!$ preserves filtered colimits (in $\MTM(\Hck)$ and therefore in $\Sat$), since this holds for $i^!$ and the truncation functors for the compactly generated t-structure.

Given that the (bounded) Satake categories are presentable (in fact compactly generated), we are in a position to apply the following general paradigm (which already appears in \refeq{DM.colim.lim} and \refeq{obj MTM colimit}): for a filtered diagram $I \r \PrL, i \mapsto C_i$, with transition functors $C_i \stackrel{L_{ij}} \r C_j$, and the corresponding diagram $I^\opp \r \Pr$ obtained by passing to right adjoints, denoted by $C_j \stackrel{R_{ij}} \r C_i$, there is an equivalence
$$\colim_{L_{ij}} C_i = \lim_{R_{ij}} C_i =: C.$$
Writing $L_i : C_i \r C$ and $R_i : C \r C_i$ for the canonical insertion and evaluation functors, the natural map $\colim L_i R_i c \r c$ is an isomorphism provided that the $L_{ij}$ are fully faithful and that the $R_{ij}$ preserve filtered colimits. This is proven in \cite[Lemma~1.3.6]{Gaitsgory:Generalities} for DG-categories; for a filtered diagram, one only needs the right adjoints to preserve filtered colimits.

Thus, an object $M \in \MTM(\Hck)$ can be written as $M = \colim i_{W!} M_W$, for $M_W \in \MTM(\Hck^W)$, and $M \in \Sat$ iff all the $M_W \in \Sat_W$. 

We now check that $\Sat$ is an abelian subcategory of $\MTM(\Hck)$: if $f : A \r B$ is a map in $\Sat$, we have to show its (co)kernel, computed in $\MTM(\Hck)$, lies in $\Sat$.
The evaluation in $\MTM(\Hck^W)$ of this (co)kernel is $\colim_{W' \ge W} \pe i^! \operatorname{(co)ker} f_{W'}$.
Here the (co)kernel is a priori in $\MTM(\Hck^W)$, but lies in $\Sat_{W'}$.
By the above, the term in the colimit therefore lies in $\Sat_W$, and hence so does the entire expression.
\xpf

\rema
The following results about $\Sat_{\redx, W}^{G,I}$ now extend to \(\Satrx^{G,I}\): For $I = \sgl$, we have $\Satrx^{G, \sgl} = \MTMrx(L^+ G \setminus \Gr_G)$. This follows by taking the colimit over the equivalences in \thref{MTM of local Grassmannian}. \thref{Sat.Int.Ext} is also true for \(\Satrx^{G,I}\), since the inclusions $i_{W,W'}$ and $\restr{i_{W,W'}}{X^{(\phi)}}$ form a cartesian square with the inclusions $\restr{j^{(\phi)}}{\Hck_{G,I}^W}$ and $\restr{j^{(\phi)}}{\Hck_{G,I}^{W'}}$. Since $\pe i_{W,W'}^!$ commutes with $j^{(\phi), *} = \pe j^{(\phi), *}$, the latter functor is fully faithful on \(\Satrx^{G,I}\) as in \thref{Sat.Restr.FF}.
Finally, by the last part of \thref{Sat.T} it follows that $\Satrx^{T,I} = \lim_{\widetilde W} \Fun(\widetilde W, \MTMrx(X^I)) =   \Fun(X_*(T)^I, \MTMrx(X^I))$.
\xrema

\coro \thlabel{coro--SatG-T}
For \(\calF \in \MTMrx(\Hck_{G, I})\), we have \(\calF \in \Satrx^{G,I}\) if and only if \(\CT_B^I(\calF)\in \Satrx^{T,I}\).

\xcoro

\pf
For an object $\calF \in \MTM(\Hck)$, we write $\calF_W := \pe i_W^! \calF \in \MTM(\Hck^W)$. 
We then have the following chain of equivalences
\begin{align*}
\calF \in \Satrx^{G,I} & \Leftrightarrow \calF_W \in \Sat^{G,I}_{\redx, W} \ \forall \ W  & \text{by definition of }\Satrx^{G,I} \\
& \Leftrightarrow \CT_B^I(\calF_W) \in \Sat^{T,I}_{\redx,W} \ \forall \ W & \text{by \thref{Sat.Props}} \\
& \Leftrightarrow \pe i_{W,W'}^! \CT_B^I(\calF_{W'}) \in \Sat^{T,I}_{\redx,W} \ \forall \ W'>W & \text{by stability of Sat under }\pe i^!\\
& \Leftrightarrow \CT_B^I(\calF) \in \Satrx^{T,I}.
\end{align*}
For the last equivalence, we first note that $\pe i_W^! \CT_B^I(\calF) = \pe i_W^! \CT_B^I(\colim_{W'} i_{W'!} \calF_{W'}) = \colim \pe i_{W,W'}^! \CT_B^I(\calF_{W'})$. Then we conclude using that the transition maps in $\colim \pe i_{W,W'}^! \CT(\calF_{W'})$ are injective, and because $\Sat^T_W$ is closed under subobjects.
\xpf

Now that we know the property of lying in the Satake category can be checked after applying \(\CT_B^I\), the following corollary is immediate.

\coro
For any parabolic \(P\subseteq G\) with Levi \(M\), the constant term functor restricts to
\[\CT_P^I\colon \Satrx^{G,I}\to \Satrx^{M,I}.\]
\xcoro

Note that $\overline{X^{\phi}} \cong X^J$. Let $i_{\overline{\phi}} \colon \Gr_J \rightarrow \Gr_I$ be the corresponding closed immersion induced by the factorization isomorphisms \cite[3.2.1]{Zhu:Introduction}. Let $d_{\phi} = |I|-|J|$. By \thref{equivariant.functoriality}, there are functors
$$i_{\overline{\phi}}^*, \: i_{\overline{\phi}}^! \colon \DTMrx(\Hck_{G, I}) \rightarrow \DTMrx(\Hck_{G, J})$$ which restrict the usual pullback functors on the non-equivariant derived categories.

\prop 
\thlabel{Sat.Tate}
For $\mathcal{F} \in \Satrx^{G,I}$ we have
$$i_{\overline{\phi}}^*\mathcal{F} [-d_{\phi}] \in \Satrx^{G,J}, \quad i_{\overline{\phi}}^! \mathcal{F}[d_{\phi}] \in \Satrx^{G,J}.$$ Furthermore, pushforward along \(\pi_G\) induces a functor
\[\pi_{G!} u^!\colon \Satrx^{G,I} \to \DTM_{\redx}(X^I),\] where Tate motives on $X^I$ are defined with respect to the trivial stratification.
\xprop

\pf
The functors $i_{\overline{\phi}}^*$ and $i_{\overline{\phi}}^!$ commute with $\CT_B^I$, so by \thref{coro--SatG-T} we can assume that $G = T$ for the first statement. By  \thref{CT.Reflect.ULA} we can also assume that $G=T$ for the second statement. 
Then since $\Satrx^{T,I} = \Fun(X_*(T)^I, \MTMrx(X^I))$, by continuity we may reduce to the case $\mathcal{F} = \IC_{\mu, L}$ for some $\mu \in X_*(T)^I$ and $L \in \MTMrx(S)$. 
Then $\IC_{\mu, L}$ is a shifted constant sheaf supported on a copy of $X^I$, so the first result follows from relative purity applied to the inclusion $X^J \r X^I$. More precisely, if $(\mu'_j)_j \in X_*(T)^J$ is defined by $\mu'_j = \sum_{\phi^{-1}(j)} \mu_i$, then $i_{\overline{\phi}}^* \IC_{\mu, L} [-d_{\phi}] \cong \IC_{\mu', L}$ and $i_{\overline{\phi}}^! \IC_{\mu, L} [d_{\phi}] \cong \IC_{\mu', L(-d_\phi)}$. For the second statement, $\pi_{T!}(\IC_{L,\mu}) = L[I]$ is Tate because $\pi_{T}$ restricts to the identity morphism on this irreducible component.
\xpf

\lemm \thlabel{TypeII.Relation.iso} 
In the notation of \thref{TypeII.Relation}, if $\calF_1, \calF_2 \in \Satrx^{G,I}$ then the map \eqref{TypeII.Relation.map} is an isomorphism, $$f_{I}^!(\calF_1 \star \calF_2) (-I)[-I] \cong i_1^!m_{\phi !}(\calF_1 \widetilde \boxtimes_S \calF_2 ).$$
\xlemm

\pf To build \eqref{TypeII.Relation.map} we 
used the compatibility of $\bx$ vs.~*-pullbacks, !-pushforwards and the adjoint of the projection formula (cf.~the paragraph before \eqref{f! f* tensor}) to construct a map
$\calF_1 \boxtimes_{X^{I}} \calF_2 (-I)[-2|I|] \rightarrow i^!(\calF_1 \boxtimes_{S} \calF_2)$. It suffices to check this map is an isomorphism. By continuity it suffices to consider bounded objects, and then by  conservativity we may apply $\CT_B^{I}$. Since hyperbolic localization commutes with $!$-pullback over the diagonal map $X^{I} \r X^{I \sqcup I}$, \thref{coro--SatG-T} allows us reduce to the case where $G=T$. By \thref{Sat.T} we may further reduce to the case where $T$ is trivial, i.e., the $\calF_i$ are \emph{unstratified} Tate motives on $X^{I}$ and we must compute $i^!(\calF_1 \boxtimes_{S} \calF_2)$ where $i \colon X^{I} \r X^{I \sqcup I}$ is the diagonal. 
Our map is the canonical map $i^* (\calF_1 \bx_S \calF_2) \t i^! \Z \r i^! (\calF_1 \bx_S \calF_2)$ (adjoint to the projection formula, i.e., the compatibility of !-pushforwards and $\bx$), where we have used relative purity (for $X^I$ and $X^{I \sqcup I} / S$) to compute $i^! \Z$. This map is an isomorphism for $\calF_1, \calF_2 = \Z(k)$, and therefore for all unstratified Tate motives on $X^I$.
\xpf

\subsubsection{Definition of the fiber functor}
\defi \thlabel{Fiber.Functor.Def}
Recall that $X^I \stackrel{\pi_T} \gets \Gr_{T,I} \stackrel u \r \Hck_{T,I}$ are the natural maps. Using the constant term functor from \thref{CT}, we define a functor \(F^I\) as the composite 
\[F^I := \pi_{T!} u^! \CT_B^I \colon \DMrx(\Hck_{G, I})\to \DMrx(X^I).\]
We denote its restriction to \(\Satrx^{G,I}\) the same way, in which case it takes values in unstratified mixed Tate motives by \thref{Fiber.Decomposition} and \thref{Sat.Tate}:
$$F^I := \pi_{T!} u^! \CT_B^I \colon \Satrx^{G,I} \r \MTMrx(X^I).$$
We call this restriction the \emph{fiber functor}.
(Since $X = \A^1_S$, we have $\MTMrx(X^I) \cong \MTMrx(S)$, but we prefer to write $X^I$ to emphasize the rôle of $I$.)
Using the natural isomorphism \(\pi_0(\Gr_{T,I})\cong X_*(T)\), we see that \(F^I\) decomposes as a direct sum, which we denote by \(F^I=\bigoplus_{\nu\in X_*(T)}F^I_{\nu}\).
\xdefi

We will mostly be interested in the restriction of \(F^I\) to \(\Satrx^{G,I}\), but the general functor will be useful when constructing adjoints in Subsection \ref{subsec:Hopf}. See also \thref{Remark:CTforGr} for an equivalent way of defining $F^I$.

\rema
\thlabel{fiber.functor.remarks}
By \thref{fiber.functor.properties}, the fiber functor \(F^I\) is exact, conservative, and faithful, and hence deserves its name.

Recall from \thref{Fiber.Decomposition} that the fiber functor \(F^I\) is isomorphic to $\bigoplus_{n \in \Z} \pH^n \pi_{G!} u^!$, hence independent of the choice of \(T\subseteq B\subseteq G\).
Note also that \(\pH^n \pi_{G!}u^! \cong \bigoplus_{\langle 2\rho,\nu\rangle=n}F^I_\nu\).

In the context of motives with rational coefficients in the case $I = \sgl$, the fiber functor appearing in \cite[Definition~5.11]{RicharzScholbach:Motivic} is the composite of \(F^I\) and taking the associated graded of the weight filtration. 
The weight filtration is less useful in the context of integral coefficients, e.g., $\Z/n$ is not pure of weight 0. Moreover, by not taking the associated graded we are able to construct a Hopf algebra in $\MTMrx(S)$, which is helpful for showing it is reduced, and thus independent of $S$, in \thref{H.Independent}.
\xrema

\rema
\thlabel{Fiber.Functor.local}
As in \thref{CT:local}, we can define a functor
\[F:=\pi_{T!}u^!\CT_B\colon \DMrx(L^+G\backslash \Gr_G) \to \DMrx(S),\]
which satisfies a similar decomposition \(F=\bigoplus_{\nu\in X_*(T)}F_{\nu}\).
\xrema

\subsection{Fusion}
The result below will be used to show that certain local convolution products are mixed Tate in the proof of \thref{fusion.Satake}, which says that the global convolution product constructed in \refsect{subsub-TypeII} preserves the Satake category. A different (but equivalent) construction of the monoidal structure on $\CT_B$ will be given in \thref{fiber.functor.monoidal}. 

\prop \thlabel{prop--convolution-CT}
For bounded objects $\mathcal{F}_1$, $\mathcal{F}_2 \in \DTM(L^+G \backslash LG / L^+G)$, there exists an isomorphism
$$\CT_B(\mathcal{F}_1 \star \mathcal{F}_2) \cong \CT_B(\mathcal{F}_1) \star \CT_B(\mathcal{F}_2).$$ 
\xprop

\pf
The proof uses the unipotent nearby cycles functor constructed in \cite{Ayoub:Six2} (and studied further in \cite{CassvdHScholbach:Central}).
For a scheme $Y \to \A_S^1$, let $Y_\eta$ be the fiber over $\Gm$ and let $Y_s$ be the fiber over $0$. The unipotent nearby cycles functor $\Upsilon_{Y} \colon \DM(Y_\eta) \rightarrow \DM(Y_s)$ 
is defined by Ayoub \cite[Définition~3.4.8]{Ayoub:Six2} as
$$\Upsilon = p_{\Delta\sharp} i^* j_* \theta_* \theta^* p_\Delta^*,$$
$$\left [ Y_\eta \stackrel[\id \x 1]{\Gamma_{f_\eta}} \rightrightarrows Y_\eta \x \Gm \dots \right ] \stackrel \theta \r \left [Y_\eta \stackrel[\id]{\id} \rightrightarrows Y_\eta \dots \right ] \stackrel {p} \r Y_\eta.$$
(The middle diagram is the constant cosimplicial diagram, for the full definition of the one on the left we refer to loc.~cit.; $\Gamma_{f_\eta}$ is the graph of $Y_\eta \r \Gm$.)
Finally, $Y_\eta \stackrel j \r Y \stackrel i \gets Y_s$ are extended to a constant diagram indexed by $\Delta$.
The functor $\Upsilon$ satisfies several properties similar to those of nearby cycles in non-motivic setups: 
\begin{enumerate}
  \item \label{item--push.pull} It is compatible with pushforward along a proper map and compatible with pullback along a smooth map (cf.~\cite[Définition~3.1.1]{Ayoub:Six2}).
  \item \label{item--trivial} For the trivial family $X = \A^1$, $\Upsilon (p^* M) = M$, for $p : \Gm \r S$ (cf.~\cite[Proposition~3.4.9]{Ayoub:Six2}),
  \item
For $\mathcal{F}_i \in \DM(Y_{i\eta})$, there is a K\"unneth map (cf.~\cite[§3.1.3]{Ayoub:Six2})
  $$\Upsilon_{Y_1} (\mathcal{F}_1) \boxtimes \Upsilon_{Y_2} (\mathcal{F}_2) \r \Upsilon_{Y_1 \x_{\A^1} Y_2} (\mathcal{F}_1  \boxtimes_{\Gm} \mathcal{F}_2).$$  
  It is an isomorphism if $Y_1 = Y_2 = \A^1$ and $\calF_1$ and $\calF_2$ arise by pullback from $S$, as in \refit{trivial}.
  \item \label{item--HypLocComp} If $Y$ is equipped with a Zariski-locally linearizable action of $\Gm / \A^1_S$, the formation of $\Upsilon$ commutes with hyperbolic localization as in \thref{hyperbolic.localization}.
To see this, note that the diagram $\theta$ is built out of smooth (projection) maps, so that $\theta^*$ commutes with all pushforward and pullback functors.
Also, $\theta_*$ commutes with !-pullback (by base change) and *-pushforward.
Likewise $p_\Delta^*$ commutes with all pullback and pushforward functors, hence $p_{\Delta\sharp}$ commutes with !-pushforward and *-pullback.
Finally these functors preserve the $\Gm$-equivariance of a sheaf on $Y_\eta$.
By \thref{hyperbolic.localization}, $i^*$ and $j_*$ commute with hyperbolic localization.
This shows that $\Upsilon$, when applied to $\Gm$-equivariant sheaves, commutes with hyperbolic localization. 
\end{enumerate}  
We apply this in the case $I_1 = \{1\}$, $I_2 = \{2\}$, $I = I_1 \sqcup I_2$ and $\phi = \id$ to the fiber of the top part of \refeq{TypeII} over $\A^1 \cong X_s \times X \subset X^I$. This gives a diagram as follows, where the top row consists of fibers over $\Gm$ and the bottom row consists of fibers over $0$.
$$\xymatrix{
\Gr_G \times \Gr_G \times \Gm \ar[d] & \ar[d] \restr{L_{I_1I_2}G}{X_s \times X_\eta} \ar[l] \ar[r] & \ar[d] \Gr_G \times \Gr_G \times \Gm \ar[r]^-{\id} & \Gr_G \times \Gr_G \times \Gm \ar[d] \\
\Gr_G \times \Gr_G \times \A^1  & \restr{L_{I_1I_2}G}{X_s \times X} \ar[l]_-p \ar[r]^-q &  \restr{\widetilde{\Gr}_{G, \phi}}{X_s \times X} \ar[r]^-m & \restr{\Gr_{G,I}}{X_s \times X} \\
\Gr_G \times \Gr_G \ar[u] & \ar[l] \ar[u] \ar[r] \restr{L_{I_1I_2}G}{X_s \times X_s} & LG \times^{L^+G} \Gr_G \ar[r]^-m \ar[u] & \Gr_G \ar[u]
}$$
We claim that $$\Upsilon_{\restr{\Gr_{G,I}}{X_s \times X} }(\mathcal{F}_1 \boxtimes \mathcal{F}_2 \boxtimes \Z_{\Gm}) \cong \mathcal{F}_1 \star \mathcal{F}_2.$$ 
Granting this claim, by applying \refit{HypLocComp} and the compatibility of hyperbolic localization with box products, we get an isomorphism
$$\Upsilon_{\restr{\Gr_{T,I}}{X_s \times X} }(\CT_B(\mathcal{F}_1) \boxtimes \CT_B(\mathcal{F}_2) \boxtimes \Z_{\Gm}) \cong \CT_B(\mathcal{F}_1 \star \mathcal{F}_2).$$ Recall that $(\Gr_G \times \Gr_G \times \Gm)^0_{\text{red}} = (\Gr_T \x \Gr_T \x \Gm)_{\text{red}} = \bigsqcup_{X_*(T)^2} \Gm$.
The closure of the copy of $\Gm$ indexed by $(\nu_1, \nu_2)$ is isomorphic to $\A^1$, and its special fiber is the point in $(\Gr_T)_{\text{red}}$ indexed by $\nu_1 + \nu_2$. Moreover, the convolution structure 
on $\DTM(L^+T \backslash LT / L^+T) \cong \Fun(X_*(T), \DTM(L^+T \backslash S))$ is induced by the abelian group structure on $X_*(T)$ and the monoidal structure on $\DTM(S)$. Thus, by \refit{trivial}, the left side of the above isomorphism is canonically identified with $\CT_B(\mathcal{F}_1) \star \CT_B(\mathcal{F}_2)$, as desired.

To prove the claim, by applying \refit{push.pull} to $m$ it suffices to check
$$\Upsilon_{ \restr{\widetilde{\Gr}_{G, \phi}}{X_s \times X}}(\mathcal{F}_1 \boxtimes \mathcal{F}_2 \boxtimes \Z_{\Gm}) \cong  \mathcal{F}_1 \widetilde \boxtimes \mathcal{F}_2.$$ 
Note that $\Upsilon_{\Gr_G \times \A^1}(\mathcal{F}_i \boxtimes \Z_{\Gm}) = \mathcal{F}_i$.
This follows from the compatibility of $\Upsilon$ with proper !-pushforward (along $\Gr_G^{\le \mu} \r \Gr_G$) and \refit{trivial} (applied to the base scheme being $\Gr_G^{\le \mu}$). Since $p$ and $q$ admit sections Zariski-locally, then by compatibility with smooth pullback it suffices to show that the $(L^+G \times L^+G)$-equivariant K\"unneth map 
$$\Upsilon_{\Gr_G \times \A^1}(\mathcal{F}_1 \boxtimes \Z_{\Gm})  \boxtimes \Upsilon_{\Gr_G \times \A^1}(\mathcal{F}_2 \boxtimes \Z_{\Gm}) \to \Upsilon_{\Gr_G \times \Gr_G \times \A^1} (\mathcal{F}_1 \boxtimes \mathcal{F}_2 \boxtimes \Z_{\Gm})$$ is an isomorphism.
Note that this map is $(L^+G \x L^+G)$-equivariant since $\Upsilon$ commutes with smooth pullback, see also \cite[Corollary~3.30]{CassvdHScholbach:Central}.
The formation of this K\"unneth map is compatible with hyperbolic localization, so by conservativity of $\CT_B$ (\thref{CT.conservative}) and the geometry of $(\Gr_T)_{\text{red}} = \bigsqcup_{X_*(T)} S$ recalled above, we are reduced to checking that the K\"unneth map is an isomorphism for Tate motives on $\A^1$, which holds by \refit{trivial} above. 
\xpf

\rema
The K\"unneth map for unipotent nearby cycles is usually not an isomorphism, but it is an isomorphism for the full nearby cycles 
(see \cite[Théorème~10.19]{Ayoub:Realisation} and \cite[Théorème~3.5.17]{Ayoub:Six1} over a field of characteristic 0). Since the nearby cycles we are considering have unipotent monodromy in the Betti and \'etale contexts, e.g. \cite[Proposition 2.4.6]{AcharRiche:Central}, it is unsurprising that we have a K\"unneth isomorphism. However, unlike \emph{loc. cit.} we never have to consider full nearby cycles.
\xrema

\coro \thlabel{coro--t-exact}
Let $\mathcal{A}_1$, $\mathcal{A}_2 \in \MTM(L^+G \backslash LG / L^+G)$. 
\begin{enumerate}
\item We have $\mathcal{A}_1 \star \mathcal{A}_2 \in \DTM(L^+G \backslash LG / L^+G)^{\leq 0}$.
\item Moreover, if $\mathcal{A}_1 \boxtimes \mathcal{A}_2 $ is mixed Tate then $\mathcal{A}_1 \star \mathcal{A}_2 \in \MTM(L^+G \backslash LG / L^+G)$.
\end{enumerate}
\xcoro

\pf
It suffices to consider bounded objects.
We then apply the conservativity of $\CT_B$ and \thref{prop--convolution-CT} to reduce to the case $G=T$, where the result follows from right exactness of $\otimes$. 
\xpf

\rema
For rational coefficients the above corollary can be proved using conservativity and t-exactness of $\ell$-adic realization as in \cite[Lemma 5.8]{RicharzScholbach:Motivic}. 

An alternative approach would be to use cellularity of the fibers of the convolution morphism as in \cite[Remark 2.5.4]{dHL:Frobenius} and \cite{Haines:Pavings}, along with the classical fact that the convolution morphism is semi-small \cite{Lusztig:Singularities, MirkovicVilonen:Geometric}.
This allows us to show convolution is right t-exact.
However, we can currently only use this to show left t-exactness when we have a suitable description of the coconnective part of the t-structure.
This includes the case of reduced motives, and by \cite[Proposition 5.3]{EberhardtScholbach:Integral} also the case \(S=\Spec \mathbf{F}_p\) and motives with \(\mathbf{F}_p\)- or \(\Q\)-coefficients.
We omit details.
\xrema

\subsubsection{Standard and costandard motives} \label{sect--(co)standard}
In the proof of \thref{fusion.Satake} we will reduce to considering certain flat objects in $\Satrx^{G,I}$.
To handle these, we adapt the method of (co)standard sheaves \cite[\S 1.11]{BaumannRiche:Satake} to a motivic context. 

In this subsection, we only consider the local case, i.e., $\Gr_G$ instead of $\Gr_{G,I}$. 
Using the isomorphism $\Gr_{G,\{*\}} = \Gr_G \times X$ and homotopy invariance, the results in this subsection immediately yield similar results for Beilinson--Drinfeld Grassmannians in the case $I = \{*\}$.

\defi \thlabel{Standard.Def}
The \emph{standard} and \emph{costandard functors} are the functors $\MTMrx(S) \r \MTMrx(L^+G \setminus \Gr_G)$ defined as
$$\mathcal{J}^\mu_! := \pe \iota^\mu_!(p_{\mu}^*(-) [\dim \Gr_G^\mu]), \quad \mathcal{J}^\mu_* := \pe \iota^\mu_*(p_{\mu}^*(-)[\dim \Gr_G^\mu]),$$ where $p_{\mu} \colon L^+G \setminus \Gr_{G}^{\mu} \r S$ and $\iota^\mu \colon  L^+G \setminus \Gr_G^\mu \r  L^+G \setminus \Gr_G$.
\xdefi

Here $p_\mu^* : \DMrx(S) \r \DMrx(L^+G \setminus \Gr_G^\mu)$ denotes the functor whose composition with the forgetful functor to $\DMrx(\Gr_G^\mu)$ is the usual functor $p_\mu^*$ (and the components in further terms of the \v Cech nerve of the $L^+G$-action, as in \refeq{equivariant.DM}, are given by !-pullbacks along the action, resp.~projection maps).

\prop \thlabel{Standard.Flat}
For $\mu \in X_*(T)^+$ and $L \in \MTMrx(S)$ there is a functorial isomorphism 
$$ \calJ^\mu_!(L)  \stackrel \cong  \r \calJ^\mu_!(\Z) \otimes L$$
(
The expression at the right denotes the action of $L \in \DTMrx(S)$ on $\DTMrx(L^+G \setminus \Gr_G)$, i.e., it is a \emph{derived} tensor product.) Furthermore, $F(\calJ^\mu_!(\Z))$ is identified with a finitely generated free graded abelian group under the faithful embedding $\grAb \r \MTMrx(S)$ (cf.~\thref{functor i}).
\xprop

\pf
The map is obtained by applying $\tau^{\geq 0}$ to $\iota^\mu_! (p_\mu^*(\Z) [\dim \Gr_G^\mu]) \t L = \iota^\mu_! (p_\mu^*(L) [\dim \Gr_G^\mu])$, which we claim lies in $\DTMrx(L^+G \setminus \Gr_G)^{\leq 0}$. Being supported on $\Gr_G^{\le \mu}$, this motive is bounded, so the result follows from \thref{fiber.functor.properties} and the computation in \thref{Standard.Free}.
\xpf

\lemm \thlabel{Standard.Free}
The composite $F_\nu \calJ^\mu_!$ is isomorphic to the endofunctor on $\MTMrx(S)$ given by
$$L \mapsto F_\nu \calJ^\mu_!(\Z) \otimes L,$$
where $F_\nu\calJ^\mu_!(\Z)$ is a free graded abelian group of rank equal to the number of irreducible components of $\Ss_\nu^+ \cap \Gr_G^\mu$.
\xlemm

\pf
By t-exactness of $\CT_B$ and base change we have
$$F_\nu \calJ^\mu_! \cong \pH^{\langle 2\rho,\nu\rangle}((q_{\nu}^+)_! (p_\nu^+)^* \iota^\mu_! p_\mu^* [\langle 2\rho,\mu\rangle]) = \pH^{\langle 2 \rho,\mu + \nu\rangle} f_! f^*,$$ where $f \colon  \Ss_\nu^+ \cap \Gr_G^\mu \rightarrow S$.
By \thref{cellularity of intersection:torus}, $\Ss_\nu^+ \cap \Gr_G^\mu$ is cellular of equidimension
$\langle \rho,\mu + \nu\rangle$ relative to $S$, so this expression computes the top-dimensional cohomology group with compact support. Now the result follows from \thref{lemm--cellular-coh}. 
\xpf

\exam
Following up on \thref{Example PGL2} we consider $G= \PGL_2$. In this case $F_\nu \calJ^\mu_! L = L (-\frac{\mu + \nu}2)$ if $|\nu| \le \mu$ and $\nu \equiv \mu \mod 2$. In all other cases, $F_\nu \calJ^\mu_! = 0$. 
\xexam

\rema \thlabel{rema--costandard}
An argument dual to the one in \thref{Standard.Flat} involving the cohomology of the dualizing complex shows that $F( \calJ^\mu_*(\Z))$ is also a finitely generated free graded abelian group. However, since $\otimes$ is not left exact, it is not immediate that $\calJ^\mu_*(L) \cong \calJ^\mu_*(\Z) \otimes L$. Cases where we can verify that this isomorphism holds are as follows:
\begin{itemize}
\item $L \in \MTMrx(S)$ is flat in the sense that $(-) \otimes L$ is t-exact.
\item $\Ss_\nu^+ \cap \Gr_G^\mu$ has a cellular stratification for all $\nu$, as a opposed to a filtrable decomposition, so that the union of top-dimensional cells is open and the excision computation in \thref{lemm--cellular-coh} simplifies considerably. (Note that there is always an isomorphism for Betti and \'etale sheaves since one is not concerned with cellularity and Tateness.)
\item For reduced motives, since $\Z$ has global dimension $1$ and $f_* f^! \Z$ is free in degrees $-2d$ and $-2d+1$, where $d = \langle \rho,\mu + \nu\rangle$ (since it is Verdier dual to $f_! f^* \Z \in \DTMr(S)^{\leq 2d}$ which is free in degree $2d$).
\end{itemize}
\xrema

\prop \thlabel{Cok.N}
For $\mu \in X_*(T)^+$ the canonical surjection $\calJ^\mu_!(\Z) \rightarrow \IC_{\mu,\Z}$ is an isomorphism. Furthermore, the formation of $\IC_{\mu,\Z}$ commutes with Betti realization and reduction.
\xprop

\pf
By \thref{Standard.Free} and \thref{rema--costandard}, the map $F(\calJ^\mu_!(\Z)) \rightarrow F(\calJ^\mu_*(\Z))$ identifies with a map of finitely generated free graded abelian groups. Thus, to check that it is injective we may apply $- \t_\Z \Q$ and restrict to reduced motives. In this case we note that the category of compact objects $\MTMr(\Gr_{G}, \Q)^\comp$ is semisimple: by \thref{DTMr.independence} (applicable by \thref{BD.WT}), we may assume $S = \Spec \Fp$ for this, and then apply \thref{DTM.Q} and \cite[Corollary~6.4]{RicharzScholbach:Motivic}. Thus, the natural morphisms of reduced motives $\calJ^\mu_!(\Q) \r \IC_{\mu,\Q}$ and $\IC_{\mu,\Q} \r \calJ^\mu_*(\Q)$ are isomorphisms, so we indeed have an injection and consequently $\calJ^\mu_!(\Z) \cong \IC_{\mu,\Z}$. 

The compatibility of $\calJ^\mu_! \Z$ with reduction and Betti realization (which is not automatic given that the reduction functor is not t-exact) can be checked after applying the family of conservative functors $F_\nu$.
Then it follows from \thref{Standard.Free}.
\xpf

\subsubsection{Construction of the fusion product} 
In this subsection we prove that the perverse truncation of the convolution product $\star \colon \DTMrx(\Hck_{G,I}) \x \DTMrx(\Hck_{G,I} ) \r \DTMrx(\Hck_{G,I} )$ in \thref{Convolution:Global:Tate} induces a symmetric monoidal structure on $\Satrx^{G,I}$, compatibly with constant terms. As in previous treatments of geometric Satake, we proceed by relating convolution to the restriction of a fusion product $* \colon \Satrx^{G,I} \x \Satrx^{G,I} \r \Satrx^{G,I \sqcup I}$ along the diagonal $X^I \r X^{I \sqcup I}$. The fusion product is first defined for general motives via a global convolution product. Over an open subscheme of $X^{I \sqcup I}$ it is simply a box product thanks to factorization  \refeq{ZhuFactor}. The key result is \thref{fusion.Satake}, which in the case of $\Satrx^{G,I}$ allows us to express fusion as an intermediate extension. 

In what follows the reader should keep in mind the example of the decomposition $I \sqcup I$, corresponding to the surjection $I \sqcup I \r \{1,2\}$ sending one copy of $I$ to $1$ and the other to $2$. More generally, let \(\phi\colon I\twoheadrightarrow J\) denote a surjection of nonempty finite index sets, where $J$ is ordered, and write \(I_j:=\phi^{-1}(j)\) for \(j\in J\). 
Recall the setup and notation from \thref{Twisted.Prod.Global} (and its generalization to $|J|$-fold convolution products): for \(\Ff_j\in \DMrx(\Hck_{G,I_j})\), we defined the twisted product \(\widetilde \boxtimes_{j\in J} \Ff_j \in \DMrx(\widetilde{\Gr}_{G,\phi})\). 
Recall also the convolution morphism \(m_\phi\colon \widetilde{\Gr}_{G,\phi}\to \Gr_{G,I}\). 
By \thref{Convolution:Global:Tate}, the functor $m_{\phi!}\widetilde \boxtimes_{j\in J}(-)$ preserves Tate motives. 
We use the same notation to denote the restriction of this functor to the product of Satake categories.

\theo
\thlabel{fusion.Satake}
The functor \(m_{\phi!}\widetilde \boxtimes_{j\in J}(-) \colon \prod_{j\in J}\Satrx^{G,I_j} \to \DTMrx(\Gr_{G,I})\) is right t-exact, and $\pH^n (m_{\phi!}\widetilde \boxtimes_{j\in J}(-))$ restricts to a functor \(\prod_{j\in J} \Sat_{\redx}^{G,I_j} \to \Sat_{\redx}^{G,I}\) for all $n \leq 0$. 
\xtheo

\pf We proceed by several reductions. \\

\emph{Reduction to bounded objects and $|J|$=2.}
Since the relevant functors commute with filtered colimits we may restrict to bounded objects.
There is a natural associativity constraint on $m_{\phi!}\widetilde \boxtimes_{j\in J}(-)$ coming from the associativity of $\boxtimes$ and proper base change. 
This can be constructed using the general version of $L_{I_1 I_2}G$ for any number of factors in analogy to the associativity constraint on the local Satake category, cf.~\cite[Lemma 3.7]{RicharzScholbach:Motivic}. 
Associativity and the closure of the Satake category under subquotients and extensions in \(\MTMrx(\Gr_{G,I})\) (\thref{Sat.Props}) allows us to reduce to the case $J= \{1,2\}$. \\

\emph{Standard objects.} We take a brief detour and consider certain standard objects in $\Satrx^{G,I}$; see \refsect{(co)standard} for more details when $I = \{*\}$.
For $\mu \in (X_*(T)^+)^I$, let $j \colon \Gr_{G,I}^{\circ, \mu} \r \restr{\ol{\Gr_{G,I}^{\circ, \mu}}}{X^\circ}$ be the inclusion (take $\phi = \id$ in \thref{BD.strata}).
Here $j$ is a product of the identity map on $X^\circ$ with the embeddings $\Gr_{G}^{\mu_i} \r \Gr_G^{\leq \mu_i}$, and $j^\circ \colon \restr{\ol{\Gr_{G,I}^{\circ, \mu}}}{X^\circ} \r \ol{\Gr_{G,I}^{\circ, \mu}}$. 
Fix $L \in \MTMrx(S)$, and let $\calJ_{\mu !}^I (L) = j^{\circ}_{!*} {}^pj_! (L[\dim \Gr_{G,I}^{\circ, \mu} ])$.
By a proof similar to \thref{CT.Sat}, $\calJ_{\mu !}^I (L) \in \Satrx^{G,I}$.
We claim that $$\restr{\calJ_{\mu !}^I (L)}{X^\circ} \cong \restr{\boxtimes_{i \in I} \calJ_{\mu_i!}^{\{i\}}(\Z) \otimes L}{X^\circ}. \eqlabel{Fusion.1}$$ Indeed, since $\Gr_{G}^I \cong \Gr_{G^I}$, then $\restr{\calJ_{\mu !}^I (L)}{X^\circ} \cong \restr{\calJ_{\mu !}^I (\Z)}{X^\circ} \otimes L$ by \thref{Standard.Flat}.
This reduces us to $L = \Z$, where the claim follows from the K\"unneth formula for $!$-pushforward 
(\refsect{motives}), right t-exactness of $\boxtimes$, and flatness of the $\CT_B^{\{i\}} \calJ_{\mu_i!}^{\{i\}}(\Z)$ (see \thref{Standard.Free} or \cite[Proposition 1.11.1]{BaumannRiche:Satake}). 
By applying $\CT_B^I$ in order to check the left side below lies in $\Satrx^{G,I}$, it follows from \refeq{Fusion.1} and \thref{Sat.Int.Ext} that $$\calJ_{\mu !}^I (\Z) \otimes L \cong \calJ_{\mu !}^I (L). \eqlabel{Fusion.2}$$ 

\emph{Reduction to standard objects.}
Suppose we have proved the theorem for the standard objects as above; we now show how to deduce the general case. Let $I = I_1 \sqcup I_2$, and let $\phi \colon I \r \{1,2\}$ map $I_1$ to $\{1\}$ and $I_2$ to $\{2\}$.
Fix a standard object $\Ff_2$.
We will prove the theorem holds for all (bounded) $\Ff_1$ by induction on the support of $\restr{\Ff_1}{X^\circ}$. 
Thus, we may assume $\Ff_1$ has a filtration with subquotients given by IC-motives.
Again, using the closure of the Satake category under subquotients and extensions, it suffices to consider the case $\Ff_1 =  \IC_{\mu, L}$ for arbitrary $\mu = (\mu_i) \in (X_*(T)^+)^{I_1}$ and $L \in \MTMrx(S)$.
The base case occurs when $\mu_i$ is minuscule for all $i \in I_1$, so $\Gr_{G}^{\mu_i} = \Gr_{G}^{\leq \mu_i}$ and we conclude since $\Ff_1$ is a standard object.
In general, there is an exact sequence $0 \r K \r \calJ_{\mu !}^{I_1}(L) \r \IC_{\mu, L} \r 0$ in $\Satrx^{G,I_1}$, so this case follows by applying induction to $K$ and the case of standard objects.
 To finish the proof we may consider a similar induction on the support of $\restr{\Ff_2}{X^\circ}$. \\
 
 \emph{Reduction to integral standard objects.}
We claim that $$\calJ_{\mu !}^{I}(\Z) \cong m_{\phi !} (\calJ_{\mu_1}^{\{1\}} (\Z) \widetilde \boxtimes \ldots \widetilde \boxtimes \calJ_{\mu_{|I|}}^{\{|I|\}} (\Z)) \in \Satrx^{G,I}.$$  Granting the claim, by associativity of  $m_{\phi!}\widetilde \boxtimes_{j\in J}(-)$ we have $\calJ_{\mu^1 \sqcup \mu^2 !}^{I_1 \sqcup I_2}(\Z) \cong m_{\phi !}(\calJ_{\mu^1 !}^{I_1}(\Z) \widetilde \boxtimes \calJ_{\mu^2 !}^{I_2}(\Z)).$ Since all functors appearing are $\DTMrx(S)$-linear, by \refeq{Fusion.2} it follows that for all $L_1, L_2  \in \MTMrx(S)$, $$m_{\phi !}(\calJ_{\mu^1 !}^{I_1}(L_1) \widetilde \boxtimes \calJ_{\mu^2 !}^{I_2}(L_2)) \cong \calJ_{\mu^1 \sqcup \mu^2 !}^{I_1 \sqcup I_2}(\Z) \otimes (L_1 \otimes L_2).$$ By the claim, applying $\pH^n$ to the right side gives an element of $\Satrx^{G,I}$ for all $n \in \Z$, and by right exactness of $\otimes$ we get $0$ if $n \geq 0$. 

It remains to prove the claim. We can assume that $\phi = \id$. Let $\mathcal{F}_i = \calJ_{\mu_i}^{\{i\}} (\Z)$. Write $\mathcal{F}_i = \mathcal{F}_i' \boxtimes \Z_X[1]$ where $\mathcal{F}_i' = \calJ_{\mu_i !}(\Z)$ is a standard object in $\MTMrx(L^+G \backslash LG / L^+G)$. By \refeq{Fusion.1}, both objects in the claimed isomorphism agree over $X^\circ$. Since $\calJ_{\mu !}^{I}(\Z)$ is an intermediate extension, the claim will follow from a computation of the cohomological degrees of the $*$- and $!$-pullbacks of $m_{\phi!}\widetilde \boxtimes_{i\in I}(\Ff_i)$ over the strata of $X^I$. These strata are indexed by surjections $\psi \colon I \twoheadrightarrow K$. For a given $\psi$, let
 $X^\psi \subset X^I$ be the corresponding stratum, and let $i \colon \Gr_{G,I}^{\psi} \rightarrow \Gr_{G,I}$ be the inclusion. Let $I = \{1, \ldots, n\}$.
Since the $\Ff_i$ are bounded, the twisted product $\widetilde \boxtimes_{i\in I} \Ff_i$ may be formed using descent with respect to $*$-pullback instead of $!$-pullback. 
By base change and factorization \refeq{ZhuFactor}, we have
\begin{equation} \label{eq--Sat1} i^* m_{\phi!}\widetilde \boxtimes_{i\in I}(\Ff_i) \cong \underset{k \in K}{\boxtimes} (\underset{i \in \psi^{-1}(k)}{\star} \Ff_i') \boxtimes \Z_{X^\psi}[n].\end{equation} For $i^!$, we note that by \thref{WT.prod} and relative purity for $\iota \colon X^\psi \rightarrow X^I$, the inclusion $(\id \times \iota) \colon \Gr^n \times X^\psi \rightarrow \Gr^n \times X^I$ has the following property (this is a straightforward generalization of \thref{TypeII.Relation.iso} to more than two factors) $$(\id \times \iota)^!((\boxtimes_{i \in I} \Ff_i') \boxtimes \Z_{X}[n]) \cong (\boxtimes_{i \in I} \Ff_i') \boxtimes \Z_{X^\psi}(|K|-n)[2|K|-n].$$ By base change it follows that \begin{equation} \label{eq--Sat2} i^! m_{\phi!}\widetilde \boxtimes_{i\in I}(\Ff_i) \cong i^* m_{\phi!}\widetilde \boxtimes_{i\in I}(\Ff_i)(|K|-n)[2|K|-2n].\end{equation} Note that the open stratum $X^\circ$ occurs when $|K| = n$, and we have $|K| < n$ for all other strata. Thus, the claim will follow if we show that $ \underset{k \in K}{\boxtimes} (\underset{i \in \psi^{-1}(k)}{\star} \Ff_i')$ is mixed Tate, since then the \eqref{eq--Sat1} lies in degree $|K|-n < 0 $ for all strata $X^\psi \neq X^\circ$, and \eqref{eq--Sat2} lies in degree $n-|K| > 0$. For this, it suffices to show that $\CT_B(\underset{i \in \psi^{-1}(k)}{\star} \Ff_i')$ is mixed Tate and has a filtration with subquotients given by free graded $\Z$-modules. This follows from induction on $|\psi^{-1}(k)|$ and \thref{prop--convolution-CT}, starting with  \thref{Standard.Free} when $|\psi^{-1}(k)| = 1$.
\xpf

\rema
The reason for using standard objects in the proof of \thref{fusion.Satake}, instead of IC-motives, is that \refeq{Fusion.1} and \refeq{Fusion.2} are in general false for IC-motives due to torsion. For coefficients in a field, we could work directly with IC-motives. 
\xrema

\rema
\thlabel{Remark:Sign:Constraint}
Recall the fully faithful functor $j^{(\phi),*} \colon \Satrx^{G,I} \rightarrow \MTMrx(\restr{L^+_IG \backslash  \Gr_{G,I}}{X^{(\phi)}}) $ from \thref{Sat.Restr.FF}. 
By \thref{Sat.Int.Ext} and \thref{fusion.Satake}, there is a natural isomorphism
$$\pe m_{\phi!}\widetilde \boxtimes_{j\in J}(-) \cong j^{(\phi)}_{!*} (\pH^0(\restr{\boxtimes_{j\in J}(-)}{X^{(\phi)}})) \colon \prod_{j\in J} \Satrx^{G, I_j} \to \Satrx^{G,I}.$$
This functor satisfies natural commutativity and associativity constraints, induced from those of the exterior products over $X^{(\phi)}$. 
However, this naive commutativity 
constraint is not compatible with that of $\prod_{j\in J} \MTM(X^{I_j}) \xrightarrow{\boxtimes} \MTM(X^I)$ under the fiber functors $F^{I_j}$, as the two will differ by some signs. 
To correct this we modify the commutativity constraint by hand as in \cite[VI.9.4 ff.]{FarguesScholze:Geometrization}.
Namely, let us decompose 
\(\Gr_{G,I} = \Gr_{G,I}^{\mathrm{even}} \coprod \Gr_{G,I}^{\mathrm{odd}}\) into open and closed subsets, where 
\(\Gr_{G,I}^{\mathrm{even}}\) is the union of the Schubert cells corresponding to \((\mu_i)_i\in (X_*(T)^+)^I\) for 
which \(\sum_{i\in I} \langle 2\rho,\mu_i\rangle\) is even, and likewise for \(\Gr_{G,I}^{\mathrm{odd}}\). This induces a similar decomposition of \(\Hck_{G,I}\).
Then, we change the commutativity constraint by adding a minus sign when commuting the exterior product of motives concentrated on \(\Gr_{G,I}^{\mathrm{odd}}\).
If we denote the resulting functor equipped with this commutativity constraint by \(\ast\), we have changed the signs such that the diagram
\[\begin{tikzcd}
  \prod_{j\in J} \Satrx^{G,I_j} \arrow[d, "\prod_{j\in J} F^{I_j}"'] \arrow[r, "\ast"] & \Satrx^{G,I} \arrow[d, "F^I"]\\
  \prod_{j\in J} \MTMrx(X^{I_j}) \arrow[r, "\boxtimes"] & \MTMrx(X^I)
\end{tikzcd}\]
is functorial in the $I_j$ and under permutations of $I_1, \ldots, I_{|J|}$; this follows from the implicit shifts appearing in the fiber functors, via \thref{Fiber.Decomposition}. The diagram is commutative because this can be checked over $X^{(\phi)}$ by full faithfulness of $j^{(\phi),*}$, where it is immediate because the twisted product in the definition of $\ast$ becomes a box product. Here we also use the \((-)_!(-)^*\)-description of \(\CT_B^I\) and the Künneth formula, as well as \thref{piT}.
\xrema

\defi
\thlabel{fusion.product}
The functor $\ast \colon \prod_{j\in J} \Satrx^{G,I_j} \to \Satrx^{G,I}$ in  \thref{Remark:Sign:Constraint} equipped with the modified commutativity constraint is called the \emph{fusion product}.
\xdefi

\defi
\thlabel{otimes.MTM}
We endow the category $\MTMrx(X^I)$ with the tensor product defined as $\pH^0( (\mhyphen ) \otimes (\mhyphen )[-I])$.
The normalization ensures that the monoidal unit is $\Z[I]$. 
We refer to this as the \emph{underived} tensor product. 
\xdefi

\prop \thlabel{fiber.functor.monoidal}
Let $$ \mhyphen \, {\pe\star} \, \mhyphen := \pH^0( \mhyphen \star \mhyphen ) \colon \MTMrx(\Hck_{G,I}) \x \MTMrx(\Hck_{G,I} ) \r \MTMrx(\Hck_{G,I} )$$ be the perverse truncation of the functor constructed in \thref{Convolution:Global:Tate} (recall it includes a shift by $[-I]$, cf.~\refeq{GrI.Conv}).
Then $(\Satrx^{G,I}, \pe \star)$ has the structure of a symmetric monoidal category coming from the fusion product (constructed in the proof), and the constant term functors $\CT_P^I$ are symmetric monoidal. Moreover, $F^I : \MTMrx(\Hck_{G,I}) \r \MTMrx(X^I)$ is a symmetric monoidal functor, where the tensor structure on the target is as in \thref{otimes.MTM}.
\xprop

\pf The construction is analogous to \cite[VI.9.4.~ff.]{FarguesScholze:Geometrization}; see also \cite[Theorem 3.24]{Richarz:New} for more details.
Briefly, for an integer $n > 0$, consider the natural surjection \(\phi:\sqcup_{i=1}^n I\to n\) 
and the diagonal embedding $i_{\ol \phi} \colon \Gr_{G,I} \r \Gr_{G,I^n}$. 
Let $d_\phi = |I|^{n-1}$. 
The fusion product and \thref{Sat.Tate} determine a functor 
\[ \Satrx^{G,I}\times \ldots \times \Satrx^{G,I} \xrightarrow{\ast} \Satrx^{G, I\sqcup \ldots \sqcup I} \xrightarrow{i_{\ol \phi}^!(d_\phi)[d_\phi]} \Satrx^{G,I},\eqlabel{Sat.fusion}\]
which makes \(\Satrx^{G,I}\) into a symmetric monoidal category by \thref{Remark:Sign:Constraint}.
By \thref{TypeII.Relation.iso} and its generalization to more than two factors, this convolution product agrees with $\pe \star$, so $(\Satrx^{G,I}, \pe \star)$ is a symmetric monoidal category.

Over \(X^{(\phi)}\subseteq X^I\), the functor $\CT_P^I$ decomposes as $\prod_{j\in J}\CT_P^{I_j}$. As in \thref{Remark:Sign:Constraint}, there is a diagram expressing the compatibility between fusion and the constant term functors, which is also functorial in the $I_j$ and under permutations of $I_1, \ldots, I_{|J|}$. 
Base change for $!$-pullback along the diagonal embedding $X^I \r X^{I\sqcup \ldots \sqcup I}$ then shows that $\CT_P^I$ and \(F^I\) are symmetric monoidal.
\xpf

\subsubsection{Dualizability}
\label{sect--dualizability}
The following \thref{lemm-dualizing}, which is similar to \cite[\S 11]{MirkovicVilonen:Geometric} and \cite[VI.8.2]{FarguesScholze:Geometrization}, will be used in \thref{Satake.Hopf} to establish inverses for the dual group. 

Let $\sw$ be the involution of $\DM(\Hck_{G,I})$ induced by $!$-pullback along the inversion map $L_IG \r L_IG, g \mapsto g^{-1}$. 
For a prestack $Z \stackrel \pi \r Y$ over a smooth $S$-scheme $Y$, we consider the functor $\Du_{Z/Y}(\calF) := \IHom(\calF, \pi^! \Z_Y[2 \dim_S Y])$, for $\calF \in \DM(Z)$.
We will refer to it as the \emph{relative Verdier duality} functor (but note that $\Du_{Z/Y}(\calF)(\dim_S Y)$ is the usual (absolute) Verdier duality functor on $Z$).
We write $\Du_Z := \Du_{Z/S}$ and we also sometimes omit the subscript in $\Du_Z$ if the choice of $Z$ is clear. 

\defi
\thlabel{Hck.lc}
Following \cite[\S12.2.3]{ArinkinGaitsgory:Singular}, the subcategory of \emph{locally compact motives} $\DTMrx(\Hck_{G,I})^\locc \subset \DTMrx(\Hck_{G,I})$ is the full subcategory of motives whose image in $\DTMrx(\Gr_{G,I})$ lies in the subcategory $\DTMrx(\Gr_{G,I})^\comp$ of compact objects.
\xdefi

\lemm \thlabel{CT-compact}
Let $\calF \in \DTM(\Hck_{G, I})$ be bounded. Then $\calF$ is locally compact if and only if $u^! \CT_B^I(\calF) \in \DTM(\Gr_{T,I})$ is compact.
\xlemm

\pf
By boundedness it suffices to consider the restriction of the top row of \refeq{hyperbolic.diagram} to actual schemes as opposed to ind-schemes.
For maps of schemes, the !-pushforward and *-pullback functors preserve compact objects, so this implies that $u^! \CT_B(\calF)$ is compact.

For the converse, we use noetherian induction on the support of $\calF$. We will also write $\calF$ for the underlying motive on $\Gr_{G,I}$. Let $\iota \colon \Gr_{G,I}^{\phi, \mu} \r \Gr_{G,I}$ be a stratum open in the support of $\calF$, where $\phi \colon I \twoheadrightarrow J$ and $\mu = (\mu_j)_{j\in J} \in (X_*(T)^+)^J$. The intersection  $\Gr_{G,I}^{\phi, \mu} \cap \Gr_{T,I}$ is a disjoint union of copies of \(X^\phi\); let us fix the copy corresponding to the anti-dominant representative \(w_0(\mu) =(w_0(\mu_j))_{j\in J}\).
Since $\Ss_{w_0(\mu_i)}^+ \cap \Gr_G^{\mu_i} = S$,
the restriction of $u^! \CT_B^I(\calF)$ to $X^\phi \subset \Gr_{T,I}$ is (up to a shift) identified with the restriction of $\iota^* \calF$ to $X^\phi$ (see also the proof of \thref{CT.conservative}). In particular, the $*$-fiber of $\iota^* \calF$ at $X^\phi$ is compact.
Since $\calF$ is $L_I^+G$-equivariant, then $\iota^* \calF$ is the constant motive on $\Gr_{G,I}^{\phi, \mu}$ given by spreading out its $*$-fiber at $X^\phi$ by \cite[Lemma 2.2.1]{RicharzScholbach:Intersection} and \thref{Stab.Split}. Thus $\iota^* \calF$ and also $\iota_! \iota^* \calF$ are compact. To conclude by induction, it suffices to observe that the cofiber of $\iota_! \iota^* \calF \r \calF$ has the following properties: it is $L_I^+G$-equivariant, has compact constant terms, and has support smaller than that of $\calF$.
\xpf

\prop
\thlabel{duality.Hck.lc}
Verdier duality (relative to $X^I$) is an anti-equivalence on $\DTMrx(\Hck_{G,I})^\locc$. 
\xprop

\pf 
We show that the natural map $\alpha_M : M \r \Du(\Du(M))$ is an isomorphism for $M \in \DTMrx(\Hck_{G,I})^\locc$. Since $M$ is locally compact it is in particular bounded.
As in \cite[Proposition~IV.6.13]{FarguesScholze:Geometrization}, $\Du$ commutes with hyperbolic localization up to taking the inverse of the $\Gm$-action.
Since  $\CT_B^I$ is conservative (\thref{CT.conservative}) and detects Tateness (\thref{CT.DTM}) when restricted to bounded objects, and preserves locally compact objects (\thref{CT-compact}), we have reduced to the case $G=T$.

We now show that $\Du$ is an involution on $\DTMrx(\Gr_{T,I})^\comp$. 
The closures of all strata $\Gr_{T,I}^{\phi, \mu}$ are smooth (in fact, they are affine spaces), so $\Du(M)$ is Tate and $M \r \Du(\Du(M))$ is an isomorphism when $M = i_* \Z(n)$ for $i \colon {\overline{\Gr_{T,I}^{\phi, \mu}} } \rightarrow \Gr_{T,I}$.
Since these $M$ generate $\DTMrx(\Gr_{T,I})^\comp$, it follows that Verdier duality is an anti-equivalence on $\DTMrx(\Gr_{T,I})^\comp$.

Now, for $M \in \DTMrx(\Hck_{T,I})^\locc$, we consider $u^! M \in \DTMrx(\Gr_{T,I})^\comp$. We then have $u^! \alpha_M = \alpha_{u^! M}$; this uses that the closure of each stratum $\overline{\Hck_{T,I}^{\phi, \mu}}$ is a quotient of $\overline{X^\phi}$ by a pro-smooth group scheme.
Since $u^!$ is conservative, we are done.
\xpf

\rema
While our convention is that $\CT_B^I$ is defined on $\DTMrx(\Hck_{G,I})$, it can also be defined on $\DTMrx(\Gr_{G,I})$, and then the above argument shows that Verdier duality is an anti-equivalence on $\DTMrx(\Gr_{G,I})^\comp$.
\xrema

\lemm \thlabel{lemm-dualizing}
The dualizing functor with respect to the (derived) convolution product $\star$ on $\DTMrx(\Hck_{G, I})^\locc$, i.e., the internal Hom functor $\IHom_\star(-, 1)$, is given by
$$\Du^- := \Du_{\Hck_{G,I} / X^I} \circ \sw.\eqlabel{internal.Hom}$$
An object $M \in \DTMrx(\Hck_{G,I})^\locc$ is dualizable iff the resulting natural map
$$M \star N \r \Du^-(\Du^- M \star \Du^- N) \eqlabel{star.duality}$$
is an isomorphism with $N = \Du^-(M)$.
\xlemm

\pf 
We have a cartesian diagram of prestacks over $X^I$ as follows.
$$\xymatrix{
L_I^+G \setminus L_IG/L^+_I G  \ar[r]^-{\text{inv}} \ar[d]^f & L^+_IG \setminus L_IG \times^{L_I^+G} L_IG / L_I^+G \ar[d]^m \ar[r]^-p & (L_I^+G \setminus L_IG/L^+_I G)^2 \\
L^+_IG \setminus L^+_IG / L^+_IG \ar[r]^{i} & L_I^+G \setminus L_IG/L^+_I G }$$
Here $f$ is the quotient by \(L^+_IG\) of the structural map $L_IG/L^+_IG \to L^+_IG/L^+_IG\cong X^I$, the map $i$ is induced by the inclusion $L^+_IG \r L_IG$, and $\text{inv}$ is induced by the identity map on the first factor and inversion on the second factor. Let $p$ be the natural quotient map (obtained by modding out two copies of $L^+_IG$ instead of one acting diagonally).
Then the composite $\inv^* p^! = \Delta^* (\id \x \sw)$, where $\id \x \sw$ is the involution of $\DM(\Hck_{G,I}^2)$ induced by inversion on the second factor. 

The first claim is a formal manipulation similar to, say, \cite[Lemma~9.10]{BoyarchenkoDrinfeld:Character}: 
$$\eqalign{
\Hom_{\DM(\Hck_{G,I})}(M, \Du^-(N)) 
& = \Hom_{\DM(\Hck_{G,I})} (\inv^* p^! (M \boxtimes N), f^! \Z[2I])  \cr
& = \Hom_{\DM(X^I / L_I^+ G)} (f_! \inv^* p^! (M \boxtimes N), \Z[2I])  \cr
& = \Hom_{\DM(X^I / L_I^+ G)} (i^* m_! p^! (M \boxtimes N), \Z[2I])  \cr
& = \Hom_{\DM(X^I / L_I^+ G)} (i^* (M \star N)[I], \Z[2I])  \cr
& =  \Hom_{\DM(\Hck_{G,I})} (M \star N, i_* \Z[I])  \cr
& =  \Hom_{\DM(\Hck_{G,I})} (M \star N, 1). }$$
Given this and \thref{duality.Hck.lc}, the category $\DTMrx(\Hck_{G,I})^\locc$ is thus an r-category in the sense of \cite[Definition~1.5]{BoyarchenkoDrinfeld:Duality}, which gives a natural morphism as in \refeq{star.duality}.
The final claim is then just \cite[Corollary~4.5]{BoyarchenkoDrinfeld:Duality}.
\xpf

\coro
\thlabel{Sat.dualiziability}
The formation of derived duals is compatible with constant terms in the sense that for $\calF \in \DTM(\Hck_{G,I})$ there is a natural isomorphism
$$\CT_B^I \IHom_\star^G(\calF, 1) = \IHom_\star^T(\CT_{B^-}^I \calF, 1).$$
\xcoro

\pf
Recall that $\CT_B^I \circ \Du = \Du \circ \CT_{B^-}^I$ by \cite[Proposition~IV.6.13]{FarguesScholze:Geometrization}.
Also, $\CT_B^I$ commutes with $\sw$ because the diagram \refeq{hyperbolic.diagram} is induced by homomorphisms of loop groups. Hence our claim holds by \thref{lemm-dualizing}.
\xpf

\section{Tannakian reconstruction}
\label{sect--Tannakian}

\subsection{The Hopf algebra object}\label{subsec:Hopf}

We continue to work with a base scheme $S$ as in \thref{nota-basescheme}; whenever we consider categories of mixed Tate motives (denoted $\MTM$), we assume it satisfies the Beilinson--Soulé vanishing (cf.~\refeq{BS.vanishing}).
For reduced mixed Tate motives ($\MTMr$), that latter condition is not needed.

The goal of this section is to construct a Hopf algebra object $H^{G,I}_{\redx} \in \MTMrx(X^I)$ such that the Satake category is equivalent to comodules over that Hopf algebra (\thref{Satake.Hopf}).
Based on the results of the previous sections, the Satake category appears for formal reasons in a monadic adjunction (\thref{comonadic.adjunction}). 
Several steps, including an analysis of standard motives, are needed to show the relevant monad is given by tensoring with a Hopf algebra.
We refer to \cite[§1.12-1.13]{BaumannRiche:Satake} 
and \cite[§VI.10]{FarguesScholze:Geometrization} for constructions of similar Hopf algebras for different sheaf theories, where the motivic difficulties below do not appear.

Throughout, $I$ denotes a nonempty finite set, 
and $W \subset (X_*(T)^+)^I$ a finite subset closed under the Bruhat order.
Let $\Gr_{G,I}^W$ be the closure of the union of the strata $\Gr_{G,I}^{\circ, \mu}$ for $\mu \in W$, and let $i_W : \Gr_{G,I}^W \r \Gr_{G,I}$ be the closed embedding. 
In addition to the Satake category $\Satrx^{G,I}$, we consider its full subcategory $\Sat_{\redx, W}^{G,I}$ consisting of motives supported on $\Gr_{G,I}^W$.

\subsubsection{Adjunctions between motives on the Hecke prestack and on the curve}

We establish a left adjoint for \(F^I:=\pi_{T!} u^! \CT_B^I\) (\thref{Fiber.Functor.Def}). Recall that the restriction of this functor to \(\Satrx^{G, I}\) is the fiber functor. It is also possible to construct a right adjoint, but it is easier to prove properties of the left adjoint since $\otimes$ is also right exact.

\lemm 
\thlabel{DM.Hck.adjoints}
The restriction of $F^I$ to 
$\DMrx(\Hck_{G,I}^W)$ admits a left adjoint given by
$$L^I_W = \coav p^-_{W!} q^{-*}_W \pi_{T,W}^* [-\deg] \eqlabel{fiber.right.DTM}$$
Here the subscript $W$ denotes the restriction of \refeq{hyperbolic.diagram} to $\Gr_{G,I}^W \subset \Gr_{G,I}$, and $\Gr_{T,I}^W = \Gr_{T,I} \cap \Gr_{G,I}^W$.
\xlemm

\pf
We will show that the restriction of the forgetful functor $u^! : \DMrx(\Hck_{G, I}) \r \DMrx(\Gr_{G,I})$ to $\DMrx(\Hck_{G,I}^W)$ admits a left adjoint $\coav_W$. 
Once this is shown, the formula in \refeq{fiber.right.DTM} defines the left adjoint by \thref{Remark:CTforGr}.
Note that the restriction of
\refeq{hyperbolic.diagram} to $\Gr_{G,I}^W \subset \Gr_{G,I}$ consists of maps of schemes (as opposed to ind-schemes), so that the functors appearing in the definitions of $L^I_W$ (notably $q_W^{-*}$) are well defined on the categories $\DMrx$.

For each $W$, the $L^+_I G$-action on $\Gr_{G,I}^W$ factors over a smooth algebraic quotient group $L^+_I G \twoheadrightarrow H_W$, such that the kernel of this quotient map is split pro-unipotent by \cite[Lemma A.3.5, Proposition A.4.9]{RicharzScholbach:Intersection}.
By the computation of equivariant motives in \cite[Proposition~3.1.27]{RicharzScholbach:Intersection}, 
$\DMrx(\Hck_{G,I}^W) = \DMrx(H_W \setminus \Gr_{G,I}^W)$.
Thus, the coaveraging functor $\coav_{H_W}$ from \thref{averaging.functor} is left adjoint to $u^!$, as claimed.
\xpf

\prop
\thlabel{DTM.Hck.X.adjunction} 
The adjunction $L^I_W \dashv F^I$ restricts to an adjunction on the categories $\DTMrx(\Hck_{G,I}^W)$ and $\DTMrx(X^I)$.
\xprop

\pf
We need to show that $L^I_W$ preserves stratified Tate motives.
The composite $p_{W !}^{-} q_W^{-*} \pi_{T,W}^*$ maps $\DTMrx(X^I, (X^I)^\dagger)$ to the category of motives on $\Gr_{G,I}$ that are Tate motives with respect to the stratification by the $\Ss_{\nu,I}^- \cap \Gr_{G,I}^{\phi,\mu}$. 
It remains to be shown that $\coav$ maps these motives to Tate motives with respect to the (coarser!) stratification by the $\Gr_{G,I}^{\phi,\mu}$.

The formation of $\coav$ commutes with $*$-pullback over the strata $X^\phi \subset X^I$. Using the factorization property \refeq{ZhuFactor} and the fact that $\Gr_G^J \cong \Gr_{G^J}$, we only need to consider $\Gr_G$. If $i \colon \Ss_\nu^- \cap \Gr_G^\mu \r \Gr_G^\mu$ is the inclusion, it suffices to show that $u^! \coav i_! \Z$ is Tate. 
Now, this motive is $L^n G$-equivariant for $n \gg 0$ (namely, such that the \(L^+G\)-action on \(\Gr_G^\mu\) factors through \(L^nG\)); let \((L^nG)_{w_0(\mu)}\subseteq L^nG\) be the stabilizer of \(t^{w_0(\mu)}\in \Gr_G^\mu\).
Then, the equivalence \(\DM(L^nG\backslash \Gr_G^\mu) \cong \DM((L^nG)_{w_0(\mu)}\backslash S)\) restricts to an equivalence of Tate motives by \cite[Propositions 1.1 and 1.3]{RicharzScholbach:IntersectionCorrigendum} (the latter of which applies by \thref{Stab.Split}). 
Thus, it suffices to show its *-restriction along the base point $t^{w_0(\mu)} \colon S \r \Gr_G^{\mu}$ is a Tate motive.
Using \thref{averaging.functor} to compute this *-restriction, it suffices to show that $f_!(\Z) \in \DM(S)$ is Tate, where $f \colon a^{-1}(t^{w_0(\mu)}) \r S$ and $a \colon L^nG \times (\Ss_\nu^- \cap \Gr_G^\mu) \r \Gr_G^\mu$ is the action map. The fiber $a^{-1}(t^{w_0(\mu)})$ has a filtrable decomposition by the preimages of the cells $X_w \subset \Ss_\nu^- \cap \Gr_G^\mu$ (as in \thref{cellularity of intersection:torus}) under the projection onto the right factor in the source of the map $a$. Thus, by excision we may replace $\Ss_\nu^- \cap \Gr_G^\mu$ by a single cell $X_w$ in the source of $a$.
Since each cell $X_w$ is contained in some $X_\delta$ as in \thref{Av.Fiber}, then the latter implies $a^{-1}(t^{w_0(\mu)}) \cong \calP_{w_0(\mu)}^n \times X_w$, where \(\calP_{w_0(\mu)}^n\subseteq L^nG\) is the stabilizer of \(t^{w_0(\mu)}\). 
Furthermore, $\calP_{w_0(\mu)}^n$ is an extension of a split reductive $\Z$-group by a split unipotent $\Z$-group (\thref{Stab.Split}). For the structural map $\pi : \calP^n_{w_0(\mu)} \r S$, we have $\pi_! \pi^* \Z \in \DTM(S)$ by virtue of the cellular Bruhat stratification, and since $X_w$ is a cell we conclude that $f_!(\Z)$ is Tate.
\xpf

\coro
\thlabel{Corollary:Adjoints}
\begin{enumerate}
\item
\label{item--adjoints.MTM}
The preceding adjunction restricts to an adjunction
$$
\pe L^I_W = \pH^0 \coav p^-_{W!} q^{-*}_W \pi_{T,W}^* [-\deg] : \MTMrx(X^I) \rightleftarrows \MTMrx(\Hck_{G,I}^W) : F^I. \eqlabel{adjunction.MTM} $$

\item
\label{item--av.reduction}
The functors $F^I$ and $L^I_W$ are compatible with the reduction functor $\rho_\red$.
\end{enumerate}
\xcoro

\pf
\refit{adjoints.MTM}: 
This holds since $F^I$ is t-exact (\thref{fiber.functor.remarks}) so that $L^I_W$ is right t-exact.

\refit{av.reduction}: The functor $\rho_\red$ commutes with $F^I$ since the latter is a composite of the standard six functors (cf.~\refsect{reduced.motives}).
It also commutes with $L_W^I$ by \thref{averaging.functor}\refit{av.DM.red}.
\xpf

\prop
\thlabel{adjoints.box}
The adjunction $L^I_W \dashv F^I$  for Tate motives is compatible with the exterior product in the following sense.
Consider the diagrams for $k=1,2$
$$X^{I_k} \stackrel{\pi_k} \gets \Gr_{T, I_k}^{W_k} \stackrel{q^+_k} \gets \Gr_{B, I_k}^{W_k} \stackrel{p^+_k} \r \Gr_{G, I_k}^{W_k} \stackrel {u_k} \r \Hck_{G, I_k}$$
as in \refeq{hyperbolic.diagram}, where $I_k = \sgl$ is a singleton. 
Let $L_k := L^{W_k}_{I_k}$ etc.~and write
$$\eqalign{
F_{12} & := (\pi_1 \x \pi_2)_! (q^+_1 \x q^+_2)_! (p^+_1 \x p^+_2)^* (u_1 \x u_2)^!, \cr
L_{12} & := \coav_{12} (p^-_1 \x p^-_2)_! (q^-_1 \x q^-_2)^* (\pi_1 \x \pi_2)^*,}$$
where $\coav_{12} : \DTMrx(\Gr_{G, I_1} \x \Gr_{G, I_2}) \r \DTMrx(\Hck_{G, I_1} \x \Hck_{G, I_2})$ is the left adjoint to $(u_1 \x u_2)^!$.
Then there are isomorphisms (of functors $\DTMrx(X^{I_1}) \x \DTMrx(X^{I_2}) \r \DTMrx(\Hck_{G, I_1} \x_S \Hck_{G, I_2})$)
$$L_{12}(- \boxtimes -) \r L_1 (-) \boxtimes L_2(-).\eqlabel{L.box}$$
$$F_1(-) \boxtimes F_2(-) \r F_{12}(- \boxtimes -).\eqlabel{F.box}$$
\xprop

\pf
First of all, the existence of $\coav_{12}$ as stated is proven exactly the same way as for the single coaveraging functors.
The isomorphism \refeq{F.box} exists  
since \(\boxtimes\) is compatible with *-pullbacks and !-pushforwards, and also with \(\coav_{12}\) by \thref{averaging.functor}\refit{av.box}.
By adjunction, this gives the map \refeq{L.box}.
The former map is an isomorphism even when evaluated on $\DMrx(\Hck_{G, I_k})$ again because $\boxtimes$ commutes with *-pullbacks and !-pushforwards.
\xpf

\subsubsection{Adjunction for the Satake category}

We now construct an adjunction involving the fiber functor $F^I\colon \Sat_{\redx, W}^{G,I}\to \MTMrx(X^I)$. 
The left adjoint $L_{\Sat, W}^I$ will computed explicitly in \thref{R.Sat.explicit}.

\prop
\thlabel{comonadic.adjunction}
There is a monadic adjunction
$$L_{\Sat, W}^I \colon  \MTMrx(X^I) \rightleftarrows \Sat_{\redx, W}^{G, I}\colon F^I.\eqlabel{adjunction.Sat}$$
In other words, there is an equivalence 
$$\Sat^{G,I}_{\redx, W} = \Mod_{T^I_W}(\MTMrx(X^I)),$$
where $T^I_W := F^I \circ L_{\Sat, W}^I$ is the monad induced by the adjunction and $\Mod$ denotes the category of modules over that monad.
Here $\MTMrx(X^I)$ denotes the category of mixed unstratified Tate motives on $X^I$.
\xprop

\pf
The restriction of $F^I$ to $\Sat^{G,I}_{\redx, W}$ takes values in unstratified Tate motives by \thref{Sat.Tate}.
On the level of $\DTM(\Gr_{G,T}^W)$, the functor $\pi_{T!} = \pi_{T*}$ is right adjoint to $\pi_T^*$, hence it preserves products.
The constant term functor also preserves products.
The functor $F^I$ is exact by \thref{fiber.functor.properties}, and therefore preserves all limits. 
Both categories are presentable so the adjoint functor theorem guarantees the existence of a left adjoint $L_{\Sat, W}^I$. 

In addition, $F^I$ is conservative again by \thref{fiber.functor.properties}.
Being exact, it also preserves finite colimits, so that the adjunction is monadic by the Barr--Beck monadicity theorem.
\xpf

Our eventual goal is to show that $T^I$ preserves limits. We start with the case $I=\{*\}$. 

To begin, fix $\nu, \lambda \in X_*(T)$ and \(\mu\in X_*(T)^+\) such that $\Ss_\nu^- \cap \Gr_G^{\mu} \neq \emptyset$ and $\Ss_\lambda^+ \cap \Gr_G^{\mu} \neq \emptyset$.
For $n \gg 0$, we have the action map $a \colon L^nG \times (\Ss_\nu^- \cap \Gr_G^{\mu}) \r \Gr_G^{\mu}$.
Let $Y(\mu, \lambda, \nu) = a^{-1}(\Ss_\lambda^+ \cap \Gr_G^{\mu})$ for one such $n$.

\prop \thlabel{Y.Strat}
We have $\dim Y(\mu, \lambda, \nu) = \dim L^nG - \langle\rho,\nu\rangle + \langle\rho,\lambda\rangle$ and this scheme has a filtrable cellular decomposition.
\xprop

\pf 
For the dimension, note that $\dim \Ss_\lambda^+ \cap \Gr_{G}^{\mu} = \langle\rho,\mu + \lambda\rangle$, $\dim \Ss_\nu^- \cap \Gr_G^{\mu} = \langle\rho,\mu - \nu\rangle$, and $\dim \calP_\mu^n = \dim L^nG - 2 \langle\rho,\mu\rangle$.

For the decomposition, it suffices to show that for every cell $X_w \subset \Ss_\lambda^+ \cap \Gr_G^{\mu}$ (as in \thref{cellularity of intersection:torus}), the fiber $a^{-1}(X_w)$ has a filtrable cellular decomposition. The map $L^nG \r \Gr_G^{\mu}$, $g \mapsto g \cdot t^\mu$ has a section \(s\colon X_w\to L^nG\) over $X_w$ by \thref{Torsor is trivial over decomposition}; note that we are currently working with positive semi-infinite orbits.
Let $F_0$ be the fiber of $a$ over $t^\mu$. Using the left action of $L^nG$ on $L^nG \times (\Ss_\nu^- \cap \Gr_G^{\mu})$, there is an isomorphism $X_w  \times F_0 \r a^{-1}(X_w)$, $(x, f) \mapsto s(x) \cdot f$. Since $X_w$ is a cell, it suffices to decompose $F_0$. 
Using the transitive \(L^nG\)-action on \(\Gr_G^\mu\), we get an isomorphism \(F_0=a^{-1}(t^\mu) \cong a^{-1}(t^{w_0(\mu)})\).
But this was shown to admit a filtrable cellular decomposition in the proof of \thref{DTM.Hck.X.adjunction}.
\xpf

\prop \thlabel{prop--L.free}
For $A \in \MTMrx(X)$ there is a canonical isomorphism $$L_W^{\{*\}}(A) \cong L_W^{\{*\}}(\Z[1]) \otimes A[-1].$$ Furthermore, $F^{\{*\}}(L_W^{\{*\}}(\Z[1]))$ is a finitely generated free graded abelian group.
\xprop

\pf
By restricting to connected components of $\Gr_{T, \{*\}}$, we get a finite direct sum decomposition $L_{W}^{\{*\}} = \oplus_{\nu \in \widetilde W} L_{W}^{\nu}$. Let $\widetilde a \colon L^nG \times (\Ss_\nu^- \cap \Gr_{G}^{W}) \r \Gr_{G}^{W}$ be the action map and let $\widetilde Y (W, \lambda, \nu) = \widetilde{a}^{-1}(\Ss_\lambda^+ \cap \Gr_G^{W})$ for $\lambda \in \widetilde{W}$. 
Then by \thref{averaging.functor} we have $$F_{\lambda}^{\{*\}} (L_{W}^{\nu}(\Z[1])) = f_! \Z(\dim L^nG)[2\dim L^nG - \langle2\rho,\nu\rangle + \langle2\rho,\lambda\rangle + 1],$$ where $f \colon \widetilde Y(W, \lambda, \nu) \times X \r X$ is the projection, cf.~\cite[Proposition 1.12.1]{BaumannRiche:Satake}.
Here we are using the isomorphism $\Gr_{G, \sgl} = \Gr_G \times X$.
The scheme $\widetilde Y(W, \lambda, \nu)$ is stratified by the $Y(\mu, \lambda, \nu)$ for $\mu \in W$, so by \thref{Y.Strat}, $\widetilde Y(W, \lambda, \nu)$ has a filtrable cellular decomposition, and it has dimension $\dim L^nG - \langle\rho,\nu\rangle + \langle\rho,\lambda\rangle$. Thus, $F_{\lambda}^{\{*\}} (\pe L_{W}^{\nu}(\Z[1]))$ computes the top cohomology of a cellular scheme, so by \thref{lemm--cellular-coh} it is identified with a finitely generated free graded abelian group. Moreover, as in the proof of \thref{Standard.Flat} it follows that $L_W^{\{*\}}(A) \cong L_W^{\{*\}}(\Z[1]) \otimes A[-1]$.
\xpf

\prop 
\thlabel{R.Sat.explicit}  Suppose that $W = \prod_i W_i$. 
\begin{enumerate}
\item 
\label{item--R.Sat.explicit(1)}
The left adjoint in \refeq{adjunction.Sat} is given by 
$$L^I_{\Sat, W} = j_{!*} j^* \pe L^I_W$$
where $j : \Hck_{G,I}^{\circ, W} \r \Hck_{G,I}^W$ is the inclusion over the open locus of pairwise distinct points.
\item 
\label{item--TGI.Cok}
For $A \in \MTMrx(X^I)$ there is a functorial isomorphism 
$$L^I_{\Sat, W}(A) \cong *_{i \in I} \pe L_{W_i}^{\{i\}}(\Z[1]) \otimes A[-I].$$
\end{enumerate}
\xprop

\pf
\refit{R.Sat.explicit(1)} 
On the larger categories of stratified mixed Tate motives, the left adjoint was computed in \thref{Corollary:Adjoints}\refit{adjoints.MTM} as $\pe L^I_{W}$. However, this functor need not send unstratified mixed Tate motives on $X^I$ to  $\Sat_{\redx, W}^{G, I}$. It therefore suffices to modify this functor to make the last property hold, in such a way that the group $\Hom(\pe L_W^I(A), \calF)$ is unchanged if $\calF \in \Satrx^{G,I}$. We claim that $j_{!*} j^* \pe L^I_W$ does the job.

First, suppose we have shown that $j_{!*} j^* \pe L_W^I$ takes values in the subcategory $\Satrx^{G,I} \subset \MTMrx(\Hck_{G,I})$. Then we may conclude using \thref{Ext.To.Ext2}: for unstratified  $A \in \MTMrx(X^I)$ and $\calF \in \Satrx^{G,I}$ we have
$$\Hom(\pe L_W^I(A), \calF) \cong \Hom(j_{!*} j^* \pe L_W^I(A), \calF).$$ 
In order to apply \thref{Ext.To.Ext2}, we must show that $\pe L^I_{W}(A)$ has no quotients supported over $X^I \setminus X^\circ$. To prove this,
let $i \colon \Hck_{G,I}^W \setminus \Hck_{G,I}^{\circ, W} \r \Hck_{G,I}^W$ be the complement of $j$. 
We claim that if $\calF \in \DTMrx^{\leq 0}(X^I)$ is an unstratified Tate motive, then $i^* L_W^I(\calF) \in \DTMrx^{\leq-1}(\Hck_{G,I}^W \setminus \Hck_{G,I}^{\circ, W})$. 
By base change applied to $\coav p^-_!$ (\thref{averaging.functor}\refit{av.base.change}), the formation of $i^* L_W^I$ commutes with *-restriction to any of the $\binom{|I|}{2}$ hyperplanes which comprise $X^I \setminus X^\circ$. 
Hence the claim that there are no quotients supported over $X^I \setminus X^\circ$ follows from the right t-exactness of $L_W^{I}$ and $i^*$. 

Now to complete the proof of \refit{R.Sat.explicit(1)}, it remains to show that $j_{!*} j^* \pe L_W^I$ takes values in the subcategory $\Satrx^{G,I} \subset \MTMrx(\Hck_{G,I})$.
Over $X^\circ$, we have a K\"unneth formula for the left adjoints by \thref{adjoints.box}. Hence $j^* L_W^I(A) \cong \restr{\boxtimes_{i \in I} L^{\{i\}}_{W_i} (\Z) \otimes A}{X^\circ}$. 
We claim that the natural maps $L^{\{i\}}_{W_i} (A) \r \pe L^{\{i\}}_{W_i} (A)$
induce an isomorphism $\boxtimes_{i \in I} \pe L^{\{i\}}_{W_i} (\Z[1]) \otimes A[-I] \cong \pH^0(\boxtimes_{i \in I}  L^{\{i\}}_{W_i} (\Z[1]) \otimes A[-I])$.
Arguing as in the proof of \thref{Standard.Flat}, this follows from the freeness of the $F^{\{i\}}(L^{\{i\}}_{W_i} (\Z[1]))$, proved in \thref{prop--L.free}.
Thus, 
$$\restr{\boxtimes_{i \in I} \pe L^{\{i\}}_{W_i} (\Z[1]) \otimes A[-I]}{X^\circ} \cong j^* ( \pe L_W^I(A)). \eqlabel{R.I.Fusion}$$
Moreover, if $|I|=1$, we have $\Satrx^{G, \sgl} = \MTMrx(\Hck_{G, \sgl})$. 
Now we conclude by using that fusion preserves the Satake categories by \thref{fusion.Satake}. 

\refit{TGI.Cok}:
By $\DTMrx(X^I)$-linearity, there is an isomorphism $L_W^I(\Z[I]) \otimes A[-I] \cong L_W^I(A)$ functorial in $W$ and $A$.
Using \refeq{R.I.Fusion}, applying $j_{!*} \pe j^*$ produces the desired isomorphism.
\xpf

\prop
\thlabel{left adjoint dualizable}
For $W = \prod_i W_i  \subset (X_*(T)^+)^I$ the object $L^I_{\Sat, W}(\Z[I])$ is dualizable in $\DTMrx(\Hck_{G,I})$ (i.e., with respect to the derived convolution product) and also in $\Satrx^{G,I}$ (i.e., with respect to the underived convolution product). 
\xprop

\pf
Consider the subcategory $\mathcal{C} \subset \DTMrx(\Hck_{G, I})^\locc \cap \Satrx^{G,I}$ consisting of objects $M$ such that $\CT_B^I(M)$ and $\CT_{B^-}^I(M)$ are given by finite free graded abelian groups on each connected component of $\Gr_{T, I}$ under the equivalence $\Satrx^{T,I} = \Fun(X_*(T)^I, \MTMrx(X^I))$. 
We claim that on $\calC$, $\star$ and its truncation $\pe \star$ agree, and turn $\calC$ into a rigid monoidal category. 
To see this, let $M, N \in \calC$. We show that the derived convolution product $M \star N \in \Satrx^{G,I}$. By \thref{TypeII.Relation.iso}, the motive underlying $M \star N$ agrees with  $i^! m_{\phi !}(M \widetilde \boxtimes N)(I)[I]$, where $i \colon \Gr_{G,I} \r \Gr_{G,I \sqcup I}$ is the diagonal embedding and $m_{\phi} \colon {\widetilde \Gr}_{G,\phi} \r \Gr_{G,I \sqcup I}$ is the convolution map associated to $\phi \colon I \sqcup I \twoheadrightarrow \{1,2\}$. 
Now by \thref{fusion.Satake} we have $m_{\phi !}(M \widetilde \boxtimes N) \in \Satrx^{G,I \sqcup I}$, since our assumption implies the constant terms of the restriction of this motive to $(X^{I \sqcup I})^\circ$ are external tensor products of finite free graded abelian groups. Thus $M \star N \in \Satrx^{G,I}$ by \thref{Sat.Tate}.

We now show that $\calC$ is preserved under $D^- := \IHom_\star(- , 1)$.
For $M \in \calC$, by \thref{Sat.dualiziability} we are reduced to showing that $D^-(\CT_{B}^I M) \in \Satrx^{T,I}$. By assumption $\CT_B^I M$ is a finite free graded abelian group on each connected component of $\Gr_{T, I}$. Then $D^-(\CT_B^I M) \in \Satrx^{G,I}$ is obtained by taking the dual of each finite free graded graded abelian group and swapping each $\mu \in X_*(T)^I$ with its inverse $-\mu$.

The fiber functor $F^I \colon \mathcal{C} \rightarrow \MTMrx(X^I)$ is monoidal with respect to $\pe \star$ by \thref{fiber.functor.monoidal}. It also intertwines the duality functor $D^-(-)$ on $\mathcal{C}$ with the usual duality functor on the subcategory of $\MTMrx(X^I)$ consisting of finite free graded abelian groups: By \thref{Sat.dualiziability} and since $F^I = \bigoplus_{n\in \Z} \pH^n \pi_{G,!} u^! = \pi_{T, !} u^! \CT_B^I $ is independent of the choice of $B$, this claim reduces to the fact that $\pi_{T,!}$ commutes with Verdier duality and is invariant under $\sw$. 

Hence $F^I$ gives a tensor functor out of the r-category $\mathcal{C}$ which intertwines dualizing functors. Since the internal Hom functor in an r-category is determined by the dualizing functor and the monoidal structure, cf. \cite[Remark 9.12(iii)]{BoyarchenkoDrinfeld:Character}, $F^I$ also preserves inner Hom's (with respect to the truncated convolution product). 
Since $F^I$ is conservative, it therefore detects dualizability.
By definition, $F^I(\calC)$  consists of dualizable objects in $\MTM(X^I)$.
Therefore, any object in $\calC$ is dualizable, i.e., $\calC$ is rigid.

We now show that $\calF := L^I_{\Sat, W}(\Z[I])$ lies in $\calC$.
(This argument would be a lot simpler if the truncation functors for the motivic t-structure preserve compact objects, which we do not presently know in sufficient generality, cf.~\thref{compact t-structure}.)
If $I = \sgl$ this follows from \thref{prop--L.free} and \thref{CT-compact}. For general $I$, let us unravel \thref{R.Sat.explicit}\refit{TGI.Cok}. If $\phi$ is the identity map of $I$, we have the associated convolution map $m_{\phi} \colon \widetilde{\Gr}_{G,\phi} \rightarrow \Gr_{G,I}$, and the motive underlying $L^I_{\Sat, W}(\Z[I])$ is $\pe m_{\phi !} \widetilde \boxtimes_{i \in I} L^{\{i\}}_{W_i}(\Z[1])$ by \thref{Remark:Sign:Constraint}. We claim that the perverse truncation is redundant, so that $L^I_{\Sat, W}(\Z[I])$ is obtained from applying the six functors to the $L^{\{i\}}_{W_i}(\Z[1])$, and hence it is compact. To see that $m_{\phi !} \widetilde \boxtimes_{i \in I} L^{\{i\}}_{W_i}(\Z[1])$ is already mixed Tate, note that every perverse cohomology sheaf of this motive lies in $\Satrx^{G,I}$ by \thref{fusion.Satake}. Thus, the claim may be checked over $X^\circ$, where it follows from the fact that the $\CT_B^{\{i\}}L^{\{i\}}_{W_i}(\Z[1])$ are finite free graded abelian groups on each connected component of $\Gr_{T, \{i\}}$, and hence their external tensor product is mixed Tate.
\xpf

From now on we assume that $W$ is of the form  $W = \prod_i W_i  \subset (X_*(T)^+)^I$.

\coro
\thlabel{object.flat}
The monad $T_W^I = F^I L^I_{\Sat, W}$ is given by tensoring with $T_W^I(\Z[I])$.
This object, which a priori is an algebra object in $\MTMrx(X^I)$, is actually a finitely generated free graded abelian group. 
\xcoro

\pf
By \thref{R.Sat.explicit}\refit{TGI.Cok} we have $T_W^I (A) = F^I L^I_{\Sat, W} (A) = F^I \left ( *_{i \in I} \pe L_{W_i}^{\{i\}}(\Z[1]) \otimes A[-I] \right )$. The $\DTMrx(X^I)$-linearity and monoidality of $F^I$ (as in \thref{Remark:Sign:Constraint}), as well as the fact that each $F^{\{i\}}L_W^{\{i\}}(\Z[1])$ a finite free graded abelian group (\thref{prop--L.free}), imply that this
is isomorphic to $\boxtimes_{i \in I} F^{\{i\}}L_W^{\{i\}}(\Z[1]) \t A[-I]$. In particular, the tensor product agrees with the underived one, i.e., it is computed in $\MTMrx(X^I)$.
\xpf

\subsubsection{The Hopf algebra} 
We now show that once we dualize the monadic bounded-level adjunctions in \thref{comonadic.adjunction}, they assemble to a comonadic adjunction for the global (unbounded) Satake category.

\theo
\thlabel{Satake.Hopf}
The fiber functor 
$$F^I: (\Satrx^{G,I}, \pe \star) \r (\MTMrx(X^I), \pe \t)$$
is comonadic.
The associated comonad on $\MTMrx(X^I)$ is given by tensoring with 
$$H_{\redx}^{G,I} := \colim_W T^I_W(\Z[I])^\dual.$$ 
This coalgebra is in fact a Hopf algebra, so that we obtain an equivalence of symmetric monoidal categories
$$\Satrx^{G,I} = \coMod_{H_{\redx}^{G,I}}(\MTMrx(X^I)).\eqlabel{Sat.comod}$$
\xtheo

\pf
In the sequel, all (co)modules will be understood to be in the category $\MTMrx(X^I)$ of unstratified mixed Tate motives (with its underived tensor product).
By \thref{object.flat}, $\Sat_W := \Sat^{G,I}_{\redx, W}$ identifies with the category of modules over the algebra $A_W = T^I_W(\Z[I])$.
Crucially, this object is finite free, and therefore dualizable.
Thus, its (underived or derived) dual $H_W = \IHom(A_W, \Z[I])$ is a coalgebra. 
We have the following diagram,
where $(i_{WW'})_!$ are, in the language of (co)modules, the canonical restriction functors,
and $(i_W)_!$ is the canonical insertion functor,
while $F^I_W := F^I|_{\Sat_W}$ is the fiber functor, which identifies with forgetting the (co)module structures:
$$\xymatrix{
\coMod_{H_W} = \Mod_{A_W} = \Sat_W \ar@<.5ex>@{^{(}->}[r]^(.45){(i_{WW'})_!} \ar[dr]_{F^I_W} &
\coMod_{H_{W'}} = \Mod_{A_{W'}} = \Sat_{W'} \ar@<.5ex>@{^{(}->}[r]^(.65){(i_{W'})_!} \ar[d]^{F^I_{W'}} \ar@<.5ex>[l]^{\pe i_{WW'}^!} &
\Sat := \colim \Sat_W \ar[dl]^{F^I = \colim F^I_W} \ar@<.5ex>[l]^(.35){\pe i_{W'}^!} \\
& \MTMrx(X^I)
}$$
Note that $F^I_W$ has a left adjoint, denoted by $L_W$, given by tensoring with $A_W$.
It also has a right adjoint, denoted by $R_W$, given by the cofree comodule, i.e., by tensoring with $H_W$.
We have $\pe i_{WW'}^! R_{W'} = R_W$.
By \refeq{eqn.Sat.colim.lim}, $F^I$ then has a right adjoint $R$, which is such that $\pe i_W^! R(c) = R_W(c)$ for each $W$.
According to the standard presentation of objects in $\Sat$, cf.~\refeq{obj MTM colimit}, this means that
$$R(c) = \colim_W i_{W!} R_W(c).$$

By \thref{fiber.functor.properties}, $F^I$ and $F^I_W$ are exact and faithful.
Thus, $(F^I, R)$ is a comonadic adjunction, and $\Sat$ identifies with comodules over the coalgebra $H_{\redx}^{G,I}$ as claimed. 

The functor $F^I$ is monoidal with respect to the (truncated) convolution product $\pe \star$ and the (underived) tensor product $\pe \otimes$ on $\MTMrx(X^I)$.
Therefore, $R$ is lax symmetric monoidal, which makes $H_{\redx}^{G,I} = F^IR(\Z[I])$ into an algebra object.
A routine, if tedious argument (e.g.,~\cite[Theorem~7.1]{Moerdijk:Monads} dualized; the multiplication map $H_{\redx}^{G,I} \t H_{\redx}^{G,I} \r H_{\redx}^{G,I}$ for the algebra structure is dual to the one in Definition~1.1 there) shows that $H_{\redx}^{G,I}$ is then a bialgebra, and hence \refeq{Sat.comod} holds on the level of symmetric monoidal categories.

As in \cite[Proposition VI.10.2]{FarguesScholze:Geometrization}, the antipode of $H_{\redx}^{G,I}$ is constructed using the dualizability (in $\Satrx^{G,I}$) of the objects $L^I_{\Sat,W}(\Z[I])$ (\thref{left adjoint dualizable}).
\xpf

\coro
\thlabel{Hopf:Alg:Exterior}
There is a natural isomorphism of Hopf algebras
$$\otimes_{i \in I} H_{\redx}^{G, \{i\}} \cong H_{\redx}^{G,I}.$$
\xcoro

\pf
This follows from \thref{R.Sat.explicit}\refit{TGI.Cok} and the fact that $F^I$ is symmetric monoidal.
\xpf

\subsubsection{Rational and modular coefficients} 

In the identification of the dual group below, we also need to work with $\Q$- and $\Fp$-coefficients.
Let $\Lambda = \Q$ or $\Fp$.
We define the category 
$$\Satrx^{G, \Lambda} := \MTMrx(\Hck_{G,\sgl}, \Lambda) \subset \MTMrx(\Hck_{G,\sgl}) = \Satrx^{G, \sgl}$$ of mixed Tate motives with rational, resp.~modular coefficients as in Subsection \ref{subsection:coefficients}.

\lemm
The full subcategory $\Satrx^{G, \Lambda} \subset \Satrx^{G,\{*\}}$ is stable under the fusion product (\thref{fusion.product}).
In addition, the underived tensor product functor 
$$\pH^0 (- \t \Lambda) : \Satrx^{G,\{*\}} \r \Satrx^{G, \Lambda}$$
is symmetric monoidal.
\xlemm

\pf
This is clear for $\Lambda = \Q$, since $\calF \t \Q = \colim (\dots \calF \stackrel n \r \calF \dots)$ and all our functors are additive and preserve filtered colimits.

We now consider the case $\Lambda = \Fp$. 
We write $\calF / p := \coker (\calF \stackrel p \r \calF)$ for an object $\calF$ in an abelian category.
Then $\Satrx^{G, \Fp} = \{\calF \in \Satrx^{G,\{*\}}\mid p \id_{\calF} = 0\} = \{\calF \mid \calF = \calF / p\}.$
Both functors in \refeq{Sat.fusion} are right exact (combine \thref{fusion.Satake} and the t-exactness of $i^!_{\ol \phi}[d_\phi]$, which follows from \thref{Sat.Tate}) so that in particular $(\calF_1 \pe \star \calF_2) / p = \calF_1 \pe \star (\calF_2 / p)$.
This shows our claims for $\Fp$-coefficients.
\xpf

The adjunction $F \dashv R$ constructed in the proof of \thref{Satake.Hopf} restricts to an adjunction
$$\Satrx^{G, \Lambda} \stackrel [R]{F} \rightleftarrows \MTMrx(X, \Lambda).$$
and $F R(\Lambda[1]) = FR(\Z[1]) \t \Lambda = \pH^0(FR(\Z[1]) \t \Lambda)$, according to the ind-freeness of $H = FR(\Z[1])$. The image of \(H_{\redx}^{G,\sgl}\) under \(\pH^0(-\otimes \Lambda)\) is again a commutative Hopf algebra object. If we denote the corresponding group by 
\(\widetilde{G}_{\redx,\Lambda} \in \MTMrx(X,\Lambda)^{\opp}\), we get a version of \refeq{Sat.comod} for \(\Lambda\)-coefficients.

\coro
\thlabel{Sat.Lambda}
There is an equivalence of symmetric monoidal categories 
$$\Satrx^{G, \Lambda} = \Rep_{\widetilde{G}_{\redx,\Lambda}} (\MTMrx(X, \Lambda)).$$
\xcoro

\subsection{The dual group}\label{subsec:identification}

In this subsection we identify the group associated to the Hopf algebra object $\H^{G,I}_{\redx} \in \MTMrx(X^I)$.
We first show in \thref{H.Independent} that the unreduced Hopf algebra $H^{G,I}$ is in fact reduced.
We then move on to compute the reduced Hopf algebra $H^{G,I}_\red$.

\subsubsection{Independence of the base and compatibility with realization} 

\theo \thlabel{H.Independent}
There is a natural isomorphism of Hopf algebras
$$H^{G, I} = i H_\red^{G, I},\eqlabel{H.red.unred}$$
where 
$$i : \MTMr(X^I) =\grAb \r \MTM(X^I), \quad \Z(k) \mapsto \Z(k)\eqlabel{functor.i}$$
is the natural symmetric monoidal functor.
\xtheo

\pf
On the level of the underlying objects in $\MTM$, this holds by \thref{object.flat} and the explicit computation in \thref{prop--L.free}.
We conclude that the (co)multiplication maps of $H^{G,I}$ and $iH_\red^{G,I}$ agree since $H_\red^{G,I}$ is ind-free, and $i$ is fully faithful when restricted to such objects (cf.~\thref{functor i}). 
\xpf

\coro
\thlabel{Betti comparison}
For $S = \Spec \Q$, the truncation of the Betti realization, $\pe \rho_\Betti := \pH^0 \rho_\Betti$ is a monoidal functor (with respect to the truncated convolution products $\pe \star$).
It fits into the following commutative diagram:
$$\xymatrix{
\MTMr(L^+ G \setminus \Gr_G) \ar@{.>}[r]^i \ar[d]^\cong &  \MTM(L^+ G \setminus \Gr_G)  \ar[d]^\cong \ar[r]^{\pe \rho_\Betti} & \Perv_{L^+G}(\Gr_G) \ar[d]^\cong \\
\Rep_{\widetilde G_\red}(\MTMr(S)) \ar[r]^i & \Rep_{\widetilde G}(\MTM(S)) \ar[r]^{\pe \rho_\Betti} & \Rep_{\rho_\Betti(\widetilde G)}(\Ab).}\eqlabel{compatibility Betti}$$
In particular, the composite $\pe \rho_\Betti \circ i$ is a realization functor for reduced motives in this situation, answering affirmatively a question posed in \cite[\S 1.6.1]{EberhardtScholbach:Integral}.
\xcoro

\pf
For an admissibly Whitney--Tate stratified (ind-)scheme $X$ over $S = \Spec \Q$, the Betti realization functor $\rho_\Betti : \DTM(X) \r \D(X^\an)$ and the reduction functor $\rho_\red : \DTM(X) \r \DTMr(X)$ are (at least) right t-exact. Their truncations $\pe \rho = \tau^{\ge 0} \circ \rho : \MTM(X) \r \Perv(X^\an)$, resp.~$\MTMr(X)$, are monoidal (with respect to the underived tensor product).

The derived convolution product $\star$ commutes with $\rho (:= \rho_\Betti, \rho_\red)$, by compatibility of $\rho$ with the six functors.
It is right t-exact by \thref{coro--t-exact}, so that $\pe \star$ commutes with $\pe \rho$.
Moreover, by the exactness of $\CT_B$, $\pe \rho$ also commutes with the fiber functors (on the level of reduced motives, motives, resp.~perverse sheaves).
Finally, $\pe \rho (\widetilde G) = \rho(\widetilde G)$, since the Hopf algebra $H^G$ is a filtered colimit of finite free graded abelian groups (\thref{prop--L.free}). This shows that the right half of \refeq{compatibility Betti} commutes, as does a similar diagram involving the truncated reduction functors $\pe \rho_\red$.
We have seen above $\widetilde G = i(\widetilde G_\red)$, giving rise to the bottom left functor $i$. 
The top-left horizontal functor $i$ is the one making the left half commutative.
\xpf

\subsubsection{Identification of the Hopf algebra for reduced motives}
We now identify the group $\widetilde{G}_{I, \redx} \in \MTMrx(X^I)^{\opp}$.
By \thref{Hopf:Alg:Exterior}, it suffices to do this for \(I=\{*\}\). 
By \thref{H.Independent}, it is enough to consider reduced motives, in which case we have \(\MTMr(X)\cong \grAb\).
In particular, we can describe \(H^G_{\red}:=H_{\red}^{G,\{*\}}\) by the affine \(\Z\)-group scheme $\widetilde{G}$ which underlies $\widetilde{G}_{\{*\}, \red}$, together with a grading of its global sections. 
This is similar to \cite[§VI.11]{FarguesScholze:Geometrization}, where we have a \(\Gm\)-action instead of a \(W_E\)-action. We will therefore follow the methods of loc.~ cit.~ in the proof of the theorem below.

Consider the Langlands dual group \(\hat{G}\), which is the pinned Chevalley group scheme with root datum dual to the root datum of \(G\). 
In particular, it comes with a fixed choice \(\hat{T}\subseteq \hat{B} \subseteq \hat{G}\) of maximal Torus and Borel.
As in \cite{FarguesScholze:Geometrization}, we need to modify it to get a canonical identification of \(\widetilde{G}\).
Namely, for any simple root \(a\) of \(\hat{G}\), instead of requiring an isomorphism \(\mathrm{Lie} (\hat{U}_a)\cong \Z\) in the data of our pinning, we choose an isomorphism 
\(\mathrm{Lie} (\hat{U}_a) \cong \Z(-1)\).
This is equivalent to choosing an isomorphism of the Hopf algebra of \(\hat{U}_a\) with the tensor algebra on \(\Z(-1)\).
(Note that $\Z(-1)$ does not have a preferred basis. 
Moreover, the appearance of a negative Tate twist differs from \cite[§VI.11]{FarguesScholze:Geometrization}, we refer to \thref{FS.difference} for an explanation why this happens.)  
In particular, this induces a \(\Gm\)-action on all the root groups of \(\hat{G}\). Letting \(\Gm\) act trivially on the root datum \((X_*(T),\Phi^\vee,X^*(T),\Phi)\), we get a \(\Gm\)-action on (the modified) \(\hat{G}\), denoted \(\psi_1\).

Now, let \(\hat{T}_\adj\subseteq \hat{G}_\adj\) be the adjoint torus of \(\hat{G}\),
and let \(2\rho_\adj\colon \Gm\to \hat{T}_\adj\) be the composition of \(2\rho\) with the projection \(\hat{T}\twoheadrightarrow \hat{T}_\adj\). Then we get a \(\Gm\)-action on \(\hat{G}\) by \[\psi_2\colon \Gm\xrightarrow{\rho_\adj} \hat{T}_\adj\to \Aut(\hat{G}), \eqlabel{grading on dual group}\]
where \(\rho_\adj\colon \Gm\to \hat{T}_\adj\) is the unique square root of \(2\rho_\adj\), and \(\hat{T}_\adj\) acts on \(\hat{G}\) by conjugation.

\rema
\thlabel{Tate twist negative degree}
To ensure that the above \(\Gm\)-action gives the correct grading in the following theorem, we consider the Tate twist \(\Z(1)\) to be in degree \(-1\) under the equivalence \(\MTMr(S)\cong \grAb\). We keep this convention for the rest of the paper. 
This agrees with the usual convention, and also with \cite{Zhu:Integral} and \cite{EberhardtScholbach:Integral}.
\xrema

\theo
\thlabel{Identification of dual group:Thm}
The \(\Gm\)-actions \(\psi_1,\psi_2\) on \(\hat{G}\) agree.
Equipped with this action, there is a canonical \(\Gm\)-equivariant isomorphism $\widetilde{G} \cong \hat{G}$. 
\xtheo

The existence of a (non-equivariant) isomorphism of $\Z$-group schemes $\widetilde{G} \cong \hat{G}$ can be deduced from \thref{H.Independent} ff. and \cite{MirkovicVilonen:Geometric}. Following \cite[§VI.11]{FarguesScholze:Geometrization}, we will prove \thref{Identification of dual group:Thm} from scratch in a way which also gives the \(\Gm\)-action, or equivalently the \(\Z\)-grading, and makes the isomorphism canonical.
We start by fixing a pinning of \(G\), which extends the choice of \((T,B)\). We will use this pinning to construct the isomorphism \(\widetilde{G}\cong \hat{G}\), and afterwards show this isomorphism does not depend on the choice of pinning.

By \thref{Sat.T}, we have \(\Satr^{T, \{*\}} = \Fun(X_*(T),\MTMr(X))\). This implies that \(H^T_{\red}\) has degree 0, and that \(\widetilde{T}\) is the torus with character group \(X_*(T)\), i.e., \(\widetilde{T}\cong \hat{T}\) canonically. 
The \(\Gm\)-actions \(\psi_1,\psi_2\) are also trivial. 

As the constant term functor \(\CT_B^{\{*\}}\colon \Satr^{G,\{*\}}\to \Satr^{T,\{*\}}\) is symmetric monoidal and commutes with the fiber functors, we get an induced morphism \(H^G_{\red}\to H^T_{\red}\) of Hopf algebras in \(\MTMr(X)\), and hence a homomorphism \(\widetilde{T}\to \widetilde{G}\). To show this is a closed immersion, it is enough to show that each \(\IC_{\nu,\Z}\in \Sat_\red^{T,\{*\}}\), for \(\nu\in X_*(T)\), is a quotient of an object lying in the image of \(\CT_B^{\{*\}}\), e.g.~ by \cite[Theorem 4.1.2 (ii)]{DuongHai:Tannakian}.
This holds since $F_\nu^\sgl (\calJ_!^\mu(\Z))$ for \(\mu\in X_*(T)^+\) is free and nonzero for each $\nu$ in the Weyl-orbit of $\mu$ by \thref{Standard.Free}.

\prop \thlabel{Generic.Split}
The generic fiber \(\widetilde{G}_\Q\) is a split connected reductive group and \(\widetilde{T}_{\Q} \subset \widetilde{G}_\Q\) is a maximal torus. 
\xprop

\pf
The proof follows \cite[\S 7]{MirkovicVilonen:Geometric}.
Restricting the equivalence of \thref{Sat.Lambda} to compact objects with $\Lambda = \Q$, we see that $\Rep_{\widetilde{G}_{\Q}}(\Vect_\Q^\fd)$ is semisimple, cf.~ the proof of \thref{Cok.N}, and the irreducible objects are parametrized by $X_*(T)^+$.
As \(X_*(T)^+\) is a finitely generated monoid and \(\IC_{\mu_1+\mu_2,\Q}\) is a subquotient of \(\IC_{\mu_1,\Q}\star \IC_{\mu_2,\Q}\) for all \(\mu_1,\mu_2\in X_*(T)^+\), we see that \(\widetilde{G}_\Q\) is algebraic by \cite[Proposition 2.20]{DeligneMilne:Tannakian}. As for any \(0\neq \mu\in X_*(T)^+\), the intersection motive \(\IC_{2\mu,\Q}\) is a subquotient of \(\IC_{\mu,\Q}
\star \IC_{\mu,\Q}\), \cite[Corollary 2.22]{DeligneMilne:Tannakian} tells us that \(\widetilde{G}_\Q\) is connected. Using this, we deduce by \cite[Proposition 2.23]{DeligneMilne:Tannakian} that \(\widetilde{G}_\Q\) is reductive. 
The rank of a reductive group over \(\Q\) agrees with the rank of its representation ring by \cite[Théorème 7.2]{Tits:Representations},
so \(\widetilde{T}_\Q\) is a maximal torus of \(\widetilde{G}_\Q\), and hence  \(\widetilde{G}_\Q\) is split. 
\xpf

We have the following generalization of \cite[Lemma VI.11.3]{FarguesScholze:Geometrization}.
\lemm
\thlabel{Iso:Criterion}
Let \(f\colon M\to N\) be a morphism of flat abelian groups. If \(M/p\to N/p\) is injective for all primes \(p\) and \(M\otimes_{\Z} \Q\to N\otimes_{\Z} \Q\) is an isomorphism, then \(f\) is an isomorphism.
\xlemm
\pf
Since \(M\) is flat, and hence torsion-free, we conclude that \(f\) is injective.
For surjectivity, consider some \(x\in N\). There exists some \(n\geq 1\) such that \(nx=f(y)\) lies in the image of \(f\); let \(n\geq 1\) be minimal with this property. Then \(y\) is a nontrivial element in the kernel of \(M/p\to N/p\) for any prime \(p\) dividing \(n\). Our assumption on the maps \(M/p\to N/p\) then tells us that \(n=1\), so that \(f\) is surjective.
\xpf

Consider the quotient \(H_{\red}^G\twoheadrightarrow K\) stabilizing the filtration \(\bigoplus_{n\geq i} \pH^n\pi_{G!}u^!\) of the fiber functor, i.e., the maximal Hopf algebra quotient \(K\) of \(H_{\red}^G\) such that the maps 
\[F^{\sgl}(\Ff)\to H_{\red}^G \otimes F^{\sgl}(\Ff) \to K\otimes F^{\sgl}(\Ff)\]
send the subobject \(\bigoplus_{n\geq i} \pH^n\pi_{G!}u^!(\Ff)\) to \(K\otimes \bigoplus_{n\geq i} \pH^n \pi_{G!}u^!(\Ff)\subseteq K\otimes F^{\sgl}(\Ff)\) for each \(i\in \Z\) and \(\Ff\in \Sat_\red^{G,\sgl}\).
This corresponds to a subgroup \(\widetilde{B}\subseteq \widetilde{G}\) containing \(\widetilde{T}\), and throughout the rest of the proof we will show that it is a Borel subgroup of \(\widetilde{G}\).
A preliminary result towards this is the following:

\lemm\thlabel{flatness of Borel}
The subgroup \(\widetilde{B}\subseteq \widetilde{G}\) is flat over \(\Spec \Z\).
\xlemm
\pf
Recall the algebras \(A_W\in \MTM_\red(X)\) from the proof of \thref{Satake.Hopf}, which gave rise to \(H_{\red}^G\) by considering their duals and taking the colimit.
These algebras are the image under \(F^{\sgl}\) of \(L_{\Sat,W}^{\sgl}\), and hence come equipped with a direct sum decomposition \(\bigoplus_{n \in \Z} A_{W,n}\), according to \(F^{\sgl}= \bigoplus_{n \in \Z} \pH^n \pi_{G!} u^!\).
The subalgebra of \(A_W\) that preserves the filtration \(\bigoplus_{n\geq i} \pH^n \pi_{G!} u^!\) on the fiber functor, is then exactly  
\(\bigoplus_{n\geq 0} A_{W,n}\).
This is finite free, as a direct summand of the finite free \(A_{W}\).
Its dual, which is a quotient of \(H_W\), is hence finite free as well.
Since the Hopf algebra of \(\widetilde{B}\) is a colimit of these quotients, we conclude that it is flat.
\xpf

Flatness of \(\widetilde{B}\) allows us to check smoothness fiberwise, although we first need to ensure \(\widetilde{B}\) is of finite type over \(\Spec \Z\). 
This will all be shown throughout the proof of \thref{Identification of dual group:Thm}.

\prop\thlabel{case pgl2}
For $G  = \PGL_2$, the actions \(\psi_1,\psi_2\colon \Gm\to \Aut(\hat{G})\) agree, and the standard pinning of $G$ induces a graded isomorphism $\widetilde{G} \cong \hat{G}$. 
\xprop

\pf
The Langlands dual group of \(\PGL_2\) is \(\SL_2\). 
Consider the minuscule dominant cocharacter \(\mu\), for which \(\Gr_{\PGL_2,\{*\}}^\mu\cong \P_S\times_S X\). Moreover, we have \(\Gr_{\PGL_2,\{*\}}^\mu\cap \Ss_{\mu,\{*\}}^+\cong \A^1_S\times_S X\), and \(\Gr_{\PGL_2,\{*\}}^\mu\cap \Ss_{-\mu,\{*\}}^+\cong X\), while the intersection of \(\Gr_{\PGL_2,\{*\}}^\mu\) with the other semi-infinite orbits is empty. 
In particular, \(F^{\{*\}}(\IC_{\mu,\Z}) \cong \Z(-1)\oplus \Z\) is an \(H_{\red}^{\PGL_2}\)-comodule, where we omit the shifts by \([1]\) for simplicity.
This induces a homomorphism \(\widetilde{G}\to \GL(\Z(-1)\oplus \Z)\). We claim this is a closed immersion, inducing an isomorphism \(\widetilde{G} \cong \SL(\Z(-1)\oplus \Z)\).

Indeed, \(\hat{T}\) acts on \(\Z(-1)\oplus \Z\) by weights \(\pm 1\),
as \(\CT_B^{\{*\}}(\IC_{\mu,\Z})\) is concentrated on the connected components corresponding to \(\pm 1\) under the isomorphism \(\pi_0(\Gr_{T,\{*\}})\cong \Z\). 
Thus, the image of \(\hat{T}\) lands in \(\SL(\Z(-1)\oplus \Z)\). 
In particular, the claim over \(\Spec \Q\) follows
as we already know \(\widetilde{G}_\Q\) is reductive with maximal torus \(\hat{T}_\Q \cong \mathbf{G}_{m,\Q}\), and \(\widetilde{G}_\Q\) is not a torus by considering its representation ring.
By flatness, the $\Z$-morphism \(\widetilde{G}\to \GL(\Z(-1)\oplus \Z)\) factors through \(\SL(\Z(-1)\oplus \Z)\),
and we get a map \(\widetilde{G}_{\Fp}\to \SL(\Z(-1)\oplus \Z)_{\Fp}\) for any prime \(p\). Let \(K_p\) denote the image of this map, so that we have a surjection \(\widetilde{G}_{\Fp}\to K_p\). 
The irreducible (ungraded) representations of \(\widetilde{G}_{\Fp}\) are parametrized by \(X_*(T)^+\).
In particular, the irreducible representations of \(K_p\) can be indexed by a subset of \(X_*(T)^+\), so that \(K_p=\SL(\Z(-1)\oplus \Z)_{\Fp}\) by \cite[Lemma VI.11.2]{FarguesScholze:Geometrization}. 
Then \thref{Iso:Criterion}, used on the level of (ungraded) Hopf algebras, tells us that \(\widetilde{G}\to \SL(\Z(-1)\oplus \Z)\) is an isomorphism.
It is moreover clear that \(\widetilde{B}\subset \widetilde{G}\) is the positive Borel in this case.

As \(\widetilde{B}\) stabilizes \(\Z(-1)\subseteq \Z(-1)\oplus \Z\), the Lie algebra of its unipotent radical \(\widetilde{U}\) can be identified with \(\IHom(\Z,\Z(-1))\cong \Z(-1)\). This gives the \(\Gm\)-equivariant isomorphism \(\widetilde{G}\cong \hat{G}\) when \(G=\PGL_2\), where \(\PGL_2\) is equipped with the action determined by \(\psi_1\).

To identify \(\psi_1\) and \(\psi_2\), we again consider $\widetilde U \subset \widetilde B \subset \widetilde G \cong \SL_2$. 
It is clear that \(\hat{T}\cong \widetilde{T}\to \widetilde{G}\) is the maximal torus appearing in the pinning of \(\SL_2\). 
Restricting the action of \(\widetilde{U}\) on \(\Z(-1)\oplus \Z\) to \(\Z\) gives a (graded) morphism
\[\Z(-1)\oplus \Z\to \Z[t]\otimes \Z  \colon (n,m)\mapsto t\otimes n,\] where \(\Z[t]\) is the Hopf algebra of \(\widetilde{U}\), equipped with a grading as a quotient of \(H_{\red}^G\). 
In particular, \(t\) lies in degree \(1\), and similarly we see that the coordinate of the unipotent radical of the opposite Borel of \(\widetilde{G}\) lies in degree $-1$. 
This shows that the grading on \(\widetilde{G}\) corresponds to the \(\Gm\)-action $\Gm\to \Aut(\SL_2)$, where an invertible element $x$ acts by
\[\begin{pmatrix}
  a&b\\
  c&d
\end{pmatrix} \mapsto 
\begin{pmatrix}
  a&xb\\
  x^{-1}c&d
\end{pmatrix} = \begin{pmatrix}
x&0\\0&1
\end{pmatrix} \begin{pmatrix}
  a&b\\
  c&d
\end{pmatrix} \begin{pmatrix}
x^{-1}&0\\0&1
\end{pmatrix}.\eqlabel{grading on dual group PGL2}\]
This is the grading \refeq{grading on dual group}, so this finishes the case \(G=\PGL_2\).
\xpf

\rema\thlabel{FS.difference}
We can now see why the Lie algebras of the simple root groups of \(\hat{G}\) are identified with a negative Tate twist, instead of with a positive one as in \cite[§VI.11]{FarguesScholze:Geometrization}. 
The positive Borel in \(\hat{G}\) can only raise the weights of a \(\hat{G}\)-representation.
So since \(\pH^n\pi_{G!}u^! \cong \bigoplus_{\langle 2\rho,\nu\rangle=n} F_\nu^{\sgl}\) by \thref{fiber.functor.remarks}, the positive Borel stabilizes the filtration \(\bigoplus_{n\geq i} \pH^n\pi_{G!}u^!\) of the fiber functor, instead of the one defined via \(\leq\) as in \cite[§VI.11]{FarguesScholze:Geometrization}.
In the particular case \(G=\PGL_2\) and the standard representation of \(\widetilde{G}\cong \SL_2\), the positive Borel stabilizes the highest weight space, which is identified with the negative Tate twist (contrary to \cite[Lemma VI.11.3 ff.]{FarguesScholze:Geometrization}), as it arises from the compactly supported cohomology of the (shifted) constant object on the affine line.
\xrema

Returning to arbitrary $G$, the adjoint quotient \(G\twoheadrightarrow G_\adj\) induces a map \(\Gr_{G,\{*\}}\to \Gr_{G_\adj,\{*\}}\) which restricts to a universal homeomorphism on each connected component. As \(\DM\) does not satisfy étale descent we cannot conclude that this map induces equivalences of motives on these connected components. However, we can show that we get equivalences for Tate motives, even for unreduced motives.

\lemm \thlabel{DTM.Homeo}
Let \(X\) be a connected component of \(\Gr_{G,\{*\}}\) and \(Y\) the connected component of \(\Gr_{G_\adj,\{*\}}\) to which it maps under the natural morphism. If we denote the induced morphism by \(\alpha\colon X\to Y\), then \(\alpha^*\colon \DTMrx(Y)\to \DTMrx(X)\) is an equivalence.
\xlemm
\pf
It suffices to prove the lemma for the stratification by Iwahori-orbits.
In this case, we will show that the unit $\id \r \alpha_* \alpha^*$ and counit $\alpha^* \alpha_* \r \id$ are equivalences.
The map \(\alpha\) is an \(LG\)-equivariant universal homeomorphism by \thref{reduction to simply connected case}, which induces an isomorphism on Iwahori-orbits.
Indeed, the map \(G\to G_{\adj}\) induces isomorphisms on affine root groups, so this follows from \cite[(4.3.10)]{RicharzScholbach:Intersection}.
Hence, the lemma is immediate if we have a single cell.
Because $\alpha$ is ind-proper, the unit and counit maps commute with $*$-pullback to any union of cells.
We thus conclude by localization and induction on the number of cells.
\xpf

\prop \thlabel{prop--ssrankOne}
If $G$ has semisimple rank $1$, then the actions \(\psi_1\) and \(\psi_2\) agree, and a pinning of $G$ induces a graded isomorphism $\widetilde G \cong \hat G$.
\xprop

\pf
The adjoint quotient $G_{\adj}$ can be identified with $\PGL_2$ via the pinning of $G$.
Since \(\pi_0(\Gr_{G,\{*\}})\cong \pi_1(G)\) canonically, it follows from \thref{DTM.Homeo} that \(\DTMr(\Gr_{G,\{*\}}) \cong \DTMr(\pi_1(G)\times_{\pi_1(\PGL_2)} \Gr_{\PGL_2,\{*\}})\). As \(\pi_1(\PGL_2)\cong \Z/2\), every object in \(\DTMr(\Gr_{\PGL_2, \{*\}})\), and hence every object in the Satake category, is equipped with a \(\Z/2\)-grading. 
Then \(\Satr^{G,\{*\}}\) is equivalent to the category of objects of \(\Satr^{\PGL_2,\{*\}}\), equipped with a \(\pi_1(G)\)-grading that refines the \(\Z/2\)-grading. In particular, we get \(\widetilde{G}\cong \widetilde{\PGL_2}\overset{\mu_2}{\times} \widetilde{Z}\), where \(\widetilde{Z}\) is the multiplicative group scheme with character group \(\pi_1(G)\); note that $\widetilde{Z}$ is a torus exactly when \(\pi_1(G)\) is torsion-free. By \cite[XII, Proposition 4.11]{SGA3:2}, the center \(\hat{Z}\) of \(\hat{G}\) is multiplicative, with character group canonically isomorphic to \(\pi_1(G)\) by \cite[Proposition 1.10]{Borovoi:Abelian}. Thus, the isomorphism \(\hat{G}\cong (\hat{G})_\sico\overset{\mu_2}{\times} \hat{Z}\)
implies that \(\hat{G}\cong \widetilde{G}\) canonically, inducing \(\hat{T}\cong \widetilde{T}\) and \(\hat{B} \cong \widetilde{B}\); in particular, \(\widetilde{B}\) is again a Borel of \(\widetilde{G}\). 
The fact that \(\psi_1=\psi_2\), and that the isomorphism above is \(\Gm\)-equivariant, follows from the case \(G=\PGL_2\) covered in \thref{case pgl2}.
\xpf

Finally, we consider a general group \(G\), still equipped with a fixed pinning. 
To any simple coroot \(a\) of \(G\), we can associate a minimal parabolic with Levi quotient strictly containing the maximal torus: \(T\subsetneq M_a\subseteq P_a\subseteq G\). We also have the symmetric monoidal constant term functor \(\CT_{P_a}^\sgl\colon \Satr^{G,\{*\}}\to \Satr^{M_a,\{*\}}\),
which commutes with the fiber functors. So it induces a morphism \(H_{\red}^G\to H_{\red}^{M_a}\) of Hopf algebras, and hence a homomorphism \(\hat{M}_a \cong \widetilde{M}_a\to \widetilde{G}\). 
By \thref{CT.composition}, the morphism \(\widetilde{M}_a \to \widetilde{G}\) commutes with the closed immersions \(\hat{T}\to \hat{M}_a \cong \widetilde{M}_a\) and \(\hat{T}\to \widetilde{G}\). 
To show this is a closed immersion on the generic fiber, we apply \cite[Proposition 2.21]{DeligneMilne:Tannakian}: consider objects of the form \(\IC_{\lambda,\Q}\in \Satr^{M_a,\{*\}}\), with \(\lambda\) dominant for \(M_a\). 
As \(\IC_{\lambda,\Q}\) is a quotient of a twist of \(\CT_{P_a}^\sgl(\IC_{\mu,\Q})\), where \(\mu\in X_*(T)\) is the unique dominant (for \(G\)) representative of \(\lambda\), we see that $\hat{M}_{a, \Q} \r \widetilde{G}_{\Q}$ is a closed immersion. 
(Recall that any cocharacter in \(X_*(T)\) has a unique dominant element in its orbit under the Weyl group action.) 
\prop \thlabel{Dual.Fiber}
The \(\Gm\)-actions \(\psi_1\) and \(\psi_2\) agree in general. 
Moreover, the closed immersions \(\hat{T}_\Q\to \widetilde{G}_\Q\) and \(\hat{M}_{a,\Q}\to \widetilde{G}_\Q\), which involve a choice of pinning of $G$, extend uniquely to a graded isomorphism \(\hat{G}_\Q \cong \widetilde{G}_\Q\).
\xprop

\pf
By \thref{Generic.Split}, \(\widetilde{G}_\Q\) is a split reductive group with maximal torus $\widetilde{T}$. As \(\widetilde{M}_{a, \Q} \to \widetilde{G}_{\Q}\) is a closed immersion, it induces an embedding on Lie algebras, so that \(a\in X_*(T) = X^*(\hat{T})\) determines a root of \(\widetilde{G}_\Q\), while the root \(a^\vee\) associated to \(a\) determines a coroot of \(\widetilde{G}_\Q\). Note that, as \(\widetilde{T}\cong \hat{T}\), passing to dual groups preserves the pairing between roots and coroots, up to reversing their roles.
Hence, the simple reflections are also contained in the Weyl group of \(\widetilde{G}_\Q\), giving an inclusion \(W=W(G,T)\subseteq \widetilde{W}:= W(\widetilde{G}_\Q,\widetilde{T}_\Q)\), as subgroups of \(\Aut(X_*(T))\).
Let us denote as before the roots of \(G\) by \(\Phi:=\Phi(G,T)\), and the coroots by \(\Phi^\vee:=\Phi^\vee(G,T)\).
Then, as all (co)roots are a \(W\)-translate of a simple (co)root, this implies 
that \(\Phi^\vee\subseteq \Phi(\widetilde{G}_\Q,\widetilde{T}_\Q)\) and \(\Phi\subseteq \Phi^\vee(\widetilde{G}_\Q,\widetilde{T}_\Q)\), as subsets of \(X_*(T) \cong X^*(\hat{T})\) and \(X^*(T) \cong X_*(\hat{T})\) respectively.

To show that these inclusions are equalities, note that for \(\lambda\in X^*(\widetilde{T})^+\), the weights of the simple \(\widetilde{G}_\Q\)-representation of highest weight \(\lambda\) are those weights \(\lambda'\) in the convex hull of the \(\widetilde{W}\)-orbit of \(\lambda\) such that \(\lambda-\lambda'\) is in the root lattice of \(\widetilde{G}_\Q\). 
Let \(a\in \Phi(\widetilde{G}_\Q,\widetilde{T}_\Q)\) be a simple root, and choose \(\lambda\) regular, so that the corresponding simple \(\widetilde{G}_\Q\)-representation has \(\lambda-a\) as a weight. 
Then the restriction of \(\CT_B^\sgl(\IC_{\lambda,\Q})\) to \(\lambda-a\in \pi_0(\Gr_{T,\sgl})\) does not vanish by \thref{examples of intersections}, \thref{Standard.Free} and \thref{Cok.N}.
This implies that \(\lambda-(\lambda-a) = a\) lies in the coroot lattice of \(G\), so that \(a\) must be a simple coroot (as \(a\in \Phi(\widetilde{G}_\Q,\widetilde{T}_\Q)\) is simple and we have \(\Phi^\vee\subseteq \Phi(\widetilde{G}_\Q,\widetilde{T}_\Q)\)).
This shows that the simple roots of \(\widetilde{G}_\Q\) and the simple coroots of \(G\) are in bijection, hence that the inclusion of Weyl groups above is an equality.
Since any (co)root is a Weyl group translate of a simple (co)root, we then conclude that \(\Phi^\vee = \Phi(\widetilde{G}_\Q,\widetilde{T}_\Q)\) and \(\Phi = \Phi^\vee(\widetilde{G}_\Q,\widetilde{T}_\Q)\).
Thus, the root data of \(\widetilde{G}_\Q\) and \(\hat{G}_\Q\) agree, so we get the desired isomorphism \(\hat{G}_\Q \cong \widetilde{G}_\Q\). Finally, the fact that \(\psi_1=\psi_2\) and the compatibility with the gradings follows from the case where \(G\) is a torus or has semisimple rank 1, as the closed immersions induced by constant term functors are compatible with the gradings.
\xpf

\prop \thlabel{Dual.Z}
The isomorphism in \thref{Dual.Fiber} extends uniquely to an integral isomorphism $\hat{G} \cong \widetilde{G}$.
\xprop

\pf
Consider a prime \(p\), and the ring of integers \(\breve{\Z}_p\subseteq \breve{\Q}_p\) of the completion of the maximal unramified extension of \(\Qp\). Viewing \(\hat{G}(\breve{\Z}_p)\) and \(\widetilde{G}(\breve{\Z}_p)\) as subsets of \(\hat{G}(\breve{\Q}_p)\cong \widetilde{G}(\breve{\Q}_p)\), we know that \(\hat{G}(\breve{\Z}_p) \subseteq \widetilde{G}(\breve{\Z}_p)\), as the former is generated by the \(\hat{M_a}(\breve{\Z}_p) \cong \widetilde{M}_a(\breve{\Z}_p)\).
Let \(\widetilde{G}\to \GL_n\) be a representation corresponding to some object in \(\Satr^{G,\{*\}}\) which is a closed immersion on the generic fiber.
Then \(\hat{G}_{\Qp}\cong \widetilde{G}_{\Qp}\to \GL_{n,\Qp}\) extends to a morphism \(\hat{G}_{\Zp}\to \GL_{n,\Zp}\) by \cite[Proposition 1.7.6]{BruhatTits:Groups2}.
Using a similar argument as in \thref{prop--ssrankOne}, we may assume \(G\) is adjoint, in which case \(\hat{G}\) is simply connected.
Since odd special orthogonal groups are not simply connected, we can apply \cite[Corollary 1.3]{PrasadYu:QuasiReductive} to see that \(\hat{G}_{\Zp}\to \GL_{n,\Zp}\) is also a closed immersion.
By flatness of \(\widetilde{G}\), this map then factors as \(\widetilde{G}_{\Zp}\to \hat{G}_{\Zp}\), which is an isomorphism on the generic fiber. But it is also surjective on the special fiber, as any point in \(\hat{G}(\overline{\mathbf{F}}_p)\) can be lifted to a point in \(\hat{G}(\breve{\Z}_p)\) by smoothness of \(\hat{G}\) and completeness of \(\breve{\Z}_p\), and we already know that \(\hat{G}(\breve{\Z}_p)\subseteq \widetilde{G}(\breve{\Z}_p)\). 
Since any surjection of algebraic groups over a field with reduced target is faithfully flat \cite[Summary 1.71]{Milne:Algebraic}, hence induces an injection on global sections, \thref{Iso:Criterion} then tells us that \(\widetilde{G}_{\Zp}\to \hat{G}_{\Zp}\) is an isomorphism. 

As the previous paragraph is valid for all primes \(p\), we see that all fibers of \(\widetilde{G}\to \Spec \Z\) are reductive, so that \(\widetilde{G}\) is reductive by \cite[Theorem 1.5]{PrasadYu:QuasiReductive}. The closed immersion \(\hat{T}\to \widetilde{G}\) determines a maximal torus over all geometric fibers by rank considerations, so \(\hat{T}\) is a maximal torus of \(\widetilde{G}\). Then \(\widetilde{G}\) is split by \cite[Example 5.1.4]{Conrad:Groups} since $\text{Pic} (\Z)$ is trivial.
The corresponding root datum can be determined on the generic fiber, so that the previously obtained identification of \(\widetilde{G}\) on the generic fiber gives us \(\hat{G}\cong \widetilde{G}\), although a priori only non-canonically. But then \cite[XXIII, Théorème 4.1]{SGA3:3} gives us a unique isomorphism \(\hat{G}\cong \widetilde{G}\) extending the isomorphisms \(\hat{G}_{\Zp}\cong \widetilde{G}_{\Zp}\) for each prime \(p\). 
\xpf

This finishes the proof of \thref{Identification of dual group:Thm} when \(G\) is equipped with a pinning; it remains to show the isomorphism \(\hat{G}\cong \widetilde{G}\) is independent of the choices.
Note that we can now conclude that the subgroup \(\widetilde{B}\subseteq \widetilde{G}\) is the Borel corresponding to \(\hat{B}\).
Indeed, as \(\widetilde{G}\) is reductive, \(\widetilde{B}\) is of finite type over \(\Spec \Z\).
Since \(\widetilde{B}\) is defined as the stabilizer of a flag (in a generalized sense), it is fiberwise a parabolic, and in particular fiberwise smooth. 
Hence \(\widetilde{B}\) is smooth by \thref{flatness of Borel}.
(In fact, the fibers of \(\widetilde{B}\) can be shown to agree with the canonical Borels constructed in \cite[\S 7]{MirkovicVilonen:Geometric}, \cite[Corollary 5.3.20]{Zhu:Introduction}, \cite[\S 1.9.2]{BaumannRiche:Satake}.)
Since \(\widetilde{B}\) clearly contains the natural Borel \(\hat{B}\subseteq \hat{G}\cong \widetilde{G}\), it must be a standard parabolic.
Now, the intersection of \(\widetilde{B}\subseteq \widetilde{G}\) with any minimal Levi \(\widetilde{M}_a\subseteq \widetilde{G}\) is a Borel of \(\widetilde{M}_a\).
Indeed, this follows from the compatibility of the constant term functors as in \thref{CT.composition}, and the similar observation in the case where \(G\) was of semisimple rank 1.
We deduce that \(\widetilde{B}\) is itself a Borel of \(\widetilde{G}\).  

\prop
The isomorphism in \thref{Dual.Z} is independent of the pinning of $G$.
\xprop

\pf
First, assume $G$ has semisimple rank one. As \(\widetilde{T}\) is the stabilizer of the cohomological grading of the fiber functor, and \(\widetilde{B}\) is the stabilizer of the corresponding filtration, they are independent of $(T, B)$ by \thref{fiber.functor.properties}. We used the rest of the pinning to identify \(G_\adj\) with \(\PGL_2\).
In the adjoint case, we had a canonical graded isomorphism 
\(\mathrm{Lie} (\widetilde{U})\cong \IHom(\Z,\Z(-1)) \cong \Z(-1)\), cf.~ the discussion preceding \refeq{grading on dual group PGL2}. 
This does not depend on how we identified \(G_\adj\) with \(\PGL_2\): automorphisms of \(\PGL_2\) induce automorphisms of the minuscule Schubert cell \(\P_S\times_S X\), but any such automorphism acts trivially on  \(F^{\{*\}}(\IC_{\mu,\Z}) \cong \Z(-1)\oplus \Z\).

To show independence for general $G$, it suffices to show that for each simple coroot \(a\), the constant term functors \(\CT_{P_a}^{\{*\}}\colon \Satr^{G,\sgl}\to \Satr^{M_a,\sgl}\) are independent of the choice of Borel, and hence of the choice of parabolic with Levi quotient \(M_a\).
Indeed, this will give a canonical embedding \(\widetilde{M}_a\to \widetilde{G}\), and since \(\widetilde{G}\) is generated by these Levi's, the independence of the isomorphism follows from the rank 1 case.
Consider the flag variety \(\FlCal\), non-canonically isomorphic to \(G/B\), parametrizing the Borels of \(G\). The quotient of the universal Borel \(B_{\FlCal}\subseteq G_{\FlCal}:=G\times_S \FlCal\) by its unipotent radical is a torus, the universal Cartan \(T_{\FlCal}\). It is defined over \(\Z\), as \(\FlCal\) is simply connected, and split. 

Consider a simple coroot \(a\). Then there is the universal parabolic \(P_{a,\FlCal}\) with Levi quotient \(M_{a,\FlCal}\). Let \(\widetilde{\FlCal}_a\) be the \(\Gm\)-torsor over \(\FlCal\) parametrizing the pinnings of \(M_{a,\FlCal}\). The group \(M_{a,\widetilde{\FlCal}_a}\) admits a pinning by construction, and is hence constant by \cite[XXIII, Théorème 4.1]{SGA3:3}.

Now, we claim that we can repeat the whole story, replacing our base \(S\) by \(\widetilde{\FlCal}_a\). In particular, we have the symmetric monoidal constant term functor
\[\CT_{P_{a,\widetilde{\FlCal}_a}}^{\sgl} \colon \DTMr(\Hck_{G_{\FlCal},\{*\},\widetilde{\FlCal}_a})\to \DTMr(\Hck_{M_{a,\FlCal},\{*\},\widetilde{\FlCal}_a}),\]
where the reduced motives are defined using the base \(\widetilde{\FlCal}_a\). Since \(P_{a,\FlCal}\) does not come from base change from $\Spec \Z$, the only thing we must check is that \(\CT_{P_{a,\widetilde{\FlCal}_a}}^{\sgl}\) preserves stratified Tate motives. 
For this we observe that the proof of \thref{cellularity of intersection:torus} shows that the intersections of the Schubert cells and semi-infinite orbits over \(\widetilde{\FlCal}_a\) admit filtrable decompositions by vector bundles and punctured vector bundles: these vector bundles arise from the root groups of the universal parabolics.
Hence, preservation of \(\DTM_\red\) follows as in \thref{CT.DTM}, using that vector bundles have Tate motives by the projective bundle formula.
Now by independence of the base, \thref{independence.star}, *-restricting to a fiber of \(\widetilde{\FlCal}_a\) and shifting by \(\dim (\widetilde{\FlCal}_a)\) induces an equivalence of Satake categories (defined using different bases), compatible with the constant term functors. This shows that the constant term functor \(\CT_{P_a}^{\{*\}}\) is independent of the choice of Borel. \xpf

\rema
Consider the semidirect product \(\GT:=\hat{G}\rtimes \Gm\), with \(\Gm\) acting on \(\hat{G}\) via \refeq{grading on dual group}.
Modulo the identification of the root groups of \(\hat{G}\) with Tate twists, \(\GT\) agrees with Deligne's modification of the dual group from \cite{FrenkelGross:Rigid}, defined as \((\hat{G}\times \Gm)/(2\rho \times \id)(\mu_2)\), cf.~ \cite[Remark 6.7]{RicharzScholbach:Motivic}. It also agrees with the group \(\hat{G}^T\) from \cite[§1]{Zhu:Integral}.
\xrema

Since graded abelian groups are equivalent to abelian groups equipped with a \(\Gm\)-action, we can also describe the Satake category as a representation category of (ungraded) abelian groups as follows.

\coro \thlabel{CG.Coro}
There is a canonical equivalence \(\Satr^{G,\{*\}}\cong \Rep_{\GT}(\Ab)\).
\xcoro

This gives further evidence for Bernstein's suggestion that the C-group from \cite{BuzzardGee:Conjectures}, which is the L-group of \(\GT\), might be more appropriate in the Langlands program than the L-group \cite[Remark 9 (2)]{Zhu:Integral}.

\subsection{The Vinberg monoid}
In this section we consider a subcategory of anti-effective stratified Tate motives for the purpose of geometrizing Hecke algebras over $\Z[\mathbf{q}]$, where $\mathbf{q}$ is an indeterminate. 
Recall that for $\ell$-adic sheaves over $S = \Fq$, the trace of the geometric Frobenius element on $\Q_{\ell}(-1)$ is $q$. Thus, it is natural to single out anti-effective motives, so that $\ell$-adic realization geometrizes the specialization map on Hecke algebras $\mathbf{q} \mapsto q$. 

By \thref{Fl.universally.WT} and \thref{Convolution.Anti}, we have the symmetric monoidal category
$$\Satrx^{G, \anti} := \Satrx^{G, \sgl, \anti} := \MTMrx(\Hck_{G, \sgl})^\anti.$$
For $G=1$, $\Satr^{G, \anti} = \MTMr(X)^\anti \subset \Satr^{G, \sgl} = \MTMr(X)$ identifies with the full subcategory consisting of those graded abelian groups concentrated in nonnegative degrees (under the monoidal isomorphism \(f^*[1]\colon \MTMr(S)\to \MTMr(X)\) for \(f\colon X\to S\), where we still consider the Tate twist \(\Z(1)\) to be negatively graded).
For general $G$, \(\Satr^{G,\anti}\) is generated by the \(\IC_{\mu,L}\) for \(\mu\in X_*(T)^+\) and \(L\in \MTMr(S)^{\anti}\subseteq \MTMr(S)\).

\prop \thlabel{Sat.Anti}
Let $M \in \Satrx^{G, \{*\}}$. Then $M \in \Satrx^{G, \anti}$ if and only if $F^{\sgl}(M) \in \MTMrx(X)^{\anti}$.
\xprop

\pf 
If $M \in \Satrx^{G, \anti}$, by excision and the filtrable decomposition in \thref{cellularity of intersection:torus} we have $F^{\sgl}(M) \in \MTMrx(X)^{\anti}$. For the converse, suppose $F^\sgl(M) \in \MTMrx(X)^{\anti}$. If $M$ is bounded, by \thref{Orbit.Eq} it admits a finite filtration with subquotients given by $\IC$-motives  $\IC_{\mu, L}$ for $L \in \MTMrx(S)^\comp$ and $\mu \in X_*(T)^+$. Each $F^\sgl(\IC_{\mu, L})$ is a subquotient of $F^\sgl(M)$, and consequently each $f^*L[1]$ is also a subquotient of $F^\sgl(M)$, as $ \Gr_{G, \sgl}^{\mu} \cap \Ss_{w_0(\mu), \sgl}^+= X$ where $w_0$ is the longest element in the Weyl group. 
Because $\MTMrx(X)^{\anti}$ is closed under subquotients (\thref{Anti.Orthogonal}), this implies each $L$ is anti-effective, and hence so is $\IC_{\mu, L}$ by \thref{t-structure.stratified} (see also \thref{defi--MTM(anti)}).
This implies that $M$ is also anti-effective. 
If $M$ is not necessarily bounded, we can present it as a filtered colimit $M = \colim M_i$ of bounded  subobjects $M_i \subset M$ in $\MTMrx(\Hck_{G,\sgl})$. Being a subobject of $F^\sgl(M)$, $F^\sgl(M_i)$ is anti-effective (\thref{Anti.Orthogonal}), hence so is $M_i$ and therefore also $M = \colim M_i$. 
\xpf

In Zhu's integral Satake isomorphism \cite[Proposition 5]{Zhu:Integral}, the Vinberg monoid appears instead of the usual dual group. 
We now explain how this monoid naturally appears from our motivic Satake equivalence. 
Afterwards, we construct a generic Satake isomorphism between the generic spherical Hecke algebra and the representation ring of the Vinberg monoid, 
which are naturally \(\Z[\qq]\)-algebras for some indeterminate \(\qq\). 
In an attempt to not further lengthen the paper, we will only recall the necessary definitions for the Vinberg monoid, 
and we refer to \cite{Vinberg:Reductive, XiaoZhu:Vector, Zhu:Integral} and the references there for more details.

Denote by \(X^*(\hat{T}_\adj)_{\text{pos}}\subseteq X^*(\hat{T})^+_{\text{pos}}\subseteq X^*(\hat{T})\) the submonoids of characters generated by the simple roots, respectively the dominant characters and the simple roots. 
Viewing \(\Z[\hat{G}]\) as a \(\hat{G}\times \hat{G}\)-module via left and right multiplication,  the global sections \(\Z[\hat{G}]\) admit an \(X^*(\hat{T})^+_{\text{pos}}\)-multi-filtration
\(\Z[\hat{G}] = \sum_{\mu\in X^*(\hat{T})^+_{\text{pos}}} \fil_\mu\Z[\hat{G}]\), where \(\fil_\mu\Z[\hat{G}]\) is the maximal \(\hat{G}\times \hat{G}\)-submodule 
such that all its weights \((\lambda,\lambda')\in X^*(\hat{T})\times X^*(\hat{T})\) satisfy \(\lambda\leq -w_0(\mu)\) and \(\lambda'\leq \mu\). 
Here \(w_0\) is the longest element of the Weyl group of \(G\).
We then define Vinberg's universal monoid \(V_{\hat{G}} = \Spec (\bigoplus_{\mu\in X^*(\hat{T})^+_{\text{pos}}} \fil_\mu \Z[\hat{G}])\), with 
the natural (co)multiplication map and monoid morphism \( d_{\rho_\adj} \colon V_{\hat{G}}\to \Spec \Z[X^*(\hat{T}_\adj)_{\text{pos}}]=:\hat{T}^+_\adj\).
The dominant cocharacter \(\rho_\adj\) extends to a monoid morphism \(\rho_\adj\colon \A^1\to \hat{T}_\adj^+\), and the monoid $V_{\hat{G},\rho_\adj}$ is defined as in the commutative diagram below, in which all squares are cartesian.
$$\xymatrix{
\hat{G} \times \hat{T} \ar@{->>}[r] & \hat{G}\overset{Z_{\hat{G}}}{\times} \hat{T} \ar[d] \ar[r] & V_{\hat{G}} \ar[d] & V_{\hat{G},\rho_\adj} \ar[d]^{d_{\rho_\adj}} \ar[l] & (\hat{G}\times \Gm)/(2\rho\times \id)(\mu_2) \ar@/_1pc/[lll]_{\id \times 2 \rho}  \ar[l] \ar[d] & \ar@{->>}[l] \hat{G}\times \Gm \\ & \hat{T}_{\adj} \ar[r] & \hat{T}^+_{\adj} & \A^1 \ar[l]_{\rho_\adj} & \Gm \ar[l] &
}$$
There is an isomorphism $(\hat{G}\times \Gm)/(2\rho\times \id)(\mu_2) \cong \hat{G} \rtimes \Gm = \GT$, $(g,t) \mapsto (g 2\rho(t)^{-1}, t^2)$.
In particular, \(\hat{G}\overset{Z_{\hat{G}}}{\times} \hat{T} \subseteq V_{\hat{G}}\) and \(\GT\subseteq V_{\hat{G},\rho_\adj}\) are the respective groups of units.

The following theorem was already explained by T. Richarz during a talk in the Harvard Number Theory Seminar in April 2021 for rational coefficients. In the same talk, he asked whether it was possible to do this integrally, and mentioned \thref{Rational:K0}. The theorem below gives an affirmative answer to this question.

\theo
\thlabel{Reps of Vinberg monoid}
The equivalence \(\Satr^{G,\{*\}}\cong \Rep_{\GT}(\Ab)\) restricts to an equivalence \(\Satr^{G,\{*\},\anti}\cong \Rep_{V_{\hat{G},\rho_\adj}}(\Ab)\).
\xtheo
\pf
The proof is a generalization of the arguments in \cite[Lemma 21]{Zhu:Integral} to integral coefficients.
First, as \(\GT\subseteq V_{\hat{G},\rho_\adj}\) is open and dense, the restriction functor \(\Rep_{V_{\hat{G},\rho_\adj}}(\Ab) \to \Rep_{\GT}(\Ab)\) is fully faithful. 
In particular, it suffices to identify the two full subcategories under the equivalence \(\Satr^{G,\{*\}}\cong \Rep_{\GT}(\Ab)\). 

The composite $\hat{G}\times \hat{T} \r V_{\hat{G}}$ in the above diagram 
corresponds to the natural inclusion of coalgebras
\[\bigoplus_{\mu\in X^*(\hat{T})^+_{\text{pos}}} \fil_\mu \Z[\hat{G}]\to \Z[\hat{G}]\otimes_\Z \Z[X^*(\hat{T})].\] 
This map is given by sending each \(\fil_\mu \Z[\hat{G}]\) into \(\Z[\hat{G}] \otimes \Z e^\mu\), where $e^\mu \in \Z[X^*(\hat{T})]$ generates the rank 1 subgroup corresponding to \(\mu\); this follows from \cite[(3.2.3) ff.]{XiaoZhu:Vector}.

We first determine when a $\hat{G}\overset{Z_{\hat{G}}}{\times} \hat{T}$-representation \(M\) extends (necessarily uniquely) to a \(V_{\hat{G}}\)-representation.
Consider the decomposition of \(M\) into its \(\hat{T}\times \hat{T}\)-weight spaces. We may assume that $1 \times \hat{T}$ acts on $M$ by a fixed weight $\mu \in X^*(\hat{T})$. Then $M$ extends to \(V_{\hat{G}}\)-representation if only if $\mu \in X^*(\hat{T})^+_{\text{pos}}$, and as a $\hat{G}$-representation, the coaction map sends $M$ into $\fil_\mu \Z[\hat{G}] \otimes_{\Z} M$. The second condition means that $\lambda \leq -w_0(\mu)$ for each weight $\lambda$ of $T \times 1$ appearing in $M$. Since these weights are symmetric under the Weyl group, it follows that \(M\) extends to a \(V_{\hat{G}}\)-representation if and only if each weight $(\lambda, \mu)$ of $\hat{T} \times \hat{T}$ satisfies $\mu + \lambda_{-} \in X^*(\hat{T}_\adj)_{\text{pos}}$, where $\lambda_-$ is the unique anti-dominant weight in the Weyl-orbit of $\lambda$.

We now specialize the previous argument to \(V_{\hat{G},\rho_\adj}\). The composite $\hat{G}\times \Gm \r V_{\hat{G}}$ sends $\fil_{\mu}\Z[\hat{G}]$ into $\fil_{\mu}\Z[\hat{G}] \otimes t^{\langle 2\rho,\mu\rangle}$, where $t$ is the coordinate of $\Gm$. Thus, a representation of $\hat{G} \times \Gm$ extends to $V_{\hat{G},\rho_\adj}$ if and only if for each weight $(\lambda, n)$ appearing, there exists some $\mu \in X^*(\hat{T})^{+}_{\text{pos}}$ such that $\langle2\rho,\mu\rangle = n$ and $\mu + \lambda_- \in X^*(\hat{T}_\adj)_{\text{pos}}$. We claim this condition is equivalent to the following two conditions on the weights $(\lambda, n)$:
$$(-1)^{\langle2 \rho,\lambda\rangle} = (-1)^n \quad \text{and} \quad  \langle 2\rho,\lambda_-\rangle \geq -n. \eqlabel{weight.condition}$$
The necessity of these conditions is straightforward; to see sufficiency take $\mu = -\lambda_- + \nu$, where $\nu \in X^*(\hat{T})_{\text{pos}}$ is any element such that $\langle 2\rho,\nu\rangle = n +\langle 2 \rho,\lambda_-\rangle$.

Under the isomorphism $(\hat{G}\times \Gm)/(2\rho\times \id)(\mu_2) \cong \GT$, $(g,t) \mapsto (g 2\rho(t)^{-1}, t^2)$, a weight $(\lambda, n)$ of $\hat{T} \rtimes \Gm \subset \GT$ pulls back to the weight $(\lambda, -\langle2\rho,\lambda\rangle+2n)$. 
We claim that this 
isomorphism identifies $\Rep_{V_{\hat{G},\rho_\adj}}(\Ab)$ with the subcategory of $\Rep_{\GT}(\Ab)$ of representations with nonnegative $\Gm$-weights. Indeed, the first condition in \refeq{weight.condition} ensures that a representation of $\hat{G}\times \Gm$ descends to the quotient by $(2\rho\times \id)(\mu_2)$. Now suppose we have a representation of $\hat{G}_1$ which extends to the Vinberg monoid. Then for a weight $(\lambda, n)$ appearing in this representation, by \refeq{weight.condition} we must have $\langle 2 \rho,\lambda_-\rangle \geq \langle 2 \rho,\lambda\rangle - 2n$. Thus $ \langle  2\rho,\lambda_- - \lambda \rangle \geq -2n$, so $n \geq \langle  \rho,\lambda - \lambda_- \rangle \geq 0$. Conversely, suppose we have a representation of $\hat{G}_1$ for which all weights $(\lambda, n)$ satisfy $n \geq 0$. Then we must show that the corresponding weight of $\hat{G} \times \Gm$, which is $(\lambda, -\langle 2 \rho,\lambda\rangle + 2n) = (\lambda, m)$, satisfies $\langle 2\rho,\lambda_-\rangle \geq -m$, or equivalently, $n \geq \langle \rho,\lambda - \lambda_- \rangle$. Using the action of the Weyl group on the given representation, and the fact that $(\hat{G}\times \Gm)/(2\rho\times \id)(\mu_2) \cong \GT$ untwists the action of $\Gm$ on the left side, it follows that the weight $(\lambda_-, n - \langle \rho,\lambda - \lambda_-\rangle)$ also occurs (as a representation of $\hat{G}_1$). Hence $n -  \langle \rho,\lambda - \lambda_- \rangle \geq 0$, and we may now conclude using \thref{Sat.Anti}. 
\xpf

\rema
The criterion we obtained for extending a representation to $V_{\hat{G}}$ or $V_{\hat{G},\rho_\adj}$ is equivalent to the condition that the representation extends to the closure of a maximal torus. For normal reductive monoids over an algebraically closed field, this condition is always sufficient, cf. \cite[Remark 5.3]{Renner:Linear}.
\xrema

We conclude with a generalization of \cite[Proposition 5]{Zhu:Integral} for generic Hecke algebras.
Recall that for \(S=\Spec \Fq\) the spectrum of a finite field, the \emph{spherical Hecke algebra} of \(G\) is the ring 
\(\calH^{\sph}_G:=C_\comp(G(\Fq\pot{t})\backslash G(\Fq\rpot{t})/G(\Fq\pot{t}),\Z)\) 
of locally constant, compactly supported, bi-\(G(\Fq\pot{t})\)-invariant, \(\Z\)-valued functions on \(G(\Fq\rpot{t})\), 
equipped with the convolution product \(\star\). This is a free \(\Z\)-module, with a basis given by the 
characteristic functions \(\mathbf{1}_\mu\) of \(G(\pot{t})\mu(t)G(\pot{t})\), for \(\mu\in X_*(T)^+\). 
The convolution is given by
\[\mathbf{1}_\mu\star \mathbf{1}_\lambda = \sum_{\nu\in X_*(T)^+} N_{\mu,\lambda,\nu}(q) \cdot \mathbf{1}_\nu,\]
for uniquely determined polynomials \(N_{\mu,\lambda,\nu}\) with integral coefficients. 
This follows by considering a second basis of \(\calH^{\sph}_G\otimes_{\Z} \Z[q^{\frac{1}{2}}]\) consisting of the \(\phi_\mu:=q^{\langle \rho,\mu \rangle} \chi_\mu\), where \(\chi_\mu\in \Z[X^*(\hat{T})]\) is the character of the irreducible complex representation of \(\hat{G}\) of highest weight \(\mu\). Indeed, the change of basis between \(\{\phi_\mu\}\) and \(\{\mathbf{1}_\mu\}\) is given by integral polynomials in \(q\) \cite[Proposition 4.4]{Gross:Satake}, while the multiplication for the basis \(\{\phi_\mu\}\) is determined by integers independent of \(q\) \cite[Corollary 8.7]{Lusztig:Singularities}. 
This suggests the following definition, generalizing \cite[Definition 6.2.2]{PepinSchmidt:Generic} for \(G=\GL_2\).

\defi
\thlabel{Defi:Generic:Hecke}
Let \(\qq\) be an indeterminate. The \emph{generic spherical Hecke algebra} \(\calH^{\sph}_G(\qq)\) of \(G\) is the free \(\Z[\qq]\)-module with basis \(\{T_\mu\mid \mu\in X_*(T)^+\}\), and multiplication
\[T_\mu\cdot T_\lambda = \sum_{\nu\in X_*(T)^+} N_{\mu,\lambda,\nu}(\qq) T_\nu.\]
\xdefi

To show that this generic spherical Hecke algebra agrees with the representation ring of the Vinberg monoid, we first show that the representation ring is not affected by rationalizing. 
Recall that the representation ring of an (algebraic) group or monoid \(M\) is defined as the Grothendieck group of the category of finitely generated representations of \(M\), and we denote it by \(R(M)\). 

\lemm
\thlabel{Rational:K0}
The rationalization functor \(\Rep_{V_{\hat{G},\rho_\adj}}(\Ab)^\comp\to \Rep_{V_{\hat{G},\rho_\adj}\otimes \Q}(\Vect_\Q)^\comp\) induces an isomorphism on Grothendieck rings. 
\xlemm
\pf Note that the source category consists of representations on finitely generated abelian groups and the target category consists of representations on finite dimensional $\Q$-vector spaces. 
As the rationalization functor is symmetric monoidal, it induces a ring homomorphism on Grothendieck rings. This ring homomorphism is surjective, as \(\Rep_{V_{\hat{G},\rho_\adj}\otimes \Q}(\Vect_\Q)^\comp\cong \Satr^{G,\Q,\anti,\comp}\) is semisimple with simple objects \(\calJ_!^\mu(\Q)(n)\) for \(\mu\in X_*(T)^+\) and \(n\leq 0\), and \(\calJ_!^{\mu}(\Z)(n)\otimes\Q\cong \calJ_!(\Q)(n)\) by \thref{Standard.Flat}.

For injectivity, we first claim that any torsion \(V_{\hat{G},\rho_\adj}\)-representation vanishes in the Grothendieck ring. 
To prove this, we work in the Satake category. 
Let \(\calF\in \Sat_\red^{G,\sgl,\anti}\) be compact and torsion, and \(\mu\in X_*(T)^+\) a maximal element in the support of \(\Ff\).
Then we can find a finite graded abelian group \(L\in \MTMr(S)^{\anti,\comp}\), and a natural map \(\Jj_!^\mu(L)\to \Ff\) whose kernel and cokernel in \(\Sat_\red^{G,\sgl,\anti}\) have support strictly smaller than \(\Ff\).
Thus, by noetherian induction we reduce to the case where $\calF = \calJ_!^\mu(L)$.
Then we can find some degreewise free \(L'\) and a presentation \(0 \to L' \to L' \to L \to 0\). This gives a short exact sequence \(0\to \calJ_!^\mu(L') \to \calJ_!^\mu(L') \to \calJ_!^\mu(L) \to 0\), so that \([\calJ_!^\mu(L)]=0\in K_0(\Rep_{V_{\hat{G},\rho_\adj}}(\Ab)^\comp)\).

Now to prove injectivity, suppose
\(M_1,M_2\in \Rep_{V_{\hat{G},\rho_\adj}}(\Ab)^\comp\) become
isomorphic after rationalizing. Then we can scale such a rational isomorphism so that it preserves the integral subrepresentations. 
The kernel and cokernel of this integral morphism are torsion representations, which are trivial in $K_0$.
\xpf

We can now construct a \emph{generic Satake isomorphism}, 
again generalizing \cite[Theorem 6.2.3]{PepinSchmidt:Generic} for \(G=\GL_2\). (The authors informed us they also knew how to generalize their proof to general split reductive groups.)
Note that the Grothendieck ring of \(\Satr^{G,\sgl,\anti, \comp}\) is naturally a \(\Z[\qq]\)-algebra, 
where multiplication by \(\qq\) corresponds to twisting by \((-1)\).

\coro
\thlabel{Generic:Hecke:Iso}
There is a unique isomorphism \(\Psi\) between \(\calH^{\sph}_G(\qq)\) and the representation ring of \(V_{\hat{G},\rho_\adj}\) such that for any prime power \(q\), the diagram 
\[
\begin{tikzcd}
  \calH_G^{\sph}(\qq) \arrow[r, "\Psi"] \arrow[d, "\qq=q"] & R(V_{\hat{G},\rho_\adj}) \arrow[d, "{[}d_{\rho_\adj}{]}=q"]\\
  \calH_G^{\sph}\otimes_{\Z} \Z[q^{\pm \frac{1}{2}}] \arrow[r, "\Psi_\cl"] & \Z[q^{\pm \frac{1}{2}}][X_*(T)]^{W_0} \cong R(\hat{G}) \otimes_{\Z} \Z[q^{\pm \frac{1}{2}}]
\end{tikzcd}
\eqlabel{Satake:Iso:Diagram}
\]
commutes. Here, \(\Psi_\cl\) denotes the classical Satake isomorphism, cf.~ \cite[Proposition 3.6]{Gross:Satake}, 
and the rightmost map is obtained by taking the character \(R(V_{\hat{G},\rho_\adj})\to R(\GT) \to \Z[X^*(\hat{T}\rtimes \Gm)]\), and then setting the character of the projection \(\hat{T}\rtimes \Gm\to \Gm\) equal to \(q\).
In particular, \(\calH^{\sph}_G(\qq)\) is commutative and unital. 
\xcoro
\pf
By \thref{Rational:K0}, it suffices construct an isomorphism between \(\calH^{\sph}_G(\qq)\) and the Grothendieck ring of 
\(\Rep_{V_{\hat{G},\rho_\adj},\Q}(\Vect_\Q)^\comp \cong \Sat^{G,\Q,\anti,\comp}\), which we denote by \(R^{\anti}\).
Recall that \(\Z[q^{\pm \frac{1}{2}}][X_*(T)]^{W_0}\) admits a natural \(\Z[q^{\pm \frac{1}{2}}]\)-basis given by the characters \(\chi_\mu\) of the simple complex algebraic representations of \(\hat{G}\), for \(\mu\in X_*(T)^+\).
Define \(f_\mu:=\Psi_{\cl}^{-1}(q^{\langle \rho,\mu\rangle} \chi_\mu)\). 
By \cite[(3.12) and Proposition 4.4]{Gross:Satake}, we have 
\(f_\mu = \mathbf{1}_\mu + \sum_{\lambda<\mu} d_{\mu,\lambda}(q) \mathbf{1}_\lambda\), for uniquely determined polynomials \(d_{\mu,\lambda}(\qq)\in \Z[\qq]\). 
Setting \(f_\mu(\qq):= T_\mu + \sum_{\lambda<\mu}d_{\mu,\lambda}(\qq) T_\lambda\), 
we get a second \(\Z[\qq]\)-basis \(\{f_\mu(\qq)\mid \mu\in X_*(T)^+\}\) of \(\calH^{\sph}_G(\qq)\). 

Now, consider the group homomorphism \(\Psi\colon \calH^{\sph}_G(\qq)\to R^{\anti}\) sending 
\(\qq^n\cdot f_\mu(\qq)\) to \([\IC_{\mu,\Q}(-n)]\), where \([-]\) denotes the class of an object in \(R^{\anti}\). 
As the simple objects in \(\Satr^{G,\Q,\anti,\comp}\) are exactly the \(\IC_{\mu,\Q}(-n)\) 
for \(\mu\in X_*(T)^+\) and \(n\in \Z_{\geq 0}\), and because the \(f_\mu(\qq)\) form a \(\Z[\qq]\)-basis of \(\calH^{\sph}_G(\qq)\), 
it follows that \(\Psi\) is both injective and surjective. 
For \(\Psi\) to be an isomorphism of rings, we need to show certain equalities of polynomials. 
But since the classical Satake isomorphism is a ring morphism, the polynomials in question 
agree for all prime powers \(q\), so that they must be equal. 
Hence \(\Psi\) gives the desired isomorphism \(\calH^{\sph}_G(\qq)\cong R^{\anti}\) of rings. 
Since we defined \(\Psi\) using \(\Psi_\cl\) and scaling by some power of \(q\), corresponding to the \(\Gm\)-action on 
\(\hat{T}\) appearing in \(\hat{T}\rtimes \Gm\), we conclude that \refeq{Satake:Iso:Diagram} commutes. 
\xpf

\rema
In \cite[(1.12)]{Zhu:Integral}, Zhu modifies the Satake isomorphism so that it is defined over $\Z$. Using the commutative diagram in \cite[Lemma 25 ff.]{Zhu:Integral}, one can also relate $\Psi$ to \cite[(1.12)]{Zhu:Integral}.
Taking the base change along \(\Z\to \Fp\) recovers the mod \(p\) 
Satake isomorphism as in \cite{Herzig:Satake, HenniartVigneras:Satake}, cf.~ also \cite[Corollary 7]{Zhu:Integral}.
\xrema

\bibliographystyle{alphaurl}
\bibliography{bib}

\end{document}